# HF=HM II: Reeb orbits and holomorphic curves for the ech/Heegaard-Floer correspondence

Cagatay Kutluhan, Yi-Jen Lee and Clifford Henry Taubes

ABSTRACT: This is the second of five papers that construct an isomorphism between the Seiberg-Witten Floer homology and the Heegaard Floer homology of a given compact, oriented 3-manifold. The isomorphism is given as a composition of three isomorphisms; the first of these relates a version of embedded contact homology on an an auxillary manifold to the Heegaard Floer homology on the original. This paper describes this auxilliary manifold, its geometry, and the relationship between the generators of the embedded contact homology chain complex and those of the Heegaard Floer chain complex. The pseudoholomorphic curves that define the differential on the embedded contact homology chain complex are also described here as a first step to relate the differential on the latter complex with that on the Heegaard Floer complex.

This paper and its sequel [KLTIII] supply the geometric and analytic results to relate Peter Ozsváth and Zóltan Szabó's Heegaard-Floer homology [OS1], [OS2] and a version of Michael Hutching's embedded contact homology [Hu1]. The precise relationship is reported in [KLTI].

By way of background, the Heegaard Floer homology of a given, compact and oriented three manifold is computed using a suitably chosen Morse function and associated pseudogradient vector field. This manifold is denoted by M. As explained in what follows, this data is used to construct geometric data on the connect sum of M with a certain number of copies of $S^1 \times S^2$. Use Y in what follows to denote this connect sum but with orientation reversed from that on M. The geometric data on Y can be used to define a stable Hamiltonian version of Michael Hutching's embedded contact homology.

The plan for this paper is to first describe how the data that is needed to compute the Heegaard-Floer homology for M is used to construct Y and the stable Hamiltonian structure needed to compute the relevant version of embedded contact homology for Y. With the respective geometric structures in place, the generators for the relevant version of embedded contact homology for Y are determined. As is explained in Proposition 2.8 below, each generator of the Heegaard Floer chain complex on M determines a set of generators for the relevant embedded contact homology complex on Y. With the generators understood, the balance of this article explains how the pseudoholomorphic curves on Y that are used to compute the embedded contact homology differential determine the data that is used on M to compute the Heegaard Floer differential. This is done using Robert Lipshitz' reformulation of Heegaard-Floer homology [L].

What follows directly is a table of contents.

1. THE GEOMETRY OF M, THE MANIFOLD Y AND THE GEOMETRY ON Y.
    The manifold Y and its stable Hamiltonian structure are constructed from the Heegaard Floer geometry on M.

2. THE CLOSED INTEGRAL CURVES OF $v$.
    The relevant generators of the embedded contact homology chain on Y complex are described.

3. PSEUDOHOLOMORPHIC SUBVARIETIES AND THE GEOMETRY OF $\mathbb{R} \times Y$.
    Certain sorts of almost complex structures on $\mathbb{R} \times Y$ are introduced. They are then used to construct the pseudoholomorphic foliations of $\mathbb{R} \times Y$ that play a central role in the subsequent sections and in [KLTIII]

4. ECH-HF SUBVARIETIES.
    This section introduces the pseudoholomorphic subvarieties that define the differential and various other endomorphisms of Y's embedded contact homology. These are the ech-HF subvarieties.

5. QUANTATIVE CONCERNS.

   This section states and then proves a crucial aprior bound on the value of the integral of the stable Hamiltonian 2-form over an ech-HF subvariety. This bound plays a central role in all of the subsequent analysis in the article and in [KLTIII].

6. HEEGAARD FLOER CURVES

   This section constitutes a digression to describe the sorts of pseudoholomorphic subvarieties that are used in Lipshitz' formulation of Heegaard Floer homology on M.

7. ECH-HF SUBVARIETIES AND HEEGAARD FLOER CURVES.

   This section defines a correspondence that maps any given ech-HF subvariety to a subvariety of the sort that is used by Lipshitz to describe the Heegaard Floer homology of M.

A. APPENDIX.

   The Appendix proves that the relevant version of embedded contact homology on Y can be defined using the rules laid out by Hutchings.

This article supplies the proofs of Theorems 2.1 and 2.2 in [KLTI]; the former is restated as Proposition 2.8 and proved in Section 2 and the latter is restated in the Appendix where it is proved.

The following notational conventions are used: What is denoted by $c_0$ is in all cases is a constant in $(1, \infty)$ whose value is independent of all relevant parameters. The value of $c_0$ can increase between subsequent appearances. A second convention concerns a function that is denoted by $\chi$. This is a fixed, non-increasing function on $\mathbb{R}$ with value 1 on $(-\infty, 0]$ and equal to 0 on $[1, \infty)$.

**Acknowledgements**

The authors owe much to Michael Hutchings for sharing his thoughts about embedded contact homology. Another debt is owed to Robert Lipshitz for sharing thoughts about his version of Heegaard Floer homology.

The authors also acknowledge the generous support of the Mathematical Sciences Research Institute where the three authors did most of the work reported here. The second and third authors were also supported in part by the National Science Foundation. The third author was also supported by a David Eisenbud Fellowship.

## 1. The geometry of M, the manifold Y and the geometry of Y

Let M denote a compact, oriented 3-manifold with a self-indexing Morse function, $f: M \to [0, 3]$ with one index 0 critical point, and one index 3 critical point. Let G denote the number of index 1 critical points; this is the same as the number of index 2 critical points. The function $f$ is used in what follows to construct what is called a *stable Hamiltonian structure* on $Y = M \#_{G+1} (S^1 \times S^2)$. This geometric data on Y consists of a pair $(a, w)$ of 1-form and 2-form with $dw = 0$, $da \subset \text{span}(w)$ and $a \wedge w$ is nowhere zero. The orientation for Y is chosen so that $a \wedge w > 0$. The final subsection summarizes the most relevant aspects of the resulting geometry.

The constructions require first fixing a 1-1 correspondence between the set of index 2 critical points of $f$ and index 1 critical points. Let $\Lambda$ denote the resulting ordered set of paired critical points. The connect sum is defined by attaching 1-handles at paired index 2 and index 1 critical points, and then at the index 0 and index 3 critical points. with the handle attached by removing small radius balls about the respective members of each pair.

**a) Attaching the handles**

This part of the section explains how the 1-handles are attached to M so as to obtain Y. There are three parts to the discussion. The first part describes the 1-handle that is labeled by a given pair from $\Lambda$. The second part describes the 1-handle that is labeled by the index 0 and 3 critical points. The final part introduces some notational conventions that are used implicitly throughout this article. The definitions that follow involve a chosen constant $\delta_* \in (0, c_0^{-1})$. An upper bound for $\delta_*$ is determined by M and $f$. Any choice below this upper bound will work.

*Part 1*: Let $p \in M$ denote an index 1 critical point. Take coordinates $(x, y, z)$ near p so that $f$ appears as

$$f = 1 + x^2 + y^2 - 2z^2.$$
(1.1)

Fix $\delta_* \in (0, \tfrac{1}{100})$ such that these coordinates are defined for $(x^2 + y^2 + z^2)^{1/2} \leq 10\delta_*$. By way of notation, if $d \in (0, 10\delta_*)$, the ball in this coordinate chart where $(x^2 + y^2 + z^2)^{1/2} < d$ is said to be the *radius d coordinate ball* centered on p. Implicit in the term 'coordinate ball' is the use of the special coordinate system where $f$ appears as above. Introduce the spherical coordinates $(r_+, (\theta_+, \varphi_+))$ of the usual sort: The Euclidean coordinates $(x, y, z)$ are given in terms of the spherical ones by the rule $x = r_+ \sin\theta_+ \cos\varphi_+$, $y = r_+ \sin\theta_+ \sin\varphi_+$, and $z = r_+ \cos\theta_+$. It proves convenient in what follows to define a function $t_+ = \ln r_+$.

Let $p'$ denote an index 2 critical point. Take coordinates $(x, y, z)$ here so that the function $f$ appears as

$$f = 2 - (x^2 + y^2 - 2z^2).$$
(1.2)

Introduce again radial coordinates $(r_-, (\theta_-, \varphi_-))$. The constant $\delta_*$ can be chosen so that these coordinates are defined where the radial coordinate $r_- \leq 10\delta_*$. By analogy with the case of an index 1 critical point, given $d \in (0, 10\delta_*)$, the ball in this coordinate chart where $r_- < d$ is said to be the radius d coordinate ball centered at p´. Use $t_-$ to denote $\ln r_-$.

Let $\mathfrak{p} = (p, p')$ denote a pair from $\Lambda$. The convention here and subsequently has the left entry, p, the index 1 critical point and the right entry, p´, the index 2 critical point. The pair $\mathfrak{p}$ labels one of the 1-handle that is used to obtain Y from M. This handle is denoted by $\mathcal{H}_\mathfrak{p}$. The attaching of the 1-handles requires the choice of a parameter, this a positive number R. The number R is constrained so that $R \gg -100\ln\delta_*$. There are no upper bound constraints. What follows explains how to attach the 1-handle $\mathcal{H}_\mathfrak{p}$.

A 1-handle is by definition diffeomorphic to $[-1, 1] \times S^2$ and it is attached by identifying the respective boundaries to constant radii spheres in the radius $8\delta_*$ coordinate balls centered on p and p´. To say more, introduce coordinates $(u, (\theta, \phi))$ for the 1-handle where $u \in [-R + \ln(7\delta_*), R - \ln(7\delta_*)]$ and where $(\theta, \phi)$ are the spherical coordinates on $S^2$. The 1-handle $\mathcal{H}_\mathfrak{p}$ is defined by making the following identifications:

- *Identify the* $t_+ \in [-2R - \ln(7\delta_*), \ln(7\delta_*)]$ *part of the radius $7\delta_*$ coordinate ball centered on p with $\mathcal{H}_\mathfrak{p}$ by writing* $(t_+ = u - R, \theta_+ = \theta, \varphi_+ = \phi)$

- *Identify the* $t_- \in [-2R - \ln(7\delta_*), \ln(7\delta_*)]$ *part of the radius $7\delta_*$ coordinate ball centered on p´ with $\mathcal{H}_\mathfrak{p}$ by writing* $(t_- = -u - R, \theta_- = \pi - \theta, \varphi_- = \phi)$

(1.3)

Adding the G handles $\{\mathcal{H}_\mathfrak{p}\}_{\mathfrak{p} \in \Lambda}$ gives a new manifold that is diffeomorphic to the connect sum of M with G copies of $S^1 \times S^2$.

*Part 2*: Fix Euclidean coordinates on respective small radius balls about the index 0 and index 3 critical points of $f$ so that $f$ appears as

$$f = x^2 + y^2 + z^2 \text{ and } f = 3 - x^2 - y^2 - z^2.$$

(1.4)

Take $\delta_*$ so that these coordinates are defined for $(x^2 + y^2 + z^2)^{1/2} \leq 10\delta_*$. Use $B_0$ and $B_3$ to denote the respective radius $10\delta_*$ coordinate balls about the index 0 and index 3 critical points. If $d \in (0, 10\delta_*)$, the subset of either ball where $(x^2 + y^2 + z^2)^{1/2} < d$ is called the radius d coordinate ball about the critical point in question. Introduce the spherical coordinates $(r_+, \theta_+, \phi_+)$ for $B_0$ and $(r_-, \theta_-, \phi_-)$ for $B_3$. Use $t_+$ and $t_-$ to denote $\ln r_+$ and $\ln r_-$.

The attaching spheres for the relevant 1-handle are the radius $7\delta_*$ coordinate balls centered on these two critical points. The handle is again parametrized using coordinates $(u, \theta, \phi)$ where $u \in [-R - \ln 7\delta_*, R + \ln 7\delta_*]$ and with $(\theta, \phi)$ the spherical angles. The handle attachment here is defined by the (1.3) with it understood that p is the index 0 critical point and p´ is the index 1 critical point.

This 1-handle in Y is denoted in what follows by $\mathcal{H}_0$.

*Part 3*: This part introduces three conventions. First, the various handles from the set $\{\mathcal{H}_0, \{\mathcal{H}_\mathfrak{p}\}_{\mathfrak{p}\in\Lambda}\}$ are identified with their images in Y. To state the second, suppose that $r \in [e^{-R}, 7\delta_*)$ The complement of the union of the $u \in (-R - \ln r, R + \ln r)$ part of each 1-handle defines a submanifold with boundary in Y. The latter is denoted by $M_r$. The identifications depicted in (1.3) identify $M_r$ with the complement in M of the union of the radius r coordinate balls about each of the critical points of $f$. This part of M is also denoted by $M_r$ and these two versions of $M_r$ are, for the most part, identified implicitly in all that follows. The third convention involves the function $f$. The latter is defined on M, and so it is defined on the incarnation of $M_r$ in M. As a consequence, the function $f$ is defined also on the incarnation of $M_r$ in Y.

**b) Geometry of Y: The pair $(a, w)$ on $\cup_{\mathfrak{p}\in\Lambda} \mathcal{H}_\mathfrak{p}$**

The construction of the desired pair $(a, w)$ of 1-form and 2-form on $\cup_{\mathfrak{p}\in\Lambda} \mathcal{H}_\mathfrak{p}$ requires specification of two additional parameters. The first is a positive number, this denoted by $\delta$. This number is constrained only by the upper bound $\delta < c_0^{-1}\delta_*$. The second parameter is also a positive number, this denoted by $x_0$. With $\delta$ chosen, the constant $x_0$ is constrained only by an upper bound $x_0 < \delta^3$. With $(\delta, x_0)$ chosen, the lower bound for the parameter R must be revised upwards to $R > -c_0 \ln x_0$. The constant $c_0$ in both cases depends only on the particular choice for the function $\chi$.

Fix $\mathfrak{p} \in \Lambda$. The specification of $(a, w)$ on $\mathcal{H}_\mathfrak{p}$ has three parts.

*Part 1*: The definitions require the introduction of three new functions of the coordinate u. All are defined using the chosen, non-increaing function $\chi: \mathbb{R} \to [0, 1]$. Recall that the latter equals 1 on $(-\infty, 0]$ and it equals 0 on $[1, \infty)$. The first of the new functions is denoted by $x$, and it is given by the rule

$$x(u) = x_0 \chi(|u| - R - \ln\delta - 12)$$

(1.5)

Thus, $x = 0$ where $|u| \geq R + \ln\delta + 11$ and $x = x_0$ where $|u| \leq R + \ln\delta + 12$. The second and third functions are denoted respectively by $\chi_-$ and $\chi_+$. These are defined by the rules

$$\chi_-(u) = \chi(u + \tfrac{1}{4} R) \quad and \quad \chi_+(u) = \chi(-u - \tfrac{1}{4} R) .$$

(1.6)

Thus, $\chi_+ = 0$ where $u \leq -\frac{1}{4}R - 1$ and $\chi_+ = 1$ where $u \geq -\frac{1}{4}R$. Note that this function is nondecreasing. Meanwhile, $\chi_- = 0$ where $u > \frac{1}{4}R + 1$ and $\chi_- = 1$ where $u \leq \frac{1}{4}R$ and so $\chi_-$ is nonincreasing. Use $\chi_+'$ and $\chi_-'$ in what follows to denote their respective u-derivatives.

*Part 2*: Define the 1-form $a$ by the rule

$$a = (x + (2\chi_+ + \chi_+')e^{2(u-R)} + (2\chi_- - \chi_-')e^{-2(u+R)})(1 - 3\cos^2\theta)\,du$$
$$- \sqrt{6}\,(x + 2(\chi_+ e^{2(u-R)} + \chi_- e^{-2(u+R)}))\cos\theta\sin^2\theta\,d\phi$$
$$+ 6(\chi_+ e^{2(u-R)} - \chi_- e^{-2(u+R)})\cos\theta\sin\theta\,d\theta,$$
(1.7)

and the 2-form $w$ on $\mathcal{H}_p$ by the rule $w = da$. The admittedly complicated formula for $a$ endows $a$ and $w$ with various desirable properties of which two are that $w$ is closed and $a \wedge w > 0$. The fact that $dw = 0$ follows directly from (1.8). The assertion that $a \wedge w > 0$ holds provided that $\delta < c_0^{-1}$, $x_0 < \delta^3$ and that $R > -c_0 \ln x_0$.

To verify this last claim, introduce on $\mathcal{H}_p$ the 1-form

$$a_* = (1 - 3\cos^2\theta)\,du - \sqrt{6}\cos\theta\sin^2\theta\,d\phi$$
(1.8)

and the function

$$f_* = (\chi_+ e^{2(u-R)} - \chi_- e^{-2(u+R)})(1 - 3\cos^2\theta)$$
(1.9)

so as to write

$$a = x a_* - 2\sqrt{6}(\chi_+ e^{2(u-R)} + \chi_- e^{-2(u+R)})\cos\theta\sin^2\theta\,d\phi + df_*.$$
(1.10)

The 1-form $a_*$ is contact 1-form. This observation with (1.7) imply that $a \wedge w > 0$ where $|u| \leq R + \ln x_0 - c_0$. This is because $x = x_0$ and so $a$ can be written as $a = x_0 a_* + \mathfrak{r}$ and $w$ as $w = x_0\,da_* + \mathfrak{e}$ where $\mathfrak{r}$ and $\mathfrak{e}$ contribute to $x_0^2 a_* \wedge da_*$ in the resulting expression for $a \wedge w$ by a term whose absolute value is at most $c_0 e^{-2(|u|-R)}(x_0 + e^{-2(|u|-R)})$. To continue, introduce $f_- = e^{2(u-R)}(1 - 3\cos^2\theta)$. Then

- $df_- \wedge d(-2\sqrt{6}e^{2(u-R)}\cos\theta\sin^2\theta\,d\phi) = 4\sqrt{6}e^{4(u-R)}(1 + 3\cos^4\theta)\,du\sin\theta\,d\theta\,d\phi$,
- $df_- \wedge da_* = 2\sqrt{6}e^{2(u-R)}(1 - 3\cos^2\theta)^2\,du\sin\theta\,d\theta\,d\phi$.

(1.11)

The top line in (1.11) implies that $a \wedge w > 0$ where $x = 0$, and it implies that $a \wedge w > 0$ where $dx \neq 0$ if $\delta < c_0^{-1}$ and $x_0 \leq \delta^3$. The two bullets together with the fact that $a_*$ is a contact form implies that $a \wedge w > 0$ where $\chi_\delta = 1$.

### c) Geometry of Y: The pair $(a, w)$ on $f^{-1}([1, 2]) \cap M_\delta$

The definition of $a$ and $w$ on $M_\delta$ requires the choice of additional data. The first of these is a class $c_{1M}$ in $H^2(M; 2\mathbb{Z})$/torsion. The second is a suitably constrained *pseudo-gradient* vector field on M for the function $f$. Having chosen such a vector field, let $\mathfrak{v}$ denote the normalized version that pairs with $df$ to give 1. The vector field $\mathfrak{v}$ is defined on the complement of the critical points of $f$. There are various constraints to the choice for $\mathfrak{v}$; these are described in what follows when they are needed. The construction has nine parts.

*Part 1*: There are constraints on $\mathfrak{v}$ near the index 1 and index 2 critical points of $f$ on M. These are described here. Let p denote an index 1 critical point. Reintroduce the coordinates $(t_+, (\theta_+, \varphi_+))$ for the radius $8\delta_*$ coordinate ball centered at p.

CONSTRAINT 1: *The vector field $\mathfrak{v}$ on the radius $8\delta_*$ coordinate ball centered at p is*

$$\mathfrak{v} = \frac{1}{2(1+3\cos^4\theta_+)} e^{-2t_+}((1 - 3\cos^2\theta_+)\tfrac{\partial}{\partial t_+} + 2\cos\theta_+ \sin\theta_+ \tfrac{\partial}{\partial \theta_+}) \, .$$

*The analogous formula for $\mathfrak{v}$ on the radius $8\delta_*$ coordinate ball centered on a given index 2 is obtained from (1.12) by first replacing $(t_+, \theta_+)$ with $(t_-, \theta_-)$ and then multiplying the resulting expression by -1.*

The fact that $\mathfrak{v}(f) = 1$ can be seen by writing $df$ on this part of M in the $(t_+, (\theta_+, \varphi_+))$ coordinates:

$$df = 2\, e^{2t_+}((1 - 3\cos^2\theta_+) dt_+ + 3\cos\theta_+ \sin\theta_+ d\theta_+) \, .$$

(1.12)

The analogous equation on the radius $8\delta_*$ coordinate ball centered on an index 2 critical point of $f$ is obtained from (1.12) by replacing $t_+$ with $t_-$ and multiplying the resulting expression by -1.

It follows from Constraint 1 that $\mathfrak{v}$ spans the kernel of the 2-form $w$ on the part of $M_\delta$ in the radius $7\delta_*$ coordinate ball centered on any given index 1 or index 2 critical point of $f$. To see why this is, let p denote an index 1 critical point. The claim follows via a straightforward computation from the fact 2-form $w$ appears where $t_+ > \ln \delta$ in the radius $7\delta_*$ coordinate ball as

$$w = -2\sqrt{6}\, d\,(e^{2t_+} \cos\theta_+ \sin^2\theta_+ \, d\varphi_+)$$

(1.13)

The corresponding formula in the case when p is an index 2 critical point verifies that $\mathfrak{v}$ also generates $w$'s kernel in the $t_- > \ln\delta$ part of the radius $4\delta_*$ coordinate ball.

*Part 2*: Introduce $\Sigma = f^{-1}(\frac{3}{2})$ to denote the Heegaard surface. The construction of $a$ and $w$ involve what Ozsváth and Szabó call a Heegaard diagram of $\Sigma$. Such a diagram is determined by the choice of $\mathfrak{v}$. In particular, the construction of Heegaard Floer homology in [OS1] requires special sorts of Heegaard diagrams, and so constrains on $\mathfrak{v}$. As is evident in what follows, the same sorts of constraints are needed here to construct $a$ and $w$.

CONSTRAINT 2: *The ascending disks from the index 1 critical points of $f$ should have transversal intersection with the descending disks from the index 2 critical points.*

Let p denote a given index 1 critical point of $f$. The ascending disk from p intersects $\Sigma$ as an embedded circle, this denoted by $C_{p+}$. Let p´ denote a given index 2 critical point. The descending disk from p´ intersects $\Sigma$ in an embedded circle, this denoted by $C_{p'-}$. These two circles intersect transversally if and only if the corresponding ascending disk and descending disk intersect transversally. By way of notation, $C_+$ is used in what follows to denote the union of the index 1 critical point versions of $C_{p+}$ and $C_-$ is used to denote the union of the index 2 critical point versions of $C_{p-}$.

The Lie transport by $\mathfrak{v}$ identifies the $\theta_+ = \frac{\pi}{2}$ equatorial circle in any given constant $t_+ \geq \ln\delta$ sphere in the radius $8\delta_*$ coordinate ball with $C_{p+}$, and it identifies the part of such a sphere where $1 - 3\cos^2\theta_+ > 0$ with an annular neighborhood of this same circle. Let $T_{p+}$ for the moment denote the annular neighborhood of $C_{p+}$ that corresponds to the annulus in the radius $7\delta_*$ coordinate ball where $1 - 3\cos^2\theta_+ > 0$. Choosing $\delta_* < c_0^{-1}$ will guarantee that distinct index 1 critical points define respective versions of $T_{p+}$ with disjoint closures in $\Sigma$.

*Part 3*: Introduce the function $\hat{h}_+ = 2e^{2t_+}\cos\theta_+\sin^2\theta_+$ on the $t_+ > \ln\delta$ part of the radius $8\delta_*$ coordinate ball centered on p. This function has no critical points, and it is such that $\mathfrak{v}(\hat{h}_+) = 0$. Likewise, $\mathfrak{v}(\varphi_+) = 0$. This understood, Lie transport by $\mathfrak{v}$ of the functions $(\varphi_+, \hat{h}_+)$ identifies the latter as coordinates on $T_{p+}$. Given that (1.13) writes the 2-form $w$ as $\sqrt{6}d\varphi_+ \wedge d\hat{h}_+$, this same Lie transport identifies $w$ with an area form on $T_{p+}$.

Meanwhile, when p´ is a given index 2 critical point of $f$, then Lie transport by $\mathfrak{v}$ defines a corresponding annular neighborhood, $T_{p´-}$ of the circle $C_{p´-}$. As with the index 1 critical point case, distinct index 2 critical points can be assumed to define respective versions of $T_{p´-}$ with disjoint closures in $\Sigma$. Given an index 2 critical point p´, Lie transport by $\mathfrak{v}$ endows $T_{p´-}$ with the corresponding coordinates $(\varphi_-, \hat{h}_-)$, and gives the latter the area form $-\sqrt{6}\,d\varphi_- \wedge d\hat{h}_-$.

Let p and p´ denote a respective pair of index 1 and index 2 critical point of $f$. If $\delta_* < c_0^{-1}$, then each component of the intersection between $T_{p+}$ and $T_{p´-}$ will appear in suitable coordinates as the intersection in $\mathbb{R}^2$ of $(-1, 1) \times \mathbb{R}$ with $\mathbb{R} \times (1, -1)$. Make such a choice for $\delta_*$.

CONSTRAINT 3: *A suitable choice for $\mathfrak{v}$ guarantees that $(d\varphi_+, d\hat{h}_+) = \pm(d\hat{h}_-, d\varphi_-)$ on the intersection of $T_{p+}$ with $T_{p´-}$ with the + sign taken when the pair of vectors $(\frac{\partial}{\partial \varphi_+}, \frac{\partial}{\partial \varphi_-})$ define an oriented basis for $T\Sigma$ at the corresponding point in $C_{p+} \cap C_{p´-}$. This constraint guarantees that the area form $\sqrt{6}\,d\varphi_+ \wedge d\hat{h}_+$ for $T_{p+}$ is equal to the area form $-\sqrt{6}\,d\varphi_- \wedge d\hat{h}_-$ for $T_{p-}$ on $T_{p+} \cap T_{p´-}$.*

There are no obstructions to choosing $\mathfrak{v}$ so as to satisfy this contraint.

*Part 4*: Introduce $T_+$ to denote the union of the various index 1 critical point versions of $T_{p+}$ and introduce $T_-$ to denote the union of the various index 2 critical point versions of $T_{p-}$. Fix a point in the complement of the closure of $T_+ \cup T_-$. This point is denoted by $z_0$.

The data consisting of $\{z_0, C_+, C_-\}$ defines what Ozsváth and Szabó call a *pointed Heegaard diagram*. A given class in $H^2(M; \mathbb{Z})$/torsion determines what they call a *strongly admissible Heegaard diagram*. Required here is a certain sort of strongly admissible diagram for the chosen class $c_{1M}$. Diagrams of the desired sort are constructed in Lemma 5.4 of [OS1].

Some terminology from [OS1] is required to say more about how the needed Heegaard diagrams are used. A fundamental domain is a component of $\Sigma - (C_+ \cup C_-)$. The fundamental domain with the base point is denoted by $\mathcal{D}_0$. The set of fundamental domains is ordered and labeled as $\{\mathcal{D}_i\}_{i=0,\ldots,N}$. Each fundamental domains is viewed in what follows as a 2-chain on $\Sigma$. A periodic domain $\mathcal{P} \subset \Sigma$ is a 2-chain whose boundary is an integer weighted sum of components of $C_+$ and $C_-$ that can be written as the formal sum $\mathcal{P} = \sum_{0 \le i \le N} z_i \mathcal{D}_i$ where $z_i \in \mathbb{Z}$ and where the coefficient $z_0$ of the distinguished domain $\mathcal{D}_0$ is zero. The vector space of periodic domains is isomorphic to $H_2(M; \mathbb{Z})$ with the isomorphism given by adding to the 2-chain $\mathcal{P}$ a corresponding weighted union of ascending and descending disks of the respective index 1 and index 2 critical points to create a closed 2-chain in M.

The lemma that follows brings the notion of a strongly admissible Heegaard diagram into the story.

**Lemma 1.1**: *There exist strongly admissible Heegaard diagrams for the class $c_{1M}$ such that $\Sigma$ admits an area form with the following properties:*
- *The area of $\Sigma$ is 2.*
- *The integral of $w_\Sigma$ over any given periodic domain in $\Sigma$ is equal to the value of $c_{1M}$ on the corresponding homology class in $H_2(M;\mathbb{Z})$.*

*Proof of Lemma 1.1*: When $\mathcal{P}$ is a periodic domain, use $c_{1M}(\mathcal{P})$ to denote the value of the class $c_{1M}$ on the class determined by $\mathcal{P}$ in $H_2(M;\mathbb{Z})$. This is an even integer. Ozsváth and Szabó introduce in Definition 5.3 of [OS1] the notion of a renormalized, $c_{1M}$-periodic domain. This is a 2-chain of the form $\mathcal{Q} = \mathcal{P} - \frac{1}{2} c_{1M}(\mathcal{P})\,[\Sigma]$ with $\mathcal{P}$ being a periodic domain. Such a domain defines a closed cycle of the form $\mathcal{Q} = \sum_{0 \leq i \leq N} z_i \mathcal{D}_i$. Their proof of Lemma 5.4 in [OS1] constructs a $c_{1M}$-strongly admissible, pointed Heegaard diagram for which each renormalized, $c_{1M}$-periodic domain has both positive and negative coefficients. They assert at the very start of their proof that this is what the proof will obtain.

Granted such a Heegaard diagram, the argument used for the proof of Lemma 4.12 in [OS1] can be repeated to find an area 2-form on $\Sigma$ such that each renormalized, $c_{1M}$-periodic domain has total signed area equal to zero. With $\mathcal{Q}$ written as above, this means that $\sum_{0 \leq i \leq N} z_i A_i = 0$ where $A_i$ here denotes the area of the corresponding fundamental domain $\mathcal{D}_i$. Here is the argument: The vector space over $\mathbb{R}$ spanned by the renormalized $c_{1M}$-periodic domains generates a vector subspace of $\mathbb{R}^{N+1}$ with dimension equal to the dimension of $H_2(M;\mathbb{R})$. This vector subspace intersects the positive quadrant only at the origin. As a consequence, there is a vector in the interior of the positive quadrant that is orthogonal to all vectors in this span. This vector can be written as $\sum_i A_i \mathcal{D}_i$ with $A_i > 0$. The desired area form can be any area form for $\Sigma$ such that the area of any given fundamental domain $\mathcal{D}_i$ is the corresponding coefficient $A_i$.

Normalize this area form so that $\sum_i A_i = 2$. This implies that the area of $\Sigma$ is 2. Let $\mathcal{P}$ denote a periodic domain. As the signed area of $\mathcal{Q} = \mathcal{P} - \frac{1}{2} c_{1M}(\mathcal{P})\,[\Sigma]$ is zero, it follows that the signed area of $\mathcal{P}$ is $c_{1M}(\mathcal{P})$. This being the case, the area form obeys the conditions asserted by the lemma.

*Part 5*: With Lemma 1.1 in hand, here is the next constraint on $\mathfrak{v}$:

CONSTRAINT 4: *The pseudogradient vector field $\mathfrak{v}$ gives a Heegaard diagram of the sort described in Lemma 1.1.*

Let $w_\Sigma$ denote an area form for $\Sigma$ of the sort described by Lemma 1.1. There is no obstruction to choosing $w_\Sigma$ so that it agrees on each index 1 critical point version of $T_{p+}$ with the area form $\sqrt{6}\, d\varphi_+ \wedge d\hat{h}_+$ if $\delta_* < c_0^{-1}$. It follows from what is said in Part 3 that there is no obstruction to also requiring that $w_\Sigma = -\sqrt{6}\, d\varphi_- \wedge d\hat{h}_-$ on each index 2 critical point version of $T_{p-}$. Fix such a form $w_\Sigma$. Note the following important point: The form $w_\Sigma$ can be assumed independent of the data $\delta$, $x_0$, and R that are used to define the geometry of the various 1-handles.

Use the integral curves of $\mathfrak{v}$ to define a diffeomorphism between $(1, 2) \times \Sigma$ and the domain $f^{-1}(1, 2) \subset M$. Let $t$ denote the Euclidean coordinate on the interal $[1, 2]$. This diffeomorphism identifies $\mathfrak{v}$ with the vector field $\frac{\partial}{\partial t}$ and $f$ with the coordinate $t$. Pull-back via the projection from $[1, 2] \times \Sigma$ to $\Sigma$ defines $w_\Sigma$ as a 2-form on the whole of $[1, 2] \times \Sigma$ and thus on the $f^{-1}([1, 2])$ part of $M_\delta$. It follows from (1.13) that this form extends $w$ from $\cup_{\mathfrak{p} \in \Lambda} \mathcal{H}_\mathfrak{p}$ to the latter's union with the $f^{-1}([1, 2])$ part of $M_\delta$. This extension is denoted by $w$ also.

*Part 6*: With $w$ extended to $f^{-1}([1, 2]) \cap M_\delta$, turn now to the task of extending the 1-form $a$. The next lemma is needed for this. The lemma uses $b_1$ to denote the dimension of $H^1(M; \mathbb{R})$, the first Betti number of M.

**Lemma 1.2**: *There exists a set of $b_1$ points in $\Sigma - (T_+ \cup T_-)$ and a 1-form, $a_\Sigma$, that is defined on the complement of $z_0$ and this set which has the following properties:*
- $da_\Sigma = w_\Sigma$.
- *The integral of $a_\Sigma$ is zero on each component of $C_+$ and on each component of $C_-$.*

*Proof of Lemma 1.2*: Let $\mathcal{I}$ denote the free $\mathbb{Z}$ module of rank 2G that is generated by the components of $C_+$ and $C_-$. The Mayer-Vietoris exact sequence for the cohomology of M defines an exact sequence

$$0 \to H^1(M; \mathbb{Z}) \to H^1(\Sigma; \mathbb{Z}) \to \mathrm{Hom}(\mathcal{I}; \mathbb{Z}) \to H^2(M; \mathbb{Z}) \to 0$$

(1.14)

There is an analogous sequence with $\mathbb{R}$ coefficients. The latter version is used for what follows.

Use the fact that $H^2_{\text{De Rham}}(\Sigma - z_0) = 0$ to choose a smooth 1-form on $\Sigma - z_0$ whose exterior derivative is $w_\Sigma$. Let $a_0$ denote the latter. Integration of $a_0$ over the various components of $C_+$ and $C_-$ define an element $\underline{a}_0 \in \text{Hom}(\mathcal{I}; \mathbb{R})$. Let F denote the homomorphism in the $\mathbb{R}$-coefficient version of (1.14) from $\text{Hom}(\mathcal{I}; \mathbb{R})$ to $H^2(M; \mathbb{R})$. If $F(\underline{a}_0) = 0$, then it follows from (1.14) that $a_0$ can be modified by adding a suitable closed 1-form on $\Sigma$ so that its integral over the various components of $C_+$ and $C_-$ are zero.

For each $i \in \{1, \ldots, N\}$, fix a point in the domain $\mathcal{D}_i$ that sits in the complement of the latter's intersection with the closure of $T_+ \cup T_-$. Let $z_i$ denote this point. Any given point $R = (R_1, \ldots, R_N) \in \mathbb{R}^N$ determines a homomorphism from $H_2(M; \mathbb{R})$ to $\mathbb{R}$ as follows: Represent the given class as a periodic domain, $\mathcal{P}$. Write $\mathcal{P}$ as $\sum_{1 \le i \le N} z_i \mathcal{D}_i$. The homomorphism in question sends $\sum_{1 \le i \le N} R_i z_i$. This then defines a surjective linear map $\psi : \mathbb{R}^N \to H^2(M; \mathbb{R})$.

Fix $b_1$ points from $\{z_i\}_{1 \le i \le N}$ so that the restriction of $\psi$ to the corresponding $b_1$-dimensional subspace in $\mathbb{R}^N$ is an isomorphism. Let $\mathbb{V}$ denote this subspace. Meanwhile, the given $a_0$ can be modified on $\Sigma - \{z_i\}_{0=1,\ldots,N}$ so that its integral over the boundary of each $i \in \{1, \ldots, N\}$ version of $\mathcal{D}_i$ is $R_i$. This modification changes $F(\underline{a}_0)$ to $F(\underline{a}_0) + \psi(R)$. This understood, and given what was said in the previous paragraph, there exists $R \in \mathbb{V}$ such that the resulting version of $a_0$ has integral zero on each component of $C_+$ and on each component of $C_-$.

*Part 7*: Let $¥ \subset \Sigma$ denote the set composed of $z_0$ and the $b_1$ points given by Lemma 1.2. Let $a_\Sigma$ denote the resulting 1-form. The constant $\delta_*$ can be chosen so that there is a collection $\{D_z\}_{z \in ¥}$ of disjoint disks in $\Sigma - (T_+ \cup T_{p-})$, each with area $\frac{\pi}{100(1+b_1)} \delta_*^2$, and with the following additional properties: Fix $z \in ¥$. Then $z \in D_z$ and there exists a real number $R_z$ and radial coordinates $(\rho, \varphi)$ for $D_z$ that make $w_\Sigma$ and $a_\Sigma$ appear as

$$w_\Sigma = \rho \, d\rho \, d\varphi \quad \text{and} \quad a_\Sigma = (\tfrac{1}{\pi} R_z + \tfrac{1}{2} \rho^2) \, d\varphi \, .$$

(1.15)

Granted the preceding, define a smooth 1-form, $\hat{a}_\Sigma$, on $\Sigma$ as follows: Set $\hat{a}_\Sigma$ to equal $a_\Sigma$ on $\Sigma - (\cup_{z \in ¥} D_z)$. Fix $z \in ¥$. To define $\hat{a}_\Sigma$ on $D_z$, first define the function $\chi_\#$ of the coordinate $\rho$ by the rule $\chi_\#(\rho) = \chi(1 - \frac{1000(1+b_1)}{\delta_*} \rho)$. This function is zero near $z$, but 1 on the complement of a compact set in $D_z$. Set

$$\hat{a}_\Sigma = \chi_\# (\tfrac{1}{\pi} R_z + \tfrac{1}{2} \rho^2) \, d\varphi \, .$$

(1.16)

Note for future reference that $d\hat{a}_\Sigma = H w_\Sigma$ where $H$ is a smooth function that equals 1 on the complement of $\cup_{z \in ¥} D_z$.

*Part 8*: Let p denote an index 1 or index 2 critical point of $f$. When $r \in (0, 7)$, use respectively $T^r_{p+} \subset T_{p+}$ or $T^r_{p-} \subset T_{p-}$ to denote the subannulus that comes from the radius $r\delta_*$ coordinate ball centered on p.

Suppose that p is an index 1 critical point of $f$. As can be seen from (1.18) and (1.12), the forms $a$ and $df$ restrict to the $r_+ > \delta$ part of the radius $7\delta_*$ coordinate ball centered on p as a 1-form that is invariant with respect to the Lie transport defined by $\mathfrak{v}$. Lie transport by $\mathfrak{v}$ of $a - df$ from this part of p's radius $7\delta_*$ coordinate ball defines a 1-form, $a_{p+}$, on $T^7_{p+}$ that obeys $da_{p+} = w_\Sigma$. If p is an index 2 critical point of $f$, then the analogous construction defines from $a - df$ a 1-form, $a_{p-}$, on $T^7_{p-}$.

Let p and p´ now denote respective index 1 and index 2 critical points of $f$ and suppose that $T \subset T_{p+} \cap T_{p´-}$ is a given component. As $T_{p+} \cap T_{p´-}$ is contractible and as both $a_{p+}$ and $a_{p-}$ are anti-derivatives of $w_\Sigma$ on T, their restriction to $T^7_{p+} \cap T^7_{p´-}$ can be written as $a_{p+} - a_{p-} = d\ell$ where $\ell$ is a smooth function on $T^7_{p+} \cap T^7_{p´-}$. Fix a compactly supported, non-negative function, q, on $T^7_{p+} \cap T^7_{p´-}$ that is equal to 1 on $T^6_{p+} \cap T^6_{p´-}$. Define a 1-form $A_{p+}$ on $T^6_{p+}$ by the rule $A_{p+} = a_{p+} + \frac{1}{2} d(q \cdot \ell)$. Note that $dA_{p+} = w_\Sigma$. Meanwhile, define a 1-form $A_{p´-}$ on $T^6_{p´-}$ by the rule $A_{p´-} = a_{p´-} - \frac{1}{2} d(q \cdot \ell)$. The forms $A_{p+}$ and $A_{p´-}$ agree on $T^6_{p+} \cap T^6_{p´-}$.

Make these modifications to all intersections between $T_+$ and $T_-$. Given $r \in (0, 7)$, use $T^r_+ \subset T_+$ and $T^r_- \subset T_-$ denote the subdomains given by the union of relevant versions of $T^r_{p\pm}$. The modified 1-forms define a smooth form, A, on $T^6_+ \cup T^6_-$ that obeys $dA = w_\Sigma$.

The 1-form A is a smooth 1-form defined on $T^6_+ \cup T^6_-$ with $dA = w_\Sigma$. Meanwhile, $\hat{a}_\Sigma$ on $T^6_+ \cup T^6_-$ also obeys $d\hat{a}_\Sigma = w_\Sigma$. It follows from (1.7) that A has integral zero over each component of $C_+$ and over each component of $C_-$. By construction, this is also true for $\hat{a}_\Sigma$. As a consequence, $A - \hat{a}_\Sigma$ is exact and so $A - \hat{a}_\Sigma = d\mathfrak{y}$ where $\mathfrak{y}$ is a smooth function on $T^6_+ \cup T^6_-$.

Fix a smooth function, $q_\Sigma$, with compact support on $T^6_+ \cup T^6_-$ that is equal to 1 on $T^5_+ \cup T^5_-$. A smooth 1-form, $\hat{a}$, is defined on $\Sigma$ by the rule

$$\hat{a} = \hat{a}_\Sigma + d(q_\Sigma \mathfrak{y}) .$$

(1.17)

The 1-form obeys $d\hat{a} = H\omega_\Sigma$ on $\Sigma$ where H is smooth such that equals 1 on $\Sigma - (\cup_{z \in ¥} D_z)$.

*Part 9*: This final part of the construction requires the introduction of yet more 'bump' functions, plus a constant $r \geq 0$ whose lower bound depends only on the chosen bump functions. There is no upper bound to the choice for $r$.

The first of the bump function, $q_*$, is a smooth, compactly supported function on $T^5_+ \cap T^5_-$ that is equal to 1 on $T^4_+ \cap T^4_-$. The second, Q, is a smooth function on $[1, \frac{3}{2}]$ that equals 1 where $t > \frac{5}{4}$ and equals 0 where $t < \frac{9}{8}$. The third, K, is a smooth, non-negative function on $[1, \frac{3}{2}]$, this equal to 1 where $t \geq \frac{9}{8}$ and equal to 0 where $t < \frac{17}{16}$. The key feature of κ are that it equals 1 where $\frac{d}{dt}Q \neq 0$ and that it vanishes where $t < 1 + 8\delta_*^2$.

With such choices in hand, the extension of $a$ to the $f \in [1, \frac{3}{2}]$ part of $M_\delta$ is defined by the 1-form

$$a = (1 + r\kappa)dt + \hat{a} + \tfrac{1}{2}d(Qq_*\mathfrak{k}).$$

(1.18)

The extension of $a$ to the $f \in [\frac{3}{2}, 2]$ part of $M_\delta$ is defined by the version of (1.18) that is obtained by replacing the function $t \to Q(t)$ by the function $t \to -Q(3-t)$ and the function κ by the function $t \to \kappa(3-t)$.

The extension just defined is such that $a \wedge w > 0$ if $r > c_0$. In addition

- $da = w$ *on the complement of* $[1, 2] \times (\cup_{z \in ¥} D_z)$;
- $da = Hw$ *on* $[1, 2] \times (\cup_{z \in ¥} D_z)$.

(1.19)

Thus, the pair $(a, w)$ defines a stable Hamiltonian structure on the union of the $f \in [1, 2]$ part of $M_\delta$ with $\cup_{p \in \Lambda} \mathcal{H}_p$.

### d) Geometry of Y: The pair $(a, w)$ on the $f \in [0, 1] \cup [2, 3]$ part of $M_\delta$ and $\mathcal{H}_0$

This section finishes the construction of the pair $(a, w)$ by extending what was constructed in Sections 1b and 1c to the remainder of Y. Various constraints on $\mathfrak{v}$ must be imposed here. All of these can be made where $f \in [0,1] \cup [2, 3]$ and thus they can be imposed without compromising those listed previously. The construction has four parts.

*Part 1*: This part of the subsection extends $a$ and $w$ to a neighborhood in Y of the union the $f \in [1, 2]$ part of $M_\delta$ with $\cup_{p \in \Lambda} \mathcal{H}_p$. To do this, let $\Sigma_1 \subset M_\delta$ denote $f^{-1}(1) \cap M_\delta$. This is an embedded surface with boundary. It is diffeomorphic to the complement in $S^2$ of G disjoint disks. Let p denote an index 1 critical point of $f$. As can be seen from Constraint 1 in Section 1c, the vector field $\mathfrak{v}$ is along the part of the boundary of $\Sigma_1$ near p is tangent to the boundary sphere of the radius δ coordinate ball centered at p. This understood, Lie transport of $\Sigma_1$ by $\mathfrak{v}$ defines an embedding of $[1 - 2\delta_*^2, 1 + \delta_*^2] \times \Sigma_1$ into $M_\delta$ with the three properties: First, the image of the $t = 0$ slice is the identity map on $\Sigma_1$. Second, the image of the embedding intersects the boundary of $M_\delta$ only on the boundary of $\Sigma_1$. Third, the Euclidean coordiante $t$ on the $[1 - 2\delta_*^2, 1 + \delta_*^2]$ factor pulls back as $f$, and the corresponding coordinate vector field $\frac{\partial}{\partial t}$ pushes forward as $\mathfrak{v}$.

Lie transport by $\mathfrak{v}$ of the image in $M_\delta$ of the $t = 1 + \delta_*^2$ boundary identifies $\Sigma_1$ with the complement in $\Sigma$ of the union of the annular neighborhoods of the index 1 critical points of $f$ where the coordinate $h_+$ obeys $|h_+| \leq \frac{2}{3\sqrt{3}} \delta_*^2$. Use this identification to view $w_\Sigma$ as an area form on $\Sigma_1$ and thus a closed 2-form $w_\Sigma$ on $[1-2\delta_*^2, 1+\delta_*^2] \times \Sigma_1$. The embedding into $M_\delta$ identifies the latter on the $f \geq 1$ part of the image with $w$. This understood, the embedding extends $w$ to the image of the whole embedding. This extension agrees on the intersection with any given $\mathfrak{p} \in \Lambda$ version of $\mathcal{H}_\mathfrak{p}$ with the latter's version of $w$ and so extends $w$ to the union $\cup_{\mathfrak{p} \in \Lambda} \mathcal{H}_\mathfrak{p}$ with the $f \in [1-2\delta_*^2, 2]$ part of $\mathcal{M}_\delta$.

Meanwhile, the 1-form $a$ depicted in (1.18) is invariant with respect to Lie transport by $\mathfrak{v}$ on the $[1, 1+\delta_*^2] \times \Sigma_1$ part of $M_\delta$, and so Lie transport by $\mathfrak{v}$ defines a smooth extension of (1.18) to the $[1-2\delta_*^2, 1+\delta_*^2] \times \Sigma_1$ part of $M_\delta$. It follows from (1.7) that this extension agrees with the 1-form $a$ on the $r_+ > \delta$ part of the radius $\delta_*$ coordinate ball centered at each index 1 critical point of $f$. This extension of $a$ to the union of $\cup_{\mathfrak{p} \in \Lambda} \mathcal{H}_\mathfrak{p}$ with the $f \in [1-2\delta_*^2, 2]$ part of $\mathcal{M}_\delta$ is such that $a$ and $w$ define a stable Hamiltonian structure.

The analogous construction defined where $f \in [2-\delta_*^2, 2+2\delta_*^2]$ extends $a$ and $w$ to the union of $\cup_{\mathfrak{p} \in \Lambda} \mathcal{H}_\mathfrak{p}$ with the $f \in [1-2\delta_*^2, 2+2\delta_*^2]$ part of $M_\delta$ as a stable Hamiltonian structure with $\mathfrak{v}$ generating the kernel of $w$ in $M_\delta$.

*Part 2*: The extension of $a$ and $w$ to the remainder of $M_\delta$ and thus to the whole of Y requires more constraints on $\mathfrak{v}$. These constraints can be met by modifying any given pseudogradient for $f$ in the radius $10\delta_*$ coordinate balls about the index 0 and index 3 critical points of $f$, and on the $f \in [100\delta_*^2, 1-3\delta_*^2]$ and $f \in [2-3\delta_*^2, 3-100\delta_*^2]$ parts of $M_\delta$. What follows states the constraint on $\mathfrak{v}$ in the coordinate balls

CONSTRAINT 1: *The pseudogradient $\mathfrak{v}$ appears in the radius $8\delta_*$ coordinate ball centered on the index 0 critical point of $f$ as the vector field*

$$\mathfrak{v} = \tfrac{1}{2} e^{-2t_+} \tfrac{\partial}{\partial t_+} .$$

*The constraint on $\mathfrak{v}$ in the radius $8\delta_*$ coordinate ball centered on the index 3 critical point of $f$ is obtained from what is written above by first replacing $t_+$ by $t_-$ on the right hand side and then reversing the sign of the resulting expression.*

To set the background for the next constraint, keep in mind that the $f = 1 - \delta_*^2$ slice of $M_\delta$ is an embedded 2-sphere. The integral of $w$ over this slice is 2. This is also the area given by the 2-form $\frac{1}{2\pi}\sin\theta_+ d\theta_+ d\varphi_+$ to the boundary of any coordinate ball of radius less the $8\delta_*$ centered on the index 0 critical point of $f$.

The gradient flow of $\mathfrak{v}$ identifies the set of points $¥ \in \Sigma$ with a set of points in this slice, and it identifies the set of disks $\{D_z\}_{z\in ¥}$ with a corresponding set of disks in this slice. This understood, view $¥$ and the set $\{D_z\}_{z\in ¥}$ as subsets of the $f = 1 - \delta_*^2$ sphere in $M_\delta$. Keep in mind that the integral of $w$ over each such disk is $\frac{\pi}{100(1+b_1)}\delta_*^2$. Choose a set of $1+b_1$ points in the boundary 2-sphere of the radius $7\delta_*$ coordinate ball about the index 0 critical point. These points should be equally spaced along the equator. Let $¥_0$ denote this set. Fix a round disk about each point in $¥_0$ with area $\frac{\pi}{100(1+b_1)}\delta_*^2$. Given $z \in ¥_0$, use $D_{0z}$ to denote the corresponding disk.

Note next that the $f = 1 - \delta_*^2$ slice of $M_\delta$ is identified via the integral curves of any given pseudogradient vector field for $f$ with the boundary of any radius $\delta$ or greater coordinate ball centered on the index 0 critical point of f.

With the last three paragraphs as background, use arguments from Moser's proof that area forms are determined up to diffeomorphism by their total area to impose the following constraint:

CONSTRAINT 2: *The pseudogradient $\mathfrak{v}$ can be chosen where $f \leq 1 - 3\delta_*^2$ subject to Constraint 1 so that the resulting identification of the $f = 1 - \delta_*^2$ slice in $M_\delta$ with the boundary of coordinate ball of any radius $r \in [\delta, 8\delta_*]$ coordinate ball centered on the index zero critical point of f is such that the following is true:*

- *The identification makes w equal to $\sin\theta_+ d\theta_+ d\varphi_+$.*
- *The identification sends $¥$ to $¥_0$ and it identifies the sets $\{D_z\}_{z\in ¥}$ and $\{D_{0z}\}_{z\in ¥_0}$.*

Assume henceforth that $\mathfrak{v}$ obeys this constraint also.

The final constraint concerns identification between the partnered disks in the sets $\{D_z\}_{z\in ¥}$ and $\{D_{0z}\}_{z\in ¥_0}$. The statement of the constraint refers to radial coordinates $(\rho_0, \varphi_0)$ for each $z \in ¥_0$ version of $D_{0z}$. These are chosen so that the form $\frac{1}{2\pi}\sin\theta_+ d\theta_+ d\varphi_+$ appears as the form $\rho_0 d\rho_0 d\varphi_0$. Meanwhile, the corresponding $z \in ¥$ version of $D_z$ has its chosen coordinates $(\rho, \varphi)$.

CONSTRAINT 3: *The identification given by pseudogradient $\mathfrak{v}$ of each $z \in ¥$ version of $D_z$ with its $z \in ¥_0$ partner $D_{0z}$ identifies the coordinates $(\rho, \varphi)$ with $(\rho_0, \varphi_0)$.*

There is no obstruction to choosing $\mathfrak{v}$ so as to obey this constraint.

*Part 3*: Use the integral curves of the vector field $\mathfrak{v}$ to define a diffeomorphism between the $f \in [\delta^2, 1 - \delta_*^2]$ part of $M_\delta$ and $[\delta^2, 1 - \delta_*^2] \times S^2$. The latter identifes $\mathfrak{v}$ with the Euclidean vector field on the first factor and it identifies the coordinate $t$ on the first factor with $f$. Use this identification to extend $w$ to the $f \in [\delta^2, 1-\delta_*^2]$ part of $M_\delta$.

The extension of $a$ requires the choice of a certain smooth 1-form on $S^2$. This form is denoted by $\hat{a}_0$. It is chosen so that

- $d\hat{a}_0 = \frac{1}{2\pi} \sin\theta_+ d\theta_+ d\varphi_+$ on $S^2 - (\cup_{z \in \yen_0} D_{0z})$.
- $\hat{a}_0 = \chi_\#(R_z + \frac{1}{2}\rho_0^2 d\varphi_0)$ on any given $z \in \yen_0$ version of $D_{0z}$.

(1.20)

The 1-form $a$ on the $f \in [1-2\delta_*^2, 1-\delta_*^2]$ part of $M_\delta$ is invariant with respect to Lie transport by $\mathfrak{v}$. Use Lie transport by $\mathfrak{v}$ to extend $a$ to the $f \in [\delta^2, 1-\delta_*^2]$ part of $M_\delta$. Let $a'$ denote this extension. The identification given by the integral curves of $\mathfrak{v}$ between the $f \in [\delta^2, 1-\delta_*^2]$ part of $M_\delta$ and $[\delta^2, 1-\delta_*^2] \times S^2$ writes $a'$ as $dt + \hat{a}_1$ where $\hat{a}_1$ is a smooth, $t$-independent 1-form on the $S^2$ factor. It follows from (1.20) that $d\hat{a}_1 = d\hat{a}_0$. There is, as a consequence, a smooth function $\mathfrak{f}_0$ on $S^2$ such that $\hat{a}_1 = \hat{a}_0 + d\mathfrak{f}_0$.

Fix a smooth function on $[0, 1]$ that equals 1 on $[\frac{3}{4}, 1]$ and equals 0 on $[0, \frac{1}{4}]$. Denote this function by $Q_0$. Choose a second smooth function on this domain, this one denoted by $\kappa_0$. Require that $\kappa_0 \geq 0$, that $\kappa_0 = 1$ on $[\frac{3}{16}, \frac{13}{16}]$, and that $\kappa_0 = 0$ on both $[0, \frac{1}{8}]$ and $[\frac{7}{8}, 1]$. Let $r$ denote the constant used in (1.21). The desired extension of $a$ to the $f \in [\delta^2, 1-2\delta_*^2]$ part of $M_\delta$ is obtained by first identifying this part of $M_\delta$ as done above as with $[\delta^2, 1-2\delta_*^2] \times S^2$; then set the extension equal to

$$a = (1 + r\kappa_0) dt + \hat{a}_0 + d(Q_0 \mathfrak{f}_0).$$

(1.21)

An analogous construction for the $f \in [2+2\delta_*^2, 3-\delta^2]$ part of $M_\delta$ extends $a$ now to the whole of $M_\delta \cup (\cup_{\mathfrak{p} \in \Lambda} \mathcal{H}_\mathfrak{p})$. The latter identifies $w$ where $t_- \in [\delta, 7\delta_*]$ with the 2-form $-\frac{1}{2\pi} \sin\theta_- d\theta_- d\varphi_-$, and it writes $a$ as in (1.21) but with $\hat{a}_0$ replaced by a 1-form $\hat{a}_3$, with $\mathfrak{f}_0$ replaced by a function $\mathfrak{f}_3$ and with suitable replacements for $Q_0$ and $\kappa_0$. In particular, the 1-form $\hat{a}_3$ is the pull-back of $\hat{a}_0$ by the involution of $S^2$ that sends $\theta$ to $\pi - \theta$.

If $r > c_0$, then the resulting extensions defines a 1-form $a$ and a closed 2-form $w$ on the whole of $M_\delta \cup (\cup_\mathfrak{p} \mathcal{H}_\mathfrak{p})$ with the property that $a \wedge w > 0$ and $da \in \text{span}(w)$.

*Part 4*: This last part of the subsection extends $a$ and $w$ over the handle $\mathcal{H}_0$. To this end, reintroduce the coordinates $(u, (\theta, \phi))$ for $\mathcal{H}_0$ as defined in Section 1a. Define this extension by setting

- $a = 2(\chi_+ e^{2(|u|-R)} + \chi_- e^{-2(|u|+R)}) du + \hat{a}_0$ ,

- $w = \sin\theta \, d\theta \, d\phi$.

(1.22)

The definition is such

$$dw = 0, \quad a \wedge w > 0 \quad and \quad da \in \mathrm{span}(w).$$

(1.23)

on $\mathcal{H}_0$. In particular, the pair $(a, w)$ defines a stable Hamiltonian structure on $\mathcal{H}_0$ that extends $a$ and $w$ from $M_\delta \cup (\cup_{\mathfrak{p} \in \Lambda} \mathcal{H}_\mathfrak{p})$ to the whole of Y.

### e) Properties of *a* and *w*

The pair $a$ and $w$ as just constructed define a stable Hamiltonian structure on Y. Five properties of this pair play key roles in what follows in this paper and in its sequel [KLTIII]. These properties are listed below. The notation uses $v$ to the vector field on Y that generates the kernel of $w$ and has pairing equal to 1 with $a$.

Here are the first two properties:

PROPERTY 1: *On $\mathcal{H}_0$*: $w = \sin\theta \, d\theta \wedge d\phi$ *and* $v = \frac{1}{2(\chi_+ e^{2(u-R)} + \chi_- e^{-2(u+R)})} \frac{\partial}{\partial u}$.

PROPERTY 2: *On $M_\delta$: The 2-form $w$ on $M_\delta$ restricts to each constant f slice as an area form. The vector field $v$ on $M_\delta$ is the pseudogradient vector field $\mathfrak{v}$.*

The third property introduces by way of notation functions f and g of the coordinate u which are given by

$$f = x + 2(\chi_+ e^{2(u-R)} + \chi_- e^{-2(u+R)}) \quad and \quad g = (\chi_+ e^{2(u-R)} - \chi_- e^{-2(u+R)}).$$

(1.24)

Their derivatives are denoted by f´ and g´. This upcoming third property uses $\alpha$ to denote a certain positive function of the coordinates u and $\theta$.

PROPERTY 3: *Fix $\mathfrak{p} \in \Lambda$. Then a, w and v on $\mathcal{H}_\mathfrak{p}$ are*
- $a = (x + g´)(1 - 3\cos^2\theta) \, du - \sqrt{6} f \cos\theta \sin^2\theta \, d\phi + 6 g \cos\theta \sin\theta \, d\theta$,
- $w = 6x \cos\theta \sin\theta \, d\theta \wedge du - \sqrt{6} \, d\{f \cos\theta \sin^2\theta \, d\phi\}$,
- $v = \alpha^{-1}\{f(1 - 3\cos^2\theta) \partial_u - \sqrt{6} x \cos\theta \, \partial_\phi + f´ \cos\theta \sin\theta \, \partial_\theta\}$.

The fourth property concerns a closed 1-form on Y that extends df from $M_\delta$. This 1 form is denoted in what follows by $\upsilon_\diamond$. The form $\upsilon_\diamond$ restricts to any given $\mathfrak{p} \in \Lambda$ version of $\mathcal{H}_\mathfrak{p}$ as $df_*$ with $f_*$ given by (1.9). Meanwhile, $\upsilon_\diamond$ on $\mathcal{H}_0$ is given by

$$\upsilon_\Diamond = 2(\chi_+ e^{2(|u|-R)} + \chi_- e^{-2(|u|+R)}) du \; .$$

(1.25)

Note that $\upsilon_\Diamond$ is exact on $M_\delta \cup (\cup_{\mathfrak{p} \in \Lambda} \mathcal{H}_\mathfrak{p})$ as the function $f_*$ extends $f$ over any given $\mathfrak{p} \in \Lambda$ version of $\mathcal{H}_\mathfrak{p}$. Even so, $\upsilon_\Diamond$ is not an exact form on the whole of Y.

PROPERTY 4: *The 1-form $\upsilon_\Diamond$ is such that $\upsilon_\Diamond \wedge w \geq 0$ with equality only on the locus in each $\mathfrak{p} \in \Lambda$ version of $\mathcal{H}_\mathfrak{p}$ where both* $u = 0$ *and* $1 - 3\cos^2\theta = 0$.

The fact that $\upsilon_\Diamond \wedge w > 0$ on $Y - (\cup_{\mathfrak{p} \in \Lambda} \mathcal{H}_\mathfrak{p})$ follows directly from the first and second properties. Meanwhile the formula for $w$ in the third property implies the following on any given $\mathfrak{p} \in \Lambda$ version of $\mathcal{H}_\mathfrak{p}$:

$$\upsilon_\Diamond \wedge w = \sqrt{6} (g'f (1 - 3\cos^2\theta)^2 + 6gf' \cos^2\theta \sin^2\theta) du \sin\theta\, d\theta\, d\phi \; .$$

(1.26)

The final property concerns the homology of Y. Note in this regard that the Mayer-Vietoris sequence defines a canonical isomorphism

$$H_2(Y; \mathbb{Z}) = H_2(M; \mathbb{Z}) \oplus H_2(\mathcal{H}_0; \mathbb{Z}) \oplus (\oplus_{\mathfrak{p} \in \Lambda} H_2(\mathcal{H}_\mathfrak{p}; \mathbb{Z}))$$

(1.27)

The summands that correspond to the various 1-handles are isomorphic to $\mathbb{Z}$; and any oriented, cross-sectional sphere is a generator. The convention in what follows is to orient these spheres with the 2-form $\sin\theta\, d\theta\, d\phi$.

PROPERTY 5: *Integration of the 2-form $w$ defines the linear map from $H_2(Y; \mathbb{Z})$ to $\mathbb{Z}$ that*
- *has value 2 on the generator of $H_2(\mathcal{H}_0; \mathbb{Z})$;*
- *has value zero on each $\mathfrak{p} \in \Lambda$ version of $H_2(\mathcal{H}_\mathfrak{p}; \mathbb{Z})$;*
- *acts on the $H_2(M; \mathbb{Z})$ summand in (1.27) as the pairing with the chosen class* $c_{1M}$.

Some additional, less central properties of $a$, $w$ and $\nu$ are pointed out as needed in subsequent sections of this paper.

**f) Heegaard-Floer homology on M and embedded contact homology on Y**
 This subsection briefly describe the purely geometric aspects of Heegaard-Floer homology and embedded contact homology that have bearing to what follows in this paper and its sequel [KLTIII]. This subsection has five parts.

*Part 1*: The definition of Heegaard-Floer homology requires the Morse function $f$, the choice of a $\text{Spin}^{\mathbb{C}}$ structure on M, and then the choice of a Heegaard diagram that is strongly admissible for the $\text{Spin}^{\mathbb{C}}$ structure's first Chern class. These sorts of Heegaard diagrams are described by Definition 4.10 in [OS1]. The definition involves only the image of the first Chern class in $H^2(M;\mathbb{Z})$/torsion. This understood, take $c_{1M}$ in Lemma 1.1 to be this image in what follows.

There are three versions of Heegaard-Floer homology, and the chain complex for these are constructed from a suitably constrained, and decorated sets of G points in $\Sigma$. The G points must lie in $C_- \cap C_+$ with no two in the same version of $C_{p+}$ nor in the same version of $C_{p-}$. Use $\mathcal{Z}_{HF}$ to denote the collection of such sets of G points. The decoration, when present, adds an integer label to a given set of points from $\mathcal{Z}_{HF}$. The differential for the chain complex involves certain sorts of pseudoholomorphic disks in the G-fold symmetric product of $\Sigma$.

It proves more convenient for the purposes at hand to use Robert Lipshitz' reformulation of Heegaard-Floer homology [L] so as to view the set $\mathcal{Z}_{HF}$ as a collection $\hat{\mathfrak{v}} = \{\mathfrak{v}_1, \ldots, \mathfrak{v}_G\}$ of G integral curves of $\mathfrak{v}$ in $f^{-1}([1, 2])$ that pair the index 1 critical points of $f$ to the index 2 critical points. The correspondence between this view and original comes by writing $f^{-1}((1, 2))$ as $(1, 2) \times \Sigma$. This done by using the integral curves of the pseudogradient vector field as explained in Part 5 of Section 1c. By way of a reminder, this identification equates $\mathfrak{v}$ with the Euclidean vector field on the $(1, 2)$ factor, and so identifies the set $\hat{\mathfrak{v}}$ with a set of G points in $C_+ \cap C_-$. The pseudoholomorphic disks in the G-fold symmetric product of $\Sigma$ appear in Lipshitz' formulation as pseudoholomorphic maps from certain sorts of surfaces with boundary into $\mathbb{R} \times [1, 2] \times \Sigma$. The upcoming Section 6 has a detailed description of the sorts of maps that arise and so no more will be said about them here.

The Lipshitz view of Heegaard-Floer geometry as geometry on $[1, 2] \times \Sigma$ is used in what follows to interpret the Heegaard-Floer geometry as geometry on $f^{-1}((1, 2)) \subset M_\delta$ and thus geometry on Y and $\mathbb{R} \times Y$. This is the bridge that connects the embedded contact homology geometry to the Heegaard-Floer geometry.

More is said about this connection in Sections 2, 6 and 7 of this article.

*Part 2*: The embedded contact homology chain complex requires the choice of a $\text{Spin}^{\mathbb{C}}$ structure on Y. The restriction of a $\text{Spin}^{\mathbb{C}}$ structure on Y to $M_\delta$ determines a unique $\text{Spin}^{\mathbb{C}}$ structure on M. Moreover, the associated $\text{Spin}^{\mathbb{C}}$ structure on M together with the pairing of its first Chern class with the generators of the 1-handle summands in (1.27) canonically characterizes the $\text{Spin}^{\mathbb{C}}$ structure on Y. With this understood, the discussion henceforth concerns the $\text{Spin}^{\mathbb{C}}$ structure on Y with the following properties:

- *It determines the given $\mathrm{Spin}^{\mathbb{C}}$ structure on M.*
- *Its first Chern class has pairing two with the generator of $H_2(\mathcal{H}_0; \mathbb{Z})$.*
- *Its first Chern class has pairing zero with the generators of $\oplus_{p \in \Lambda} H_2(\mathcal{H}_p; \mathbb{Z})$.*

(1.28)

As explained in [Hu1], any given generator of embedded contact homology chain complex can be viewed as a finite set, $\Theta$, whose typical element is a pair $(\gamma, m)$ where $\gamma$ is a closed, integral curve of $v$ and where m is a non-negative integer. The chosen $\mathrm{Spin}^{\mathbb{C}}$ structure further constrains the pairs that comprise $\Theta$. It suffices for the present to describe the constraint that comes from the first Chern class of the associated $\mathrm{Spin}^{\mathbb{C}}$ structure. To do so, note that the kernel of $a$ defines a 2-plane subbundle in TY that is oriented by the 2-form $w$. Let $K^{-1}$ denote this oriented, real 2-plane bundle. Use $e_{K^{-1}}$ to denote the Euler class of $K^{-1}$, this a class in $H^2(Y; \mathbb{Z})$. Suppose that $(\gamma, m) \in \Theta$. The vector field $v$ orients $\gamma$, and with this orientation understood, the loop $\gamma$ defines a generator of $H_1(Y; \mathbb{Z})$. Use $[\gamma]^{Pd}$ to denote the Poincaré dual class in $H^2(Y; \mathbb{Z})$. Let $c_1$ denote the first Chern class of the $\mathrm{Spin}^{\mathbb{C}}$-structure. The pairs that comprise $\Theta$ are constrained so that

$$c_1 = e_{K^{-1}} + 2\sum_{(\gamma,m) \in \Theta} m [\gamma]^{Pd}.$$

(1.29)

Note for future reference that the class $e_{K^{-1}}$ has the following properties:

- *It has pairing 2 with the generator of $H_2(\mathcal{H}_0; \mathbb{Z})$.*
- *It has pairing -2 with the generators of $\cup_{p \in \Lambda} H_2(\mathcal{H}_p; \mathbb{Z})$.*

(1.30)

The first bullet follows from what is said in PROPERTY 1 of Section 1e, and the second follows from what is said in PROPERTY 3 in Section 1e.

*Part 3*: Hutchings [Hu1] puts an additional constraint on the integers that can appear in pairs from $\Theta$. This part of the subsection describes the latter constraint. It is important to keep in mind with regards to this constraint that the generating set for embedded contact homology can be defined only if certain closed integral curves of $v$ are suitably non-degenerate. More is said about this momentarily. Assume for now that such is the case. What follows directly introduces some background that is needed to define the integer constraint.

Let $\gamma$ denote a closed integral curve of $v$. Lie transport by $v$ on a circumnavigation of $\gamma$ starting at p defines an automorphism $U_\gamma: K^{-1}|_p \to K^{-1}|_p$ which preserves the 2-form $w$ since $w$ is annihilated by $v$'s Lie derivative. Thus, $U_\gamma$ acts symplectically on $K^{-1}|_p$ and so it

has determinant 1. The eigenvalues of $U_\gamma$ are independent of p because a change in p changes $U_\gamma$ by a conjugation.

The curve $\gamma$ is said to be *non-degenerate* when $U_\gamma$ has no eigenvalue equal to 1 or -1. The curve is said to be *totally non-degenerate* when no power of $U_\gamma$ has eigenvalue 1 or -1. The curve $\gamma$ is said to be *hyperbolic* when $U_\gamma$ has real eigenvalues, and it is said to be *elliptic* otherwise. The eigenvalues are on the unit circle when $\gamma$ is elliptic.

Use $\mathcal{Z}_{ech}$ to denote the collection of sets of pairs that obey the first Chern class constraint. The corresponding embedded contact homology chain complex can be defined when all closed integral curves that appear in sets from $\mathcal{Z}_{ech}$ are totally non-degenerate. The upcoming Section 2 explains why this is so for the $Spin^{\mathbb{C}}$ structure that is specified in (1.28). Granted that such is the case, the embedded contact homology chain complex uses only those elements from $\mathcal{Z}_{ech}$ that pair hyperbolic closed orbits with the integer 1. The latter subset is denoted in what follows by $\mathcal{Z}_{ech,M}$.

*Part 4*: The differential for the embedded contact homology chain complex and various other important endomorphisms require the choice of a suitable almost complex structure on $\mathbb{R} \times Y$. The almost complex structure is constrained to be invariant with respect to translations along the $\mathbb{R}$ factor of $\mathbb{R} \times Y$. Use *s* in what follows to denote the Euclidean coordinate for this $\mathbb{R}$ factor. The almost complex structure must map $\frac{\partial}{\partial s}$ to $v$ and it must tame *w* in the sense that it define with *w* a positive definite, bilinear form on $K^{-1} \subset \mathbb{R} \times Y$. Let J denote such an almost complex structure.

A subset $C \subset \mathbb{R} \times Y$ is said here to be J-holomorphic subvariety in $\mathbb{R} \times Y$ when the following conditions are met:

- *C is a closed, Hausdorff dimension 2 subset with no-totally disconnected components.*
- *The complement in C of a finite set of points is a submanifold with J-invariant tangent space.*
- *The integral of w over C is finite.*
- *There is an $s \in R$ independent bound for the integral of $ds \wedge a$ over C's intersection with $[s, s+1] \times Y$.*

(1.31)

With regards to the last two bullets, note that the constraint given in the second bullet allows J to orient the tangent space to the smooth part of C. With this orientation understood, *w* restricts as non-negative 2-forms on this tangent space. The almost complex structure J is not required to preserve the kernel of *a* and so there is no assumption that d$s \wedge a$ is positive everywhere on the tangent space of C. Even so, there exists $r_* \geq 1$ such that the 2-form

$$\omega = ds \wedge a + r_* w$$

(1.32)

tames J in the sense that bilinear form $\omega(\cdot, J\cdot)$ is positive definite on $T(\mathbb{R} \times Y)$. Indeed, this follows because $w(\cdot, J\cdot)$ is positive definite on the kernel of $a$ and the pairing of $a$ with $v$ is 1.

The conditions in (1.32) with the observation that $\omega$ tames J imply the following: There exists $s_C \geq 1$ such that the $|s| > s_C$ part of C is a disjoint union of embedded cylinders on which $s$ has no critical points. Moreover, the large $|s|$ slices of any such cylinder appear in Y as a braid in a fixed radius tubular neighborhood of some closed integral curve of $v$, and the resulting $s$-parametrized family of braids converges pointwise to this integral curve as $|s| \to \infty$. This is explained by Hofer [HWZ] under some assumptions, and by Siefring [S] in general. See also Section 2 of [HT] for a discussion of the cases that arise in the present context.

A cylinder of the sort just described is said to be an *end* of C. A given end $\mathcal{E} \subset C$ is said to be *positive* when $s$ is unbounded from above on $\mathcal{E}$ and negative when $s$ is unbounded from below. Associated to the given end $\mathcal{E}$ is a pair $(\gamma_\mathcal{E}, m_\mathcal{E})$ of the following sort: First, $\gamma_\mathcal{E}$ is the closed integral curve of $v$ that gives the limit of the large $|s|$ slices of $\mathcal{E}$. Meanwhile, $m_\mathcal{E}$ is the positive integer that describes the degree of the associated projection from any given large $|s|$ slice of $\mathcal{E}$ to $\gamma_\mathcal{E}$. Note that the 1-form $a$ orients the large $|s|$ slices and so this degree is in all cases positive.

A J-holomorphic subvariety is said to be *irreducible* when the complement of any finite set of points is connected.

*Part 5*: The differential and some other relevant endomomorphism of the embedded contact homology chain complex are defined using sets of following sort: A suitable set $\vartheta$ consists of a finite collection of pairs, with each of the form (C, n) where C is an irreducible embedded, J-holomorphic submanifold and where n is a positive integer. The integer n is constrained to equal 1 unless C is $\mathbb{R}$-invariant. The collection $\vartheta$ is further constrained in various ways. What follows are the three constraints that are germain to what follows here and in [KLTIII]. First, no two pair share the same J-holomorphic subvariety. The second constraint involves the positive ends of the subvarieties from $\vartheta$ and the third involves the negative ends of the subvarieties from $\vartheta$.

To say more about the second constraint, let $\Theta_{\vartheta+}$ denote the collection of pairs of the form $(\gamma, m)$ where $\gamma = \gamma_{\mathcal{E}}$ with $\mathcal{E}$ being a positive end of some subvariety from $\vartheta$. Meanwhile, m is sum of positive integers. A given term in the sum is labeled by a pair $((C, n), \mathcal{E})$ where $(C, n) \in \vartheta$ and $\mathcal{E}$ is a positive end of C with $\gamma_{\mathcal{E}} = \gamma$. The corresponding contribution to m is $n m_{\mathcal{E}}$. The set $\Theta_{\vartheta+}$ is defined so that no two pair share the same closed integral curve. The second constraint requires that $\Theta_{\vartheta+}$ define an element in $\mathcal{Z}_{ech,M}$.

There is a negative end analog of $\Theta_{\vartheta+}$, this denoted by $\Theta_{\vartheta-}$, and the third constraint requires that $\Theta_{\vartheta-}$ to likewise define an element in $\mathcal{Z}_{ech,M}$.

## 2. The closed integral curves of $v$

This section first describes various properties of the closed integral curves of $v$. These properties are used to characterization and then list the elements that comprise the generating set $\mathcal{Z}_{ech,M}$. Proposition 2.8 gives the list.

### a) The closed integral curves that intersect $\mathcal{H}_0$.

Write $\mathcal{H}_0$ as $[-R - \ln 7\delta_*, R + \ln 7\delta_*] \times S^2$ as in Section 1a. PROPERTY 1 in Section 1e asserts that the vector field $v$ on $\mathcal{H}_0$ is proportional to the vector field $\partial_u$. This the case, all integral curves of $v$ that intersect $\mathcal{H}_0$ do so as $[-R - \ln 7\delta_*, R + \ln 7\delta_*] \times$ point. All such curves intersect any given constant $u \in [-R - \ln 7\delta_*, R + \ln 7\delta_*]$ cross-sectional sphere with the same intersection sign, this being positive when these 2-spheres are oriented by $w$. In particular, any closed integral curve of $v$ that intersects $\mathcal{H}_0$ defines a non torsion class in $H_1(Y; \mathbb{Z})$ with the following property: Its Poincaré dual in $H^2(Y; \mathbb{Z})$ has positive pairing with the generator of the $H_2(\mathcal{H}_0; \mathbb{Z})$ summand in (1.27). This last fact with (1.28) and (1.29) imply the following:

*The integral curves of $v$ from an element in $\mathcal{Z}_{ech,M}$ do not intersect $\mathcal{H}_0$.*

(2.1)

Granted this, then no more will be said here about the integral curves of $v$ that intersect $\mathcal{H}_0$ except for the following observation: The construction in Section 1d provides $1 + b_1$ of curves of this sort. These curves are labeled by the points in the set ¥ that is defined in Parts 6 and 7 of Section 1c. If $z \in$ ¥, the corresponding curve intersects the surface $\Sigma$ at z. This curve is denoted by $\gamma^{(z)}$.

### b) Closed integral curves that are disjoint from $\mathcal{H}_0$

An integral curve of $v$ that intersect $M_\delta$ where $f \in (0, 1]$ or where $f \in [2, 3)$ must intersect $\mathcal{H}_0$. This follows from PROPERTY 2 in Section 1e: The vector field $v$ on this part of Y is a pseudogradient vector field for $f$. This last fact with (2.1) lead to the following:

*The integral curves of $v$ from an element in $\mathcal{Z}_{ech,M}$ sit entirely in $f^{-1}(1,2) \cup (\cup_\mathfrak{p} \mathcal{H}_\mathfrak{p})$.*

(2.2)

The six parts that follow describe the closed integral curves of $v$ that obey (2.2).

*Part 1*: This part of the subsection describes the parts of the integral curves of $v$ that sit entirely in any given $\mathfrak{p} \in \Lambda$ version of $\mathcal{H}_\mathfrak{p}$. The following lemma summarizes. It uses the coordinates $(u, (\theta, \phi))$ for $\mathcal{H}_\mathfrak{p}$ from Section 1a.

**Lemma 2.1**: *Fix $\mathfrak{p} \in \Lambda$. The two circles in $\mathcal{H}_\mathfrak{p}$ where $\{u = 0, \cos(\theta) = \pm \frac{1}{\sqrt{3}}\}$ are closed integral curves of $v$. Both closed integral curves are hyperbolic and such that the associated linear return map has postive eigenvalues. These are the only integral curves of $v$ that lie entirely in $\mathcal{H}_\mathfrak{p}$.*

*Proof of Lemma 2.1*: The fact that the circles in question are integral curves of $v$ follows from the formula for $v$ given in PROPERTY 3 of Section 1e. Note in this regard that the function f is positive at $u = 0$, but f´ is zero and so $v$ is proportional to the vector field $\partial_\phi$ where both $u = 0$ and $1 - 3\cos^2\theta = 0$. The fact that these curves are hyperbolic follows from the form of the linearization of $v$ along the $\{u = 0, \cos(\theta) = \pm \frac{1}{\sqrt{3}}\}$ locus. This is because the linearization is proportional to

$$\partial_\phi + \beta((x_0 + 4e^{-2R})((\theta - \theta_*)\partial_u + 16e^{-2R} u \partial_\theta)$$

(2.3)

where $\beta$ is a suitable constant and where $\theta_*$ is the relevant angle.

What follows explains why these are the only two integral curves of $v$ that lie entirely in $\mathcal{H}_\mathfrak{p}$. To start, note that the $\theta = \frac{\pi}{2}$ slice is foliated by the integral curves of $v$. Since the function f is positive, the function u increases along any integral curve of $v$ in this slice, so there can be no such integral curve that stays entirely in $\mathcal{H}_\mathfrak{p}$. The $\theta = 0, \pi$ arcs in $\mathcal{H}_\mathfrak{p}$ are also integral curves of $v$, but neither is closed in $\mathcal{H}_\mathfrak{p}$. Consider next a point in $\mathcal{H}_\mathfrak{p}$ where $\theta \in (0, \frac{\pi}{2})$. Let $\theta_*$ denote the angle with $\cos(\theta_*) = \frac{1}{\sqrt{3}}$ If both $u > 0$ and $\theta \geq \theta_*$, then u increases as does $\theta$ along the integral curve of $v$ from the point, and so the integral curve must exit $\mathcal{H}_\mathfrak{p}$ at the $u = R + \ln 7\delta_*$ end. If $u > 0$ and $\theta < \theta_*$, then u decreases initially and θ increases. If the curve does not limit to the $(u = 0, θ = θ_*)$ circle, then either θ becomes greater than $θ_*$ at some positive u value, in which case both u and θ increase from that point on and the curve exits the $u = R + \ln 7δ_*$ end of $\mathcal{H}_p$. On the other hand, if u becomes negative with $θ < θ_*$, then θ begins to decrease, and u continues to decrease, so the curve exits $\mathcal{H}_p$ at $u = -R - \ln 7δ_*$ end. A similar analysis applies when u is negative at the given point.

*Part 2*: The upcoming lemma states an important feature of the integral curves of $v$ in a given $\mathfrak{p} \in \Lambda$ version of $\mathcal{H}_p$ that cross from one boundary two sphere to the other.

**Lemma 2.2**: *There exists $κ ≥ 100$ depends only on the Heegaard Floer data and has the following significance: Assume that $δ < κ^{-1}$ and that $x_0 ≤ κ^{-1}δ^2$. Fix $\mathfrak{p} \in \Lambda$. Suppose that $γ$ is a connected component in the $|u| ≤ \ln R + δ_*$ part of $\mathcal{H}_p$ of an integral curve of $v$. Let $(θ_-, φ_-)$ and $(θ_+, φ_+)$ denote the respective coordinates of $γ$'s starting and ending points on the $|u| = R + \ln δ_*$ locus.*
- *Either $\cos(θ) = 0$ or $|\cos θ| > 0$ along the whole of $γ$.*
- *The following assertions are equivalent:*
  a) *The end points of $γ$ are on different spheres.*
  b) *$|1 - 3\cos^2 θ| > 0$ along the whole of $γ$.*
  c) *$θ_+ = θ_-$.*
  d) *$|\cos θ_±|\sin^2 θ_± < \frac{2}{3\sqrt{3}} δ_*^{-2}(x_0 + 4e^{-R})$*
- *The start and end points of $γ$ are on the respective sphere where $u = -R - \ln δ_*$ and the sphere where $u = R + \ln δ_*$ if and only if $1 - 3\cos^2 θ > 0$ on $γ$.*
- *Assume that $u = -R - \ln δ_*$ on the start point of $γ$ and that $u = R + \ln δ_*$ on the end point. Lift the coordinate $φ$ along $γ$ to an $\mathbb{R}$-valued function and let $Δφ$ denote the net change in this function from starting point to ending point along $γ$. Then*

$$Δφ = -\sqrt{6} \int_{[-R-\ln δ_*, R+\ln δ_*]} \frac{x(u)}{f(u)} \frac{\cos(θ(u))}{1 - 3\cos^2(θ(u))} du$$

*where the function $u \to θ(u)$ is such that $θ(u)$ is the unique solution to the equation*

$$\cos θ \sin^2 θ = 2δ_*^2 \cos θ_± \sin^2 θ_± \frac{1}{f(u)}.$$

*with $1 - 3\cos^2 θ > 0$.*
- *Assume that $u = -R - \ln δ_*$ on the start point of $γ$ and that $u = R + \ln δ_*$ on the end point. The integral of the 1-form $a$ over the $u \in [-R - \ln δ_*, R + \ln δ_*]$ portion of $γ$ is equal to*

$$2R\, x_0 (1 + \mathfrak{e}) \int_{[-R-\ln\delta_*,\, R+\ln\delta_*]} \frac{1 + 3\cos^4\theta(u)}{|1 - 3\cos^2\theta(u)|}\, du + \mathfrak{e}'$$

where $\mathfrak{e}$ and $\mathfrak{e}'$ are such that $|\mathfrak{e}| \leq \kappa^{-1}$ and $|\mathfrak{e}'| \leq \kappa$.

This lemma is proved momentarily. What follows directly is a lemma to elaborate on the assertion in the fourth bullet regarding the function $u \to \theta(u)$

**Lemma 2.3**: *Let $\mathcal{J} \subset [0, \pi]$ denote the open set in $[0, \pi]$ where $|\cos\theta|\sin^2\theta < \frac{2}{3\sqrt{3}}$. The set $\mathcal{J}$ has three components, these denoted respectively by $\mathcal{J}_+, \mathcal{J}_-$ and $\mathcal{J}_0$. The components $\mathcal{J}_+$ and $\mathcal{J}_-$ are the respective $\cos\theta > 0$ and $\cos\theta < 0$ parts of the set where $1 - 3\cos^2\theta < 0$. The component $\mathcal{J}_0$ is the region where $1 - 3\cos^2\theta > 0$. The maps from $\mathcal{J}_+$ to $[0, \frac{2}{3\sqrt{3}})$, from $\mathcal{J}_-$ to $(-\frac{2}{3\sqrt{3}}, 0]$, and from $\mathcal{J}_0$ to $(-\frac{2}{3\sqrt{3}}, \frac{2}{3\sqrt{3}})$ given by $\theta \to \cos\theta\sin^2\theta$ are homeomorphisms. Moreover, the map from $\mathcal{J}_0$ is a diffeomorphism, while those from $\mathcal{J}_+$ and $\mathcal{J}_-$ are diffeomorphisms on their interiors.*

*Proof of Lemma 2.3*: This follows from the inverse function theorem given that the differential of the map $\theta \to \cos\theta\sin^2\theta$ is equal to $-(1 - 3\cos^2\theta)\sin\theta$.

*Proof of lemma 2.2*: The first bullet's assertion follows from what is said in the proof of Lemma 2.1.

What is said in the proof of Lemma 2.1 also proves that Item a) of the second bullet implies Item b) and that Item c) implies Item a). The fact that Item b) implies Item a) follows from the formula for $v$ given in PROPERTY 3 of Section 1e as the latter implies that the coefficient in front of $\partial_u$ is nowhere zero. To obtain Item c) from Item b), use the fact that the coefficient in front of $\partial_u$ is nowhere zero when Item b) holds. Then u can be used as an affine parameter along $\gamma$ and the fact that $f'(-u) = -f'(u)$ implies Item c). To obtain Item d) from Item b), view $\theta$ as a function of u. Then the equation for its u-derivative can be integrated to see that

$$\frac{\cos(\theta(u))\sin^2(\theta(u))}{\cos(\theta_\pm)\sin^2(\theta_\pm)} = \frac{f(R+\ln\delta_*)}{f(u)}.$$

(2.4)

To obtain the constraint on $\theta_\pm$, note first that the function $\theta \to |\cos\theta|\sin^2\theta$ is maximized where $\cos^2\theta = \frac{1}{3}$. Meanwhile, $f(R + \ln\delta_*) = 2\delta_*^2$ and the minimum of the function $f(\cdot)$ is at $u = 0$, where $f(0) = x_0 + 4e^{-R}$. These observations with (2.4) give the constraint on $|\cos\theta_\pm| \sin^2\theta_\pm$ in Item d). Similar considerations prove that the other items follow from Item d).

The assertion that is made by the third bullet follows from the fact that the coefficient in front of $\partial_u$ in the formula for $v$ in PROPERTY 3 of Section 1e is positive if and only if $1 - 3\cos^2\theta > 0$.

To prove the fourth bullet, view an $\mathbb{R}$-values lift of the coordinate $\phi$ along $\gamma$ as a function of u also. This done, then the formula for $v$ in PROPERTY 3 of Section 1e implies that the function $u \to \phi(u)$ obeys

$$\frac{d\phi}{du} = -\sqrt{6} \, \frac{x(u)}{f(u)} \, \frac{\cos\theta(u)}{1 - 3\cos^2\theta(u)}$$

(2.5)

The formula given in the fourth bullet follows by integrating (2.5).

The last bullet follows using the formula for $a$ in PROPERTY 3 of Section 1e with (2.4) and (2.5).

*Part 3*: This part of the subsection constitutes a digression to relate certain integral curves of $v$ with those of the pseudogradient vector field $\mathfrak{v}$ for $f$. To this end, keep in mind that $v = \mathfrak{v}$ on $M_\delta$. The next lemma states a consequence of this fact.

**Lemma 2.4**: *There exists $\kappa \geq 1$ with the following significance: Assume that $\delta < \kappa^{-1}\delta_*^3$. Let p and p´ denote respective index 1 and index 2 critical points of $f$. Suppose that $\gamma$ is a closed, connected segment of an integral curve of $v$ that starts on the boundary of the radius $\delta$ coordinate ball centered at p and ends on the boundary of the radius $\delta$ coordinate ball centered at p´. Assume, in addition, that the interior of $\gamma$ is disjoint from $\mathcal{H}_0 \cup (\cup_{p \in \Lambda} \mathcal{H}_p)$. Then $\gamma$ sits in the radius $\kappa\delta$ tubular neighborhood of an integral curve of $\mathfrak{v}$ in M that starts at p and ends at p´.*

*Proof of Lemma 2.4*: The segment $\gamma$ is an integral curve of $\mathfrak{v}$. Parametrize $\gamma$ as an embedding $s \to \gamma(s)$ from [0, 1] such that $\gamma(0)$ sits on the boundary of the radius $\delta$ coordinate ball centered at p. The vector field $\mathfrak{v}$ in the radius $\delta_*$ ball about p is depicted in (1.14) and $v = \mathfrak{v}$ where the distance from p is at least $\delta$. It follows as a consequence that $(1 - 3\cos^2\theta) \geq 0$ on $\gamma(0)$. This implies that $f \geq 1$ on $\gamma$. Note that $f$ on $\gamma$ can not be greater than 2 for if so, then the either the interior of $\gamma$ must intersect some radius $\delta$ coordinate ball about an index 2 critical point, or else $\gamma$ would intersect $\mathcal{H}_0$. This understood, the whole of $\gamma$ must be a segment of an integral curve of $\mathfrak{v}$ lying in $f^{-1}[1, 2]$ and $\gamma(1)$ must lie on the $f \leq 2$ part of the boundary of the radius $\delta$ coordinate ball centered on p´. This is to say that there is an integral curve of $\mathfrak{v}$ that connects the radius $\delta$ coordinate ball centered on p with that centered on p´. It then follows (by arguing to the contrary to obtain a contradiction) that there exists $c_0 \geq 1$ with the following significance: If $\delta \leq c_0^{-1}$, then γ must lie in the radius $c_0 \delta$ tubular neighborhood of an integral curve of $\mathfrak{v}$ that runs from p to p´.

*Part 4*: Lemma 2.4 has the following consequence:

**Lemma 2.5**: *There exists $\kappa \geq 1$ with the following significance: Assume that $\delta < \kappa^{-1} \delta_*^3$. Suppose that γ is a closed, integral curve of v that does not intersect $\mathcal{H}_0$ and does not lie entirely in some $\mathfrak{p} \in \Lambda$ version of $\mathcal{H}_\mathfrak{p}$. Then*
- *The intersection of γ with $M_\delta$ must lie where $f \in [1, 2]$ and each component of this intersection lies in a radius $\kappa \delta$ tubular neighborhood of an integral curve of $\mathfrak{v}$ that runs from an index 1 critical point of f to an index 2 critical point.*
- *A properly embedded, connected segment in γ must cross at least one $\mathfrak{p} \in \Lambda$ version of $\mathcal{H}_\mathfrak{p}$ starting on the boundary 2-sphere at the index 2 critical point end of $\mathcal{H}_\mathfrak{p}$ and ending on the boundary 2-sphere at the index 1 critical point end.*
- *Let $\mathfrak{p} = (p, p´)$ denote a given point in $\Lambda$. If γ intersects the boundary of the radius δ coordinate ball about p´, then it does so where $f \leq 2$ and the continuation of γ from this intersection point crosses $\mathcal{H}_\mathfrak{p}$ to an intersection point with the boundary of the radius δ coordinate ball about p.*
- *Let $\mathfrak{p} = (p, p´)$ again denote a given point in $\Lambda$. If γ intersects the boundary of the radius δ coordinate ball about p, then it does so where $f \geq 1$ and the continuation of γ backwards from this intersection point crosses $\mathcal{H}_\mathfrak{p}$ in reverse to an intersection point with the boundary of the radius δ coordinate ball p´.*

*Proof of Lemma 2.5*: Let $\gamma_1 \subset \gamma$ denote the closure of a connected component of γ's intersection with $M_\delta$. Since $v = \mathfrak{v}$ on $\gamma_1$, this can not be the whole of γ. As a consequence, $\gamma_1$ must have a boundary. One boundary point must lie on the $f \geq 1$ part of the boundary of the radius δ coordinate ball centered on some index 1 critical point, and the other on the $f \leq 2$ part of the boundary of the radius δ coordinate ball centered on an index 2 critical point. This follows from the fact that $\mathfrak{v}$ near such critical points obeys CONSTRAINT 1 in Section 1c. This understood, then the claim made by the first bullet follows directly from Lemma 2.4.

Let p´ denote the endpoint of $\gamma_1$ on the boundary of the radius δ coordinate ball centered on the index 2 critical point. Let $\gamma_2$ denote the connected component of γ in the corresponding version of $\mathcal{H}_\mathfrak{p}$ that has p´ as its starting point. The other endpoint of this segment must lie on the boundary 2 sphere at the index 1 critical point end of $\mathcal{H}_\mathfrak{p}$. Indeed, if this were not the case, then it follows from CONSTRAINT 1 in Section 1c that the other endpoint lies on the $f > 2$ part of the boundary of the radius $\delta$ ball centered at p´, and so continuing along $\gamma$ from this end point would lead to an intersection with $\mathcal{H}_0$. These last observations imply what is claimed by the second and third bullets of the lemma. A very much analogous argument proves the claim made by the fourth bullet.

The following is a direct consequence of Lemma 2.5:

**Corollary 2.6**: *There exists $\kappa \geq 1$ such that when $\delta < \kappa^{-1}\delta_*^3$, then what follows is true. Suppose that $\gamma$ is a closed, integral curve of $v$ that does not intersect $\mathcal{H}_0$ and does not lie entirely in some $\mathfrak{p} \in \Lambda$ version of $\mathcal{H}_\mathfrak{p}$. Then there exists a positive integer N and a cyclically ordered set $\{\mathfrak{p}^1, \ldots, \mathfrak{p}^N\}$ of not necessarily distinct elements from $\Lambda$ with the following property: The curve $\gamma$ is an end-to-end concatenation of 2N closed, connected segments $\{\gamma^{\mathfrak{p}^1}, \gamma_1, \gamma^{\mathfrak{p}^2}, \gamma_2, \ldots, \gamma^{\mathfrak{p}^N}, \gamma_N\}$ where*

- *For each $k \in \{1, \ldots, N\}$, the segment $\gamma^{\mathfrak{p}^k}$ is a component of $\gamma \cap \mathcal{H}_{\mathfrak{p}^k}$. This segment starts on the $f \leq 2$ part of the boundary 2-sphere of the index 2 critical point end of the handle and ends on the $f \geq 1$ part of the boundary of the index 1 critical point end of the handle.*
- *For each $k \in \{1, \ldots, N\}$, the segment $\gamma_k$ lies in $M_\delta$. It starts at the endpoint of $\gamma^{\mathfrak{p}^k}$ and, for $k < N$, it ends at the starting point of $\gamma^{\mathfrak{p}^{k+1}}$. The segment $\gamma_N$ ends at the starting point of $\gamma^{\mathfrak{p}^1}$.*
- *For each $k \in \{1, \ldots, N\}$, the segment $\gamma_k$ lies in the radius $\kappa\delta$ tubular neighborhood of an integral curve of $\mathfrak{v}$ that runs from the index 1 critical point component of $\mathfrak{p}^k$ to the index 2 critical point component of $\mathfrak{p}^{k+1}$ when $k < N$, and to the index 2 critical point component of $\mathfrak{p}^1$ when $k = N$.*

This corollary motivates the following definition: An *index 1-2 cycle* is a non-empty, finite, cyclically ordered set $\{\upsilon_1, \ldots, \upsilon_N\}$ whose components are the closures of integral curves of $\mathfrak{v}$ that connect index 1 critical points of $f$ to index 2 critical points of $f$ and are such that the following is true: For $k < N$, the ending point of $\upsilon_k$ and the starting point of $\upsilon_{k+1}$ define a pair from $\Lambda$. This is also true for the ending point of $\upsilon_N$ and the starting point of $\upsilon_1$. The integer N is said to be the length of the index 1-2 cycle.

*Part 5*: The preceding corollary asserts that any closed integral curve of $v$ that avoids $\mathcal{H}_0$ and does not lie entirely in some $\mathfrak{p} \in \Lambda$ version of $\mathcal{H}_\mathfrak{p}$ defines an index 1-2 cycle. The upcoming Proposition 2.7 describes the set closed integral curves of $v$ that correspond to a given index 1-2 cycle.

The proposition refers to an *orientation sign*, either + or -, that is associated to a given integral curve of $\mathfrak{v}$ that runs from an index 1 critical point of $f$ to an index 2 critical point. What follows directly explains how this sign is determined. To start, remark that the pseudogradient $\mathfrak{v}$ is such that the ascending disks from the index 1 critical points of $f$ have transversal intersection with the descending disks from the index 2 critical points. These intersections are the integral curves of $\mathfrak{v}$ that run from the index 1 to the index 2 critical points. Let p denote a given index 1 critical point. The intersection of the ascending disk from p with $\Sigma = f^{-1}(\frac{3}{2})$ is the circle $C_{p+}$ from Part 2 of Section 1c. By way of a reminder, the integral curves of $\mathfrak{v}$ identify $C_{p+}$ with the $\theta_+ = \frac{\pi}{2}$ circle in any constant radius sphere in the radius $8\delta_*$ coordinate ball centered at p. Introduce from Part 3 of Section 1c the functions $(\varphi_+, \hat{h}_+)$ for the tubular neighorhood $T_{p+}$ of $C_{p+}$ and then orient $C_{p+}$ using the 1-form $d\varphi_+$.

Let p´ denote a given index 2 critical point of $f$ and let $C_{p´-}$ denote the intersection between the descending disk from p´ with $\Sigma$. Introduce the coordinates $(\varphi_-, \hat{h}_-)$ from Part 3 of Section 1c and orient $C_{p-}$ using the 1-form $d\varphi_-$. Meanwhile, the level set $\Sigma$ is oriented using the area form $w_\Sigma$. Granted these orientation, any given intersection point in $\Sigma$ between $C_{p+}$ and $C_{p-}$ has an associated sign; this sign is obtained by comparing the orientation of $\Sigma$ with that defined by the pair $(\frac{\partial}{\partial \varphi_+}, \frac{\partial}{\partial \varphi_-})$ of respective oriented tangent vectors to $C_{p+}$ and $C_{p-}$ at the intersection point. Each such intersection point corresponds to a unique integral curve of $\mathfrak{v}$ from p to p´. This orientation sign is the sign associated to the given integral curve.

**Proposition 2.7**: *There exists $\kappa \geq 1$, and given $\delta < \kappa^{-1}\delta_*^3$, there exists $\kappa_\delta > \kappa |\ln \delta|$ with the following significance: Take $R \geq \kappa|\ln\delta|$. Let $\{\mathfrak{v}_1, \ldots, \mathfrak{v}_N\}$ denote an index 1-2 cycle. The set of closed orbits of $v$ that correspond to this cycle enjoys a 1-1 correspondence with the set of N-tuples of integers. To elaborate, let $\{\mathfrak{k}_1, \ldots, \mathfrak{k}_N\}$ denote a given N-tuple of integers. The corresponding closed integral curve of $v$ has a decomposition into segments as described in Corollary 2.6 such that for each $k \in \{1, \ldots, N\}$, the change in the angle $\phi$ for the segment $\gamma^{p^k}$ can be written as $\sigma_k + 2\pi m_k$ with $\sigma_k \in [0, 2\pi)$. In all cases, the resulting closed integral curve is hyperbolic; and the sign of the eigenvalues is $(-1)^N$ times the product of the orientation signs of the integral curves $\{\mathfrak{v}_1, \ldots, \mathfrak{v}_N\}$. The integral, $\ell$, of the 1-form a over this integral curve obeys*

$$\kappa^{-1}\sum_{1\leq k\leq N}(x_0 R + |\mathfrak{k}_k|) - \kappa \leq \ell \leq \kappa\sum_{1\leq k\leq N}(x_0 R + |\mathfrak{k}_k|)) + \kappa.$$

*Two paired sets of length $N$ index 1-2 cycle and $N$ integers label the same closed integral curve of $v$ if one is obtained from the other by the action of $\mathbb{Z}/(N\mathbb{Z})$ that relabels both the*

set $\{\upsilon_1, \ldots, \upsilon_N\}$ *and the set* $\{\mathfrak{k}_1, \ldots, \mathfrak{k}_N\}$ *by a simultaneous cyclic permutation of the labels* $\{1, \ldots, N\}$.

To elaborate on this business of $\mathbb{Z}/(N\mathbb{Z})$ permutations, suppose for the moment that $\{(\upsilon_1, \ldots, \upsilon_N), (\mathfrak{k}_1, \ldots, \mathfrak{k}_N)\}$ is a data set as described in the proposition such that $\mathbb{Z}/(N\mathbb{Z})$ acts with non-trivial stabilizer. Let $n \in \{2, \ldots, N-1\}$ denote the smallest integer such that the cycle of length n acts trivially on this set. The number n is a divisor of N. A suitable cyclic permutation of $(\upsilon_1, \ldots, \upsilon_N)$ writes this set as N/n consecutive copies of $(\upsilon_1, \ldots, \upsilon_n)$ and the same cyclic permutation simultaneously writes the set $(\mathfrak{k}_1, \ldots, \mathfrak{k}_N)$ as N/n consecutive copies of $(\mathfrak{k}_1, \ldots, \mathfrak{k}_n)$. The former is an index 1-2 cycle also; and the data $\{(\upsilon_1, \ldots, \upsilon_n), (\mathfrak{k}_1, \ldots, \mathfrak{k}_n)\}$ defines the same integral curve of $\nu$ as does the original, larger set. (The larger set can be viewed as labeling an N/n times multiple cover of the smaller set.) This understood, no generality is lost by restricting attention to sets where $\mathbb{Z}/(N\mathbb{Z})$ acts freely. This restriction is implicit in what follows.

*Proof of Proposition 2.7*: The proof is given in four steps.

Step 1: Keep in mind for what follows that the pseudogradient $\mathfrak{v}$ that gives $\nu$ on $M_\delta$ was chosen in part to guarantee that the ascending disks from the index 1 critical points of $f$ intersect the descending disks from the index 2 critical points in a transversal fashion on $\Sigma$ as described in CONSTRAINT 3 of Section 1c.

Let p denote an index 1 critical point of $f$ and let $\upsilon$ denote an integral curve of $\mathfrak{v}$ that runs from p to some index 2 critical point, p´. Reintroduce the coordinates $(\varphi_+, \hat{h}_+)$ for the annular neighborhood $T_{p+} \subset \Sigma$ of $C_{p+}$. Keep in mind that $\varphi_+$ is $\mathbb{R}/2\pi\mathbb{Z}$ valued and that $|\hat{h}_+| < \frac{128}{3\sqrt{3}} \delta_*^2$. The functions $(\varphi_+, \hat{h}_+)$ are annihilated by $\mathfrak{v}$; and this has the following implication: An integral curve of $\mathfrak{v}$ that intersects $T_{p+}$ with coordinates $\varphi_+ = \varphi$ and $\hat{h}_+ = \hat{h}$ intersects the radius $8\delta_*$ coordinate ball with the coordinates $(r_+, \theta_+, \varphi_+)$ with $r_+$ and $\theta_+$ such that $2r_+^2 \cos\theta_+ \sin^2\theta_+ = \hat{h}_+$ and with $\varphi_+ = \varphi$.

Reintroduce the corresponding coordinates $(\varphi_-, \hat{h}_-)$ for the annular neighborhood $T_{p´-}$ of $C_{p´-}$. As noted in CONSTRAINT 3 of Section 1c, the respective ± pairs of coordinate differentials on $\upsilon$'s component of $T_{p+} \cap T_{p´-}$ are such that $(d\varphi_+, d\hat{h}_+) = \varepsilon_\upsilon (d\hat{h}_-, d\varphi_-)$ where $\varepsilon_\upsilon \in \{-1, 1\}$ is the orientation sign for $\upsilon$. The intersection point of $\Sigma$ with $\upsilon$ is a point in $C_{p+} \cap C_{p´-}$. This point occurs where $\hat{h}_+ = 0$ and $\hat{h}_- = 0$. Write the respective $\varphi_+$ and $\varphi_-$ coordinates of this point as $\phi_{\upsilon+}$ and $\phi_{\upsilon-}$. The aforementioned identity between the ± pairs of coordinate differentials implies that the functions $(\varphi_+, \hat{h}_+)$ and $(\varphi_-, \hat{h}_-)$ are related on $\upsilon$'s component of $T_{p+} \cap T_{p´-}$ as follows: Lift $\phi_{\upsilon+}$ to $[0, 2\pi)$ and likewise lift $\phi_{\upsilon-}$. These lifts determine lifts of $\varphi_+$ and $\varphi_-$ to $\mathbb{R}$-valued functions. So interpreted, these $\mathbb{R}$-valued versions with $\hat{h}_+$ and $\hat{h}_-$ obey

$$\varphi_- - \phi_{\upsilon-} = \varepsilon_\upsilon \hat{h}_+ \quad \text{and} \quad \hat{h}_- = \varepsilon_\upsilon (\varphi_+ - \phi_{\upsilon+})$$

(2.6)

on $\upsilon$'s component of $T_{p+} \cap T_{p'-}$.

Step 2: Let $\gamma$ denote a given integral curve of the vector field $v$ that leaves the radius $\delta_*$ coordinate ball centered p_ with given values for $\varphi_+$ and $\hat{h}_+$ on the boundary 2-sphere. Write the $\hat{h}_+$ value as $\hat{h}_{+0}$. Let $\varphi_{+0}$ denote a lift to $\mathbb{R}$ of the $\varphi_+$ value. Assume that this lift with $\hat{h}_+$ are such that

$$|\varphi_{+0} - \phi_{\upsilon+}| + |\hat{h}_{+0}| < \tfrac{4}{3\sqrt{3}} \delta_*^2.$$

(2.7)

This guarantees that $\gamma$ intersects $\upsilon$'s component of $T_{p+} \cap T_{p-}$ and thus the boundary of the radius $\delta_*$ coordinate ball centered on p´.

Let $\mathfrak{p}$ denote the pair from $\Lambda$ with index 2 critical point p´ and let $p_1$ denote the index 1 critical point that is paired with p´ in $\mathfrak{p}$. Given (2.7), Lemma 2.2 asserts the following: The integral curve $\gamma$ passes through the 1-handle $\mathcal{H}_\mathfrak{p}$ so as to exit the radius $\delta_*$ ball centered on $p_1$ if and only if this intersection point occurs where $(\varphi_-, \hat{h}_-)$ are such that

$$|\hat{h}_-| < \tfrac{2}{3\sqrt{3}} (x_0 + 4e^{-R}).$$

(2.8)

Assume that (2.8) also holds for the chosen integral curve. What with (2.6), this requires that

$$|\varphi_{+0} - \phi_{\upsilon+}| < \tfrac{2}{3\sqrt{3}} (x_0 + 4e^{-R}).$$

(2.9)

Let $(\varphi_{+1}, \hat{h}_{+1})$ denote the the corresponding pair that are defined by the values of $(\varphi_+, \hat{h}_+)$ where the integral curve $\gamma$ intersects the boundary of the radius $\delta_*$ ball centered on $p_1$. It follows from Lemma 2.2 and (2.6) that $(\varphi_{+1}, \hat{h}_{+1})$ are given in terms of $(\varphi_{+0}, \hat{h}_{+0})$ by the rule

- $\varphi_{+1} - \phi_{\upsilon-} = \varepsilon_\upsilon \hat{h}_{+0} + \Delta\phi_\mathfrak{p}$ ,
- $\hat{h}_{+1} = -\varepsilon_\upsilon (\varphi_{+0} - \phi_{\upsilon+})$ ,

(2.10)

where $\Delta\phi_\mathfrak{p}$ denotes the integral that is written in the fourth bullet of Lemma 2.2.

*Step 3*: This step proves Proposition 2.6 for the simplest case, this when the index 1-2 cycle has the single element $\mathfrak{p} = \mathfrak{p}^{(1)}$. Write $\mathfrak{p} = (p, p')$. Apply what is said in the previous steps with $p_1 = p$. The segment, $\gamma$, is a closed integral curve of $v$ if and only if (2.10) holds with

$$h_{+1} = h_{+0} \quad \text{and} \quad \varphi_{+1} = \varphi_{+0} + 2\pi \mathfrak{k} \text{ for some } \mathfrak{k} \in \mathbb{Z}_-.$$
(2.11)

To determine the solutions, use the bottom equation in (2.10) to write $\varphi_{+0} - \phi_{\upsilon+} = -\varepsilon_\upsilon h_{+1}$. Having done so, (2.11) follows from the top equation in (2.10) if and only if

$$2\pi \mathfrak{k} = \phi_{\upsilon-} - \phi_{\upsilon+} + 2\varepsilon_\upsilon h_{+1} + \Delta\phi_\mathfrak{p}.$$
(2.12)

To see about solving (2.12), first use Lemma 2.3 to find $c_0 > 1$ such that if $\delta_*^{-2} x_0 \leq c_0^{-1}$, then the set in $[0, \pi]$ where $|\cos\theta| \sin^2\theta < \frac{2}{3\sqrt{3}} \delta_*^{-2}(x_0 + 4e^{-2R})$ has three components with one an interval centered on $\frac{\pi}{2}$. The function $h = \delta_*^2 \cos\theta \sin^2\theta$ maps this component diffeomorphically to the open interval $(-\frac{2}{3\sqrt{3}}(x_0 + 4e^{-2R}), \frac{2}{3\sqrt{3}}(x_0 + 4e^{-2R}))$. Let $\mathcal{I}$ denote this interval. As is explained in the next paragraph, there is a unique $h_{+1} \in \mathcal{I}$ that obeys (2.12) for any given integer m.

Write $h_{+1}$ as $h$ so as to view the right hand side of (2.12) as a map, $h \to F(h)$, from $\mathcal{I}$ to $\mathbb{R}$. It follows from the fourth bullet of Lemma 2.2 that F is a proper map from $\mathcal{I}$ to $\mathbb{R}$ if $R\delta_*^2 x_0^{-1} \geq c_0$. Indeed, the chain rule with the definition of $\theta(u)$ in the fourth bullet of Lemma 2.2 can be used to see that the derivative of F on $\mathcal{I}$ obeys

$$\frac{dF}{dh} = -\sqrt{6} \int_{[-R-\ln\delta, R+\ln\delta]} \frac{x(u)}{f(u)^2} \frac{1 + 3\cos^2(\theta(u))}{(1 - 3\cos^2(\theta(u)))^3} du + 2\varepsilon_\upsilon.$$
(2.13)

The integral that appears on the right hand side of (2.13) is no less than $c_0^{-1} R x_0^{-1}$. This fact about (2.13) implies that F is a diffeomorphism from $\mathcal{I}$ to $\mathbb{R}$ if $R x_0^{-1} > c_0$.

Assuming the bound $R x_0^{-1} > c_0$, it then follows that for each $m \in \mathbb{Z}$, there is a unique $h_{1+} = h_m \in \mathcal{I}$ such that solves (2.12). Note that (2.9) follows now automatically from the top bullet in (2.10), and thus (2.7) as well.

To continue with the proof of Proposition 2.7 when $N = 1$, let $\gamma$ denote the closed integral curve of $v$ that is defined as above for a given $\mathfrak{k} \in \mathbb{Z}$. It follows from (2.6), (2.10), (2.12) and (2.13) that the linearized return map defines a matrix in $SL(2; \mathbb{R})$ that has the form

$$\begin{pmatrix} -\varepsilon_\upsilon \sigma R & \varepsilon_\upsilon \\ -\varepsilon_\upsilon & 0 \end{pmatrix}$$

(2.14)

with $\sigma > c_0^{-1}$. The eigenvalues of this matrix are real; and one has distance $\mathcal{O}(R)$ from the unit circle when R is large. The sign of the large eigenvalue is that of $-\varepsilon_\upsilon$, thus (-1) times the orientation sign given to $\upsilon$.

The final assertion of Proposition 2.7 gives a bound from above and below for the integral of a along $\gamma$. This bound follows from the last bullet of Lemma 2.2.

Step 4: This step proves Proposition 2.7 for the N > 1 cases. To set notation, the index k in what follows is from the set $\{1, \ldots, N\}$, and if k = N, then k +1 denotes the integer 1.

Write the pairs in $\Lambda$ that are associated to the index 1 starting points of the curves from the set $\{\upsilon_1, \ldots, \upsilon_N\}$ as $\{\mathfrak{p}^1, \ldots, \mathfrak{p}^N\}$ and write any given $k \in \{1, \ldots, N\}$ version of $\mathfrak{p}^k$ as $(\mathfrak{p}^{(k)}, \mathfrak{p}'^{(k)})$. The closed integral curve $\upsilon_k$ from the cycle runs from $\mathfrak{p}^{(k)}$ to $\mathfrak{p}'^{(k+1)}$. For each $k \in \{1, \ldots, N\}$, let $\phi_{k+}$ denote the $[0, 2\pi)$ lift of the $\varphi_+$ coordinate of the intersection point between $\upsilon_k$ and the boundary sphere of the radius $\delta_*$ coordinate ball centered on $\mathfrak{p}^{(k)}$, and use $\phi_{k-}$ to denote the $[0, 2\pi)$ lift of the $\varphi_-$ coordinate of the intersection point between $\upsilon_k$ and the boundary sphere in the radius $\delta_*$ coordinate ball centered on $\mathfrak{p}'^{(k)}$. Each $k \in \{1, \ldots, N\}$ version of $\upsilon_k$ intersects $\Sigma$ in a distinct component of $T_- \cap T_+$. Each such component has the two sets of coordinates $(\varphi_+, \hat{h}_+)$ and $(\varphi_-, \hat{h}_-)$, and these are related by (2.6). Use $\varepsilon_{\upsilon_k}$ to the sign $\varepsilon_\upsilon$ that appears in $\upsilon_k$'s version of (2.6).

Suppose that $\{\gamma_{1*}, \ldots, \gamma_{N*}\}$ are segments of integral curves of $v$ in $M_\delta$ such that each $k \in \{1, \ldots, N\}$ version of $\gamma_k$ intersects $\upsilon_k$'s component of $T_+ \cap T_-$. Use $(\varphi_{+k}, \hat{h}_{+k})$ to denote the $(\varphi_+, \hat{h}_+)$ coordinates of this intersection point with $\varphi_{+k}$ a lift to $\mathbb{R}$ of $\varphi_+$ such that the pair $(\varphi_{+k} - \phi_{k+}, \hat{h}_{+k})$ obey (2.7) Meanwhile, use $(\varphi_{-k}, \hat{h}_{-k})$ to denote the $(\varphi_-, \hat{h}_-)$ coordinates with $\varphi_{-k}$ denoting a lift to $\mathbb{R}$ of $\varphi_-$ such that $(\varphi_{-k} - \phi_{k-}, \hat{h}_{-k})$ also obey (2.7). Assume that $\hat{h}_{-k}$ obeys (2.8) and so $\varphi_{+k} - \phi_{k+}$ obeys (2.9). This condition is necessary and sufficient to conclude that $\gamma_k$ extends across $\mathcal{H}_{\mathfrak{p}^{k+1}}$ as a segment of an integral curve of $v$.

The segment $\gamma_k$ enters the handle $\mathcal{H}_{\mathfrak{p}^{k+1}}$ with the value $-\hat{h}_{-k}$ for $\hat{h} = f\cos\theta \sin^2\theta$ and with $\phi$ coordinate that lifts to $\mathbb{R}$ as $\varphi_{-k}$. The extended segment exits this handle with the lifted $\phi$ coordinate and $\hat{h}$ equal to $\varphi_{-k} + \Delta\phi_{\mathfrak{p}^{k+1}}$ and $-\hat{h}_{-k}$. Here, $\Delta\phi_{\mathfrak{p}^{(k+1)}}$ is the relevant version of the integral on the right hand side of the third bullet in Lemma 2.2. This understood, the extended $\gamma_k$ concatenates with $\gamma_{k+1}$ to form a segment of an integral curve of $v$ if and only if there exists $\mathfrak{k}_{k+1} \in \mathbb{Z}$ such that

$$\varphi_{-k} + \Delta\varphi_{p^{(k+1)}} = \varphi_{+(k+1)} + 2\pi \mathfrak{k}_{k+1} \quad and \quad -\hat{h}_{-k} = \hat{h}_{k+1} .$$

(2.15)

Meanwhile, $\varphi_{-k} - \phi_{k-} = \varepsilon_{\upsilon_k} \hat{h}_{+k}$ and $\hat{h}_{-k} = \varepsilon_{\upsilon_k} (\varphi_{+k} - \phi_{k+})$. Thus, the N extended segments concatenate to define a closed integral curve of $v$ if and only if the following two conditions are met for each $k \in \{1, \ldots, N\}$: There exists an integer $\mathfrak{k}_{k+1}$ such that

$$2\pi \mathfrak{k}_{k+1} = \phi_{k-} - \phi_{(k+1)+} + \varepsilon_k \hat{h}_{+k} + \Delta\varphi_{p^{k+1}} .$$

(2.16)

What follows explains why there is a unique solution to the system of equations in (2.16) for any given set $\{\mathfrak{k}_1, \ldots, \mathfrak{k}_N\} \in \times_N \mathbb{Z}$ if $\delta_*^{-2} x_0 < c_0^{-1}$.

Reintroduce from Step 3 the interval $\mathcal{I} = (-\frac{2}{3\sqrt{3}}(x_0 + 4e^{-2R}), \frac{2}{3\sqrt{3}}(x_0 + 4e^{-2R}))$. The right hand side of (2.16) defines a proper, smooth map from $\times_N \mathcal{I}$ to $\times_N \mathbb{R}$. It follows from (2.13) that this map defines a diffeomorphism if $R\delta_*^2 x_0^{-1} > c_0$. Assume this. Then there exists, for each N-tuple $\{\mathfrak{k}_1, \ldots, \mathfrak{k}_N\} \in \times_N \mathbb{Z}$, a unique N-tuple $(\hat{h}_{+1}, \ldots \hat{h}_{+N}) \in \times_N \mathcal{I}$ that solves (2.16). As in the N = 1 case, (2.9) and (2.7) follow automatically.

The preceding gives the desired 1-1 correspondence between $\times_N \mathbb{Z}$ and the set of closed trajectories of $v$ that define the given index 1-2 cycle. The fact that each such closed trajectory is hyperbolic can be seen as follows: View the right hand side of (2.15) as defining a map from $\mathbb{R}^2$ to itself. The differential of this map at $(\varphi_{-k}, \hat{h}_{-k})$ is a matrix in $SL(2; \mathbb{R})$. The linearized return map is the appropriately ordered product of these N matrices. It follows from (2.13) that this ordered product has trace that can be written as

$$\sigma_N (-1)^{N_+} R^N \prod_{1 \leq k \leq n} \varepsilon_{\upsilon_k}$$

(2.17)

where $\sigma_N > c_N^{-1}$ with $c_N \geq 1$ determined by N. Here, $N_+$ is the number of segment of the index 1-2 cycle with positive orientation sign. Granted (2.17), it follows that the closed integral curve in question is hyperbolic, and its eigenvalues have the sign asserted by Proposition 2.7. The assertion about the integral of $a$ over this closed integral curve follows directly using the fifth bullet of Lemma 2.2.

The final assertion about the $\mathbb{Z}/(N\mathbb{Z})$ action is a consequence of the fact that the map from $\times_N \mathcal{I}_N$ to $\times_N \mathbb{R}$ defined by the right hand side of (2.16) is a diffeomorphism when $R\delta_*^2 x_0^{-1}$ is large.

### c) The set $\mathcal{Z}_{ech,M}$

What is said in Sections 2a and 2b allow for complete description of the set $\mathcal{Z}_{ech,M}$. This is provided by the upcoming Proposition 2.8. The four parts that follow directly supply the necessary background.

*Part 1*: The correspondence is canonical with the choice of the following data: Choose a lift to $[0, 2\pi)$ for the $\varphi_+$ coordinate function on the radius $8\delta_*$ coordinate ball about each index 1 critical point of $f$. Assume that such a lift has been chosen. Note that these G lifts supply via the handle identifications lifts to $[0, 2\pi)$ of the $\varphi_-$ coordinate function on the radius $8\delta_*$ ball about any given index 2 critical point of $f$. The correspondence also requires a labeling of the elements in $\Lambda$ as $\{\mathfrak{p}^1, \ldots, \mathfrak{p}^G\}$.

*Part 2*: Reintroduce from Part 1 of Section 1f the set $\mathcal{Z}_{HF}$. As noted in Section 1f, this is the set that is used to define the generators for the Heegaard-Floer chain complex. In any event, view $\mathcal{Z}_{HF}$ as the set of unordered G-tuples whose elements are integral curves of $\mathfrak{v}$ that start at an index 1 critical point of $f$ and end at an index 2 critical point with the G-tuple constrained as follows: No two elements start at the same index 1 critical point and no two end at the same index 2 critical point.

*Part 3*: Fix $\hat{\upsilon} = \{\upsilon_1, \ldots, \upsilon_G\} \in \mathcal{Z}_{HF}$. The set $\{\upsilon_1, \ldots, \upsilon_G\}$ has a unique partition modulo ordering into a set of at most G irreducible index 1-2 cycles. Let $\hat{O}$ denote this unordered set of index 1-2 cycles. A set of this sort is called an *HF-cycle*. Any given HF-cycle defines a corresponding partition of $\Lambda$. Let $\hat{o}$ denote a given index 1-2 cycle from $\hat{O}$, and let $\Lambda_{\hat{o}} \subset \Lambda$ denote the corresponding subset. Let $\{\mathfrak{k}_\mathfrak{p}\}_{\mathfrak{p} \in \Lambda_{\hat{o}}}$ denote a given set of integers. Proposition 2.7 associates to the data $(\hat{o}, \{\mathfrak{k}_\mathfrak{p}\}_{\mathfrak{p} \in \Lambda_{\hat{o}}})$ a unique, closed integral curve of $v$. Use $\gamma_{\hat{o}}$ to denote the latter curve.

*Part 4*: Fix $\mathfrak{p} \in \Lambda$. Use $\hat{\gamma}_\mathfrak{p}^+$ denote the integral curve of $v$ in $\mathcal{H}_\mathfrak{p}$ that sits in the locus where $u = 0$ and $\cos\theta = \frac{1}{\sqrt{3}}$. Use $\hat{\gamma}_\mathfrak{p}^-$ to denote the integral curve of $v$ in $\mathcal{H}_\mathfrak{p}$ that sits where $u = 0$ and $\cos\theta = -\frac{1}{\sqrt{3}}$.

Introduce O to denote the 4-element set $\{0, 1, -1, (1, -1)\}$. Let $o_\mathfrak{p}$ denote a given element in O. What follows associates an element in $\mathcal{Z}$ to $o_\mathfrak{p}$. The latter is denoted by $\Theta_\mathfrak{p}$.

- If $o_\mathfrak{p} = 0$, *then* $\Theta_\mathfrak{p} = \emptyset$.
- If $o_\mathfrak{p} = 1$ *or* $-1$, *then* $\Theta_\mathfrak{p} = (\hat{\gamma}_\mathfrak{p}^+, 1)$ *or* $(\hat{\gamma}_\mathfrak{p}^-, 1)$ *as the case may be*.
- If $o_\mathfrak{p} = (1, -1)$, *then* $\Theta_\mathfrak{p} = \{((\hat{\gamma}_\mathfrak{p}^+, 1), (\hat{\gamma}_\mathfrak{p}^-, 1))$.

(2.18)

What follows is the promised description of $\mathcal{Z}_{ech,M}$.

**Proposition 2.8**: *There exists $\kappa > 1$ such that if $\delta < \kappa^{-1}$ and $x_0 < \kappa^{-1}\delta^2$, then the elements in $\mathcal{Z}_{ech,M}$ are in 1-1 correspondence with $\mathcal{Z}_{HF} \times (\times_{\mathfrak{p} \in \Lambda}(\mathbb{Z} \times O))$. This correspondence is such that a given element $(\hat{\upsilon}, \{\mathfrak{k}_\mathfrak{p}, o_\mathfrak{p}\}_{\mathfrak{p} \in \Lambda}) \in \mathcal{Z}_{HF} \times (\times_{\mathfrak{p} \in \Lambda}(\mathbb{Z} \times O))$ corresponds to the element $\Theta = \{(\gamma_{\hat{o}}, 1)\}_{\hat{o} \in \hat{O}} \cup \{\cup_{\mathfrak{p} \in \Lambda} \Theta_\mathfrak{p}\} \in \mathcal{Z}_{ech,M}$*

*Proof of Proposition 2.8*: The proof must establish that any collection $\Theta$ of the sort described has the $\text{Spin}^{\mathbb{C}}$ structure that is described in (1.28), and it must establish that its list of such collections is complete.

Assume for the moment that the elements define the desired $\text{Spin}^{\mathbb{C}}$ structure to focus on the second aspect of the proof. According to (2.2), each closed integral curve of $v$ from a given element in $\mathcal{Z}_{ech,M}$ sits entirely in $f^{-1}(1, 2) \cup (\cup_{\mathfrak{p} \in \Lambda} \mathcal{H}_\mathfrak{p})$. Let $\Theta \in \mathcal{Z}_{ech,M}$. Fix $\mathfrak{p} \in \Lambda$. It follows from Lemma 2.1 that $\Theta$ can contain either 0, 1 or two pairs of the form $(\gamma, 1)$ where $\gamma$ is an integral curve of $v$ that lies entirely in $\mathcal{H}_\mathfrak{p}$. Meanwhile, Lemma 2.2 with (1.29) and (1.30) has the following implication: There is one and only one segment of an integral curve of $v$ that crosses $\mathcal{H}_\mathfrak{p}$ from the $u = -R - \ln 7\delta_*$ end to the $u = R + \ln 7\delta_*$ end. This last observation, Proposition 2.7, and what was just said about the integral curves entirely inside $\cup_{\mathfrak{p} \in \Lambda} \mathcal{H}_\mathfrak{p}$ lead directly to Proposition 2.8.

Consider now the $\text{Spin}^{\mathbb{C}}$ structure issue. Given what is said in (1.30), the construction guarantees that the second and third bullets of (1.28) are obeyed. Note in this regard that the integral curves $\{\hat{\gamma}_\mathfrak{p}^+, \hat{\gamma}_\mathfrak{p}^-\}_{\mathfrak{p} \in \Lambda}$ are null-homologous as each lies in the $u = 0$ sphere of its labeling handle.

The verification of the top bullet in (1.28) is straightforward, but lengthy and so is omitted but for a proof of the following claim:

*Any two sets $\Theta$ and $\Theta'$ that are described by Proposition 2.8 determine the same $\text{Spin}^{\mathbb{C}}$ structure on M.*

To prove this, note first that set of $\text{Spin}^{\mathbb{C}}$-structures on Y is an affine space modeled on $H^2(Y; \mathbb{Z})$. In particular, the difference between any two $\text{Spin}^{\mathbb{C}}$ structures determines and is determined by a class in $H^2(Y; \mathbb{Z})$. With this fact in mind, suppose $\Theta$ and $\Theta'$ are each a collection of pairs of closed orbit of $v$ and positive integer. The class given by the 2-cycle $\sum_{(\gamma,m) \in \Theta} m[\gamma]^{Pd} - \sum_{(\gamma,m) \in \Theta'} m[\gamma]^{Pd}$ gives the class in $H^2(Y; \mathbb{Z})$ that determines the difference between the respective $\text{Spin}^{\mathbb{C}}$ structures that label $\Theta$ and $\Theta'$. If $\Theta$ and $\Theta'$ are both described by Proposition 2.8, then this difference class has pairing zero with the cross sectional spheres in all of the handles. As a consequence, its Poincaré dual is represented by a 1-cycle with support in the $f^{-1}(1, 2)$ part of M. What follows explains how to construct the desired 1-cycle from the respective HF-cycles that label $\Theta$ and $\Theta'$. View these HF-cycles as collections of G arcs in M that begin and end at the respective sets of index 1 and index 2 critical points of $f$. Concatenate the ends of $\Theta$'s segments with those of $\Theta'$. The result is an oriented 1-cycle with the orientation defined by segments from $\Theta$. A slight isotopy of this cycle sits where $f \in (1, 2)$. This cycle represents the Poincaré dual of the difference class in question.

## 3. Pseudoholomorphic subvarieties and the geometry of $\mathbb{R} \times Y$

This section takes a preliminary step towards describing the pseudoholomorphic subvarieties on $\mathbb{R} \times Y$ that are used to compute the differential and certain other endomorphisms of the $\mathcal{Z}_{ech,M}$-generated embedded contact homology chain complex. The first subsection describes the sorts of almost complex structures that are considered here. The remaining subsections describe a foliation of most of $\mathbb{R} \times Y$ by pseudoholomorphic submanifolds. The missing parts are the closed orbits of $\nu$ that sit where $u = 0$ and $1 - 3\cos^2\theta = 0$ in the various $\mathfrak{p} \in \Lambda$ versions of $\mathcal{H}_\mathfrak{p}$. This foliation plays an absolutely central role in the subsequent analysis of the curves that are relevant to $\mathcal{Z}_{ech,M}$.

### a) Almost complex structures and holomorphic subvarieties

The notion of pseudoholomorphic requires first the specification of an almost complex structure on $\mathbb{R} \times Y$. The almost complex structures that are considered in this section are constrained in certain ways. The constraints are described in Part 1 of this subsection. Additional constraints are imposed in Section 7 so as to make contact with the Heegaard-Floer story. Part 2 of this subsection describes two particularly important features of the almost complex structures that obey Part 1's constraints. Part 3 describes a topology for the set of pseudoholomorphic subvarieties.

*Part 1*: Let J denote a given almost complex structure. What follows lists the constraints that are imposed in this section on J.

CONSTRAINT 1: J *maps the Euclidean tangent vector $\partial_s$ to the $\mathbb{R}$ factor of $\mathbb{R} \times Y$ to the vector field $\nu$.*

CONSTRAINT 2: J *is not changed by constant translations of the coordinate $s$ on the $\mathbb{R}$ factor of $\mathbb{R} \times Y$.*

The statement of the third constraint refers to a certain 1-form, $\hat{a}$, on Y. The definion of $\hat{a}$ involves the 1-form $\upsilon_\diamond$ that is defined in Section 1e. The definition also requires a certain function of the coordinate u on any given $\mathfrak{p} \in \Lambda$ version of $\mathcal{H}_\mathfrak{p}$. This function is denoted by $\chi_\delta$ and it is defined by the rule

$$u \to \chi_\delta(u) = \chi(|u| - R - \ln\delta - 10).$$

(3.1)

The function $\chi_\delta$ equals 0 where $|u| \geq R + \ln\delta + 11$ and it equals 1 where $|u| < R + \ln\delta + 10$. In particular, $\chi_\delta = 1$ where the function $x$ is non-zero. The 1-form $\hat{a}$ is given by the rule

$$\hat{a} = \upsilon_\diamond \ \text{on} \ M_\delta \cup \mathcal{H}_0 \quad \text{and} \quad \hat{a} = \chi_\delta a + (1-\chi_\delta)\upsilon_\diamond \ \text{on} \ \cup_{\mathfrak{p}\in\Lambda} \mathcal{H}_\mathfrak{p}.$$

(3.2)

What was said above about $\chi_\delta$ guarantees that $\hat{a}$ is smooth, that $\hat{a} \wedge w > 0$ and that $\hat{a}$ has pairing 1 with $v$. The fact that $\hat{a} \wedge w > 0$ implies that the kernel of $\hat{a}$ defines a 2-dimensional subbundle in TY and that $w$ restricts to this subbundle as a symplectic form.

CONSTRAINT 3: J *preserves the kernel of the 1-form $\hat{a}$ and as such, it defines the orientation given by w.*

The remaining constraints concern the behavior of J on each $\mathfrak{p} \in \Lambda$ version of $\mathbb{R} \times \mathcal{H}_\mathfrak{p}$.

CONSTRAINT 4: J *on $\mathbb{R} \times \mathcal{H}_\mathfrak{p}$ is unchanged by $\mathbb{R}/2\pi\mathbb{Z}$ translations of the coordinate $\phi$.*

To set the background for the next constraint, use PROPERTY 3 and PROPERTY 4 of Section 1e to see that the kernel of $\hat{a}$ is spanned where u and $(1 - 3\cos^2\theta)\sin\theta$ are not both zero by the vector fields

- $e_1 = -6\,g\cos\theta\sin\theta\,\partial_u + (x + g')(1 - 3\cos^2\theta)\,\partial_\theta$.
- $e_2 = \partial_\phi + \sqrt{6}\chi_\delta f\cos\theta\sin^2\theta\,v$.

(3.3)

With the preceding understood, what follows is the last of the constraints that are imposed in this section.

CONSTRAINT 5: J *on $\mathbb{R} \times \mathcal{H}_\mathfrak{p}$ is such that $J e_1 = \sigma^{-1} e_2$ with $\sigma$ a positive function of* u *and* $\theta$.

Note that the function $\sigma$ is bounded near the locus where $u = 0$ and $1 - \cos^2\theta = 0$, and is further constrained there so as to insure that J is smooth. (This last constraint is not crucial to what follows here and in [KLTIII], but it does simplify some arguments.)

Assume henceforth that J obeys these constraints.

*Part 2*: The almost complex structure J has the following two properties:

- *The bilinear form $w(\cdot, J\cdot)$ is positive definite on the kernel of the 1-form a.*
- *There exists $r \geq 1$ such that the 2-form $\omega = ds \wedge a + rw$ tames J.*

(3.4)

To see about the top bullet, note first that CONSTRAINT 3 in Part 1 means that the bilinear form $w(\cdot, J\cdot)$ restricts to the kernel of $\hat{a}$ as a positive definite inner product. With this in mind, write a given vector $e \in \text{kernel}(a)$ as $e = x\nu + \hat{e}$ where $x \in \mathbb{R}$ and $\hat{e} \in \text{kernel}(\hat{a})$. The fact that $a \wedge w > 0$ implies that $|x| \leq c_0 w(\hat{e}, J\hat{e})^{1/2}$. Meanwhile, $w(e, Je)$ is equal to $w(\hat{e}, J\hat{e})$ because $\nu$ generates the kernel of $w$. To see about the second bullet, write a given tangent vector as $q = x_1 \partial_s + x_2 \nu + \hat{e}$ where $\hat{e} \in \text{kernel}(\hat{a})$. Then

$$\omega(q, Jq) = x_1^2 + x_2^2 + x_1 a(J\hat{e}) - x_2 a(\hat{e}) + r w(\hat{e}, J\hat{e}).$$

(3.5)

This is positive if and only if $q \neq 0$ when $r \geq c_0$ because $w(\cdot, J\cdot)$ is positive definite on the kernel of $\hat{a}$.

Note for reference below that J is compatible with the 2-form

$$\hat{\omega} = ds \wedge \hat{a} + w$$

(3.6)

in the sense that the bilinear form $\hat{\omega}(\cdot, J\cdot)$ defines a Riemannian metric on $T(\mathbb{R} \times Y)$. This metric is used implicitly in what follows to measure distances and to define the norms of vectors, covectors and other tensors.

*Part 3*: Given that J maps $\partial_s$ to $\nu$, the properties in (3.4) imply that the almost complex structure J obeys the conditions that are assumed in Part 4 of Section 1f. This being the case, (1.35) defines what is meant by a J-holomorphic subvariety. The set of J-holomorphic subvarieties is given the topology whereby a basis of open neighborhoods of a given subvariety C consists of sets that are indexed by a pair of positive number and smooth, compactly supported 2-form on $\mathbb{R} \times Y$. If $(\varepsilon, \varpi)$ is such a pair, then a J-holomorphic subvariety C´ is in the corresponding open neighborhood when

- $\sup_{z \in C} \text{dist}(z, C') + \sup_{z \in C'} \text{dist}(z, C) < \varepsilon$.
- $|\int_C \varpi - \int_{C'} \varpi| < \varepsilon$.

(3.7)

The resulting topological space is called the *moduli space* of J-holomorphic subvarieties.

The group $\mathbb{R}$ acts as a continuous group of homeomorphisms on the moduli space of J-holomorphic subvarieties. The action of the generator $1 \in \mathbb{R}$ on a given subvariety C translates C by 1 along the $\mathbb{R}$ factor in $\mathbb{R} \times Y$. Let Z denote a given space with an action of $\mathbb{R}$. A given map from a moduli space component to Z or from Z to a component is said to *$\mathbb{R}$-equivariant* when it intertwines the given $\mathbb{R}$ action on Z with this particular $\mathbb{R}$ action on the moduli space of J-holomorphic subvarieties.

**b) Level sets**

The almost complex structure J preserves the kernel of the 1-form $\upsilon_\diamond$ on $M_\delta \cup \mathcal{H}_0$. Meanwhile, the kernel of $\upsilon_\diamond$ is tangent to a surface either with or without boundary. These surface are the level sets in $M_\delta$ and the level sets of the coordinate u in $\mathcal{H}_0$. Those with boundary are the level sets of $f$ with $f \in [1 - 2\delta_*^2, 1 + \delta_*^2]$ or $f \in [2 - \delta_*^2, 2 + \delta_*^2]$. The surfaces without boundary are J-holomorphic submanifolds as defined by (1.31). As noted in the propositions that follow, those with boundary are subsets of J-holomorphic submanifolds that are described by (1.31). These propositions also indicate that there are important distinctions to be made depending on the following three alternatives: Either the level set in question lies in the union of $\mathcal{H}_0$ with the $f \notin [1, 2]$ part of $M_\delta$; or it is a level set in $M_\delta$ where $f \in (1, 2)$; or it is either the $f = 1$ or $f = 2$ level set in $M_\delta$. The proofs of the propositions are given at the end of the subsection.

The first proposition talks about the $f \notin [1, 2]$ case.

**Proposition 3.1**: *There exists a component, $\mathcal{M}_0$, of the moduli space of J-holomorphic subvarieties which has the properties listed below.*
- *Each subvariety from $\mathcal{M}_0$ is an embedded 2-sphere in $\mathbb{R} \times Y$.*
- *$\mathcal{M}_0$ has the structure of a smooth manifold which is $\mathbb{R}$-equivariantly diffeomorphic to $\mathbb{R} \times (-1, 1)$ with the action on the latter space given by the constant translations along the $\mathbb{R}$ factor. Such a diffeomorphism can be chosen so that:*
    a) *Suppose that $(x, y) \in \mathbb{R} \times [-\frac{1}{2}, \frac{1}{2}]$. The corresponding submanifold is $\{x\} \times S_y$ where $S_y \subset \mathcal{H}_0$ is the level set $u = 2(R + \ln\delta)y$.*
    b) *Suppose that $(x, y) \in \mathbb{R} \times [\frac{1}{2}, 1)$. The corresponding submanifold intersects the subset $\mathbb{R} \times M_\delta$ as $\{x\} \times S_y$ where $S_y$ is the $f = \delta^2 + (1 - \delta^2)(2y - 1)$.*
    c) *Suppose that $(x, y) \in \mathbb{R} \times (-1, -\frac{1}{2}]$. The corresponding submanifold intersects the*

subset $\mathbb{R} \times M_\delta$ as $\{x\} \times S_y$ where $S_y$ is the $f = 3 - \delta^2 - (1 - \delta^2)(2y - 1)$.

What follows are some parenthetical remarks for the cognicenti concerning $\mathcal{M}_0$. A certain Fredholm operator can be canonically associated to any given J-holomorphic subvariety. See, e.g. Definition 4.12 in [Hu2]. The $C \in \mathcal{M}_0$ version of this operator has Fredholm index equal to 2 and trivial cokernel. A subvariety is said to be *obstructed* when this operator has non-trivial cokernel. The subvariety is *unobstructed* when the cokernel is trivial. All of the subvarieties in $\mathcal{M}_0$ are unobstructed. Hutching in Definition 2.14 of [Hu2] has associated an integer to each component of the moduli space of J-holomorphic subvarieties. This integer is called the *ech index*. Each component of $\mathcal{M}_0$ has ech index equal to 2.

The next proposition talks about the $f \in (1, 2)$ level sets.

**Proposition 3.2**: *There exists a component, $\mathcal{M}_\Sigma$, of the moduli space of J-holomorphic subvarieties which has the properties listed below.*
- *Each subvariety from $\mathcal{M}_\Sigma$ is an embedded surface of genus G in $\mathbb{R} \times Y$.*
- *$\mathcal{M}_\Sigma$ has the structure of a smooth manifold which is $\mathbb{R}$-equivariantly diffeomorphic to $\mathbb{R} \times (1, 2)$ with the action on the latter space defined to be the translations along the $\mathbb{R}$ factor. Such a diffeomorphism can be chosen so that if $(x, y) \in \mathbb{R}$, then the corresponding manifold intersects $\mathbb{R} \times M_\delta$ as $\{x\} \times S_y$ where $S_y$ is the $f = y$ level set.*

By way of a parenthetical remark, note that each of the subvarieties from $\mathcal{M}_\Sigma$ is obstructed. The relevant Fredholm operator has index $2 - 2G$ and cokernel dimension equal to $2G$.

The final proposition in this subsection concerns the $f = 1$ and $f = 2$ level sets. These subvarieties have ends. When talking about the ends of these subvarieties, the proposition refers to the integral curve of $v$ and positive integer that is associated to a given end of a given J-holomorphic subvariety.

**Proposition 3.3**: *There exist components, $\mathcal{M}_1$ and $\mathcal{M}_2$, of the moduli space of J-holomorphic subvarieties which have the properties listed below.*
- *Each subvariety from either $\mathcal{M}_1$ or $\mathcal{M}_2$ is a properly embedded, 2G-times punctured sphere in $\mathbb{R} \times Y$.*
- *Each such submanifold from $\mathcal{M}_1$ or $\mathcal{M}_2$ has 2G negative ends and no positive ends.*
- *The ends of each such subvariety are labeled by the set $\{(\mathfrak{p}, \varepsilon)\}_{\mathfrak{p} \in \Lambda, \varepsilon \in \{-, +\}}$. This labeling is such that the integral curve of $v$ that is associated to the $(\mathfrak{p}, \varepsilon)$ labeled end is the curve $\hat{\gamma}_\mathfrak{p}^\varepsilon$, this the locus in $\mathcal{H}_\mathfrak{p}$ where $u = 0$ and $\cos\theta = \varepsilon \frac{1}{\sqrt{3}}$. The associated integer in each case is equal to 1.*

- $\mathcal{M}_1$ and $\mathcal{M}_2$ *have the structure of a smooth manifold that is $\mathbb{R}$-equivariantly diffeomorphic to the group $\mathbb{R}$. These diffeomorphisms can be chosen so that:*
  a) *The $x \in \mathbb{R}$ labeled submanifold in $\mathcal{M}_1$ intersects $\mathbb{R} \times M_\delta$ as $\{x\} \times f^{-1}(1)$.*
  b) *The $x \in \mathbb{R}$ labeled submanifold in $\mathcal{M}_2$ intersects $\mathbb{R} \times M_\delta$ as $\{x\} \times f^{-1}(2)$.*

The submanifolds in $\mathcal{M}_1$ and in $\mathcal{M}_2$ are obstructed; and both have ech-index 2 - 2G.

The remainder of this subsection is occupied with the

***Proofs of Propositions 3.1-3.3***: The proof of these propositions is given in six steps. The first five steps describe the J-holomorphic subvarieties that restrict to $\mathbb{R} \times M_\delta$ as $\{x\} \times S$ with $S \subset M_\delta$ a level set of $f$ with $f$ in $[1 - 2\delta^2, 1+\delta^2]$ or with $f$ in $[2 - \delta^2, 2+\delta^2]$. Note that the arguments are very similar to those used in Section 4a of [T1]

<u>Step 1</u>: Fix $p \in \Lambda$. Let $S \subset \mathcal{H}_p$ denote a closed, connected submanifold which is invariant under rotations of the angle $\phi$. Suppose in addition that the tangent space of S is mapped to itself by J. The fact that S is $\phi$-invariant implies that $\partial_\phi$ is tangent to S. Use (3.3) to write

$$\partial_\phi = e_2 - \sqrt{6}\chi_\delta f \cos\theta \sin^2\theta \, \nu \;,$$
(3.8)

and then use the fifth constraint in Part 1 to see that

$$J \cdot \partial_\phi = -\sigma e_1 + \sqrt{6}\chi_\delta f \cos\theta \sin^2\theta \, \partial_s \;.$$
(3.9)

As J maps TS to itself, the vector field $J \cdot \partial_\phi$ is also tangent to S. Use (3.3) and (3.9) to restate this fact as the follows: The vector field $J \cdot \partial_\phi$ is $-\sigma$ times

$$6g \cos\theta \sin\theta \, \partial_u - (x + g')(1 - 3\cos^2\theta) \partial_\theta + \sigma^{-1}\sqrt{6}\chi_\delta f \cos\theta \sin^2\theta \, \partial_s \;.$$
(3.10)

As such, the integral curves of this vector field foliate the complement in S of the fixed points of the $S^1$ action that rotates the angle $\phi$; and any two such curves are related by a rotation of the angle $\phi$.

<u>Step 2</u>: Let S be as described in Step 1 and let $\upsilon \subset S$ denote an integral curve of the vector field in (3.10). The function $u \to g(u)$ vanishes where $u = 0$ and this implies that u is either strictly positive, strictly negative or identically zero on S. The case where $u = 0$ on S is relevant to what is said in the next subsection and is not discussed further here.

The coefficient of $\partial_\theta$ in (3.10) vanishes where $1 - 3\cos^2\theta = 0$. This implies that the function $1 - 3\cos^2\theta$ is either strictly negative, or strictly positive, or identically zero on S. These different cases are considered respectively in Steps 3, 4 and 5. What is said in these steps explains why these case respectively concern Propositions 3.1, 3.2 and 3.3.

Step 3: This step assumes that $1 - 3\cos^2\theta < 0$ on S. It follows as a consequence that $\cos\theta$ is nowhere zero on S. Assume here that $\cos\theta > 0$. Suppose that g(u) and thus u is positive on S. If this is the case, then $\theta$, u and $s$ all decrease when traveling backwards along the curve $\upsilon$. As a consequence, the curve $\upsilon$ must end where $\cos\theta = 1$, which is by necessity, a fixed point of $\phi$ on S. This fixed point is the point on S where both u and $s$ take on their minimum values.

Moving in the oriented direction on $\upsilon$ out of this point increases u, $\theta$ and $s$. The function u increases faster and faster relative to $\theta$ as $1 - 3\cos^2\theta$ increases towards zero. As a consequence, u is eventually greater than $R + \ln\delta$ on $\upsilon$ and this part of $\upsilon$ is in the radius $7\delta_*$ coordinate ball centered on $\mathfrak{p}$'s index 1 critical point. Noting that S is an orbit orbit of $\upsilon$ with respect to the $S^1$ action generated by constant translations of $\phi$, it follows that S is an embedded disk in $\mathbb{R} \times \mathcal{H}_\mathfrak{p}$.

The function $x$ is zero where $u > R + \ln\delta$ on $\upsilon$, and therefore $s$ is constant on this part of $\upsilon$. As $J \cdot \partial_s = v$, such a slice must be a level set of $f$ and as $1 - 3\cos^2\theta < 0$, this level set is one where $f < 1$. This implies that S seemlessly extends some submanifold in $\mathbb{R} \times M_\delta$ of the form $\{x\} \times S_y$ as described in Item b) of Proposition 3.1 by attaching a disk to its $\cos\theta > 0$ boundary component on the boundary of the radius $\delta$ coordinate ball centered on $\mathfrak{p}$'s index 1 critical point. The case where $\cos\theta < 0$ on S has the analogous conclusion and is proved by essentially identical arguments.

Let $\{x\} \times S_y \subset \mathbb{R} \times M_\delta$ denote a surface with $S_y$ being an $f \in [1 - 2\delta^2, 1)$ level set. Integrate (3.10) from any point on its boundary components on the radius $\delta$ coordinate ball centered on $\mathfrak{p}$'s index 1 critical point to obtain two disks as just described that extends $\{x\} \times S_y$ into $\mathbb{R} \times \mathcal{H}_\mathfrak{p}$.

Much the same argument applies to the case where $u < 0$ on S. In this case, S extends some submanifold in $\mathbb{R} \times M_\delta$ of the form $\{x\} \times S_y$ with $S_Y$ as described in Item c) of Proposition 3.1.

Apply the preceding construction for each $\mathfrak{p} \in \Lambda$ and for a given surface in $\mathbb{R} \times Y$ of the form $\{x\} \times S_y$ with $S_y \subset M_\delta$ an $f \in [1 - 2\delta^2, 1)$ or $f \in (2, 2 - 2\delta^2]$ level set. The result is a compact, J-holomorphic submanifold in $\mathbb{R} \times Y$ that is diffeomorphic to $S^2$.

Step 4: This step assumes that $1 - 3\cos^2\theta > 0$ on S. Consider first the case where $u > 0$ and so $g > 0$. Start on $\upsilon$ at a point where $\cos\theta < 0$. The coefficient that multiplies $\partial_\theta$ in (3.10) is negative and so $\theta$ decreases along $\upsilon$. Meanwhile, both $s$ and u also decrease. Eventually, the $\theta = \frac{\pi}{2}$ point is reached, and then $\cos\theta$ turns positive. Both u and s start to increase along $\upsilon$ and do so from then on. Thus, the minimum of both u and s occur at the $\theta = \frac{\pi}{2}$ point on $\upsilon$. As $\cos\theta$ approaches $\frac{1}{\sqrt{3}}$, the rate of increase of u relative to $\theta$ becomes very large, and so eventually $u \geq R + \ln\delta$. Here $x$ is zero, so s is constant as is $f$. Rotating $\upsilon$ by changing $\phi$ shows that the $u \geq R + \ln\delta$ part of S where $\cos\theta > 0$ is an annular neighborhood of a boundary component of a surface of the form $\{x\} \times S_y \subset \mathbb{R} \times M_\delta$ where $S_y$ is some $y \in (1, 1+\delta^2)$ level set of $f$. Going in the opposite direction on $\upsilon$ from the $\cos\theta < 0$ starting point has u increasing and s also increasing, so the same story applies to the part of S where $u \geq R + \ln\delta$ and $\cos\theta < 0$. Note that the same surface $S_y$ appears here because changing $\theta$ to $\pi - \theta$ changes the sign of all terms in (3.10). Granted this, then attaching the $u < R + \ln\delta$ part of S to $\{x\} \times S_y$ seemlessly extends the latter by attaching a cylinder to its two boundary components on the radius $\delta$ coordinate ball centered on the index 1 critical point from $\mathfrak{p}$.

Let $(x, y) \in \mathbb{R} \times (1, 1+\delta^2]$ and let $\{x\} \times S_y$ denote the surface in $\mathbb{R} \times M_\delta$ with $S_y$ the level set of $f = y$. Integrate (3.10) from any point on its $\cos\theta < 0$ boundary on the radius $\delta$ coordinate ball centered on $\mathfrak{p}$'s index 1 critical point to obtain a cylinder S as just described that extends $\{x\} \times S_y$ into $\mathbb{R} \times \mathcal{H}_\mathfrak{p}$.

Much the same argument applies in the case where $u < 0$ on $\upsilon$. Here, the cylinder S extends a surface of the form $\{x\} \times S_y$ with $y \in [2-\delta^2, 2)$ and with $S_y$ again denoting the $f = y$ level set.

Apply the preceding construction for each $\mathfrak{p} \in \Lambda$ and for a given surface in $\mathbb{R} \times Y$ of the form $\{x\} \times S_y$ with $S_y \subset M_\delta$ an $f \in (1, 1+\delta^2]$ or $f \in [2-\delta^2, 2)$ level set. The result is a compact, J-holomorphic submanifold in $\mathbb{R} \times Y$ that is diffeomorphic to $\Sigma$.

Step 5: Consider here the case where $1 - 3\cos^2\theta = 0$ on S. Suppose that $u > 0$ and that $\cos\theta > 0$. Both u and s decrease as $\upsilon$ is traveled in reverse. The rate of decrease of u as it nears zero is proportional to u, while the rate of change of s remains bounded away from zero. As a consequence, s diverges towards $-\infty$ as $u \to 0$ on $\upsilon$. Moving in the forward direction on $\upsilon$ increases u and eventually $u \geq R + \ln\delta$ on $\upsilon$. Here, s is constant as is $f$. As S is invariant under translations of the coordinate $\phi$, this implies that the part of S where $u \geq R + \ln\delta$ is an annular neighborhood of a boundary component of a surface in $\mathbb{R} \times M_\delta$ of the form $\{x\} \times S_y$ where $S_y$ is the $f = 1$ level set in $M_\delta$. In particular, the part of S where $u \leq R + \ln\delta$ seemlessly attaches to this boundary component to extend the surface $\{x\} \times S_Y$. Integrating (3.10) will extend any surface of this sort into the $\cos\theta > 0$ part of $\mathbb{R} \times \mathcal{H}_\mathfrak{p}$. A very similar construction extends the $\cos\theta = -\frac{1}{\sqrt{3}}$ part of the intersection of such a surface with the radius $\delta$ coordinate ball centered on the index 1 critical point of $f$ from $\mathfrak{p}$.

But for cosmetics, the same construction applies in the case where $u < 0$ on $\upsilon$. Here, the cylinder S extends a surface of the form $\{x\} \times S_y$ with $S_y$ the $f = 2$ level set, and any surface of this sort can be so extended by integrating (3.10).

Apply the preceding construction for each $\mathfrak{p} \in \Lambda$ and for a given surface in $\mathbb{R} \times Y$ of the form $\{x\} \times S_y$ with $S_y \subset M_\delta$ an $f = 1$ or $f = 2$ level set. The result is a compact, J-holomorphic submanifold in $\mathbb{R} \times Y$ that is diffeomorphic to a sphere with 2G punctures. Such a submanifold has 2G ends, all negative ends, and with the $s \ll -1$ asymptotics as described by Proposition 3.3.

Step 6: This step completes the proof of Proposition 3.1. But for notation, the arguments given below complete the proofs of Propositions 3.2 and 3.3. This being the case, the final arguments for the latter two propositions are not supplied. These submanifolds together with those desribed in the first bullet of Proposition 3.1 form a set that is parametrized by $\mathbb{R} \times (-1, 1)$ as in Proposition 3.1. It is a straightforward task to verify that this parametrization is continuous with respect to the topology defined in Part 3 of Section 3a. The parametrization is $\mathbb{R}$-equivariant because the coefficients in (3.10) do not depend on the $\mathbb{R}$ factor of $\mathbb{R} \times Y$. As explained momentarily, this homeomorphism of $\mathbb{R} \times (-1, 1)$ into the moduli space of J-holomorphic subvarieties is onto a component of this moduli space. Granted that such is the case, the smooth structure on this component is defined by that on $\mathbb{R} \times (-1, 1)$.

To see about completeness, suppose that $\{C_k\}_{k=1,2,\ldots}$ is a sequence of J-holomorphic subvarieties in $\mathbb{R} \times Y$ that converges to a subvariety $C$ in $\mathcal{M}_0$. Since all submanifolds from $\mathcal{M}_0$ are compact, this must also be the case of all large k versions of $C_k$. If $C_k$ is not in $\mathcal{M}_0$, it must either intersects some subvariety on $\mathcal{M}_0$, or else it lies entirely in the union of the $f \in [1, 2]$ part of $\mathbb{R} \times M_\delta$ with the $1 - 3\cos^2\theta \le 0$ part of each $\mathfrak{p} \in \Lambda$ version of $\mathcal{H}_\mathfrak{p}$. If it intersects some subvariety in $\mathcal{M}_0$, then it must intersect all because J-holomorphic subvarieties intersect with strictly positive local intersection number. This can't happen if $C_k$ is compact. If it intersects the $f \in (1, 2)$ part of $\mathbb{R} \times M_\delta$, then it must intersect all of the subvarieties that are constructed in Step 4, and so it again can not be compact. If it interesects $\mathbb{R} \times M_\delta$ as the $f = 1$ or $f = 2$ level set, then it is one of the subvarieties constructed in Step 5 and so it can not be compact.

By way of a converse, suppose that $\{C_k\}_{k=1,2,\ldots} \subset \mathcal{M}_0$ is a sequence that converges to some J-holomorphic subvariety $C$. If $C$ is not in $\mathcal{M}_0$, then $C$ must intersect $\mathbb{R} \times M_\delta$ as either the $f = 1$ or $f = 2$ level set. Thus, $C$ is one of the subvarieties that are constructed in Step 5. Since $C$ is non-compact, the sort of convergence required by the first bullet in (3.7) is not possible.

### c) Subvarieties that lie entirely in $\mathbb{R} \times \mathcal{H}_{\mathfrak{p}}$.

The proposition that follows describes the J-holomorphic subvarieties that sit entirely in some $\mathfrak{p} \in \Lambda$ version of $\mathbb{R} \times \mathcal{H}_{\mathfrak{p}}$.

**Proposition 3.4**: *Fix $\mathfrak{p} \in \Lambda$. A J-holomorphic subvarieties that lies entirely in $\mathbb{R} \times \mathcal{H}_{\mathfrak{p}}$ sits in the $u = 0$ slice of $\mathbb{R} \times \mathcal{H}_{\mathfrak{p}}$. These subvarieties comprise three components of the moduli space of J-holomorphic subvarieties. The latter are denoted by $\mathcal{M}_{\mathfrak{p}+}$, $\mathcal{M}_{\mathfrak{p}-}$ and $\mathcal{M}_{\mathfrak{p}0}$. These components and their constituent subvarieties are described below.*

- *The component $\mathcal{M}_{\mathfrak{p}+}$: This moduli space component is smooth manifold that is equivariantly diffeomorphic to $\mathbb{R}$.*
  a) *Each element in $\mathcal{M}_{\mathfrak{p}+}$ is an embedded, open disk that lies where $\cos\theta > \frac{1}{\sqrt{3}}$.*
  b) *The end of each element is positive. The associated closed integral curve of $v$ is the curve $\hat{\gamma}_{\mathfrak{p}}^{+}$, this the $(u = 0, \cos\theta = \frac{1}{\sqrt{3}})$ locus in $\mathcal{H}_{\mathfrak{p}}$. The associated integer is 1.*
  c) *The function $s$ has a single critical point on each element; this is a non-degenerate critical point where $s$ is smallest. It occurs where $\theta = 0$.*
  d) *The parametrization by $\mathbb{R}$ can be chosen so that the parameter gives the minimum value of $s$.*

- *The component $\mathcal{M}_{\mathfrak{p}-}$: This moduli space component is smooth manifold that is equivariantly diffeomorphic to $\mathbb{R}$.*
  a) *Each element in $\mathcal{M}_{\mathfrak{p}+}$ is an embedded, open disk that lies where $\cos\theta < -\frac{1}{\sqrt{3}}$.*
  b) *The end of each element is positive. The associated closed integral curve of $v$ is the curve $\hat{\gamma}_{\mathfrak{p}}^{-}$, this the $(u = 0, \cos\theta = -\frac{1}{\sqrt{3}})$ locus in $\mathcal{H}_{\mathfrak{p}}$. The associated integer is 1.*
  c) *The function $s$ has a single critical point on each element; this is a non-degenerate critical point where $s$ is smallest. It occurs where $\theta = \pi$.*
  d) *The parametrization by $\mathbb{R}$ can be chosen so that the parameter gives the minimum value of $s$.*

- *The component $\mathcal{M}_{\mathfrak{p}0}$: This moduli space component is a smooth manifold that is equivariantly diffeomorphic to $\mathbb{R}$.*
  a) *Each element is an embedded cylinder that lies where $1 - 3\cos^2\theta > 0$.*
  b) *The ends of each element are positive. The curves $\hat{\gamma}_{\mathfrak{p}}^{+}$ and $\hat{\gamma}_{\mathfrak{p}}^{-}$ are their respective associated closed integral curves of $v$. The associated integer is 1 for both ends.*
  c) *The function $s$ has a single critical value on $s$, this the circle where $\theta = \frac{\pi}{2}$.*
  d) *The parametrization by $\mathbb{R}$ can be chosen so that the parameter gives the minimum value of $s$.*

As a parenthetical remark, note that $\mathcal{M}_{\mathfrak{p}_-}$ and $\mathcal{M}_{\mathfrak{p}_+}$ have ech index equal to 1, and each element is unobstructed. Meanwhile, $\mathcal{M}_{\mathfrak{p}0}$ has ech index 0 and is obstructed.

*Proof of Proposition 3.4*: Except for some minor changes, the arguments are the same as those used in the previous subsection to prove Propositions 3.1-3.3. In fact, the analysis of the J-holomorphic subvarieties that sit where $u = 0$ is simpler by virtue of the fact that the vector field in (3.10) can be replaced where $u = 0$ by the vector field

$$-(1-3\cos^2\theta)\,\partial_\theta + \sigma^{-1}\sqrt{6}\cos\theta\sin^2\theta\,\partial_s\,.$$

(3.11)

The reader can also see a very similar argument used in Example 2 in Section 4a of [T1] for the argument.

### d) Compact J-holomorphic subvarieties

The following lemma summarizes the contents of this last subsection.

**Lemma 3.5**: *The moduli space components $\mathcal{M}_0$ and $\mathcal{M}_\Sigma$ contain the only irreducible, compact J-holomorphic subvarieties in $\mathbb{R} \times Y$.*

*Proof of Lemma 3.5*: To see why this is, suppose to the contrary that C is an irreducible, compact J-holomorphic subvariety that is not in $\mathcal{M}_0$ nor $\mathcal{M}_\Sigma$. This subvariety can not be in $\mathcal{M}_1$ nor can it be in $\mathcal{M}_0$, nor any $\mathfrak{p} \in \Lambda$ version $\mathcal{M}_{\mathfrak{p}\pm,0}$ as all of the latter contain only non-compact subvarieties. This being the case, then C must intersect some submanifold from $\mathcal{M}_0$ or from $\mathcal{M}_\Sigma$. Given that J-holomorphic subvarieties intersect with strictly positive local intersection number, and given that all manifolds involved are compact, C must intersect all subvarieties from $\mathcal{M}_0$ if it intersects just one. Likewise, it must intersect all subvarieties from $\mathcal{M}_\Sigma$ if it intersects just one. This is can happen only if $|s|$ is unbounded on C.

### 4. Ech-HF subvarieties

Let $\vartheta$ denote as in Part 5 of Section 1f a finite collection of sets of the form (C, n) such that C is an irreducible, J-holomorphic subvariety and n is a positive integer. Part 5 of Section 1f explains how to associate to $\vartheta$ two elements, $\Theta_{\vartheta+}$ and $\Theta_{\vartheta-}$, in $\mathcal{Z}_{\text{ech}}$. Of primary interest in what follows are those $\vartheta$ where $\Theta_{\vartheta+}$ and $\Theta_{\vartheta-}$ are in $\mathcal{Z}_{\text{ech},M}$. Note that all pairs in such a set $\vartheta$ whose subvariety component is non-compact have integer component equal to 1. This is because all pairs in $\mathcal{Z}_{\text{ech},M}$ have integer component 1. This is also why distinct closed integral curves of $v$ are associated to distinct ends from the union of the subvarieties in $\vartheta$, and why each end from this union has associated integer 1.

If all integer components of $\vartheta$ are 1, then $\vartheta$ can be viewed as the set of irreducible components of $\cup_{(C,1)\in\vartheta} C$. This view is used in what follows. This view also motivates the following definition: A J-holomorphic subvariety whose irreducible components define a set $\vartheta$ with $\Theta_{\vartheta+}$ and $\Theta_{\vartheta-}$ in $\mathcal{Z}_{ech,M}$ is said to be an *ech-subvariety*. When C denotes an ech-subvariety, the corresponding sets in $\mathcal{Z}_{ech,M}$ are henceforth labeled as $\Theta_{C+}$ and $\Theta_{C-}$. An ech-subvariety is said to be an *ech-HF subvariety* if it lacks irreducible components that are described by Propositions 3.1-3.4.

This section begins the task of describing the ech-HF subvarieties.

**a) Projections to $\mathcal{H}_0$, the $f \notin (1, 2)$ part of $M_\delta$ and the $1 - 3\cos^2\theta \le 0$ part of $\mathcal{H}_p$**

The following proposition summarizes the conclusions of this subsection.

**Proposition 4.1**: *Let $C \subset \mathbb{R} \times Y$ denote an ech-HF subvariety. Then*

- *C has no intersection with $\mathbb{R} \times \mathcal{H}_0$ and with the parts of $\mathbb{R} \times M_\delta$ where $f \le 1$ or $f \ge 2$.*
- *C has no intersection with the $1 - 3\cos^2\theta < 0$ part of any $\mathfrak{p} \in \Lambda$ version of $\mathbb{R} \times \mathcal{H}_p$.*
- *An irreducible component of C that intersects the locus in $\mathbb{R} \times \mathcal{H}_p$ where $1 - 3\cos^2\theta = 0$ does so at $u = 0$. This component is the cylinder $\mathbb{R} \times \gamma$ where $\gamma \subset \mathcal{H}_p$ is the relevant $u = 0$ and $1 - 3\cos^2\theta = 0$ closed integral curve of $v$; thus either $\hat{\gamma}_\mathfrak{p}^+$ or $\hat{\gamma}_\mathfrak{p}^-$.*

*Proof of Proposition 4.1*: If the projection of C to Y contains a point in $\mathcal{H}_0$ or in $M_\delta$ where $f < 1$ or $f > 2$, then C must intersect a J-holomorphic subvariety in $\mathcal{M}_0$. Such an intersection has locally positive intersection number, and this means that C must intersect every subvariety in $\mathcal{M}_0$. This is true despite the fact that C is non-compact because the closed integral curves of $v$ that are associated to C's ends are disjoint from $\mathcal{H}_0$, from the part of $M_\delta$ where $f \notin (1, 2)$, and from the part $1 - 3\cos^2\theta < 0$ part of any $\mathfrak{p} \in \Lambda$ version of $\mathcal{H}_p$. However, this same observation about C's ends implies the following: Fix $\varepsilon > 0$ and there exists $x_\varepsilon > 1$ such that C has no intersections with the subvarieties from $\mathcal{M}_0$ that are parametrized as in Proposition 3.1 by a values of $(x, y)$ in $\mathbb{R} \times (-1+\varepsilon, 1-\varepsilon)$ with $|x| \gg x_\varepsilon$.

Suppose that C's projection to Y intersects $M_\delta$ where $f = 1$ or where $f = 2$, or that it intersects some $\mathfrak{p} \in \Lambda$ version of $\mathcal{H}_p$ where $1 - 3\cos^2\theta = 0$. Two cases can be distinguished. The first occurs when such an intersection occurs at some point that is not on the locus in some $\mathfrak{p} \in \Lambda$ version of $\mathcal{H}_p$ where both $u = 0$ and $1 - 3\cos^2\theta = 0$. Since C is has no component in either $\mathcal{M}_0$ or $\mathcal{M}_1$, so C must intersect some submanifold from $\mathcal{M}_0$ also. This is ruled out by what is said in the previous paragraph.

In the second case, the projection of C to Y intersects the relevant set only in some $\mathfrak{p} \in \Lambda$ version of $\mathcal{H}_\mathfrak{p}$ where $u = 0$ and $1 - 3\cos^2\theta = 0$. For the sake of argument, suppose this occurs where $\cos\theta = \frac{1}{\sqrt{3}}$. The subvariety C therefore intersects the J-holomorphic subvariety $\mathbb{R} \times \hat{\gamma}_\mathfrak{p}^+$. Suppose that the irreducible component involved is not this same $\mathbb{R} \times \hat{\gamma}_\mathfrak{p}^+$. Fix $\varepsilon \in (0, 1 - \frac{1}{\sqrt{3}}]$ and an angle $\phi_* \in \mathbb{R}/2\pi\mathbb{Z}$. It follows from (2.4) and Lemma 2.2 that there exists a segment of an integral curve of $v$ in $\mathcal{H}_\mathfrak{p}$ with the following properties: It enters $\mathcal{H}_\mathfrak{p}$ from the index 1 critical point end; it exits from the index 2 critical point end; it lies where the function $1 - 3\cos^2\theta$ is negative; and intersects the $u = 0$ locus at the point where $\cos\theta = \frac{1}{\sqrt{3}} + \varepsilon$ and $\phi = \phi_*$. If $\varepsilon$ is sufficiently small, this curve must intersect C. In particular, there are points on C where $1 - 3\cos^2\theta$ is negative, and so the previous analysis applies.

**b) Projection to the $f \in (1, 2)$ part of $M_\delta$**

This subsection makes some first observations about the intersection between a given ech-HF subvariety and the $f \in (1, 2)$ part of $M_\delta$. The discussion here has two parts.

*Part 1*: The lemma that follows states a fact that underpins all of the subsequent analysis.

**Lemma 4.2**: *Let* C *denote and ech-HF subvariety. Then* C *has intersection number* G *with each subvariety from* $\mathcal{M}_\Sigma$.

*Proof of Lemma 4.2*: The positive ends of C determine a set $\Theta_{C+} \in \mathcal{Z}_{ech,M}$. The integral curves of $v$ that comprise $\Theta_{C+}$ are pairwise disjoint in Y because their restrictions to $M_\delta$ are pairwise disjoint. Let $Z_+ \subset Y$ denote the union of the integral curves of $v$ from $\Theta_{C+}$. The submanifold $\mathbb{R} \times Z_+ \subset \mathbb{R} \times Y$ is a J-holomorphic subvariety.

The intersection of $Z_+$ with $M_\delta$ is a disjoint union of segments with each starting on the boundary of the radius $\delta$ coordinate ball centered on an index 1 critical point of $f$ and ending on the boundary of the radius $\delta$ coordinate ball centered on an index 2 critical point of $f$. Let $\{\gamma_1, \ldots, \gamma_G\}$ denote these segments. The submanifold with boundary $\mathbb{R} \times (\cup_{1 \le k \le G} \gamma_k) \subset \mathbb{R} \times Z$ is a properly embedded cylinder with boundary in $\mathbb{R} \times M_\delta$.

Fix a given point $y \in (1+\delta^2, 2-\delta^2)$ and let $S_y \subset M_\delta$ denote the corresponding level set of $f$. If $x \in \mathbb{R}$, then $\{x\} \times S_y \in \mathcal{M}_\Sigma$. This subvariety has intersection number G with $\mathbb{R} \times (\cup_{1 \le k \le G} \gamma_k)$ and thus with $\mathbb{R} \times Z_+$ since $S_y$ is disjoint from $Z_+ \cap (\cup_{\mathfrak{p} \in \Lambda} \mathcal{H}_\mathfrak{p})$. This implies that C has intersection number G with all $x \gg 1$ versions of $\{x\} \times S_y$.

Given that intersection numbers are locally constant under small perturbations, and given that $|s|$ is bounded on any given submanifold from $\mathcal{M}_\Sigma$, what is said in the previous paragraph implies that C has intersection number G with every submanifold from $\mathcal{M}_\Sigma$.

*Part 2*: This part uses Lemma 4.2 with Proposition 4.1 to say something about the intersections between a given ech-HF subvariety and the surfaces in $\mathcal{M}_\Sigma$ that are parametrized by pairs (x, y) with y near 1 or with y near 2. The relevant remarks are supplied by the next lemma.

**Lemma 4.3**: *Let C denote an ech-HF subvariety. Fix $\varepsilon \in (\delta, \delta_*)$ and there exists $z > 0$ with the following significance:*
- *Fix $(x, y) \in \mathbb{R} \times (1, 1+z]$ and let $S \in \mathcal{M}_\Sigma$ denote the corresponding submanifold. Then projection to Y of $C \cap S$ has precisely one point in the $|u| < R + \ln\varepsilon$ part of each $\mathfrak{p} \in \Lambda$ version of $\mathcal{H}_\mathfrak{p}$.*
- *Fix $(x, y) \in \mathbb{R} \times [2-z, 2)$ and let $S \in \mathcal{M}_\Sigma$ denote the corresponding submanifold. Then projection to Y of $C \cap S$ has precisely one point in the $|u| < R + \ln\varepsilon$ part of each $\mathfrak{p} \in \Lambda$ version of $\mathcal{H}_\mathfrak{p}$.*

**Proof of Lemma 4.3**: The submanifold C has G intersections with $\mathbb{R} \times Z_+$. Suppose that there exists a sequence $\{(x_k, z_k)\}_{k=1,2,\ldots} \subset \mathbb{R} \times (0, \delta_*^2]$ which does not converge and is such that for each $k \in \{1, 2, \ldots\}$, there exists at least one intersection point between C and the $(x = x_k, z = z_k)$ version of S whose projection to Y lies in $M_\varepsilon$. The sequence $\{x_k\}$ is bounded from above. To see why, note first that any constant $s \gg 1$ slice of C is very close to $Z_+$. Meanwhile, the curves that comprise $Z_+$ do not intersect an $f = 1 + z$ level set of $f$ in $M_\delta$ if $z < c_0^{-1}$ since this is true for the integral curves of $\mathfrak{v}$ that run from the index 1 critical points of $f$ to the index 2 critical points of $f$. A similar argument using the integral curves of $\nu$ from $\Theta_{C-}$ proves that $\{x_k\}_{k=1,2,\ldots}$ is bounded from below. As a consequence, no generality is lost by assuming that the latter sequence converges. This the case, it follows that the sequence $\{z_k\}_{k=1,2,\ldots}$ has limit zero. But then C must intersect some subvariety from $\mathcal{M}_1$ since $M_\delta$ is compact. Such an intersection is ruled out by what is said in Lemma 4.2.

Suppose now that there is a sequence $\{(x_k, z_k)\}_{k=1,2,\ldots} \subset \mathbb{R} \times (0, \delta_*^2]$ which does not converge and is such that for each $k \in \{1, 2, \ldots\}$, there exists at least two intersection points between C and the $(x = x_k, z = z_k)$ version of S whose projection to Y lies in some $\mathfrak{p} \in \Lambda$ version of $\mathcal{H}_\mathfrak{p}$ where $|u| < R + \ln\varepsilon$. An argument much like that given in the preceding paragraph proves that $\{x_k\}_{k=1,2,\ldots}$ must be bounded and so can be assumed convergent. Let $x_*$ denote the limit point. As before, the sequence $\{z_k\}_{k=1,2,\ldots}$ can be assumed to converge to zero. Fix $x_0 \gg 1$ so that the very large $s$ part of C is very close to the correspondingly large $s$ part of $\mathbb{R} \times Z_+$. If $x > x_0$ and k is large, then C has 1 intersection point with each $(x > x_0, z = z_k)$ version of S where $|u| < R + \ln \varepsilon$ in $\mathcal{H}_\mathfrak{p}$. As x is now decreased from $x_0$ towards $x_*$, there must exist for each large x, some point $x'_k$ between $x_0$ and $x_k$ with the property that C intersects the $(x = x'_k, z = z_k)$ version of S where $|u| > R + 2\ln \varepsilon$ in $\mathcal{H}_\mathfrak{p}$. But this is ruled out by what is said in the previous paragraph.

These last two lemmas provide a picture of C's intersection with $f^{-1}(1, 2)$ that plays a central role in Section 7 where the Heegaard-Floer geometry enters the story. To set the stage now, use the integral curves of $\mathfrak{v}$ as in Part 5 of Section 1c to define a diffeomorphism from $f^{-1}(1, 2) \subset M$ to $(1, 2) \times \Sigma$ which identifies $\mathfrak{v}$ with the Euclidean vector field $\partial_t$ along the $(1, 2)$ factor and which identifies $f$ with the coordinate $t$ on this same factor. The next lemma uses this diffeomorphism to view the $\mathbb{R} \times M_\delta$ part of C as a submanifold in $\mathbb{R} \times (1, 2) \times \Sigma$. The lemma also refers to the respective sets $C_+$ and $C_-$ in $\Sigma$ and their annular neighborhoods $T_+$ and $T_-$. These are introduced in Parts 2 and 4 of Section 1b.

**Lemma 4.4**: *There exists $\kappa > 1$ that is independent of C and the parameters $(\delta, x_0, R)$ and which has the following significance Suppose that C is an ech-HF subvariety. View the $\mathbb{R} \times M_\delta$ part of C as a subset of $\mathbb{R} \times (1, 2) \times \Sigma$ as explained above. There exists $\varepsilon > 0$ such that*
- *The $\mathbb{R} \times M_\delta$ part of C intersects $\mathbb{R} \times (1, 1+\varepsilon) \times \Sigma$ in $\mathbb{R} \times (1, 1+\varepsilon) \times T_+$.*
- *$\mathbb{R} \times M_\delta$ part of C intersects $\mathbb{R} \times (2-\varepsilon, 2) \times \Sigma$ in $\mathbb{R} \times (2, 2-\varepsilon) \times T_-$.*
- *The projection to $\Sigma$ of the $|s| > \varepsilon^{-1}$ part of C's intersection with $\mathbb{R} \times M_\delta$ lies in the union of the radius $\kappa\delta$ disks centered at the points in $C_- \cap C_+$.*

*Proof of Lemma 4.4*: The assertions follow directly from Lemmas 4.2 and 4.3; and from Corollary 2.6.

**c) The projection to $\cup_{\mathfrak{p} \in \Lambda} \mathcal{H}_\mathfrak{p}$**

Fix $\mathfrak{p} \in \Lambda$. This subsection talks about the intersection of a given ech-HF subvariety with $\mathbb{R} \times \mathcal{H}_\mathfrak{p}$. The section has five parts. Proposition 4.5 in Part 2 summarizes what is said. Part 1 supplies the background for the statement of this proposition. Parts 3 and 4 contain the proof of Proposition 4.5. Part 5 concerns the values of the function $\theta$ on $\mathbb{R} \times \mathcal{H}_\mathfrak{p}$ part of an ech-HF subvariety

This and subsequent sections use $\mathcal{H}^+_p$ to denote the part of $\mathcal{H}_p$ where the coordinates u and $\theta$ obey $|u| < R + \ln\delta_*$ and $1 - 3\cos^2\theta > 0$.

*Part 1*: What follows directly describes an extremely useful parametrization of $\mathcal{H}^+_p$. This parametrization involves the function $\hat{h}$ on $\mathcal{H}_p$ given by

$$\hat{h}(u, \theta) = (x + 2(\chi_+ e^{2(u-R)} + \chi_- e^{-2(u+R)}))\cos\theta \sin^2\theta .$$

(4.1)

In the notation of (1.24), the function $\hat{h} = f\cos\theta\sin^2\theta$. This function has no critical points on $\mathcal{H}^+_p$. In fact, $d\hat{h}$ has strictly negative pairing with the vector field $e_1$ in (3.3). More is true: It follows from (3.10) that $d\hat{h}$ pulls back without critical points to the $\mathcal{H}^+_p$ part of any submanifold from $\mathcal{M}_\Sigma$; and it follows from (3.11) that $d\hat{h}$ pulls back without critical points to any submanifold from $\mathcal{M}_0$. As the vector field $\partial_\phi$ annihilates $\hat{h}$, what was just said implies the following: The pair $(\hat{h}, \phi)$ define coordinates on the $\mathcal{H}^+_p$ part of any submanifold from $\mathcal{M}_\Sigma$ and on any submanifold from $\mathcal{H}_0$.

Note for future reference that the vector field $v$ also annihilates $\hat{h}$. This implies that $\hat{h}$ is constant on the segments of the integral curves of $v$ in $\mathcal{H}_p$. In particular, the level sets of $\hat{h}$ are foliated by the integral curves of $v$; and this implies that $\hat{h}$, if not constant, lacks local extrema on any submanifold in $\mathcal{H}^+_p$ whose tangent space is preserved by J.

Define an embedding from $\mathbb{R} \times \mathcal{H}^+_p$ into

$$\mathbb{R} \times (-R - \ln(\delta_*), R + \ln(\delta_*)) \times (\mathbb{R}/2\pi\mathbb{Z}) \times (-\tfrac{4}{3\sqrt{3}}\delta_*^2, \tfrac{4}{3\sqrt{3}}\delta_*^2)$$

(4.2)

as follows: Let q denote a given point in $\mathbb{R} \times \mathcal{H}^+_{p*}$. The image of q via this embedding has coordinates $(x, \hat{u}, \hat{\phi}, h)$ where $\hat{\phi} = \phi(q)$, $h = \hat{h}(q)$, and where $(x, \hat{u})$ are the respective values of $s \in \mathbb{R}$ and u at the $\theta = \tfrac{\pi}{2}$ point on the integral curve of the vector field in (3.10) that contains the point q.

Let $\Psi_p$ denote the inverse of the map just defined. The domain of $\Psi_p$ sits in (4.2) as a product $\mathbb{R} \times \mathcal{X}$ with $\mathcal{X} \subset (-R - \ln(\delta_*), R + \ln(\delta_*)) \times (\mathbb{R}/2\pi\mathbb{Z}) \times (-\tfrac{4}{3\sqrt{3}}\delta_*^2, \tfrac{4}{3\sqrt{3}}\delta_*^2)$ a certain p-independent open set.

To obtain a sense of this parameterization, let $(x, \hat{u}, \hat{\phi}, h)$ with $\hat{u} \neq 0$ denote a point in $\mathbb{R} \times \mathcal{X}$. Then $\Psi_p(x, \hat{u}, \hat{\phi}, h)$ is the point where $(\hat{h} = h, \phi = \hat{\phi})$ on the surface from $\mathcal{M}_\Sigma$ that also contains the point in $\mathbb{R} \times \mathcal{H}_p$ where $s = x$, $u = \hat{u}$, $\theta = \tfrac{\pi}{2}$, and $\phi = \hat{\phi}$. Meanwhile, a point $(x, 0, \hat{\phi}, h)$ is sent by $\Psi_p$ to the point in the surface from $\mathcal{M}_{p0}$ that also contains the point where $s = x$, $u = 0$, $\theta = \tfrac{\pi}{2}$ and $\phi = \hat{\phi}$.

The map $\Psi_{\mathfrak{p}}$ has the following very important properties:

- *The constant (x, û) surfaces in $\mathbb{R} \times X$ are mapped by $\Psi_{\mathfrak{p}}$ to J-holomorphic submanifolds.*
- *The map $\Psi_{\mathfrak{p}}$ is equivariant with respect to the $\mathbb{R}$ action on $\mathbb{R} \times X$ that translates along the first factor, and the corresponding action on $\mathbb{R} \times \mathcal{H}^+_{\mathfrak{p}}$.*
- *The map $\Psi_{\mathfrak{p}}$ is equivariant with respect to the $\mathbb{R}/2\pi\mathbb{Z}$ action that translates the coordinate $\hat{\phi}$ on $\mathbb{R} \times X$ and translates the coordinate $\phi$ on $\mathbb{R} \times \mathcal{H}^+_{\mathfrak{p}}$.*

(4.3)

These properties, and especially the first, motivate the introduction of $\Psi_{\mathfrak{p}}$. Note in particular that the constant (x, û ≠ 0) surfaces in $\mathbb{R} \times X$ are mapped to the intersections with $\mathbb{R} \times \mathcal{H}^+_{\mathfrak{p}}$ of the surfaces from $\mathcal{M}_\Sigma$. Meanwhile, the constant (x, û = 0) surfaces are mapped by $\Psi_{\mathfrak{p}}$ to the surfaces in $\mathcal{M}_{\mathfrak{p}0}$.

The formulas for the $\Psi_{\mathfrak{p}}$-pushforwards of the coordinate vector fields $\partial_x$, $\partial_{\hat{u}}$, $\partial_{\hat{\phi}}$ and $\partial_h$ are also worth keeping in mind. These are written using the notation from PROPERTY 3 in Section 1e and (3.3) as

- $\Psi_{\mathfrak{p}*}\partial_x = \partial_s$ ,
- $\Psi_{\mathfrak{p}*}\partial_{\hat{u}} = \upsilon(\nu + \alpha^{-1}\sqrt{6}x\cos\theta\,\partial_\phi + \varpi\,\partial_s)$ ,
- $\Psi_{\mathfrak{p}*}\partial_{\hat{\phi}} = \partial_\phi$ ,
- $\Psi_{\mathfrak{p}*}\partial_h = -\beta^{-1}(e_1 - \sigma^{-1}\sqrt{6}\chi_\delta f \cos\theta \sin^2\theta\,\partial_s)$ ,

(4.4)

where $\upsilon$, $\varpi$ and $\beta$ are certain functions of the pair (u, θ) with $\upsilon$ and $\beta$ positive.

*Part 2*: This part of the subsection states the promised proposition that describes the intersection of an ech-HF subvariety with a given $\mathfrak{p} \in \Lambda$ version of $\mathcal{H}^+_{\mathfrak{p}}$. To set the stage, fix $z_* \in (0, \frac{1}{4}\delta_*^2)$. The proposition refers to the subspace $\mathcal{M}_{\Sigma,*} \subset \mathcal{M}_\Sigma$ that is parametrized by (x, y) $\in \mathbb{R} \times (1, 2)$ with either $y \leq 1 + z_*$ or with $y \geq 2 - z_*$. It also refers to $\nu$'s integral curves $\hat{\gamma}^+_{\mathfrak{p}}$ and $\hat{\gamma}^-_{\mathfrak{p}}$ in $\mathcal{H}_{\mathfrak{p}}$ where (u = 0, cosθ = $\frac{1}{\sqrt{3}}$) and (u = 0, cosθ = $-\frac{1}{\sqrt{3}}$). By way of notation, the proposition introduces $I_*$ to denote $[-R - \frac{1}{2}\ln z_*, R + \frac{1}{2}\ln z_*]$.

The proposition also introduces the following terminology with regards to any given $\mathfrak{p} \in \Lambda$ version of the circles $\hat{\gamma}^+_{\mathfrak{p}}$ and $\hat{\gamma}^-_{\mathfrak{p}}$. Let $\hat{\gamma}$ denote one or the other of these integral curves of $\nu$. The curve $\hat{\gamma}$ is said to be *associated to an end of* C *in* $\mathcal{H}^+_{\mathfrak{p}}$ when the following occurs: The large |s| part of an end of C lies entirely in $\mathcal{H}^+_{\mathfrak{p}}$, and its constant |s| slices converge as $|s| \to \infty$ to $\hat{\gamma}$. Note in this regard that $\hat{\gamma}$ is *not* associated to an end of C in $\mathcal{H}^+_{\mathfrak{p}}$ if $\mathbb{R} \times \hat{\gamma}$ is a component of C.

**Proposition 4.5**: *Fix $\mathfrak{p} \in \Lambda$ and $z_* \in (0, \frac{1}{4}\delta_*^2)$. Suppose that C is an ech-HF subvariety which has intersection number 1 in $\mathcal{H}^+_\mathfrak{p}$ with each submanifold from $\mathcal{M}_{\Sigma,*}$. Then the $\Psi_\mathfrak{p}$-inverse image of C in $\mathbb{R} \times \mathcal{X}$ intersects the $\hat{u} \in I_*$ part of $\mathbb{R} \times \mathcal{X}$ as a graph of the form*

$$(x, \hat{u}) \to (x, \hat{u}, \hat{\phi} = \varphi^\mathfrak{p}(x, \hat{u}), h = \varsigma^\mathfrak{p}(x, \hat{u}))$$

*where the pair $(\varphi^\mathfrak{p}, \varsigma^\mathfrak{p})$ denotes a smooth map from certain domain in $\mathbb{R} \times I_*$ to $(\mathbb{R}/2\pi\mathbb{Z}) \times (-\frac{4}{3\sqrt{3}}\delta_*^2, \frac{4}{3\sqrt{3}}\delta_*^2)$. This domain is as follows:*

- *If neither $\hat{\gamma}^+_\mathfrak{p}$ or $\hat{\gamma}^-_\mathfrak{p}$ is associated to an end of C in $\mathbb{R} \times \mathcal{H}^+_\mathfrak{p}$, then the domain of $(\varphi^\mathfrak{p}, \varsigma^\mathfrak{p})$ is the whole strip $\mathbb{R} \times I_*$,*
- *If just one of $\hat{\gamma}^+_\mathfrak{p}$ or $\hat{\gamma}^-_\mathfrak{p}$ is associated to an end of C in $\mathbb{R} \times \mathcal{H}^+_\mathfrak{p}$, then the domain of $(\varphi^\mathfrak{p}, \varsigma^\mathfrak{p})$ is the complement in $\mathbb{R} \times I_*$ of a single $\hat{u} = 0$ point.*
- *If both $\hat{\gamma}^+_\mathfrak{p}$ and $\hat{\gamma}^-_\mathfrak{p}$ are associated to ends of C in $\mathbb{R} \times \mathcal{H}^+_\mathfrak{p}$, then the domain of $(\varphi^\mathfrak{p}, \varsigma^\mathfrak{p})$ if the complement in $\mathbb{R} \times I_*$ of two $\hat{u} = 0$ points.*

This proposition implies that C's intersection with the image of $\Psi_\mathfrak{p}$ is a smooth, properly embedded submanifold with boundary in this image.

The proof of this proposition is given in the next part of this subsection.

*Part 3*: The proof of Proposition 4.5 invokes two auxilliary lemmas. Both are proved in Part 4. The first of the lemmas is very specific to the situation at hand.

**Lemma 4.6**: *Let C denote an ech-HF subvariety.*
- *If neither $\hat{\gamma}^+_\mathfrak{p}$ nor $\hat{\gamma}^-_\mathfrak{p}$ is associated to an end of C in $\mathbb{R} \times \mathcal{H}^+_\mathfrak{p}$, then C has intersection number 1 with each subvariety from $\mathcal{M}_{\mathfrak{p}0}$.*
- *If just one of $\hat{\gamma}^+_\mathfrak{p}$ or $\hat{\gamma}^-_\mathfrak{p}$ is associated to an end of C in $\mathbb{R} \times \mathcal{H}^+_\mathfrak{p}$, then C has intersection number zero with one subvariety from $\mathcal{M}_{\mathfrak{p}0}$ and intersection number 1 with each of the others.*
- *If both $\hat{\gamma}^+_\mathfrak{p}$ and $\hat{\gamma}^-_\mathfrak{p}$ are associated to ends of C in $\mathbb{R} \times \mathcal{H}^+_\mathfrak{p}$, then C has intersection number zero with two subvarieties from $\mathcal{M}_{\mathfrak{p}0}$ and intersection number 1 with each of the others.*

The next lemma contains a fundamental observation about pseudoholomorphic curves and fibrations with pseudoholomorphic fiber. To set the stage for this lemma, suppose that U is an open set with compact closure in a smooth 4-manifold with an almost complex structure. A subset $C \subset U$ is said to be *a pseudoholomorphic subvariety* *in* U *when the following conditions are met.* First, the closure of U is contained in an open set and this open set contains a closed subset that lacks disconnected components and intersects U as C. Second, this larger closed subset contains a finite set whose complement is a smooth 2-dimensional submanifold. Third, the tangent space of this submanifold is mapped to itself by the almost complex structure.

**Lemma 4.7**: *Let* U *denote an open set in a 4-manifold with almost complex structure. Let* S *denote a 2-dimensional surface and suppose that the closure of* U *is contained in an open set with a proper fibration to* S *with pseudoholomorphic fibers. Let* C *denote a pseudoholomorphic subvariety in* U *with no open set in a fiber. Assume in addition that* C *has intersection number at most 1 with each fiber. Then* C *is a smooth submanifold and the fibering map to* S *restricts to* C *as a diffeomorphism onto an open set in* S.

Granted for the moment these two lemmas, here is the promised

***Proof of Proposition 4.5***: The composition of first the inverse of $\Psi_p$ and then the projection from $\mathbb{R} \times \mathcal{X}$ to $\mathbb{R} \times (-R - \ln(\delta_*), R + \ln(\delta_*))$ is a submersion with J-holomorphic fibers. Granted that this is the case, then the depiction by the proposition of the $\Psi_p$-inverse image of C follows from Lemmas 4.5 and 4.6 with the fact that the pair $(\phi, \hat{h})$ restrict as coordinates to the $\mathbb{R} \times \mathcal{H}^+_p$ part of surfaces from $\mathcal{M}_\Sigma$ and they also restrict as coordinates to each subvariety from $\mathcal{M}_{p0}$.

*Part 4*: This part of the subsection contains the proofs of Lemmas 4.6 and 4.7.

***Proof of Lemma 4.6***: The proof has seven steps.

   Step 1: The set $\Theta_{C+}$ contains one integral curve of $v$ whose intersection with $\mathcal{H}_p$ runs from the $u = -R - \ln(7\delta_*)$ end to the $u = R + \ln(7\delta_*)$ end. This curve intersects the $u = 0$ slice of $\mathcal{H}_p$, and this implies that each $x \gg 1$ labeled submanifolds in $\mathcal{M}_{p0}$ intersect C precisely once. The analogous argument using $\Theta_{C-}$ proves that the $x \ll -1$ labeled submanifolds in $\mathcal{M}_{p0}$ also intersect C precisely once. If C is described by the first bullet of Lemma 4.6, then the continuity of the intersections as x varies in $\mathbb{R}$ implies that each submanifold in $\mathcal{M}_{p0}$ must intersect C exactly once.

   Step 2: Before considering the second and third bullets, use (3.11) with Lemma 2.3 to see that the restriction of $\Psi_p^* s$ to the $\hat{u} = 0$ slice of the domain of $\Psi_p$ has the form

$$\Psi_p^* s = x - \mu \ln(1 - 3\cos^2 \hat{\theta}) + w(\hat{\theta})$$

(4.5)

where $\mu > 0$ is a constant, w(·) is a smooth function on $\mathbb{R}$ and $\hat{\theta}$ is the implicit function of the coordinate h on $\mathbb{R} \times X$ that is given by the unique $\theta \in (0, \pi)$ with $1 - 3\cos^2\theta > 0$ and is such that $h = (x_0 + 4e^{-2R}) \cos\theta \sin^2\theta$. The pull-back of s to where û is non-zero has the same schematic form of what is written on the right hand side of (4.5), but with $\mu$ a function of û, with $\hat{\theta}$ a function of h and û, and with w now a function of û and $\hat{\theta}$. When using (4.5) in what follows, keep in mind that the very small |u| part of $\mathcal{H}_p$ has inverse image via $\Psi_p$ in the very small |û| part of $\Psi_p$'s domain (but not vice-versa).

Step 3: This step explains why neither $\hat{\gamma}_p^+$ nor $\hat{\gamma}_p^-$ are associated to negative ends of C. Suppose to the contrary that $1 - 3c^2 \to 0$ as $s \to -\infty$ on an end $\mathcal{E} \subset C$ in $\mathbb{R} \times \mathcal{H}_p$. This can occur only if $u \to 0$ on the constant s slices of $\mathcal{E}$. This requires that û limit to zero on $\Psi_p^{-1}(\mathcal{E})$. This observation implies that there exists a path $\upsilon \subset \mathbb{R} \times I_*$ of the following sort: First, $\upsilon$ is in the domain of the functions $(\varphi^p, \varsigma^p)$ so $\upsilon$ defines a corresponding path in C. Denote this path by $\upsilon_\mathcal{E}$. This path lies in $\mathcal{E}$, and s is unbounded from below on $\upsilon_\mathcal{E}$. Granted these properties, use (4.5) to see that $s \to -\infty$ on $\upsilon_E$ only if $x \to -\infty$ on $\upsilon$. Since $u \to 0$ on $\upsilon_\mathcal{E}$ as $s \to -\infty$, so $\hat{u} \to 0$ on $\upsilon$ as $x \to -\infty$ on $\upsilon$. The existence of such a path is not possible given that the $s \ll -1$ part of C in $\mathbb{R} \times \mathcal{H}_p$ has a second $s \to -\infty$ limit of its constant s slices. The latter is a segment of an integral curve of $v$ from $\Theta_{C^-}$. This segment crosses the $u = 0$ slice transversely, and so its inverse image via $\Psi_p$ has transverse intersection with the $\hat{u} = 0$ slice of $\Psi_p$'s domain. It follows as a consequence that the $x \ll -1$ part of a path such as $\upsilon$ must have more than 1 inverse image in $\Psi_p^{-1}(C)$. This nonsense proves that neither $\hat{\gamma}_p^+$ nor $\hat{\gamma}_p^-$ are associated a negative end of C.

Step 4: Suppose that $\hat{\gamma}_p^+$ or $\hat{\gamma}_p^-$ are associated to a positive end of C. Let $\mathcal{E}$ denote the latter. The $s \gg 1$ part of $\mathcal{E}$ can be parametrized by the correspondingly large s part of $\mathbb{R} \times (\mathbb{R}/(2\pi\mathbb{Z}))$ as

$$(s, \phi) \to (s, u = b(s, \phi), \theta = \theta_* + a(s, \phi), \phi)$$

(4.6)

where (a, b) are smooth functions of their arguments that limit uniformly to 0 as $s \to \infty$.

Because the integral curves involved are non-degenerate, the functions (a, b) must behave at large s as

$$a = e^{-\lambda s}(a_*(\phi) + \mathfrak{e}_a) \quad and \quad b = e^{-\lambda s}(b_*(\phi) + \mathfrak{e}_b)$$

(4.7)

where λ is a positive constant and $(a_*, b_*)$ are smooth functions on $\mathbb{R}/2\pi\mathbb{Z}$ that are not simultaneously zero. Meanwhile and $\mathfrak{e}_a$ and $\mathfrak{e}_b$ are smooth functions of $(s, \phi)$ with limit zero as $s \to \infty$. This can be proved using the arguments used to prove Proposition 2.4 in [HT]. The arguments in the latter reference identify $(a_*, b_*)$ as eigenvectors of a symmetric, first order differential operator on $C^\infty(S^1; \mathbb{R}^2)$ that is defined by the Lie transport by $v$ on the normal bundle to the relevant integral curve and the metric $w(\cdot, J \cdot)$ on this bundle. The constant λ is the corresponding eigenvalue.

In the cases at hand, the integral curves are orbits of the $\mathbb{R}/2\pi\mathbb{Z}$ action that translates φ and both $v$ and J are invariant with respect to this action. This implies that any given eigenfunction of the relevant operator must have the form

$$a_*(\phi) = a_n \cos(n\phi - \phi_n) \quad and \quad b_*(\phi) = b_n \sin(n\phi - \phi_n) \tag{4.8}$$

where $a_n, b_n$ are $\mathbb{R}$-valued constants and $\phi_n \in \mathbb{R}/2\pi\mathbb{Z}$.

Step 5: In the case at hand, only the n = 0 version of (4.8) can occur because $\theta - \theta_*$ does not change sign on $\mathcal{E}$. In this case, $(a_*, b_*)$ are constants that are not both zero. Note that $a_* \geq 0$ if $\gamma_p^+$ is associated to $\mathcal{E}$ and $a_* \leq 0$ if $\gamma_p^-$ is associated to $\mathcal{E}$. As explained next, the only case that can occur here is that where $a_* \neq 0$ and $b_* = 0$.

To see why this is, note first that the argument given in Step 3 can be repeated no matter the value of $b_*$ to see that x stays bounded on $\Psi_p^{-1}(\mathcal{E})$. What follows is a consequence. There exists a bounded set $\mathrm{I} \subset \mathbb{R}$ with the following property: Given $T \geq 1$ and $\varepsilon > 0$, there exists $(x, \hat{u})$ with $x \in \mathrm{I}$ and $|\hat{u}| < \varepsilon$ that are in the domain of $(\varphi^p, \varsigma^p)$ and parametrize a point on $\mathcal{E}$ where $s > T$.

Suppose for the sake of argument that $b_* > 0$. This implies via (4.6)-(4.8) that u is positive on the $s \gg 1$ part of $\mathcal{E}$. Therefore, the $s \gg 1$ part of $\mathcal{E}$ is parameterized via $\Psi_p$ only by points where $\hat{u}$ is positive. Let $\{x_n\}_{n=1,2,\ldots}$ denote a sequence of points in $\mathbb{R}$ such that there exists $\hat{u} \in \mathrm{I}_*$ with $\hat{u} \in (0, \frac{1}{n})$ and such that $(x_n, \hat{u})$ is a point in the domain of $(\varphi^p, \varsigma^p)$ that parametrizes and $s > n$ point in $\mathcal{E}$. Let $x_* \in \mathrm{I}$ denote a limit point of the $\{x_n\}_{n=1,2,\ldots}$. The interval $\{x_*\} \times [-1, 0)$ in $\mathbb{R} \times \mathrm{I}_*$ lies in the domain of $(\varphi^p, \varsigma^p)$ and so corresponds via $\Psi_p$ to an embedded arc in C. The function s is bounded on this arc if $b_*$ is positive, and so the arc has a limit on the u = 0 slice of C. A disk about the limit point is parameterized via $(\varphi^p, \varsigma^p)$ by an open neighborhood of $(x_*, 0)$ in $\mathbb{R} \times \mathrm{I}_*$. But if this is the case, then C will have two or more intersections in $\mathbb{R} \times \mathcal{H}_p^+$ with some elements in $\mathcal{M}_{\Sigma,*}$.

Step 6: Assume here that only one of $\hat{\gamma}_p^+$ and $\gamma_p^-$ is associated to an end of C. Let $\theta_* \in (0, \pi)$ denote the angle of the component of the $1 - 3c^2 = 0$ locus in question.

Use $x_+ \in \mathbb{R}$ to denote the maximum value of x on the complement of the domain of $(\varphi^p, \varsigma^p)$, and let $x_- \in \mathbb{R}$ denote the minimum value. Use the fact that $a_* \neq 0$ and $b_* = 0$ in (4.6)-(4.8) to write (4.5) at points $(x, 0) \in \mathbb{R} \times I_*$ with $x > x_+$ as

$$(1 - \lambda\mu)\Psi_p * s = x - \mu \ln(2\sqrt{2}|a_*|) + w(\theta_*) + \cdots \tag{4.9}$$

where the unwritten terms have limit zero as $s \to \infty$. Since $x > x_+$, this equation requires that $\lambda\mu = 1$ and it identifies $x_+ = \mu \ln(2\sqrt{2}|a_*|) + w(\theta_*)$. The analogous argument finds this same value for $x_-$. It follows that the domain of $(\varphi^p, \varsigma^p)$ is the complement in $\mathbb{R} \times I_*$ of a single $\hat{u} = 0$ point. This proves what is asserted by the second bullet of Lemma 4.6

Note that the argument in Step 5 proves the following: Given $T > 1$ there exists positive $\varepsilon_T$ such that the complement of the origin of the radius $\varepsilon_T$ disk in $\mathbb{R} \times I_*$ about this missing point parametrizes only points in C from the $s > T$ part of $E$.

Step 7: Assume in this last step that both $\hat{\gamma}_p^+$ and $\hat{\gamma}_p^-$ are associated to ends of C. Let $E_+$ and $E_-$ denote the associated domains. The domain of $(\varphi^p, \varsigma^p)$ in this case is the complement in $\mathbb{R} \times I_*$ of some bounded set $\mathrm{I} \subset \mathbb{R}$ on the $\hat{u} = 0$ slice. The third bullet asserts that $\mathrm{I}$ consist of two points. To see why this is, note first that the function s on the image via $\Psi_p$ of the domain of $(\varphi^p, \varsigma^p)$ must have an index 1 critical point that corresponds to the saddle between $E_+$ and $E_-$. To elaborate, let $Y$ denote the space of compact, embedded paths in C of the following sort: Both ends are parametrized by $\hat{u} = 0$ points in the domain of $(\varphi^p, \varsigma^p)$. In addition, one end of the path lies in $E_+$ and the other lies in $E_-$. The function s achieves its minimum on each path in $Y$, and the set of these minima is bounded from above. Let $s_*$ denote the supremum.

Fix a path, $\upsilon_C \subset Y$. This path is parametrized by a corresponding path, $\upsilon \subset \mathbb{R} \times I_*$ in the domain of $(\varphi^p, \varsigma^p)$. This path must cross the $\hat{u} = 0$ locus at a point with components from $\mathrm{I}$ on both sides. For if not, the path can be homotoped rel end points in the domain of $(\varphi^p, \varsigma^p)$ to a path whose image is in $Y$ has s everywhere greater than $s_*$.

Granted that $\mathrm{I}$ has at least two components, the argument used in Step 6 can be repeated to prove that there are precisely two components and each is a single point.

*Proof of Lemma 4.7*: Use $\pi: U \to S$ to denote the fibering map. The map $\pi$ is 1-1 and continuous on C because pseudoholomorphic subvarieties have locally positive intersection number. The differential of π at each smooth point of C must be invertible. To see why, note that this differential is zero where C is tangent to a fiber. But there can be no tangencies. This can be proved by using the local structure of intersecting pseudoholomorphic curves to see that C must have intersection at least two with the fibers of π if it has a tangency See, e.g. Proposition 2.6 in [Mc], or Lemma A.7 in [T2]. The local picture given in these references also imply that C must have intersection number greater than 1 with a fiber through any singular point.

*Part 5*: The last part of the subsection first states and then proves a lemma about the values of function θ on pseudoholomorphic subvarieties in $\mathbb{R} \times \mathcal{H}_\mathfrak{p}$.

**Lemma 4.8**: *Let* $U \subset \mathbb{R} \times \mathcal{H}_\mathfrak{p}$ *denote an open set with compact closure and let* $C \subset U$ *denote an embedded, pseudoholomorphic subvariety in* U. *If the function* θ *is not constant on* C, *then any local minima and maxima of* θ *on* C *lies where* u = 0.

***Proof of Lemma 4.8***: The function cosθ on C can not take a local maximum or local minimum where $1 - 3\cos^2 = 0$, where cosθ = 0 nor where cosθ = ±1, as the product of $\mathbb{R}$ with these loci is foliated by J-holomorphic submanifolds. It remains to prove that cosθ has no local extremum on C in the complement of C's intersection with the three loci where u = 0, where $\cos\theta \sin^2\theta = 0$ and where $(1 - 3\cos^2\theta) = 0$. This is done by applying the maximum principle to a second order equation for θ's restriction to this part of C.

To start the argument, introduce $\{e^1, e^2\}$ to denote the dual 1-forms to the vector fields $e_1$ and $e_2$ that are depicted in (3.3). These are given by

$$e^1 = \beta^{-1}(f(1-3\cos^2\theta)\,d\theta - f'\cos\theta\sin\theta\,du) \quad \text{and} \quad e^2 = \beta_*^{-1}(d\phi + \sqrt{6}x\alpha^{-1}\cos\theta\,\hat{a}).$$

(4.10)

where β is the positive function of (u, θ) that appears in (4.4) and $\beta_*$ is a second, strictly positive function of (u, θ). It follows from (4.3) and from the formulae in PROPERTY 3 of Section 1e that the 1-form dθ can be written as

$$d\theta = (x + g')(1 - 3\cos^2\theta)\,e^1 + \sqrt{6}\alpha^{-1}\chi_\delta f'f\cos^2\theta\sin^3\theta\,e^2 + \alpha^{-1}f'\cos\theta\sin\theta\,\hat{a}$$

(4.11)

Let $J^T: T^*(\mathbb{R} \times Y) \to T^*(\mathbb{R} \times Y)$ denote the adjoint of J. Then

$$J^T \cdot d\theta = -\sigma(x + g')(1 - 3\cos^2\theta)\,e^2 + \sqrt{6}\sigma^{-1}\alpha^{-1}f'f\cos^2\theta\sin^3\theta\,e^1 + \alpha^{-1}f'\cos\theta\sin\theta\,ds.$$

(4.12)

This can be written as

$$J^T \cdot d\theta = \sigma(x+g')(1-3\cos^2\theta)\,d\phi + b_1 d\theta + b_2 du + \alpha^{-1} f' \cos\theta \sin\theta\, ds$$

(4.13)

where $b_1$ and $b_2$ are functions of u and θ. It follows from this that $d(J^T \cdot d\theta)$ can be written schematically as

$$d(J^T \cdot d\theta) = \mathfrak{z}_1 \wedge d\theta + (1-3\cos^2\theta)\, b_3\, du \wedge d\phi + \cos\theta \sin\theta\, b_4\, du \wedge ds,$$

(4.14)

where $\mathfrak{z}_1$ is a 1-form and $b_3$ and $b_4$ are functions of θ and u.

To continue, let λ denote a given function of u and θ that is defined and nowhere zero on the complement of the three loci where $u = 0$, $(1 - 3\cos\theta) = 0$ and $\cos\theta \sin^2\theta = 0$. Use (4.13) and (4.14) to see that

$$\lambda^{-1} d(\lambda J^T \cdot d\theta) = \mathfrak{z}_2 \wedge d\theta + (b_3 + \lambda^{-1}\lambda_u \sigma(x+g'))(1-3\cos^2\theta)\,du \wedge d\phi$$
$$+ (b_4 + \lambda^{-1}\lambda_u \alpha^{-1} f' \cos\theta \sin\theta)\,du \wedge ds.$$

(4.15)

Here, $\lambda_u$ is shorthand for $\partial_u \lambda$ and $\mathfrak{z}_2$ is another 1-form.

At the points where $f' \cos\theta \sin^2\theta \neq 0$, use (4.10) to write

$$du = -(f' \cos\theta \sin\theta)^{-1} \beta e^1 + b_5\, d\theta$$

(4.16)

where $b_5$ is a function of u and θ. This done, then (4.15) can be written as

$$\lambda d(\lambda J^T \cdot d\theta) = \mathfrak{z}_3 \wedge d\theta + (b_3 + \lambda^{-1}\lambda_u \sigma(x+g'))(1-3\cos^2\theta)\,du \wedge d\phi - (b_6 + \lambda^{-1}\lambda_u \beta\alpha^{-1})\, e^1 \wedge ds.$$

(4.17)

Here $\mathfrak{z}_3$ is a 1-form and $b_6$ is a function of u and θ. The restriction to TC of $e^1 \wedge ds$ is that of $\sigma e^2 \wedge \hat{a}$. Modulo forms with a factor of $d\theta$, the latter is

$$\beta_*^{-1} \sigma(x+g')(1-3\cos^2\theta)\,d\phi \wedge du.$$

(4.18)

Noting the function multiplying $d\phi \wedge du$, it follows from (4.17) that up to terms with $d\theta$, the right hand side of (4.17) restricts to TC as

$$(1 + \beta_*^{-1}\beta\alpha^{-1})(b_7 + \lambda^{-1}\lambda_u)\sigma(x+g')(1-3\cos^2\theta)\,du \wedge d\phi.$$

(4.19)

where $b_7$ is yet another function of u and θ. Granted (4.19), note that λ can be chosen on the complement of the three loci where $u = 0$, $\cos\theta \sin^2\theta = 0$ and $(1 - 3\cos^2\theta) = 0$ so as to make (4.19) equal to zero. This done, then the restriction of (4.17) to TC on the complement of these loci has the form

$$d(J^T d\theta) + \mathfrak{R} \cdot d\theta$$
(4.20)

where $\mathfrak{R}$ is a homomorphism from T*C to $\wedge^2$T*C. What with (4.10), the claim made in the lemma follows using the strong maximum principle.

## 5. Quantative concerns

This section is concerned with various cohomological and numerical invariants that can be associate to any given ech-HF subvariety. These invariants are used in the last subsection to give refined versions of Lemmas 4.3, 4.4 and Proposition 4.5.

### a) Homological considerations

This subsection introduces some homological considerations that play a central role in subsequent discussions. The subsection has six parts.

*Part 1*: This part sets the stage. Introduce by way of notation $Y_\Diamond$ to denote the open subset of Y given by the union of the following:

- *The $f^{-1}(1, 2)$ part of* $M_\delta$
- $\cup_{\mathfrak{p} \in \Lambda} \mathcal{H}^+_\mathfrak{p}$
- *Each* $\mathfrak{p} \in \Lambda$ *version of* $\hat{\gamma}^+_\mathfrak{p} \cup \hat{\gamma}^-_\mathfrak{p}$, *this the* $\{u = 0, 1 - 3\cos^2\theta = 0\}$ *locus in* $\mathcal{H}_\mathfrak{p}$.

(5.1)

An ech-HF subvariety is a subset of $\mathbb{R} \times Y_\Diamond$. It intersects a set described by the third bullet if and only if it contains some $\mathfrak{p} \in \Lambda$ version of the $\mathbb{R}$-invariant cylinder $\mathbb{R} \times \hat{\gamma}^+_\mathfrak{p}$ or the corresponding cylinder $\mathbb{R} \times \hat{\gamma}^-_\mathfrak{p}$. Of interest here is the homology of $Y_\Diamond$.

Let $\Theta$ denote an element in $\mathcal{Z}_{ech,M}$. The set $\Theta$ defines a class in $H_1(Y_\Diamond; \mathbb{Z})$ as follows: Orient each of its constituent closed integral curves using $v$ so as to view each as a closed 1-cycle. If $\gamma$ denotes the given curve, then $[\gamma]$ is used to denote this 1-cycle. The class defined by $\Theta$ is that of $\sum_{(\gamma,1) \in \Theta} [\gamma]$. This class is denoted by $[\Theta]$ in what follows.

Let C denote an ech-HF subvariety. The corresponding classes $[\Theta_{C_-}]$ and $[\Theta_{C_+}]$ must define the same class in $H_1(Y_\Diamond; \mathbb{Z})$ because the image of C in $Y_\Diamond$ via the projection defines a 2-chain whose boundary is the difference between these two classes. This observation places a constraint on the possible ech-HF subvarieties. Meanwhile, any two ech-HF with the same limit sets $\Theta_{C_-}$ and $\Theta_{C_+}$ define by their difference a closed 2-cycle in $Y_\Diamond$ and thus a class in $H_2(Y_\Diamond; \mathbb{Z})$.

*Part 2*: The homology of $Y_\diamond$ can be determined using the Mayer-Vietoris sequence for the decomposition in (5.1). The relevant portions of this sequence read as follows:

$$0 \to \mathbb{Z} \to H_2(Y_\diamond) \to \oplus_{\mathfrak{p} \in \Lambda}(\mathbb{Z} \oplus \mathbb{Z}) \xrightarrow{\hat{O}} H_1(\Sigma) \oplus (\oplus_{\mathfrak{p} \in \Lambda} \mathbb{Z}) \xrightarrow{P} H_1(Y_\diamond) \xrightarrow{Q} \oplus_{\mathfrak{p} \in \Lambda} \mathbb{Z} \to 0.$$
(5.2)

Here and in what follows, $\mathbb{Z}$ coefficients are implicit. By way of explanation, note that $H_2(\mathcal{H}^+_{\mathfrak{p}} \cup \hat{\gamma}^+_{\mathfrak{p}} \cup \hat{\gamma}^-_{\mathfrak{p}})$ and $H_2(Y_\diamond \cap (\mathcal{H}^+_{\mathfrak{p}*} \cup \hat{\gamma}^+_{\mathfrak{p}} \cup \hat{\gamma}^-_{\mathfrak{p}}))$ are zero. Meanwhile, $f^{-1}(1, 2) \cap M_\delta$ deformation retracts onto $\Sigma = f^{-1}(\frac{3}{2})$ and so $H_2(f^{-1}(1, 2) \cap M_\delta) = \mathbb{Z}$. This $\mathbb{Z}$ is the left most term in (5.2). Meanwhile, each $\mathfrak{p} \in \Lambda$ version of $H_1(Y_\diamond \cap (\mathcal{H}^+_{\mathfrak{p}} \cup \hat{\gamma}^+_{\mathfrak{p}} \cup \hat{\gamma}^-_{\mathfrak{p}}))$ is isomorphic to $\mathbb{Z} \oplus \mathbb{Z}$. The direct sum of the latter is the $\oplus_{\mathfrak{p} \in \Lambda}(\mathbb{Z} \oplus \mathbb{Z})$ term preceding the map $\hat{O}$. The afore-mentioned deformation retract identifies $H_1(\Sigma)$ with $H_1(f^{-1}(1, 2) \cap M_\delta)$ which is why the latter appears in the summand preceding the map P. The $\oplus_{\mathfrak{p} \in \Lambda} \mathbb{Z}$ summand in this same term is $H_1(\mathcal{H}^+_{\mathfrak{p}} \cup \hat{\gamma}^+_{\mathfrak{p}} \cup \hat{\gamma}^-_{\mathfrak{p}})$. Note for future reference that the generator of $H_1(\mathcal{H}^+_{\mathfrak{p}} \cup \hat{\gamma}^+_{\mathfrak{p}} \cup \hat{\gamma}^-_{\mathfrak{p}})$ and of a given component of $H_1(Y_\diamond \cap \mathcal{H}^+_{\mathfrak{p}})$ can be taken to be the $\theta = \frac{\pi}{2}$ circle in a constant u slice of $\mathcal{H}^+_{\mathfrak{p}}$ with the orientation defined by $\partial_\phi$.

To say something about the homomophisms in (5.2), consider first that on the right from $H_1(Y_\diamond)$ that is labeled by Q. Any oriented loop in $H_1(Y_\diamond)$ that crosses a given $\mathfrak{p} \in \Lambda$ version of $\mathcal{H}^+_{\mathfrak{p}}$ so as to have intersection number 1 with the u = 0 slice defines an element in $H_1(Y_\diamond)$ that maps via Q to the generator of $\mathfrak{p}$'s summand of the right most term in (5.2). By way of an example, suppose that $\Theta \in \mathcal{Z}_{\text{ech},M}$. Then the image of [$\Theta$] via Q is the class (1, 1, ...., 1) in this $\oplus_{\mathfrak{p} \in \Lambda} \mathbb{Z}$ term on the far right in (5.2).

Granted what was just said, the difference between the classes in $H_1(Y_\diamond)$ that come from any two elements in $\mathcal{Z}_{\text{ech},M}$ define a class in $H_1(\Sigma) \oplus (\oplus_{\mathfrak{p} \in \Lambda} \mathbb{Z})$ via the homomorphism P. These classes are determined by the homomorphism

$$\hat{O} : \oplus_{\mathfrak{p} \in \Lambda} (\mathbb{Z} \oplus \mathbb{Z}) \to H_1(\Sigma) \oplus (\oplus_{\mathfrak{p} \in \Lambda} \mathbb{Z}) \,.$$
(5.3)

To say something about $\hat{O}$, introduce for each $\mathfrak{p} \in \Lambda$, the classes $z_{\mathfrak{p}1}$ and $z_{\mathfrak{p}2}$ in $H_1(\Sigma)$ that are generated by the intersection of $\Sigma$ with the respective ascending disk from the index 1 component of $\mathfrak{p}$ and descending disk from the index 2 component. Then $\hat{O}$ sends an element $(n_{\mathfrak{p}1}, n_{\mathfrak{p}2})_{\mathfrak{p} \in \Lambda} \in \oplus_{\mathfrak{p} \in \Lambda}(\mathbb{Z} \oplus \mathbb{Z})$ to the element whose respective $H_1(\Sigma)$ and $\oplus_{\mathfrak{p} \in \Lambda} \mathbb{Z}$ summands are

$$\sum_{\mathfrak{p} \in \Lambda} (n_{\mathfrak{p}1} z_{\mathfrak{p}1} + n_{\mathfrak{p}2} z_{\mathfrak{p}2}) \quad \text{and} \quad \{(n_{\mathfrak{p}1} - n_{\mathfrak{p}2})\}_{\mathfrak{p} \in \Lambda} \,.$$
(5.4)

What follows are some additional remarks. To start, let $\hat{O}_\Sigma$ denote the homomorphism from $\oplus_{\mathfrak{p} \in \Lambda} (\mathbb{Z} \oplus \mathbb{Z})$ to $H_1(\Sigma)$ that is given by the left most term in (5.4). The cokernel of $\hat{O}_\Sigma$ is isomorphic to $H_1(M)$ and its kernel is isomorphic to $H_2(M)$. This can be seen by using the homology Mayer-Vietoris sequence for M that comes from splitting M along $\Sigma$'s incarnation as the $f = \frac{3}{2}$ level set. It follows from this that the image of $\hat{O}$ is isomorphic to $\ker(H_1(\Sigma) \to H_1(M)) \oplus H_2(M)$. Note however that this isomorphism is not canonical. It follows from this last observation that the image of P is isomorphic to $H_1(M) \oplus ((\oplus_{\mathfrak{p} \in \Lambda} \mathbb{Z})/H_2(M))$.

*Part 3*: Suppose that C is an ech-HF subvariety. Because the class $[\Theta_{C+}] - [\Theta_{C-}]$ in $H_1(Y_\diamond)$ is zero, there are classes in $H_1(\Sigma) \oplus (\oplus_{\mathfrak{p} \in \Lambda} \mathbb{Z})$ that map to the class $[\Theta_{C+}] - [\Theta_{C-}]$ via the homomorphism P in (5.2). As explained here, C determines a canonical choice for such a class with $H_1(\Sigma)$ summand equal to 0. The entry in any given $\mathfrak{p} \in \Lambda$ summand is denoted here by $\mathfrak{m}^C_\mathfrak{p}$; the latter is either 0, 1 or -1. The explanation follows.

Fix $\mathfrak{p} \in \Lambda$ and introduce $\gamma_{\mathfrak{p}+}$ to denote the segment in $\mathcal{H}^+_\mathfrak{p}$ of the integral curve from $\Theta_{C+}$ that crosses $\mathcal{H}_\mathfrak{p}$. Likewise define $\gamma_{\mathfrak{p}-}$ using $\Theta_{C-}$. Let $p_1$ and $p_2$ denote here the respective index 1 and index 2 critical points from $\mathfrak{p}$. Given $r \in (\delta_*, 7\delta_*)$, let $S_{1,r}$ and $S_{2,r}$ denote the respective boundaries of the radius r coordinate balls centered on $p_1$ and $p_2$. Choose r so that C has transversal intersection with both $\mathbb{R} \times S_{1,r}$ and $\mathbb{R} \times S_{2,r}$. If r is sufficiently generic, then both $C \cap (\mathbb{R} \times S_{1,r})$ and $C \cap (\mathbb{R} \times S_{2,r})$ will be disjoint from both $\mathbb{R} \times \gamma_{\mathfrak{p}+}$ and $\mathbb{R} \times \gamma_{\mathfrak{p}-}$. Assume that this is the case also.

Note that the intersection of C with $\mathbb{R} \times S_{1,r}$ has precisely one unbounded component; this is a smooth arc with $s$ unbounded from above and below. The function $s$ restricts to this arc as an affine parameter where $|s| \gg 1$. The constant $s$ points on the arc converge in $S_{1,r}$ as $s \to \infty$ to the point where the integral curve segment $\gamma_{\mathfrak{p}+}$ intersects $S_{1,r}$; and they converge as $s \to -\infty$ to the point where $\gamma_{\mathfrak{p}-}$ intersects $S_{1,r}$. The intersection of C with $\mathbb{R} \times S_{2,r}$ has a similar arc component.

With the preceding as background, introduce $c_1$ and $c_2$ to denote the respective images in $\mathcal{H}^+_\mathfrak{p}$ of $C \cap (\mathbb{R} \times S_{1,r})$ and $C \cap (\mathbb{R} \times S_{2,r})$. Neither $c_1$ nor $c_2$ intersects $\gamma_{\mathfrak{p}+}$, and neither intersects $\gamma_{\mathfrak{p}-}$. As a consequence, the closure of $c_1$ adds two points to $c_1$; these are the points where $\gamma_{\mathfrak{p}+}$ and $\gamma_{\mathfrak{p}-}$ intersect the radius r coordinate ball centered on $p_1$. Likewise, the closure of $c_2$ adds the points where $\gamma_{\mathfrak{p}+}$ and $\gamma_{\mathfrak{p}-}$ intersect the radius r coordinate ball centered on $p_2$.

Define a cycle in $\mathcal{H}^+_\mathfrak{p}$ as follows: Start where $\gamma_{\mathfrak{p}+}$ intersects $S_{2,r}$ and proceed along $\gamma_{\mathfrak{p}+}$ in $\mathcal{H}^+_\mathfrak{p}$ until it intersects the boundary of $S_{1,r}$. Then proceed along $c_1$ until it ends where $\gamma_{\mathfrak{p}-}$ intersects $S_{1,r}$. Travel in the reverse direction along $\gamma_{\mathfrak{p}-}$ until it intersects $S_{2,r}$. Close up the cycle by proceeding along $c_2$ to its endpoint where $\gamma_{\mathfrak{p}+}$ intersects $S_{2,r}$. This closed cycle is denoted by $\mathfrak{v}^C_p$. The cycle $\mathfrak{v}^C_p$ defines a class in $H_1(\mathcal{H}^+_p \cup \hat{\gamma}^+_p \cup \hat{\gamma}^-_p)$, and thus an integer. This is the integer $\mathfrak{m}^C_p$. The latter does not depend on the precise choice for r.

Given the definition above of $\mathfrak{m}^C_p$, what is said in Proposition 4.5 has the following consequence:

- $\mathfrak{m}^C_p = 0$ *if either none or both of $\{\hat{\gamma}^+_p, \hat{\gamma}^-_p\}$ are associated to ends of C in $\mathcal{H}^+_p$.*
- $\mathfrak{m}^C_p = -1$ *if $\hat{\gamma}^+_p$ is associated to an end of C in $\mathcal{H}^+_p$.*
- $\mathfrak{m}^C_p = 1$ *if $\hat{\gamma}^-_p$ is associated to an end of C in $\mathcal{H}^+_p$.*

(5.5)

To see that the class $(0, (\mathfrak{m}^C_p)_{p \in \Lambda}) \in H_1(\Sigma) \oplus (\oplus_{p \in \Lambda} \mathbb{Z})$ maps via P to $[\Theta_{C+}] - [\Theta_{C-}]$, define a class $M^C_\Sigma \subset H_1(\Sigma)$ as follows: Use $[\gamma^C_p]$ to denote the $\mathfrak{p} \in \Lambda$ version of the cycle that was described in the preceding paragraph. Then $[\Theta_{C+}] - [\Theta_{C-}] - \sum_{p \in \Lambda} [\gamma^C_p]$ is a closed cycle in the $f^{-1}(1, 2)$ portion of $M_\delta$ and so defines a class in $H_1(\Sigma)$. The latter is $M^C_\Sigma$, and it follows from the Mayer-Vietoris definition of P that $(M^C_\Sigma, (\mathfrak{m}^C_p)_{p \in \Lambda})$ is sent by P to the class $[\Theta_{C+}] - [\Theta_{C-}]$. However, the chain $[\Theta_{C+}] - [\Theta_{C-}] - \sum_{p \in \Lambda} [\gamma^C_p]$ in $f^{-1}(1, 2) \cap M_r$ is the boundary of $C \cap (f^{-1}(1,2) \cap M_r)$ and so represents the zero class in $H_1(\Sigma)$.

It is a straightforward task to verify that the assignment $C \to (\mathfrak{m}^C_p)_{p \in \Lambda}$ is constant as C varies in its moduli space component.

*Part 4*: Let C denote an ech-HF subvariety. This part of the subsection associates to C a set of $1 + b_1$ non-negative integer that characterizes in part of C's moduli space component. To define these integers, reintroduce the set $¥ \subset \Sigma - (T_+ \cup T_-)$ from Parts 6 and 7 of Section 1c. By way of a reminder, $¥$ contains the fiducial point $z_0$ and $b_1$ additional points. As noted in Section 2a, each point is the intersection point between $\Sigma$ and a closed integral curve of $v$. The curve through any given $z \in ¥$ is denoted by $\gamma^{(z)}$. Use $\mathfrak{n}^{(z)}_C$ to denote the intersection number between the subvariety C and the J-holomorphic cylinder $\mathbb{R} \times \gamma^{(z)}$. This is a non-negative integer, and the collection $\{\mathfrak{n}^{(z)}_C\}_{z \in ¥}$ is the desired set of integer. The $z_0$ version plays a distinguished role at times, and is denoted by $\mathfrak{n}_C$.

*Part 5*: Suppose as before that C is an ech-HF subvariety. Let $\hat{\mathfrak{v}}_+$ and $\hat{\mathfrak{v}}_-$ denote the respective HF-cycles that are used to define $\Theta_{C+}$ and $\Theta_{C-}$. View these as collections of G oriented arcs in the $f \in [1, 2]$ part of M that run from the index 1 critical points of $f$ to the index 2 critical points of $f$. Let $[\hat{\mathfrak{v}}_+]$ and $[\hat{\mathfrak{v}}_-]$ denote the corresponding 1-cycles. The next paragraph defines from the subvariety C a 2-chain in M with boundary the closed 1-cycle $[\hat{\mathfrak{v}}_+] - [\hat{\mathfrak{v}}_-]$. This 2-chain defines a canonical element in the $\mathbb{Z}$-module $H_2(M; [\hat{\mathfrak{v}}_+] - [\hat{\mathfrak{v}}_-])$. Here and in what follows, a 'chain' in M or Y of a given dimension indicates that the given subset can be viewed as a linear combination of symplexes of the indicated dimension, and as such is viewed as a chain in the appropriate singular homology chain complex.

To define the 2-chain, fix $\mathfrak{p} \in \Lambda$ and reintroduce the spheres $S_{1,r}$ and $S_{2,r}$ from Part 3. Remove the $|u| \leq R + \ln r$ portion of C's intersection with $\mathbb{R} \times \mathcal{H}^+_\mathfrak{p}$ so as to obtain a subvariety with boundary on $\mathbb{R} \times S_{1,r}$ and $\mathbb{R} \times S_{2,r}$. Now view $S_{1,r}$ as a sphere in M centered on the index 1 critical point of $f$. Cone the boundary curve $c_1 \subset S_{1,r}$ to the index 1 critical point in question by a cone that lies where $f \geq 1$ in the radius r coordinate ball centered on this index 1 critical point. Do the same for the boundary curve $c_2 \subset S_{2,r}$. Do this for all $\mathfrak{p} \in \Lambda$. The resulting chain with C's image in the complement of the radius r coordinate balls about the index 1 and index 2 critical points of $f$ defines a 2-chain whose boundary is a closed 1-chain which is homologous in the radius $c_0 \delta$ tubular neighborhood of $\cup_{\upsilon \in (\Theta_{C+} \cup \Theta_{C-})} \upsilon$ to the cycle $[\hat{\upsilon}_+] - [\hat{\upsilon}_-]$. Choose a 2-chain that gives such a homology and add it to the chain that was just defined by the coning of C's image in $M_\delta$. Use $[C]_M$ to denote the image of this chain in $H_2(M; [\hat{\upsilon}_+] - [\hat{\upsilon}_-])$.

The intersection number of this chain with the closed integral curve of $v$ through the point $z_0$ in $\Sigma$ is the integer $\mathfrak{n}_C$ that is defined in Part 4.

Suppose that C and C´ are two ech-HF subvarieties such that $\Theta_{C+}$ and $\Theta_{C'+}$ are defined using the same HF-cycle, and such that $\Theta_{C-}$ and $\Theta_{C'-}$ are also defined using the same HF-cycles. The difference $[C]_M - [C']_{M'}$ is a closed 2-cycle in M.

**b) Integral bounds**

Let C denote an ech-HF subvariety. This subsection supplies bounds on the integral over C of $w$, and the integral of $ds \wedge \hat{a}$ over compact subsets of C. These are stated in the next proposition. The proposition refers to an integer that can be associated to any given ech-HF moduli space component. This integer is called the *ech index*. Definition 2.14 in [Hu2] supplies a formula for the index. This index is denoted here by $I_{ech}(C)$. Note in this regard that the image of C in Y via the projection from $\mathbb{R} \times Y$ defines a class in the relative homology group $H_2(Y; [\Theta_{C+}] - [\Theta_{C-}])$.

**Proposition 5.1**: *There exists $\kappa > 1$ that is independent of the choice for J and the data R, $x_0$, and $\delta$ that is used to construct Y, and which has the following significance: Let C denote an ech-HF subvariety and let $I_{ech}(C)$ denote the ech index of its modui space component. Then*

- $0 \leq \int_C w \leq \kappa(1 + |I_{ech}(C)|)$.

- $(1 - 2\delta^2) G < \int_{C \cap ([s,s+1] \times M_\delta)} ds \wedge \hat{a} < G$.

*There exists $\kappa_J \geq 1$ that depends on J and the data R, $x_0$, $\delta$ and which has the following significance: Let C be as described above. Then*

- $0 \leq \int_{C \cap ([s,s+1] \times Y)} ds \wedge \hat{a} \leq \kappa_J (1 + |I_{ech}(C)|) + \Sigma_{(\gamma,1) \in \Theta_{C-}} \int_\gamma a$ .

With regards to the proof, remark that no generality is lost by considering only that case where C has no components from the set $\{\mathbb{R} \times \hat{\gamma}_p^+, \mathbb{R} \times \hat{\gamma}_p^-\}_{p \in \Lambda}$. This assumption is made implicitly in all that follows, lemmas included.

Part 1 of what follows supplies a key input to the proof, this Lemma 5.2. Part 2 of what follows contains the proof Proposition 5.1. Parts 3 and 4 of the subsection contain the proof of Lemma 5.2. Note that the lower bounds by zero in the first and third bullets follow from the fact that J is compatible with $ds \wedge \hat{a} + w$ and so no more will be said about them.

*Part 1*: The upcoming Lemma 5.2 supplies a key input to the proof of Proposition 5.1. This lemma refers to the integer $\mathfrak{n}_C$ that is defined in Part 4 of Section 5a. The lemma also refers to an element in the $\oplus_{p \in \Lambda}(\mathbb{Z} \oplus \mathbb{Z})$ summand of (5.2) that is associated to C's moduli space component. What follows directly gives the definition.

To start things off, reintroduce the respective 1-chains $[\hat{\upsilon}_+]$ and $[\hat{\upsilon}_-]$ that are defined by the HF-cycles which are used in Proposition 2.8 to label $\Theta_{C+}$ and $\Theta_{C-}$. As noted in Part 5 of Section 5a, the existence of C requires that $[\hat{\upsilon}_+] - [\hat{\upsilon}_-]$ is zero in $H_1(M; \mathbb{Z})$. The vanishing of this has the following implication: There exists a 2-chain in the $f^{-1}([1, 2])$ part of M whose boundary is $[\hat{\upsilon}_+] - [\hat{\upsilon}_-]$ that depends only on the Heegaard Floer data. Moreover, 2-chains with certain desirable properties can be constructed using the fact that the components of the 1-cycle $[\hat{\upsilon}_+] - [\hat{\upsilon}_-]$ are smoothly embedded arcs that intersect in pairs at their end points, these being index 1 and index 2 critical points of $f$. In particular, there exists a smooth map from an oriented surface with boundary into M whose image represents the 2-cycle with the properties listed in (5.6). This equation uses $\mathfrak{s}$ to denote the surface and $\mathfrak{f}: \mathfrak{s} \to M$ to denote the map.

- *The boundary of $\mathfrak{s}$ has 2G marked points. These are mapped in a 1-1 fashion to the index 1 and index 2 critical points of $f$.*
- *The complement in $\mathfrak{s}$ of the 2G marked boundary points is mapped to the $f^{-1}(1,2)$ part of M.*
- *The complement in the boundary of these 2G marked point is identified by $\mathfrak{f}$ with the interiors of the segments that comprise $\hat{\upsilon}_+ \cup \hat{\upsilon}_-$.*

- *The map f is transversal to the boundary of each radius $r \in [\delta_*, 8\delta_*]$ coordinate ball about the index 1 and index 2 critical points of f; and the image of f intersects each such boundary sphere as an immersed arc.*
- *Let $\iota \subset \partial S$ denote a component of the complement of the 2G marked boundary points. The map f restricts to a neighborhood in S of $\iota$ as an embedding which has the vector field $\mathfrak{v}$ tangent its image.*

(5.6)

Let [S] denote the 2-chain in M given by f(S)

Fix $\mathfrak{p} \in \Lambda$ and write $\mathfrak{p}$ as $(p_1, p_2)$ as in Part 3 of Section 5a. The Mayer-Vietoris construction of (5.2) identifies the left hand $\mathbb{Z}$ in the $\mathfrak{p}$-labeled summand of $\oplus_{\mathfrak{p} \in \Lambda}(\mathbb{Z} \oplus \mathbb{Z})$ with $H_1(T_{p_1+})$. Likewise, it identifies the right most $\mathbb{Z}$ with $H_1(T_{p_2-})$. The specification of an integer in each index 1 critical point version of $H_1(T_{p+})$ and in each index 2 critical point version of $H_1(T_{p-})$ gives an element in the $\oplus_{\mathfrak{p} \in \Lambda}(\mathbb{Z} \oplus \mathbb{Z})$ summand. The 2-chain [S] is used to define these integers.

Reintroduce the notation from Section 5a, and in particular the chain $[C]_M$, the constant $r \in (\delta_*, 7\delta_*)$ and the spheres $S_{1,r}$ and $S_{2,r}$. The 2-chain $[C]_M - [S]$ is a closed 2-cycle whose intersection with the sphere $S_{1,r}$ lies in its $1 - 3\cos^2\theta$ part. As such, it defines a closed 1-cycle in $T_{p_1+}$ and thus a class in $H_1(T_{p_1+})$. Use $n^C_{p_1}$ to denote the corresponding integer. The intersection of $[C]_M - [S]$ with $S_{2,r}$ defines an analogous integer which is denoted by $n^C_{p_2}$. The desired element $(n^C_{p_1}, n^C_{p_2})_{\mathfrak{p} \in \Lambda} \in \oplus_{\mathfrak{p} \in \Lambda}(\mathbb{Z} \oplus \mathbb{Z})$ is the desired element.

**Lemma 5.2**: *Fix $I \geq 1$ and there exists $\kappa > 1$ that depends only on the Heegaard Floer data and which has the following significance: Suppose that C is an ech-HF subvariety with $I_{ech}(C) \leq I$. Then $\sum_{z \in ¥} \mathfrak{n}^{(z)}_C + \sum_{\mathfrak{p} \in \Lambda}(|n^C_{p_1}| + |n^C_{p_2}|) \leq \kappa$.*

Lemma 5.2 is proved in Parts 3 and 4 of this subsection. Assume this lemma for now.

*Part 2*: This part of the subsection contains the

***Proof of Proposition 5.1***: The proof has four steps. Steps 1-3 prove the assertion made by the first bullet of Proposition 5.1. Step 4 proves the assertions that are made by the second and third bullets. Appeals to Lemma 5.2 enter only in Steps 1 and 3.

Step 1: As noted in (1.19), the 2-form $w$ can be written as d$a$ on the complement of the intersection of C with the part of $M_\delta$ that corresponds to the $(1, 2) \times (\cup_{z \in ¥} D_z)$ part of $(1, 2) \times \Sigma$. Granted this, and given that C has strictly positive intersection number with each $z \in ¥$ version of $\mathbb{R} \times \gamma^{(z)}$, an integration by parts with (1.16) writes

$$\int_C w = \sum_{\gamma \in \Theta_{C+}} \int_\gamma a - \sum_{\gamma \in \Theta_{C-}} \int_\gamma a + 2\sum_{z \in \Lambda} \mathfrak{n}^{(z)}{}_C R_z \ .$$

(5.7)

Each $z \in \mathcal{Y}$ version of $\mathfrak{n}^{(z)}{}_C$ is non-negative, and Lemma 5.2 says that each is bounded by a constant that depends only on $I_{ech}(C)$ and the Heegaard Floer data. Thus, the right most term in (5.7) is also bounded by a constant that depends only on the Heegaard Floer data.

The integrals in (5.7) over the various closed orbits of $v$ can be broken into their contributions from the parts in the 1-handles and the parts in $M_\delta$. It is a consequence of what is said in Corollary 2.6 that the $M_\delta$ contribution to the integrals in (5.7) is no greater than a constant the depends solely on the Heegaard Floer data.

Fix $\mathfrak{p} \in \Lambda$. The contribution to the integrals on the right side of (5.7) from $\mathcal{H}_\mathfrak{p}$ can be written as a sum of various terms. There is, first of all, a contribution from each integral curve in $\Theta_{C+}$ or $\Theta_{C-}$ from the set $\{\hat{\gamma}^+_\mathfrak{p}, \hat{\gamma}^-_\mathfrak{p}\}$. A look at PROPERTY 3 in Section 1e finds that each such curve in $\Theta_{C+}$ contributes $\frac{4\sqrt{2}}{3}\pi(x_0 + 4e^{-2R})$ to the right hand side of (5.7), and each in $\Theta_{C-}$ contributes minus this amount to the right hand side of (5.7).

Let $\gamma_{\mathfrak{p}+}$ and $\gamma_{\mathfrak{p}-}$ denote the segments of the integral curves in $\Theta_{C-}$ and $\Theta_{C+}$ that cross $\mathcal{H}^+_\mathfrak{p}$. The remaining contribution to the right hand side of (5.7) from $\mathcal{H}_\mathfrak{p}$ is given by

$$\int_{\gamma_{\mathfrak{p}+}} a - \int_{\gamma_{\mathfrak{p}-}} a \ .$$

(5.8)

It follows from (1.28) that the integral $a$ over $\gamma_{\mathfrak{p}+}$ is equal to

$$\int_{\gamma_{\mathfrak{p}+}} x(1 - 3\cos^2\theta)du - \sqrt{6}\hat{h}_+\Delta_+ \ ,$$

(5.9)

where $\hat{h}_+$ here denotes the constant value of the function $\hat{h}$ on $\gamma_{\mathfrak{p}+}$, and where $\Delta_+$ denotes the change in the angle $\phi$ along $\gamma_{\mathfrak{p}+}$. The latter is given by integral in the fourth bullet of Lemma 2.2. The integral of $a$ over $\gamma_{\mathfrak{p}-}$ is given by the analogous expression in terms of an integral and the corresponding $\hat{h}_-$ and $\Delta_-$.

<u>Step 2</u>: This step rewrites (5.8). This is done by first writing (5.8) as the difference between the $\gamma_{\mathfrak{p}+}$ and $\gamma_{\mathfrak{p}-}$ versions of (5.9). To start the task of rewriting the latter expression, let $\gamma$ denote for the moment a segment of an integral curve of $v$ in $\mathcal{H}^+_\mathfrak{p}$ with $|\hat{h}| < (x_0 + 2e^{-2R})\frac{2}{3\sqrt{3}}$. Recall from Lemma 2.2 that this last bound is the necessary and sufficient condition for $\gamma$ to cross $\mathcal{H}^+_\mathfrak{p}$ from its $u = -R - \ln(7\delta_*)$ end to its end where $u = R + \ln(7\delta_*)$. View $\theta$ along any segment of this sort as a function of $u$ that depends implicitly on the constant value of $\hat{h}$ along the segment. The assignment to $\gamma$ of the integral $\int_\gamma x(1 - 3\cos^2\theta)du$ defines a function of $\hbar$. Let $\hbar \to \mathfrak{x}(\hbar)$ denote this function. Use the chain rule with (4.1) to see that

$$\frac{d}{d\hbar}\mathfrak{x} = \sqrt{6}\Delta,$$

(5.10)

where $\Delta = \Delta(\hbar)$ is given by the integral in the fourth bullet of Lemma 2.2.

Use (5.10) to write the difference between the $\gamma_{\mathfrak{p}+}$ and $\gamma_{\mathfrak{p}-}$ versions of (5.9) as

$$\sqrt{6}\left(\int_{\hbar_-}^{\hbar_+} \Delta\, d\hbar - \hbar_+\Delta_+ + \hbar_-\Delta_-\right).$$

(5.11)

With (5.11) in hand, use the mean value theorem to write the integral that appears in (5.11) as $\Delta(\hbar_0)(\hbar_+ - \hbar_-)$ where $\hbar_0$ is a number between $\hbar_-$ and $\hbar_+$. This done, the (5.11) can be written as

$$\sqrt{6}(\Delta(\hbar_0) - \Delta_+)\hbar_+ - \sqrt{6}(\Delta(\hbar_0) - \Delta_-)\hbar_-.$$

(5.12)

This last expression is the desired rewrite of (5.8).

Step 3: Use the expression in the fourth bullet of Lemma 2.2 to see that the derivative of the function $\hbar \to \Delta(\hbar)$ is given by $-\sqrt{6}$ times the integral that appears on the right hand side of (2.13). This implies in particular that the function $\hbar \to \Delta(\hbar)$ is a monotonic function of $\hbar$ and so the number $\Delta(\hbar_0)$ that appears in (5.12) is between $\Delta_-$ and $\Delta_+$. Thus, both $|\Delta(\hbar_0) - \Delta_+|$ and $|\Delta(\hbar_0) - \Delta_-|$ are bounded $|\Delta_+ - \Delta_-|$. It follows as a consequence that the absolute value of the expression in (5.12) is no greater than $c_0 x_0 |\Delta_+ - \Delta_-|$.

Use Proposition 2.8 to label $\Theta_{C+}$ and $\Theta_{C-}$. This done, introduce $\mathfrak{k}_{\mathfrak{p}+}$ and $\mathfrak{k}_{\mathfrak{p}-}$ to denote the respective integers that are used in these labels for $\mathfrak{p}$'s factor of $\mathbb{Z} \times \mathfrak{o}$. Given that $|\mathfrak{m}_\mathfrak{p}| \leq 1$, the bound given by Lemma 5.2 implies that $|\mathfrak{k}_{\mathfrak{p}+} - \mathfrak{k}_{\mathfrak{p}-}|$ is bounded by a constant that depends only on $I_{\text{ech}}(C)$ and the Heegaard Floer data. This bound implies that $|\Delta_+ - \Delta_-|$ is also bounded by a constant that depends only on $I_{\text{ech}}(C)$ and the Heegaard Floer data. This fact with what is said in the preceding paragraph implies that the absolute value of what is written in (5.8) is bounded by a constant that depends solely on $I_{\text{ech}}(C)$ and the Heegaard Floer data. This fact with what is said in the first step prove the assertion made by the first bullet of Proposition 5.1.

Step 4: To prove the second bullet, note first that C intersects $\mathbb{R} \times M_\delta$ only where $f \in (1,2)$. This understood, write this part of $M_\delta$ as $(1,2) \times \Sigma$. The 1-form $\hat{a}$ appears in this guise as $dt$ where $t$ is the Euclidean coordinate on the $(1, 2)$ factor. The integral depicted in the second bullet is that of the pull-back of $ds \wedge dt$ via the projection map to $[s, s+1] \times (1,2)$ of C's intersection with the $[s, s+1] \times M_\delta$ part of $([s, s+1] \times (1,2) \times \Sigma)$. Each constant $(s, t)$ slice of $\mathbb{R} \times (1, 2) \times \Sigma$ is pseudoholomorphic, and so it follows from Lemma 4.2 that C has at most G intersection points (counting multiplicities) with each slice. This implies that the integral in question is at most G.

By way of comparison, the projection induced map to $[s, s+1] \times (1+\delta^2, 2-\delta^2)$ from C's intersection with the $[s, s+1] \times M_\delta$ part of $[s, s+1] \times (1+\delta^2, 2-\delta^2) \times \Sigma$ is a proper map of degree G and so the integral of $ds \wedge \hat{a}$ over $C \cap ([s, s+1] \times (1+\delta^2, 2-\delta^2) \times \Sigma)$ is exactly $(1-2\delta^2)G$. The latter integral is strictly less than than the integral of $ds \wedge \hat{a}$ over pull-back of $ds \wedge dt$ to $C \cap ([s, s+1] \times M_\delta)$ because $ds \wedge \hat{a}$ is non-negative on TC.

To prove the final assertion of Proposition 5.1, use the second bullet of (3.4) to bound

$$\int_{C \cap ([s,s+1] \times Y)} ds \wedge \hat{a} = r \int_{C \cap ([s,s+1] \times Y)} (ds \wedge a + w) ,$$

(5.13)

To exploit (5.13), write

$$\int_{C \cap ([s,s+1] \times Y)} ds \wedge a = \int_{[s,s+1]} \left( \int_{C \cap (\{s'\} \times Y)} a \right) ds'$$

(5.14)

and integrate by parts to write

$$\int_{C \cap (\{s'\} \times Y)} a = \int_{C \cap ((-\infty,s'] \times Y)} da + \Sigma_{(\gamma,1) \in \Theta_{C-}} \int_\gamma a .$$

(5.15)

By construction, $da = Hw$ where H is a smooth function on Y. This understood, what is written on the right hand side of (5.13) is no greater than

$$r_* \int_C w + \Sigma_{(\gamma,1) \in \Theta_{C-}} \int_\gamma a ,$$

(5.16)

where $r_* \geq 1$. An appeal to the first bullet in Proposition 5.1 finishes the proof.

*Part 3*: This part of the subsection proves Lemma 5.2 modulo an auxilliary result, this being Lemma 5.3. Lemma 5.3 is proved in Part 4 of this subsection. Lemma 5.3 is one of two key inputs to the proof. The argument for Lemma 4.14 of [OS] is the second key input.

*Proof of Lemma 5.2*: The proof has five steps.

Step 1: Let C denote a given ech-HF subvariety and let Z denote a 2-chain in $M_\delta \cup (\cup_{\mathfrak{p} \in \Lambda} \mathcal{H}_\mathfrak{p})$ that defines an element in $H_2(Y; [\Theta_{C_+}] - [\Theta_{C_-}])$. A chain of this sort has an associated ech index. Let $[C]_Y$ denote the 2-chain that is defined by C's image in Y via the projection from $\mathbb{R} \times Y$. The chain defined by $[C]_Y - Z$ is a closed 2-cycle in Y. Hutchings explains in Section 2 of [Hu2] why

$$I(C) - I(Z) = \langle c_1, [C]_Y - Z \rangle \tag{5.17}$$

where the notation is such that $\langle\,,\,\rangle \colon H^2(Y; \mathbb{Z}) \otimes H_2(Y; \mathbb{Z}) \to \mathbb{Z}$ denotes the canonical pairing and $c_1$ denotes the first Chern class of the chosen $\text{Spin}^\mathbb{C}$ structure on Y, this being the class depicted in (1.29).

**Lemma 5.3**: *There is a constant $\kappa$ that depends only on the Heegaard Floer data, and a 2-chain $Z \in H_2(Y; [\Theta_{C_+}] - [\Theta_{C_-}])$ with the following properties:*
- *Z intersects only the union of $f^{-1}(1, 2)$ part of $M_\delta$ and the $\mathfrak{p} \in \Lambda$ versions of $\mathcal{H}^+_\mathfrak{p}$.*
- *The part of Z that lies in $M_\delta$ is obtained from the chain [s] by adding a chain with support in the radius $\kappa\delta$ tubular neighborhood of the integral curves of $\mathfrak{v}$ from the HF-cycles $\hat{\mathfrak{v}}_+$ and $\hat{\mathfrak{v}}_-$.*
- *$|I_{ech}(Z)| \leq \kappa$.*

This lemma is proved in Part 4 of the subsection. The subsequent steps refer to a chosen 2-chain Z that is described by Lemma 5.3.

Step 2: It follows from (1.29) that the first Chern class $c_1$ can be represented by a closed form with support disjoint from $\cup_{\mathfrak{p} \in \Lambda} \mathcal{H}_\mathfrak{p}$. Observations from Parts 4 and 6 of Section 1b are used to construct such a form whose intersection with the $f^{-1}((1, 2))$ portion of $M_\delta$ is a 2-form that is obtained from a form on $\Sigma$ via the Lie transport defined by $\mathfrak{v}$. This form is denoted in what follows by $w_{\Sigma 0}$. It is non-negative on $\Sigma$ and it has the properties listed in the upcoming (5.18). This equation uses terminology from Parts 4 and 6 of Section 1b.

- *The form $w_{\Sigma 0}$ has support in the complement of $T_+ \cup T_-$.*
- *The integral of $w_{\Sigma 0}$ over any fundamental domain is strictly positive.*
- *The integral of of $w_{\Sigma 0}$ over $\Sigma$ is 2.*

- *The integral of $w_{\Sigma 0}$ any given periodic domain in $\Sigma$ is equal to the value of $c_{1M}$ on the corresponding homology class in $H_2(M;\mathbb{Z})$*

(5.18)

What is said in Lemma 5.3 about $\mathbb{Z}$ and what was just said about the representative form for the class $c_1$ have the following consequence: The pairing on the right hand side of (5.17) can be computed as described next. Use Lie transport by the pseudogradient vector field $\mathfrak{v}$ to map $[C]_M$ and $[\mathfrak{s}]$ into $\Sigma$. Their respective images define 2-chains in $\Sigma$; these are denoted in what follows by $C_\Sigma$ and $S_\Sigma$. The right hand side of (5.17) is the integral of the 2-form $w_{\Sigma 0}$ over $C_\Sigma - S_\Sigma$. Note in this regard that neither $C_\Sigma$ nor $S_\Sigma$ is a closed 2-chain. Even so, their respective boundaries lie in $T_+ \cup T_-$ which is disjoint from the support of $w_{\Sigma 0}$. The fact that their respective boundaries lie in $T_+ \cup T_-$ implies that both have a well defined intersection number with each point in ¥. The intersection number between $C_\Sigma$ and any given $z \in$ ¥ is the integer $\mathfrak{n}^{(z)}{}_C$; and the absolute value of that between $S_\Sigma$ and $z$ is bounded solely by the Heegaard Floer data.

The difference $C_\Sigma - S_\Sigma$ is not closed either because $[C]_M - [\mathfrak{s}]$ intersects the critical index 1 and index 2 critical points of $f$. More to the point, $C_\Sigma$ and $S_\Sigma$ can be modified in $T_+ \cup T_-$ so that the resulting chains, $C_{\Sigma 0}$ and $S_{\Sigma 0}$, are such that

- *The intersection number between $C_{\Sigma 0} - S_{\Sigma 0}$ and any given $z \in$ ¥ differs from $\mathfrak{n}^{(z)}{}_C$ by a constant that depends only on the Heegaard Floer data.*
- $\partial(C_{\Sigma 0} - S_{\Sigma 0}) = \sum_{\mathfrak{p} \in \Lambda}(\mathrm{n}^C{}_{\mathfrak{p}1} Z_{\mathfrak{p}1} + \mathrm{n}^C{}_{\mathfrak{p}2} Z_{\mathfrak{p}2})$.

(5.19)

The integral of $w_{\Sigma 0}$ over $C_{\Sigma 0} - S_{\Sigma 0}$ also computes the right hand side of (5.17) because the support of $w_{\Sigma 0}$ is disjoint from $T_+ \cup T_-$.

Lemma 5.2 asserts a particular relationship between $I_{ech}(C)$ on the one hand and the integers that comprise the sets $\{\mathfrak{n}^{(z)}{}_C\}_{z \in ¥}$ and $\{(\mathrm{n}^C{}_{\mathfrak{p}1}, \mathrm{n}^C{}_{\mathfrak{p}2})\}_{\mathfrak{p} \in \Lambda}$ on the other. The chain $C_{\Sigma 0} - S_{\Sigma 0}$ is the bridge that brings these two seemingly disparate notions together by virtue of two facts: The boundary of $C_{\Sigma 0} - S_{\Sigma 0}$ defines the two sets of integers via (5.19); and the integral of $w_{\Sigma 0}$ over this 2-cycle computes (5.17). This bridge leads to a proof of Lemma 5.2 using what are essentially the same arguments that are used to prove Lemma 4.14 in [OS]. The remaining steps supply these arguments.

<u>Step 3</u>: If z is a point of ¥, let $\gamma^{(z)}$ again denote the corresponding closed integral curve of $v$ through z. Use $[\gamma^{(z)}]$ to again denote the corresponding closed 1-cycle and $[\gamma^{(z)}]^{Pd}$ to denote the Poincaré dual in $H^2(Y;\mathbb{Z})$ of the resulting homology class. The span of the set $\{[\gamma^{(z)}]^{Pd}\}_{z \in ¥} \in H^2(Y;\mathbb{Z})$ is a $1+b_1$ dimensional subspace that is dual via the pairing $\langle\,,\,\rangle$ to the $H_2(\mathcal{H}_0;\mathbb{Z}) \oplus H_2(M;\mathbb{Z})$ summand in (1.27). The set ¥ also has its fiducial point, $z_0$; and the span of the set $\{[\gamma^{(z)}]^{Pd} - [\gamma^{(z_0)}]^{Pd}\}_{z \in \yen - z_0}$ is dual to the $H_2(M; \mathbb{Z})$ summand in (1.27).

It proves useful at this point to introduce some specific cycles to represent the homology in the summand $H_2(M; \mathbb{Z}) \oplus H_2(\mathcal{H}_0; \mathbb{Z})$. In particular, the classes in the $H_2(M; \mathbb{Z})$ summand are represented by periodic domains in $\Sigma$ as done in the proof of Lemma 1.1. Let $\mathcal{P}$ denote such a domain, and let $h(\mathcal{P}) \in H_2(M; \mathbb{Z})$ denote the corresponding homology class. The cycle $[\Sigma]$ of the Heegaard surface completes the desired set of cycles.

Each cycle from the set $\{[\gamma^{(z)}]^{Pd}\}_{z \in \yen}$ has pairing 1 with $[\Sigma]$. The elements from $\{[\gamma^{(z)}]^{Pd} - [\gamma^{(z_0)}]^{Pd}\}_{z \in \yen - z_0}$ have pairing 0 with $[\Sigma]$, and so the span of this set is dual to the span of the classes defined by the periodic domains. This understood, fix a basis, $\{\mathcal{P}_z\}_{z \in \yen - z_0}$, of periodic domains such that the corresponding set $\{h(\mathcal{P}_z)\}_{z \in \yen - z_0}$ is dual (over $\mathbb{Q}$) to the elements in the set $\{[\gamma^{(z)}]^{Pd} - [\gamma^{(z_0)}]^{Pd}\}_{z \in \yen - z_0}$.

Let $\mathcal{P}$ again denote a periodic domain. The boundary of $\mathcal{P}$ is a closed 1-cycle in $\Sigma$ that can be written as

$$\partial \mathcal{P} = \sum_{p \in \Lambda} (n^{\mathcal{P}}_{p_1} z_{p1} + n^{\mathcal{P}}_{p_2} z_{p2})$$

(5.20)

The resulting element $(n^{\mathcal{P}}_{p1}, n^{\mathcal{P}}_{p2})_{p \in \Lambda} \in \oplus_{p \in \Lambda}(\mathbb{Z} \oplus \mathbb{Z})$ is mapped by the homomorphism $\hat{O}$ in (5.3) to an element with $H_1(\Sigma)$ entry zero. In fact, the elements in $\oplus_{p \in \Lambda}(\mathbb{Z} \oplus \mathbb{Z})$ that come from the boundaries of periodic domains span the linear space of such elements.

<u>Step 4</u>: As (5.20) indicates, the class $(n^C_{p1}, n^C_{p2})_{p \in \Lambda}$ in $\oplus_{p \in \Lambda}(\mathbb{Z} \oplus \mathbb{Z})$ is mapped by (5.3)'s homomorphism $\hat{O}$ to an element in $H_1(\Sigma) \oplus (\oplus_{p \in \Lambda} \mathbb{Z})$ with $H_1(\Sigma)$ entry zero. What follows is a consequence of this and what said in Step 3: The 2-chain $C_{\Sigma 0} - S_{\Sigma 0}$ is homologous rel $C_+ \cup C_-$ to

$$(\mathfrak{n}_C + K^s_0)[\Sigma] + \sum_{z \in \yen - z_0} q^C_z \mathcal{P}_z ,$$

(5.21)

where $K^s_0$ is the intersection number between $S_{\Sigma 0}$ and $\gamma^{(z_0)}$ and where each $z \in \yen - z_0$ version of $q^C_z$ is an integer. Granted this, depiction, use (5.18) to see that the right hand side of (5.17) is equal to

$$2(\mathfrak{n}_C + K^s_0) + \sum_{z \in \yen - z_0} q^C_z (\sum_{1 \leq i \leq N} K^z_i A_i) ,$$

(5.22)

where the notation is as follows: First, $\mathcal{P}_z$ is written in terms of fundamental domains as $\mathcal{P}_z = \sum_{1 \leq i \leq N} \kappa^z{}_i \mathcal{D}_i$. Meanwhile, any given $i \in \{1, \ldots, N\}$ version of $A_i$ is the integral of $w_{\Sigma 0}$ over the domain $\mathcal{D}_i$.

To exploit this formula, note that both $C_{\Sigma 0}$ and $S_{\Sigma 0}$ are homologous rel $C_+ \cup C_-$ to integer weighted sums of periodic domains; for the purposes at hand, $C_{\Sigma 0}$ can be replaced by a 2-chain of the form $\mathfrak{n}_C \mathcal{D}_0 + \sum_{1 \leq i \leq N} \kappa^C{}_i \mathcal{D}_i$ and $S_{\Sigma 0}$ by $\sum_{0 \leq i \leq N} \kappa^S{}_i \mathcal{D}_i$. This done, then (5.17) and (5.22) imply that

$$I_{ech}(C) = 2\mathfrak{n}_C A_0 + \sum_{1 \leq i \leq n} \kappa^C{}_i A_i + \mathfrak{r},$$
(5.23)

where $|\mathfrak{r}| \leq c_0$ with $c_0$ determined solely by the Heegaard Floer data. This is so because $S$ and so $S_{\Sigma 0}$ is determined by the Heegaard Floer data, and because of what is said about $I_{ech}(Z)$ by Lemma 5.3.

Step 5: The 2-chain $C_{\Sigma 0}$ is an example of what Osváth and Szabó call a *positive* 2-chain, which is to say that the integer coefficients $\{\kappa^C{}_i\}_{1 \leq i \leq N}$ and $\mathfrak{n}_C$ are non-negative. Here is why: Any given $\kappa^C{}_i$ is the intersection number between C and a J-holomorphic subvariety defined in $\mathbb{R} \times M_\delta$. To elaborate, let $z_i \in \mathcal{D}_i$ denote the point described in Part 6 of Section 1b. Let $\gamma^{(i)} \in Y$ denote the segment of the integral curve of $v$ that contains $z_i$ and sits in $M_\delta$. Then $\kappa^C{}_i$ is the intersection number between C and $\mathbb{R} \times \gamma^{(i)}$. As noted previously, the integer $\mathfrak{n}_C$ is non-negative for the same reason.

The preceding observation, the fact that each $i \in \{0, \ldots, N\}$ version of $A_i$ is positive, and (5.23) have the following consequence: There exists $c_0 \geq 1$ that is determined solely by the Heegaard Floer data such that

$$|\mathfrak{n}_C| + \sum_{1 \leq i \leq n} \kappa^C{}_i \leq c_0 (1 + I_{ech}(C)).$$
(5.24)

What with (5.19), this bound on the elements in the set $\{\kappa^C{}_i\}_{1 \leq i \leq N}$ and on $\mathfrak{n}_C$ implies what is asserted by Lemma 5.2 because it gives an apriori bound on the coefficients that appear when writing the 2-chain $c_{\Sigma 0}$ as $\mathfrak{n}_C \mathcal{D}_0 + \sum_{1 \leq i \leq N} \kappa^C{}_i \mathcal{D}_i$.

*Part 4*: This part of the subsection contains the proof of Lemma 5.3. The proof starts momentarily. What follows directly is meant to give an idea of the strategy.

The argument for Lemma 5.3 starts by constructing certain elements in $\mathcal{Z}_{ech,M}$ from the HF-cycles $\hat{\upsilon}_+$ and $\hat{\upsilon}_-$ and the pair $(s, \mathfrak{f})$. These are denoted by $\Theta_+$ and $\Theta_-$. To say more, first write $\Theta_{C_+}$ using the notation from Proposition 2.8 as $(\hat{\upsilon}_+, (\mathfrak{k}_{p+}, O_{p+})_{p \in \Lambda})$, and likewise write $\Theta_{C_-}$ as $(\hat{\upsilon}_-, (\mathfrak{k}_{p-}, O_{p-})_{p \in \Lambda})$. The respective elements $\Theta_+$ and $\Theta_-$ are defined as follows:

$$\Theta_+ = (\hat{\upsilon}_+, (\mathfrak{k}_{\mathfrak{p}0}, O_{\mathfrak{p}+})_{\mathfrak{p}\in\Lambda}) \quad and \quad \Theta_- = (\hat{\upsilon}_-, (0, O_{\mathfrak{p}-})_{\mathfrak{p}\in\Lambda})$$

(5.25)

where each $\mathfrak{p} \in \Lambda$ version of $\mathfrak{k}_{\mathfrak{p}0}$ is chosen with the help of $\mathfrak{s}$ and $\mathfrak{f}$. The construction is such that $|\mathfrak{k}_{\mathfrak{p}0}| \leq c_0$ where $c_0$ depends only on the Heegaard Floer data.

The preceding bound on the integers $\{\mathfrak{k}_{\mathfrak{p}0}\}$ implies that the integral curve of $v$ from either $\Theta_+$ and $\Theta_-$ that crosses any given $\mathfrak{p} \in \Lambda$ version of $\mathcal{H}_\mathfrak{p}$ does so with a $(\delta, x_0, R)$-independent bound on the change in the $\phi$ coordinate. This fact is used to construct a 2-chain $Z_0 \subset Y$ with the following properties:

- $Z_0 \in H_2(Y; [\Theta_+] - [\Theta_-])$
- $|I_{ech}(Z_0)| \leq c_0$.
- $Z_0$ *lies in the union of the* $f^{-1}(1,2)$ *part of* $M_\delta$ *and the* $1 - 3\cos^2\theta \geq 0$ *part of each* $\mathfrak{p} \in \Lambda$ *version of* $\mathcal{H}_\mathfrak{p}$.
- $Z_0 \cap M_\delta$ *differs from* $[\mathfrak{s}] \cap M_\delta$ *only in the radius* $c_0\delta$ *tubular neighborhood of the elements from* $\hat{\upsilon}_+$ *and* $\hat{\upsilon}_-$.

(5.26)

Moreover, the constant $c_0$ that appears here depends only on the Heegaard Floer data.

With $Z_0$ in hand, a pair of 2-chains are constructed, one in $H_2(Y; [\Theta_{C+}] - [\Theta_+])$ and the other in $H_2(Y; [\Theta_-] - [\Theta_{C-}])$. These are denoted by $Z_+$ and $Z_-$. Their properties are described momentarily. In the meantime, note that the chain $Z = Z_+ + Z_0 + Z_-$ defines an element in $H_2(Y; [\Theta_{C+}] - [\Theta_{C-}])$. As noted in Section 2 of [Hu2],

$$I_{ech}(Z) = I_{ech}(Z_+) + I_{ech}(Z_0) + I_{ech}(Z_-) \ .$$

(5.27)

What follows says something about $Z_+$ and $Z_-$.

- *Both lie in the union of the* $f^{-1}(1,2)$ *part of* $M_\delta$ *and with* $\cup_{\mathfrak{p}\in\Lambda}\mathcal{H}_\mathfrak{p}$.
- *The parts of* $Z_+$ *and* $Z_-$ *in* $M_\delta$ *lie in the respective radius* $c_0$ *tubular neighborhoods of the curves from* $\hat{\upsilon}_+$ *and* $\hat{\upsilon}_-$.
- $|I_{ech}(Z_+)| + |I_{ech}(Z_-)| \leq c_0$.

(5.28)

As in (5.25), the constant $c_0$ depends only on the Heegaard Floer data. Where as the calculation of $I_{ech}(Z_0)$ is straightforward, those of $I_{ech}(Z_+)$ and $I_{ech}(Z_-)$ are quite delicate when $\sup_{\mathfrak{p}\in\Lambda}|\mathfrak{k}_{\mathfrak{p}+}|$ and $\sup_{\mathfrak{p}\in\Lambda}|\mathfrak{k}_{\mathfrak{p}-}|$ are large because they involve a cancellation of two terms that are on the order these numbers.

The chain $Z$ as just described satisfies Lemma 5.3's requirements.

*Proof of Lemma 5.3*: The proof has eight steps.

<u>Step 1</u>: The element $\Theta_-$ is given in (5.14). The element $\Theta_+ \in \mathcal{Z}_{ech,M}$ is given by the formula in (5.25) modulo the definition of the integers $\{\mathfrak{k}_{\mathfrak{p}0}\}_{\mathfrak{p}\in\Lambda}$. This step defines these integers. To start, suppose that $\Theta = (\hat{\upsilon}_+, (\mathfrak{k}_{\mathfrak{p}0}, O_{\mathfrak{p}+})_{\mathfrak{p}\in\Lambda})$ is a given element in $\mathcal{Z}_{ech,M}$. Fix $r \in (\delta_*, 7\delta_*)$ and fix $\mathfrak{p} \in \Lambda$. Write $\mathfrak{p}$ as $(p_1, p_2)$ and reintroduce $S_{1,r}$ and $S_{2,r}$ to denote the respective boundary spheres of the radius r coordinate balls centered around $p_1$ and $p_2$. Let $c^s_1$ and $c^s_2$ denote the respective intersections between $\mathfrak{f}(S)$ and the spheres $S_{1,r}$ and $S_{2,r}$. Let $\gamma^\Theta_{\mathfrak{p}+}$ denote the segment of the integral curve of $v$ from $\Theta$ that crosses $\mathcal{H}^+_\mathfrak{p}$ and let $\gamma^\Theta_{\mathfrak{p}-}$ denote the corresponding integral curve of $v$ from $\Theta_-$. The intersection point of $\gamma^\Theta_{\mathfrak{p}+}$ with $S_{1,r}$ has distance at most $c_0\delta$ from one endpoint of $c^s_1$ and that of $\gamma^\Theta_{\mathfrak{p}-}$ has distance at most $c_0\delta$ from the other endpoint. The intersection points of $\gamma^\Theta_{\mathfrak{p}+}$ and $\gamma^\Theta_{\mathfrak{p}-}$ with $S_{2,r}$ are similarly close to the respective endpoints of $c^s_2$.

Granted the preceding, define a 1-cycle in $\mathcal{H}^+_\mathfrak{p}$ as follows: Start where $\gamma^\Theta_{\mathfrak{p}+}$ intersects $S_{2,r}$ and proceed along $\gamma^\Theta_{\mathfrak{p}+}$ in $\mathcal{H}^+_\mathfrak{p}$ until it intersects the boundary of $S_{1,r}$. Then proceed along an arc of length at most $c_0\delta$ in this sphere to the nearby endpoint of $c^s_1$. Continue along $c^s_1$ until it ends, and then proceed along an arc of length at most $c_0\delta$ to the point where $\gamma^\Theta_{\mathfrak{p}-}$ intersects $S_{1,r}$. Travel in the reverse direction along $\gamma^\Theta_{\mathfrak{p}-}$ until it intersects $S_{2,r}$. Proceed next along an arc of length at most $c_0\delta$ to the nearby endpoint of $c^s_2$; and then along $c^s_2$ to its endpoint near $\gamma^\Theta_{\mathfrak{p}+} \cap S_{2,r}$. Close up the cycle by taking an arc of length at most $c_0\delta$ to $\gamma^\Theta_{\mathfrak{p}+} \cap S_{2,r}$. This closed 1-cycle is denoted by $\upsilon^s_\mathfrak{p}$.

The 1-cycle $\upsilon^s_\mathfrak{p}$ determines a class in $H_1(\mathcal{H}^+_\mathfrak{p} \cup \hat{\gamma}^+_\mathfrak{p} \cup \hat{\gamma}^-_\mathfrak{p})$, and thus an integer. In particular, there is a unique choice for $\mathfrak{k}_{\mathfrak{p}0}$ that makes this integer equal to the integer $\mathfrak{m}^C_\mathfrak{p} \in \{-1, 0, 1\}$ that is defined in (5.5). The corresponding set of such choices is used in (5.25) to define the element $\Theta_+$. Note in this regard that the set $\{\mathfrak{k}_{\mathfrak{p}0}\}_{\mathfrak{p}\in\Lambda}$ is determined up to an error of size at most ±1 solely by the pair S and $\mathfrak{f}$, and thus solely by the Heegaard Floer data.

<u>Step 2</u>: This step describes the chain $Z_0$. Consider first its intersection with $M_r$, this being the complement in M of the radius r coordinate spheres centered at the index 1 and index 2 critical points of $f$. The chain here is the image of S by a map that is constructed by modifying $\mathfrak{f}$ near the boundary of S.

To say more, let $\upsilon$ denote for the moment a segment $M_r \cap \hat{\upsilon}_+$. What follows is a consequence of what is said in Corollary 2.6. There exists an isotopy of $M_r$ with the following four properties: The isotopy moves only points in the radius $c_0\delta$ tubular neighborhood of $\upsilon$; and these points are moved at most distance $c_0\delta$. The isotopy preserves the relevant versions of $S_{1,r}$ and $S_{2,r}$ that contain the endpoints of $\upsilon$. The end member of the isotopy maps $\upsilon$ to the nearby segment of $M_r$'s intersection with the integral curves from $\Theta_+$. The derivatives of the isotopy to order 5 are bounded by a constant that depends only on the Heegaard Floer data. Choose such an isotopy; and likewise choose an isotopy of this sort for the other components $M_r \cap (\hat{\upsilon}_+ \cup \hat{\upsilon}_-)$. Compose $\mathfrak{f}$ with the resulting set of 2G isotopies and use $\mathfrak{f}_r$ to denote the resulting map from $\mathfrak{f}^{-1}(M_r)$ to $M_r$. The 2-chain $Z_0$ intersects $M_r$ as the $\mathfrak{f}_r$-image of $\mathfrak{f}^{-1}(M_r)$.

Fix $\mathfrak{p} \in \Lambda$ so as to consider the part of $Z_0$ in $\mathcal{H}_\mathfrak{p}$. The simplest case to consider is that where neither $\hat{\gamma}_\mathfrak{p}^+$ nor $\hat{\gamma}_\mathfrak{p}^-$ is associated to an end of C in $\mathcal{H}^+{}_\mathfrak{p}$. This implies that the integer $\mathfrak{m}^C{}_\mathfrak{p}$ is equal to zero. Here is a consequence: Use the arcs $\gamma^\Theta{}_{\mathfrak{p}+}$, $\gamma^\Theta{}_{\mathfrak{p}-}$ and the arcs that comprise the intersections of $\mathfrak{f}_r(\mathfrak{f}^{-1}(M_r))$ with $S_{1,r}$ and $S_{2,r}$ to define a version of the 1-cycle $\upsilon^s{}_\mathfrak{p}$ from Step 1. There is a smooth map from $[-1, 1] \times [-R-\ln r, R+\ln r]$ into $\mathcal{H}^+{}_\mathfrak{p}$ with the properties listed in the next equation. This equation uses $\mathfrak{f}_\mathfrak{p}$ to denote the map.

- $\mathfrak{f}_\mathfrak{p}$ *sends* $[-1, 1] \times (-R-\ln r, R+\ln r)$ *to the* $|u| < R+\ln r$ *part of* $\mathcal{H}^+{}_\mathfrak{p}$.
- $\mathfrak{f}_\mathfrak{p}$ *restricts to the respective boundary components* $\{-1, 1\} \times [-R-\ln r, R+\ln r]$ *as diffeomorphisms onto* $\gamma^\Theta{}_{\mathfrak{p}-}$ *and* $\gamma^\Theta{}_{\mathfrak{p}+}$.
- $\mathfrak{f}_\mathfrak{p}$ *restricts to the respective boundary components* $[-1, 1] \times \{-R-\ln r, R+\ln r\}$ *as immersions onto* $\mathfrak{f}_r(\mathfrak{f}^{-1}(S_{2,r}))$ *and* $\mathfrak{f}_r(\mathfrak{f}^{-1}(S_{1,r}))$.
- $\mathfrak{f}_\mathfrak{p}$ *restricts to a neighborhood of* $\{-1, 1\} \times [-R-\ln r, R+\ln r]$ *as an embedding which has the vector field* $v$ *tangent to its image*.

(5.28)

The intersection of $Z_0$ with $\mathcal{H}^+{}_\mathfrak{p}$ is defined to be the $f_\mathfrak{p}$ image of $[-1, 1] \times [-R-\ln r, R+\ln r]$.

Consider next the case when the closed integral curve $\hat{\gamma}_\mathfrak{p}^+$ is associated to an end of C in $\mathcal{H}^+{}_\mathfrak{p}$. In this case, $\mathfrak{m}^C{}_\mathfrak{p} = -1$ and the cycle $\upsilon^s{}_\mathfrak{p}$ is homotopic in $\mathcal{H}^+{}_\mathfrak{p}$ to -1 times the generator of its first homology. This implies that it is homotopic in $\mathcal{H}^+{}_\mathfrak{p} \cup \hat{\gamma}_\mathfrak{p}^+ \cup \hat{\gamma}_\mathfrak{p}^-$ to $\hat{\gamma}_\mathfrak{p}^+$. With understood, what follows describes the intersection of $Z_0$ with $\mathcal{H}^+{}_\mathfrak{p}$.

The 2-chain $Z_0$ is the image of a map into $\mathcal{H}^+{}_\mathfrak{p} \cup \hat{\gamma}_\mathfrak{p}^+$ from the complement in the strip $[-1, 1] \times [-R-\ln r, R+\ln r]$ of a disk centered at $(\frac{1}{2}, 0)$ with radius bounded by $\frac{1}{1000}$. Let $D \subset [-1, 1] \times [-R-\ln r, R+\ln r]$ denote such a disk. The map in this case is also denoted by $\mathfrak{f}_\mathfrak{p}$. This map obeys the conditions in (5.28) with it understood that the domain in the first bullet is the complement of D in $[-1, 1] \times [-R-\ln r, R+\ln r]$. The boundary of the closure of D is mapped by $\mathfrak{f}_\mathfrak{p}$ to the cycle $\hat{\gamma}_\mathfrak{p}^+$ and this boundary is the only part of the domain that lies where $1 - 3\cos^2\theta = 0$. The behavior of $\mathfrak{f}_\mathfrak{p}$ near the boundary of D is described by the next equation. The equation parametrizes a neighborhood of the boundary of D by coordinates $(\rho, \varphi)$ with $\rho$ the Euclidean coordinate on a neighborhood of 0 in $[0, 1)$ and $\varphi \in \mathbb{R}/2\pi\mathbb{Z}$. The boundary is the locus where $\rho = 0$. The angle in $[0, \pi]$ where $\cos\theta = \frac{1}{\sqrt{3}}$ is denoted by $\theta_*$.

$\mathfrak{f}_\mathfrak{p}$ *embeds a neighborhood of the boundary of* D *in* $\mathcal{H}_\mathfrak{p}$ *as the map that sends* $(\rho, \varphi)$ *to* $(u = 0, \theta = \theta_* - \frac{1}{1000}\rho, \phi = \varphi)$

(5.29)

The intersection of $Z_0$ with $\mathcal{H}_p$ in the case when only $\hat{\gamma}_p^-$ is associated to an end of C in $\mathcal{H}^+_p$ is the image of a smooth map from the complement in $[-1, 1] \times [-R - \ln r, R + \ln r]$ of a disk centered at $(-\frac{1}{2}, 0)$. The intersection of $Z_0$ with $\mathcal{H}_p$ in the case when both $\hat{\gamma}_p^+$ and $\hat{\gamma}_p^-$ are associated to ends of C in $\mathcal{H}^+_p$ is the image of a smooth map from the complement in $[-1, 1] \times [-R - \ln r, R + \ln r]$ of a disk centered at $(\frac{1}{2}, 0)$ and a disk centered at $(-\frac{1}{2}, 0)$. The map in either case is denoted by $\mathfrak{f}_p$; it obeys (5.28) with the domain in the first bullet suitably modified, and it obeys a suitable version of (5.29) near the boundary of a deleted disk.

The map $\mathfrak{f}_p$ in all cases depends on the parameters $(\delta, x_0, R)$ because $\mathcal{H}^+_p$, $v$, $\gamma^\Theta_{p-}$ and $\gamma^\Theta_{p+}$ depend on them. Even so, $\mathfrak{f}_p$ can be chosen so that this dependence is irrelevant to the purposes at hand. This is said formally by the next lemma. To set the stage, suppose that $(\delta', x_0', R')$ is a data set that can be used as in Section 1 to define Y and its stable Hamiltonian structure. Assume that $\delta' \leq \delta$, $x_0' \leq x_0$ and $R' \leq R$. The lemma uses $\mathcal{H}^+_p{}'$ to denote the $(\delta', x_0', R')$ version of the handle $\mathcal{H}^+_p$.

**Lemma 5.4**: *Suppose that $\mathfrak{f}_p$ is as described by the $(\delta, x_0, R)$ version of (5.29). There exists a diffeomorphism from $\mathcal{H}^+_p{}'$ to $\mathcal{H}^+_p$ whose composition with $\mathfrak{f}_p$ is described by the $(\delta', x_0', R')$ version of (5.28).*

*Proof of Lemma 5.4*: This follows in a straight forward manner from two facts: First, the arcs $\mathfrak{f}_r(\mathfrak{f}^{-1}(S_{1,r}))$ and $\mathfrak{f}_r^{-1}(\mathfrak{f}^{-1}(S_{2,r}))$ are essentially independent of $(\delta, x_0, R)$. Of course they depend on these parameters where the distance is bounded by $c_0\delta$ from their endpoints, this where the isotopies that give $\mathfrak{f}_r$ from $\mathfrak{f}$ are involved. However, given the bounds on the derivatives of these isotopies, this dependence is of no consequence. The second relevant fact is that the integer $\mathfrak{k}_p$ is independent of $(\delta, x_0, R)$. With the preceding understood, fix $(\delta', x_0', R')$ and identify the $|u| \leq R' + \ln r$ part of $\mathcal{H}_p$ with $[-1, 1] \times S^2$ via the map that sends the coordinates $(u, \theta, \phi)$ to $(\frac{1}{R' + \ln r} u, \theta, \phi)$. The $(\delta', x_0', R')$ versions of the arcs $\mathfrak{f}_r(\mathfrak{f}^{-1}(S_{1,r}))$ and $\mathfrak{f}_r^{-1}(\mathfrak{f}^{-1}(S_{2,r}))$ appear now as arcs on the respective $u = 1$ and $u = -1$ boundary spheres. Meanwhile, the $(\delta', x_0', R')$ version of the vector field $v$ appears as a smooth vector field in $[-1, 1] \times S^2$ and the arcs $\gamma^\Theta_{p+}$ and $\gamma^\Theta_{p-}$ appear as arcs that cross from the $u = -1$ boundary sphere to the $u = 1$ boundary sphere.

The formulas in PROPERTY 3 of Section 1e have the following implication: The vector field $v$ varies smoothly in $[-1, 1] \times S^2$ as the parameters $(\delta', x_0', R')$ vary. Meanwhile, Lemma 2.2 implies that the arcs $\gamma^\Theta_{p+}$ and $\gamma^\Theta_{p-}$ vary via ambient, compactly supported isotopies in the $1 - 3\cos^2\theta > 0$ part of $[-1, 1] \times S^2$ as the parameters $(\delta', x_0', R')$ vary. This fact, and what is said about the isotopies that define $\mathfrak{f}_r$ from $\mathfrak{f}$ imply what is asserted by the lemma.

Step 3: This step describes $Z_+$. The chain $Z_-$ has, but for notation, the same sort of description. This being the case, only $Z_+$ is discussed.

Consider first the part of $Z_+$ in $M_r$. This part lies in the radius $c_0\delta$ tubular neighborhood of the $M_r$ parts of the integral curves that comprise $\hat{\upsilon}_+$. To say more, let $\upsilon$ denote the $M_r$ part of one of these curves. The segment $\upsilon$ has distance at most $c_0\delta$ from a component of the $M_r$ part of an integral curve of $v$ from $\Theta_{C+}$, and it has distance at most $c_0\delta$ from an integral curve of $v$ from $\Theta_+$. Let $\gamma^C_{M+}$ and $\gamma^\Theta_{M+}$ denote these respective segments. These are also integral curves of $\upsilon$. Both start on the boundary of the radius r coordinate ball centered at an index 1 critical point of $f$ and both end on the boundary of a similar ball centered at one of $f$'s index 2 critical points. With $f^{-1}(1,2)$ viewed in the usual way as $(1,2) \times \Sigma$, both $\gamma^C_{M+}$ and $\gamma^\Theta_{M+}$ project to $\Sigma$ as points with distance $c_0\delta$ or less from the point in $C_+ \cap C_-$ that is defined by $\upsilon$. Fix an arc of length $c_0\delta$ or less in $\Sigma$ from the $\gamma^\Theta_{M+}$ point to the $\gamma^C_{M+}$ point. Denote the latter by $\iota$. The rectangle $[1,2] \times \iota$ when viewed in M intersects $M_r$ as an embedded rectangle that is foliated by integral curves of $\upsilon$ with one edge on $\gamma^C_{M+}$, the opposing edge on $\gamma^\Theta_{M+}$, a third edge on the boundary of the radius r coordinate ball centered on the relevant index 1 critical point of $f$, and the latter's opposing edge on the boundary of the radius r coordinate ball centered on the relevant index 2 critical point of $f$. There is an analogous rectangle for each component of $\hat{\upsilon}_+$. The union of these rectangles comprise the intersection of $Z_+$ with the closure of $M_r$.

Fix $\mathfrak{p} \in \Lambda$ so as to consider the intersection between $Z_+$ and $\mathcal{H}_\mathfrak{p}$. If $\hat{\gamma}^+_\mathfrak{p}$ is in $\Theta_{C+}$, then it is also in $\Theta_+$. In this case, $Z_+$ contains as a component the degenerate 2-chain that is given by the projection to $\mathcal{H}_\mathfrak{p}$ from $[-1,1] \times \mathcal{H}_\mathfrak{p}$ of the cylinder $[-1,1] \times \hat{\gamma}^+_\mathfrak{p}$. There is an analogous component if $\hat{\gamma}^-_\mathfrak{p}$ is in $\Theta_{C+}$ and $\Theta_+$. There is, in all cases, another component of $Z_+$ in the $|u| \leq R + \ln r$ part of $\mathcal{H}_\mathfrak{p}$, this denoted by $Z_{\mathfrak{p}+}$. What follows describes the latter.

Let $S_{1,r}$ and $S_{2,r}$ again denote the respective boundary spheres of the radius r coordinate balls that are centered on the index 1 and index 2 critical points from $\mathfrak{p}$. Let $c_{+1} \subset S_{1,r}$ and $c_{+2} \subset S_{2,r}$ denote the short arcs that comprise the intersection of these spheres with the $M_r$ part of $Z_+$. The concatenation of these two arcs with $\gamma_{C+}$ and $\gamma_+$ defines a closed cycle which is denoted by $\upsilon_{\mathfrak{p}+}$ in $\mathcal{H}^+_\mathfrak{p}$. This cycle will not be homotopically trivial in $\mathcal{H}^+_\mathfrak{p}$ unless $\gamma_{\mathfrak{p}+}$ and $\gamma^\Theta_{\mathfrak{p}+}$ have the same Proposition 2.8 integer label, this being $\mathfrak{k}_{\mathfrak{p}0}$.

The simplest case by far is that here $\upsilon_{\mathfrak{p}+}$ is null-homotopic. In this event, there is an embedded rectangle in the $|u| \leq R + \ln r$ part of $\mathcal{H}^+_\mathfrak{p}$ that is foliated by integral curves of $v$ with one edge $\gamma_{\mathfrak{p}+}$, the opposite edge $\gamma^\Theta_{\mathfrak{p}+}$, a third edge $c_{+1}$ and the opposite edge $c_{+2}$. This embedded rectangle is $Z_{\mathfrak{p}+}$.

Suppose next that $\upsilon_{\mathfrak{p}+}$ is not null-homotopic in $\mathcal{H}^+_\mathfrak{p}$. It proves convenient in this case to represent $Z_{\mathfrak{p}+}$ as the image of a submanifold in the $|u| \leq R + \ln r$ part of $[-1,1] \times \mathcal{H}^+_\mathfrak{p}$ via the projection to the $\mathcal{H}^+_\mathfrak{p}$ factor. The submanifold in question is the graph of a map from $[-1, 1] \times [-R-\ln r, R+\ln r]$ to $S^2$. To describe the latter map, first note that the cycle $\upsilon_{p+}$ can be viewed as a graph of a map from the boundary of $[-1, 1] \times [-R-\ln r, R+\ln r]$ to $S^2$. The image of the latter map on the $u = -R-\ln r$ part of the boundary is the arc $c_{+2}$ and its image on the $u = R+\ln r$ part of the boundary is the arc $c_{+1}$. The map on the boundary component $\{1\} \times [-R-\ln r, R+\ln r]$ is the map $u \to (\theta(u), \phi(u))$ that parametrizes the $v$-integral curve $\gamma_{p+}$. The map on the component $\{-1\} \times [-R-\ln r, R+\ln r]$ is the corresponding parametrizing map for the $v$-integral curve $\gamma_{p+}$. Use $q_p$ to denote the map. The extension to the interior of the rectangle $[-1,1] \times [-R-\ln r, R+\ln r]$ can be made so as to have the properties listed below in (5.30). The list uses $\mathfrak{k}_{p+}$ to denote the Proposition 2.8 label of $\gamma_{p+}$. Equation (5.30) also refers to the intersection number between the image of $q_p$ and $\theta = 0$ point in $S^2$. This intersection number is well defined by by virtue of the fact that the boundary of $[-1, 1] \times [-R-\ln r, R+\ln r]$ is mapped by $q_p$ to the $1 - 3\cos^2\theta > 0$ portion of $S^2$.

- $q_p$ *maps to the* $\cos\theta > -\frac{1}{\sqrt{3}}$ *portion of* $S^2$.
- $q_p$ *restricts to the respective boundary components* $[-1, 1] \times [-R-\ln r, R+\ln r]$ *as described above.*
- $q_p$ *restricts to a neighborhood of the boundary of* $[-1, 1] \times [-R-\ln r, R+\ln r]$ *so that variation of* u *with constant* $[-1, 1]$ *factor parametrizes part of an integral curve of* $v$.
- *The image via* $q_p$ *of* $[-1, 1] \times \{0\}$ *lies where* $1 - 3\cos^2\theta > 0$.
- *The image of* $q_p$ *has* $|\mathfrak{k}_{p+} - \mathfrak{k}_{p0}|$ *intersections with the* $\theta = 0$ *point in* $S^2$. *Each such intersection is transverse; and all have the same intersection number.*

(5.30)

The construction of a map of this sort is straightforward and so the task is left to the reader.

Step 4: This step summarizes Hutching's definition of $I_{ech}$ as a prelude to the calculations that are done in the upcoming Step 5. To set the background, suppose that $\Theta$ is a finite set whose typical element is a pair $(\gamma, m)$ with $\gamma$ a closed integral curve of $v$ and m a positive integer. Assume that no two pair from $\Theta$ share the same closed integral curve. Let $\Theta'$ denote another such set, and suppose that $[\Theta]-[\Theta']$ bounds in Y. Let Z denote a 2-chain with boundary $[\Theta]-[\Theta']$. Definition 2.14 in [echir] supplies a formula for computing $I_{ech}(Z)$ that is written as

$$I_{ech}(Z) = c_\tau(Z) + Q_\tau(Z) + \mu_\tau(\Theta) - \mu_\tau(\Theta')$$

(5.31)

What follows describes the meaning of the various terms in this equation.

To start, introduce the 2-plane bundle $K^{-1} \subset TY$, this being the kernel of $a$ with orientation given by $w$. The subscript $\tau$ in (5.31) refers to a chosen homotopy class of trivialization of $K^{-1}$ along the integral curves that comprise $\Theta$ and $\Theta'$. The individual terms in (5.32) depend on $\tau$ but their sum does not. As explained next, there is an almost canonical choice for $\tau$ in the present circumstances.

Consider first $\tau$ along a given $\mathfrak{p} \in \Lambda$ version of either $\hat{\gamma}_\mathfrak{p}^+$ or $\hat{\gamma}_\mathfrak{p}^-$. The 1-forms $\{d\theta, du\}$ map $K^{-1}$ isomorphically to the trivial bundle bundle on such a curve; and so defines the desired version of $\tau$.

Consider next the case when $\gamma \subset M_\delta \cup (\cup_{\mathfrak{p} \in \Lambda} \mathcal{H}^+_\mathfrak{p})$ is a closed integral curve of $v$ from an element in $\mathcal{Z}_{ech,M}$. Fix $\mathfrak{p} \in \Lambda$ so as to describe a trivialization of $K^{-1}$ along $\gamma$'s intersection with $\mathcal{H}_\mathfrak{p}$. To this end, reintroduce the function $\hat{h}$ that is depicted in (4.1). The 1-forms $\{d\phi, d\hat{h}\}$ map $K^{-1}|_\gamma$ isomorphically to the trivial bundle along $\gamma \cap \mathcal{H}^+_\mathfrak{p}$. Let $\tau_\mathfrak{p}$ denote this trivialization. It follows from what is said by CONSTRAINT 3 in Part 3 of Section 1b that the trivializations that comprise the $\{\tau_\mathfrak{p}\}_{\mathfrak{p} \in \Lambda}$ extend as trivializations of $K^{-1}$ along $\gamma \cap ((M_\delta - \Sigma) \cup (\cup_{\mathfrak{p} \in \Lambda} \mathcal{H}^+_\mathfrak{p}))$. Note that each component of this intersection is a segment that passes through precisely one $\mathfrak{p} \in \Lambda$ version of $\mathcal{H}^+_\mathfrak{p}$. It also follows from this same CONSTRAINT 3 that these trivializations do not agree when compared along $\Sigma$. Agreement can be obtained by suitably rotating the the fiber of the trivial bundle $\gamma \times \mathbb{R}^2$ along $\gamma$'s intersection with the $f \in (\frac{5}{4}, \frac{3}{2})$ part of $\gamma \cap M_\delta$. In particular, the rotation on any given segment of this intersection can be chosen to have angle either $\frac{\pi}{2}$ or $-\frac{\pi}{2}$ depending on whether the + or - sign appears in CONSTRAINT 3 for the relevant component of $T_+ \cap T_-$. This convention defines a homotopy class of trivialization of $K^{-1}$ along the whole of $\gamma$; and this gives the desired version of $\tau$.

With $\tau$ defined, consider next the term in (5.31) that is denoted by $\mu_\tau(\Theta)$. The definition in the case when $\Theta \in \mathcal{Z}_{ech,M}$ is straightforward by virtue of the fact that all of the closed integral curves from $\Theta$ are hyperbolic. In this case, these terms are defined as follows: Let $\gamma$ denote one of the integral curves of $v$ from $\Theta$. Lie transport of $K^{-1}|_\gamma$ by $v$ along a circumnavigation of $\gamma$ rotates a basis vector in the chosen trivialization by an angle $n\pi$ with $n \in \mathbb{Z}$. This integer n is $\mu_\tau(\Theta)$. It follows from the formula in PROPERTY 3 in Section 1e that n = 0 when $\gamma$ is some $\mathfrak{p} \in \Lambda$ version of either $\hat{\gamma}_\mathfrak{p}^+$ or $\hat{\gamma}_\mathfrak{p}^-$. Suppose next that $\gamma$ is a closed orbit from $\Theta$ that corresponds to an irreducible, HF-cycle. It follows from PROPERTY 3 in Section 1e, from what is said in Lemma 2.2, and from what is said about $\tau$ in the preceding paragraph that $|\mu_\tau(\Theta)|$ is no greater than G.

The terms $c_\tau(Z)$ and $Q_\tau(Z)$ depend on the chain Z. What follows describes these terms for the case when $\Theta$ and $\Theta'$ are from $\mathcal{Z}_{ech,M}$. To start, it is necessary to first choose a properly embedded, oriented surface with boundary in $[-1,1] \times Y$ with certain properties. Use $S$ to denote this surface.

- *The components of $S \cap (\{-1\} \times Y)$ are the integral curves of $v$ from $\Theta'$, and those from $S \cap (\{1\} \times Y)$ are the integral curves of $v$ from $\Theta$.*
- *The surface $S$ is transversal to $\{-1, 1\} \times Y$.*
- *The image of $S$ in $Y$ via the projection defines a 2-chain that gives the same class as $Z$ in $H_2(Y; [\Theta] - [\Theta'])$.*

(5.32)

To define $c_\tau(Z)$, choose a section along $S$ of $K^{-1}$ that restricts to $S \cap (\{-1, 1\} \times Y)$ as a basis vector for a trivialization defined by $\tau$. A sufficiently generic choice will vanish in a transversal fashion, and so will have a well defined Euler number. This number is independent of the chosen section; and it is independent of $S$ except to the extent that the constraints in (5.32) are imposed. This Euler number is $c_\tau(Z)$.

To define $Q_\tau(Z)$, let $N \to S$ denote the latter's normal bundle. The constraints given by the first and second bullets in (5.32) imply that $N$ along $S \cap (\{-1,1\} \times Y)$ is canonically isomorphic to $K^{-1}$. Choose a section of $N$ that restricts to $S \cap (\{-1, 1\} \times Y)$ as a basis vector for a trivialization defined by $\tau$. A sufficiently generic choice will vanish transversaly and so have a well defined Euler number. The latter is independent of the chosen section; and it is also independent of $S$ except to the extent that (5.32) is imposed. This Euler number is $Q_\tau(Z)$.

Step 5: This step, Step 6 and Step 7 derive bounds for the absolute values of $I_{ech}(Z_0), I_{ech}(Z_+)$ and $I_{ech}(Z_-)$. To this end, the simplest case to consider is $I_{ech}(Z_0)$. It follows directly from what is in Step 2 about $Z_0 \cap M_r$ and from Lemma 5.4 that the absolute value of each term that appears in (5.31) is bounded by $c_0$ where $c_0$ depends only on the Heegaard Floer data. Thus $|I_{ech}(Z_0)|$ is also bounded by $c_0$ with $c_0$ dependent only on the Heegaard Floer data.

Consider next the case of $I_{ech}(Z_+)$. As noted in Step 4, the respective absolute value of the relevant version of $\mu_\tau(\Theta)$ and $\mu_\tau(\Theta')$ is bounded by $c_0$ where $c_0$ depends only on the Heegaard Floer data. This is not true for either $c_\tau(Z)$ or $Q_\tau(Z)$.

To compute $c_\tau(Z_+)$ and $Q_\tau(Z_+)$, it is necessary to first choose a surface $S$ that is described by (5.32). If $\gamma \in \Theta_{C+}$ is from the set $\{\hat{\gamma}_p^+ \cup \hat{\gamma}_p^-\}_{p \in \Lambda}$, then $[-1, 1] \times \gamma$ is a component of $S$. The boundaries of the remaining components account for the remaining integral curves in $\Theta_{C+}$ and $\Theta_{C-}$.

The remaining components of $S$ lie in $[-1, 1] \times (M_r \cup (\cup_{p \in \Lambda} \mathcal{H}_p))$. Each such compnent is a closed annulus; thus each is diffeomorphic to $[-1, 1] \times S^1$. These annuli are in 1-1 correspondence with the integral curves from $\Theta_{C+}$ that come from the HF-cycle $\hat{\upsilon}_+$. Let $\gamma^C \in \Theta_C$ denote one of these integral curves, and let $\gamma^\Theta \subset \Theta_+$ denote the corresponding curve. Use $S_\gamma \subset S$ to denote the associated annular component of $S$. The respective boundary components $\{-1\} \times S^1$ and $\{1\} \times S^1$ of $S_\gamma$ appear in $\{-1\} \times Y$ and in $\{1\} \times Y$ as $\{-1\} \times \gamma^\Theta$ and $\{1\} \times \gamma^C$.

An annulus can be viewed as a union of rectangles with disjoint interior such that each rectangle has one edge on one of the boundary circles and the other edge on the other boundary circle. The two other edges intersect the boundary circle only at the corners. A description of $S_\gamma$ of this sort is presented in what follows.

The part of $S_\gamma$ in $[-1, 1] \times M_r$ is a disjoint union of rectangles. To describe a given such rectangle, recall first from Step 4 that $Z_+$ intersects $M_r$ as a disjoint union of embedded rectangles. A subset of these rectangles that comprise $Z_+ \cap M_r$ have one edge on a component of $\gamma^C \cap M_r$ and the opposing edge on the nearby component of $\gamma^\Theta \cap M_r$. A given rectangle is foliated by integral curves of $\mathfrak{v}$; these integral curves intersect $\Sigma$ as an arc, $\iota$, of length bounded by $c_0\delta$ in a component of $T_+ \cup T_-$. This rectangle is in the $f \in (1, 2)$ part of $M_r$, and so it can be viewed as the graph from a rectangle in $[-1, 1] \times (1, 2)$ to $T_+ \cup T_-$. The latter rectangle has one edge on $\{-1\} \times (1, 2)$ and the opposite edge on $\{1\} \times (1, 2)$. These edges map to $M_r$ as $\gamma^\Theta$ and $\gamma^C$. The other two edges map to the boundary of $M_r$. Any given slice of the $[-1, 1]$ factor maps to an integral curve of $\mathfrak{v}$. Identify $(1, 2) \times \Sigma$ with the $f^{-1}(1, 2)$ part of $M$ and this rectangle appears as a submanifold in $[-1, 1] \times M_r$. The latter is a component of the intersection of $S_\gamma$ with $[-1, 1] \times M_r$.

Each of the remaining rectangles that comprise $S_\gamma$ is the $|u| \leq R + \ln r$ part of the intersection between $S_\gamma$ and a $\mathfrak{p} \in \Lambda$ version of $[-1, 1] \times \mathcal{H}_\mathfrak{p}$. Each rectangle of this sort is a suitable parameterization of the relevant $\mathfrak{p} \in \Lambda$ version of the graph in (5.33).

Step 6: With $S_\gamma$ understood, consider now the latter's contribution to the $Z_+$ version of what is denoted by $c_\tau(\cdot)$ in (5.31). It proves useful for this purpose to first choose a particular section of the dual to $K^{-1}$ over $S_\gamma$. The pair of sections of $T^*Y$ that identify $K^{-1}$ with the trivial bundle over $\gamma^C$ and $\gamma^\Theta$ so as to extend $\tau$ without a zero over the rectangles that comprise $Z_+ \cap M_r$. One of these sections can be taken so that it restricts as $d\hat{h}$ to the $|u| \geq R + \ln r$ part of any given $\mathfrak{p} \in \Lambda$ version of $\mathcal{H}_\mathfrak{p}$. Let $\eta_r$ denote the latter. The 1-form $d\hat{h}$ is defined over the whole of $\mathcal{H}_\mathfrak{p}$ and so defines by restriction a section of $(K^{-1})^*$ over $\mathcal{H}_\mathfrak{p}$. This understood, the pull-back to $S_\gamma$ via the projection to Y of the 1-form $\eta_r$ over the $[-1, 1] \times M_r$ part of $S_\gamma$, and the pull-back to $S_\gamma$ of $d\hat{h}$ over the part of $S_\gamma$ in any given $\mathfrak{p} \in \Lambda$ version of $[-1, 1] \times \mathcal{H}_\mathfrak{p}$ defines a section over $S_\gamma$ of $(K^{-1})^*$. This section is denoted in what follows by $\eta$.

The section $\eta$ of $(K^{-1})^*$ has no zeros on the $[-1, 1] \times M_r$ part of $S_\gamma$. Fix $\mathfrak{p} \in \Lambda$ so as to consider the zero locus of $\eta$ on the $[-1, 1] \times \mathcal{H}_\mathfrak{p}$ part of $S_\gamma$. It is assumed in what follows that this part of $S_\gamma$ is not empty. The 1-form $d\hat{h}$ annihilates $K^{-1}|_{\mathcal{H}_\mathfrak{p}}$ at the points where it is proportional to $a$. This locus consists of the $\theta \in \{0, \pi\}$ arcs and the circles

where both $u = 0$ and $1 - 3\cos^2\theta = 0$. All of these zeros of $d\hat{h}|_{K^{-1}}$ are transversal. Note that the projection to $\mathcal{H}_\mathfrak{p}$ of the graph defined by (5.30) intersects only the $\theta = 0$ arc.

The fifth bullet in (5.30) and what is said in the preceding paragraph have the following implication: The bundle $K^{-1}$ over $S_\gamma$ has a smooth section with the large $|s|$ asympotics required to compute $c_\tau(\cdot)$ and which vanishes only at the points in the $\mathfrak{p} \in \Lambda$ versions of $S_\gamma \cap ([-1, 1] \times \mathcal{H}_\mathfrak{p})$ where $\theta = 0$. Moreover, the section vanishes transversaly at each such point, and all of these points make the same contribution to $c_\tau$. This section also generates the kernel of $d\hat{h}$ on $S_\gamma \cap ([-1, 1] \times \mathcal{H}_\mathfrak{p})$. The section in question is denoted by $\mathfrak{z}$.

To compute the contribution of the zeros of $\mathfrak{z}$ to $c_\tau$, introduce Euclidean coordinates $z_1 = \sin\theta \cos\phi$ and $z_2 = \sin\theta \sin\phi$ for a neighborhood in $S^2$ of the $\theta = 0$ point. The bundle $K^{-1}$ on $\mathcal{H}_\mathfrak{p}$ is trivialized near the $\theta = 0$ arc by the pair of 1-forms $\{dz_1, dz_2\}$. It follows from the formula for $w$ in PROPERTY 3 of Section 1e that the pair $(-dz_1, dz_2)$ define an orientation *perserving* identification of $K^{-1}$ near the $\theta = 0$ point with the product $\mathbb{R}^2$ bundle. Here, $\mathbb{R}^2$ is oriented in the standard manner. This identification of $K^{-1}$ with the product bundle near the $\theta = 0$ arc is assumed implicitly in what follows.

Let $(s_0, u_0) \in [-1, 1] \times [-R - \ln r, R + \ln r]$ denote a point that is sent to the $\theta = 0$ point by the map $q_\mathfrak{p}$ from (5.31). For the purposes at hand, no generality is lost by assuming that the $q_\mathfrak{p}$ near $(s_0, u_0)$ is the map

$$(s, u) \to (z_1 = s - s_0, z_2 = u - u_0) \,, \tag{5.33}$$

and that $\mathfrak{z}$ appears on a neighborhood of $(x_0, u_0)$ as the map to $\mathbb{R}^2$ given by

$$(s, u) \to (u - u_0, s - s_0) \,. \tag{5.34}$$

These formulae imply that each zero of $\mathfrak{z}$ on $S^{-1} \cap ([-1, 1] \times \mathcal{H}_\mathfrak{p})$ contributes $-1$ to $c_\tau(Z_+)$.

<u>Step 7</u>: Consider next the $Z = Z_+$ version of the term $Q_\tau$ in (5.31). The section $\mathfrak{z}$ that was just used to compute $c_\tau(Z_+)$ will be used to compute $Q_\tau(Z_+)$ as well. Its use in this regard exploits the following observation: Suppose that $X_1$ and $X_2$ are smooth manifolds, $g: X_1 \to X_2$ a smooth map and $\mathcal{V} \subset X_1 \times X_2$ the graph of $g$. Then the differential of the projection to $X_2$ identifies the normal bundle to $\mathcal{V}$ with $TX_2|_\mathcal{V}$.

As explained above, a given rectangle that comprises $S_\gamma \cap ([-1, 1] \times M_r)$ can be viewed as a graph of a map from a rectangle in $[-1, 1] \times (1, 2)$ to $T_+ \cup T_- \subset \Sigma$. Doing so identifies the normal bundle to this part of $S_\gamma$ with $T\Sigma$. Meanwhile, $K^{-1}$ is also identified as $T\Sigma$ here. This understood, the section $\mathfrak{z}$ defines a section of the normal bundle to $S_\gamma$ on this rectangle.

Fix $\mathfrak{p} \in \Lambda$. If non-empty, then the $|u| \leq R + \ln r$ part of $S_\gamma \cap ([-1, 1] \times \mathcal{H}_\mathfrak{p})$ is the graph of the map $q_\mathfrak{p}$ that is depicted in (5.32). This being the case, the differential of the projection map to $S^2$ from $[-1, 1] \times [-R-\ln r, R+\ln r] \times S^2$ identifies the normal bundle to this part of $S_\gamma$ with $TS^2$. The differential of this projection also defines a homomorphism from $K^{-1}|_{\mathcal{H}_\mathfrak{p}}$ to $TS^2$. This homomorphism is not an isomorphism only where $v$ is tangent to the $S^2$ factor in $\mathcal{H}_\mathfrak{p}$. A look at PROPERTY 3 in Section 1e finds that such is the case where $1 - 3\cos^2\theta = 0$. The bundle $K^{-1}$ on this locus is spanned by $\partial_u$ and $3f\partial_\theta - g\partial_\phi$, and so the kernel of the homomorphism to $TS^2$ is spanned by $\partial_u$.

The image of $\mathfrak{z}$ by this homomorphism to $TS^2$ defines a section of the normal bundle of $S_\gamma$ over its intersection with $[-1, 1] \times \mathcal{H}_\mathfrak{p}$. This normal bundle section can be modified slightly with no added zeros near the $|u| = R + \ln r$ boundary of its domain so that the result extends the normal bundle section that is defined by $\mathfrak{z}$ over the $|u| \leq R + \ln r$ part of $S_\gamma \cap ([-1, 1] \times M_r)$. Here is why: The intersection of $S_\gamma$ with the $|u| = R + \ln r$ spheres in $[-1, 1] \times \mathcal{H}_\mathfrak{p}$ lies where the angle $\theta$ differs from $\frac{\pi}{2}$ by at most $c_0\delta$; and the 1-form $a$ on this part of $\mathcal{H}_\mathfrak{p}$ differs $2e^{2(|u|-R)} du$ by a term with absolute value bounded by $c_0\delta$.

It follows as a consequence of what was just said that there exists a section of the normal bundle of $S_\gamma$ that has the requisite large $|s|$ asymptotics for computing $Q_\tau$ and whose zero locus consists solely of the various $\mathfrak{p} \in \Lambda$ versions of the $|u| < R + \ln r$ part of $S_\gamma \cap ([-1, 1] \times \mathcal{H}_\mathfrak{p})$ where $\theta = 0$. Note that this can be done despite the fact that the homomorphism from $K^{-1}$ to $TS^2$ has rank 1 on the $1 - 3\cos^2\theta = 0$ locus.

The assertion about the $1 - 3\cos^2\theta = 0$ locus can be justified by straightforward arguments using five facts noted previously. The first is that $S_\gamma$ intersects each sphere where $|u| = R + \ln r$ only where $\theta$ has distance at most $c_0\delta$ from $\frac{\pi}{2}$. The second fact is that each such intersection is an arc with length bounded by $c_0\delta$. The third fact is that the bundle $K^{-1}$ where $|u| = R + \ln r$ and $\theta$ is close to $\frac{\pi}{2}$ is very nearly identical to $TS^2$. In particular, it is very nearly spanned by the pair $\{\partial_\theta, \partial_\phi\}$. The fourth fact is that the vector field $3f\partial_\theta - g\partial_\phi$ is in $K^{-1}$ along the locus where $1 - 3\cos^2\theta = 0$. The final fact is that this vector field has non-zero image in $TS^2$ via the homomorphism from $K^{-1}$.

To see about the contribution to $Q_\tau(Z_+)$ from a $\theta = 0$ point in the $|u| < R + \ln r$ part of $S_\gamma \cap ([-1, 1] \times \mathcal{H}_\mathfrak{p})$, parametrize $S^2$ near the $\theta = 0$ point as before using the coordinates $(z_1, z_2)$. The normal bundle to $TS^2$ near the $\theta = 0$ point has the oriented basis $\{\frac{\partial}{\partial z_1}, \frac{\partial}{\partial z_2}\}$, and this basis is used implicitly in what follows to identify $TS^2$ near the $\theta = 0$ point with the product $\mathbb{R}^2$ bundle. With the preceding kept in mind, parameterize the $|u| < R + \ln r$ part of $S_\gamma \cap ([-1, 1] \times \mathcal{H}_\mathfrak{p})$ as the graph of the map $q_\mathfrak{p}$ in (5.33) and let $(s_0, u_0)$ again denote a $\theta = 0$ point. It follows from (5.34) that the normal bundle section defined by $\mathfrak{z}$ here appears as the map to $\mathbb{R}^2$ given by

$$(s, u) \to (-u + u_0, s - s_0) .$$
(5.35)

The determinant of the differential of this map is +1. As a consequence, each $\theta = 0$ point in the $|u| < R + \ln r$ part of $S_\gamma \cap ([-1, 1] \times \mathcal{H}_p)$ contributes +1 to $Q_\tau(Z_+)$.

Step 8: As noted in Step 5, the absolute value of the contribution of $\mu_\tau(Z_+)$ to the $Z = Z_+$ version is bounded by a constant that depends only on the Heegaard Floer data. As seen in Step 6, each $\theta = 0$ point in the $|u| < R + \ln r$ part of $S_\gamma \cap ([-1, 1] \times \mathcal{H}_p)$ contributes -1 to $c_\tau(Z_+)$, and as seen in Step 7, each such point contributes +1 to $Q_\tau(Z_+)$. This implies that the the $|u| < R + \ln r$ part of $S_\gamma \cap ([-1, 1] \times \mathcal{H}_p)$ contributes *zero* to the sum $c_\tau(Z_+) + Q_\tau(Z_+)$ that appears in the $Z = Z_+$ version of (5.31). Granted this, it then follows that $|I_{ech}(Z_+)|$ is bounded by a constant that depends solely on the Heegaard Floer data.

### c) Applications of Gromov-compactness

This subsection has two part. Part 1 states the version of Gromov's compactness theorem for J-holomorphic submanifolds that is used in Part 2. Part 2 gives the applications. Of particular interest for [KLTIII] is the upcoming Proposition 5.8 from Part 2. This lemma states conclusions that are stronger versions of those stated by Lemmas 4.3 and 4.4.

*Part 1*: To set the notation for the statement of the compactness theorem, suppose that $U \subset Y$ is a given open set with compact closure. Given $\varepsilon > 0$, introduce $U_\varepsilon$ to denote the set of points in Y with distance $\varepsilon$ or less from U. The statement of the compactness theorem also introduces the following terminology: Let $I \subset \mathbb{R}$ denote a given bounded, open interval, and let U be as just described. A subset $C \subset I \times U$ is said to be a J-holomorphic subvariety if there is an open neighborhood in $\mathbb{R} \times U_\varepsilon$ of the closure of $I \times U$ and a set in this neighborhood with the following properties: The intersection of this set with $I \times U$ is C; the set has finite, 2-dimensional Hausdorff measure; the set has no isolated components; and the complement of a finite set of points from this set is a submanifold with J-invariant tangent space. The term *weighted J-holomorphic subvariety* in $I \times U$ denotes a finite set of pairs consisting of an irreducible, J-holomorphic subvariety in $I \times U$ and positive integer with it understood that no two pair share the same subvariety.

**Proposition 5.5**: *Fix an open set* $U \subset Y$, *an open, bounded interval* $I \subset \mathbb{R}$ *and* $\mathcal{L} \geq 1$. *Given this data, there exists* $\kappa > 1$ *with the following significance: Let* $\{C\}_{n=1,2,...}$ *denote a* sequence of J-holomorphic subvarieties in $\mathbb{R} \times U$. Suppose this sequence has the following property: There exists $x \subset \mathbb{R}$ and $\varepsilon > 0$ such that

$$\int_{C_n \cap ([x-2, x+2] \times U_\varepsilon)} (ds \wedge \hat{a} + w) \le \mathcal{L}$$

for each n. Then there exists a weighted J-holomorphic subvariety $\vartheta$ and a subsequence of the original (hence renumbered consecutively from 1) with the following properties:

- $\lim_{n \to \infty} (\sup_{z \in C_n} \text{dist}(z, \cup_{(C,m) \in \vartheta} (C \cap (I \times U)) + \sup_{z \in C \cap (I \times U)} \text{dist}(z, C_n)) = 0$.
- Let $\mu$ denote a smooth, bounded 2-form defined on a neighborhood of $I \times U$ in $\mathbb{R} \times Y$. Then $\lim_{n \to \infty} (\int_{C_n} \mu - \Sigma_{(C,m) \in \vartheta} \, m \int_C \mu) = 0$.

*Proof of Proposition 5.5*: The existence of a weighted J-holomorphic subvariety and subsequence with the required properties follows from the standard sort of compactness theorems for sequences of J-holomorphic subvarieites in 4-dimensional manifolds. See for example, one of [W], [Y], Section 6 of [T3], or [Hum].

The next lemma states an important corollary to Proposition 5.5.

**Lemma 5.6**: *Fix an open set $U \subset Y$, $\mathcal{L} \ge 1$ and $\varepsilon > 0$. Given this data, there exists $\kappa > 1$ with the following significance: Let $C \subset \mathbb{R} \times Y$ denote a J-holomorphic subvariety. Suppose that $x \subset \mathbb{R}$ is such that*
- $\int_{C \cap ([x-2, x+2] \times U_1)} w < \kappa^{-1}$
- $\int_{C \cap ([x-2, x+2] \times U_1)} ds \wedge a \le \mathcal{L}$.

*Let $C' \subset C \cap ([-x-1, x+1] \times U)$ denote a connected component. Then there is an integral curve of $v$ such that each point of $(C'|_s \cap U)$ for each $s \in [x-1, x+1]$ has distance no greater than $\varepsilon$ from the intersection of this curve with $U$, and vice versa. Let $\gamma$ denote this integral curve. There exists a positive integer m which is less than $\mathcal{L}$ and such that if $\mu$ is a smooth 2-form on $[x-1, x+1] \times U$ with $\|\mu\|_\infty = 1$ and $\|\nabla \mu\|_\infty \le \varepsilon^{-1}$, then*

$$|\int_{C'} \mu - m \int_{[x-1, x+1] \times \gamma} \mu | \le \varepsilon.$$

*Proof of Lemma 5.6* The conditions given by the two bullets bound the area of C's intersection with $[x-2, x+2] \times U_1$. This understood, the assertion can be proved by assuming the contrary and invoking Proposition 5.5 so as to obtain a non-sense conclusion: There exists a sequence of J-holomorphic subvarieties whose n'th member obeys the assumptions with $\kappa = n$, $\varepsilon = n^{-1}$ that converges as described in Proposition 5.1 to a weighted J-holomorphic subvariety with the following property: The subvariety component from each constituent pair is an $\mathbb{R}$-invariant cylinder. For more details, see for example the proofs of Lemma 3.6 and Corollary 4.7 in [T4] as these state an analogous assertion in the context where $a$ is a contact 1-form and $w = da$.

*Part 2*: The next lemma states a crucial result for the proof of the upcoming Proposition 5.8 and for use in [KLTIII].

**Lemma 5.7**: *Given $\mathcal{I} > 0$ and $\rho \in (0, \delta_*]$, there exists $\kappa > 1$ with the following significance: Fix data $\delta, x_0, R$ and J suitable for defining the geometry of Y and $\mathbb{R} \times Y$ with $\delta < \frac{1}{2}\rho$. Let $C \subset \mathbb{R} \times Y$ denote an ech-HF subvariety with $I_{ech}(C) \leq \mathcal{I}$ or $\int_C w \leq \mathcal{I}$. Let $\mathfrak{c}_C$ denote maximum value of $\cos^2(\theta)$ on C's intersection with the portion of the radius $\delta_*$ coordinate ball centered on each of the index 1 and index 2 critical points of f where the radius is greater than $\rho$. Then $1 - 3\mathfrak{c}_C^2 > \kappa^{-1}$.*

*Proof of Lemma 5.7*: Were there no such $\kappa$, one could find an infinite sequence of J-holomorphic submanifolds that obeyed the assumptions of the lemma, but whose n'th member violated the conclusions with $\kappa = n$. Use the bounds from the two bullets in Proposition 5.1 to invoke Proposition 5.5 with $U = M_{\rho/2}$. Proposition 5.5 finds a subsequence that converges on compact subsets of the $f \in (1, 2)$ part of $\mathbb{R} \times M_{\rho/2}$ in a pointwise manner to a J-holomorphic subvariety in $\mathbb{R} \times M_\rho$. The limit subvariety must have an irreducilble component that intersects the $f = 1$ or $f = 2$ locus. Such a component would be the intersection of $\mathbb{R} \times M_{\rho/2}$ with a submanifold from either $\mathcal{M}_{f=1}$ or $\mathcal{M}_{f=2}$. Let $C_f$ denote this component of the limit. Either $C_f$ is disjoint from the remainder of the limit or not. If not, then all sufficiently large n submanifolds from the chosen subsequence must intersect $\mathbb{R} \times M_\delta$ at points where $f \neq (1, 2)$. This is ruled out by the assumptions. On the other hand, if $C_f$ is disjoint from the rest of the limit, then a winnowing of the chosen subsequence and renumbering will result in a new subsequence with the following property: There exists some $s \in \mathbb{R}$ such that any given large n member has a component that intersects $\mathbb{R} \times M_\delta$ with all points having distance at most $n^{-1}$ from $\{s\} \times C_f$. Moreover, all points in $\{s\} \times C_f$ will have distance at most $n^{-1}$ to this component of the n'th element from the subsequence. However, such a component of an admissible J-holomorphic subvariety will intersect at least one submanifold from $\mathcal{M}_\Sigma$, and so it will intersect all submanifolds from $\mathcal{M}_\Sigma$ unless it is an actual element in $\mathcal{M}_\Sigma$. This last option is also ruled out by the assumptions.

The next proposition subsumes some of what is asserted by Lemmas 4.3 and 4.4. This proposition identifies Proposition 3.2's moduli space $\mathcal{M}_\Sigma$ with $\mathbb{R} \times (1, 2)$ using the diffeomorphism described by Proposition 3.2.

**Proposition 5.8**: *Given $\mathcal{I} > 0$ and $\rho \in (0, \delta_*]$, there exists $\kappa > 1$ and $z_\diamond \in (0, \rho^2)$ with the following significance: Fix data $\delta, x_0, R$ and $J$ that are suitable for defining the geometry of $Y$ and $\mathbb{R} \times Y$ with $\delta < \frac{1}{100}\rho$. Let $C \subset \mathbb{R} \times Y$ denote an ech-HF subvariety with either $I_{ech}(C) \leq \mathcal{I}$ or $\int_C w \leq \mathcal{I}$.*

- *$C$ intersects the $f \leq 1 + z_\diamond$ part of $\mathbb{R} \times M_\delta$ in the union of the radius $\rho_\diamond$ coordinate balls centered at the index 1 critical points of $f$. In addition, $C$ intersects each subvariety from the $\mathbb{R} \times (1, 1+z_\diamond)$ part of $\mathcal{M}_\Sigma$ precisely once in each $\mathfrak{p} \in \Lambda$ version of $\mathbb{R} \times \mathcal{H}_\mathfrak{p}$ and this intersection point lies where $u \in (0, R + \ln\rho)$*
- *$C$ intersects the $f \geq 2 - z_\diamond$ part of $\mathbb{R} \times M_\delta$ in the union of the radius $\rho$ coordinate balls centered at the index 2 critical points of $f$. In addition, $C$ intersects each subvariety from the $\mathbb{R} \times (1, 1+z_\diamond)$ part of $\mathcal{M}_\Sigma$ precisely once in each $\mathfrak{p} \in \Lambda$ version of $\mathcal{H}_\mathfrak{p}$ and this intersection point lies where $u \in (-R - \ln\rho, 0)$.*

*Proof of Proposition 5.8*: Given $\rho \in (0, \delta_*)$, fix $\delta, x_0, R$ and $J$ so as to define the geometry of $Y$ and $\mathbb{R} \times Y$, and so that $\delta < \frac{1}{100}\rho$. Suppose $C$ is an ech-HF subvariety. Let $\mathfrak{p}$ denote an index 1 critical point of $f$. Suppose that $\varepsilon > 0$ is such that the following is true: Let $\mathfrak{p}$ denote either an index 1 or index 2 critical point of $f$. Then $1 - 3\cos^2\theta$ is greater than $\varepsilon$ on the part of $C$ in the radius $\delta_*$ coordinate ball centered on $\mathfrak{p}$ where the radius is greater that $\rho$. If $\mathfrak{p}$ is an index 1 critical point then $f > 1 + \rho^2\varepsilon$ on $C$'s intersection with this part of the radius $\delta_*$ coordinate ball centered at $\mathfrak{p}$; and if $\mathfrak{p}$ is an index 2 critical point of f, then $f < 2 - \rho^2\varepsilon$ on $C$'s intersection with this part of the radius $\delta_*$ coordinate ball centered at $\mathfrak{p}$. The next paragraph describes some consequences of these observations.

Suppose that $(x, y) \in (1, 1+\rho^2\varepsilon)$ and let $S_{(x,y)} \subset \mathcal{M}_\Sigma$ denote the corresponding element from Proposition 3.2's moduli space. If $y$ is very close to 1, then $C$ has precisely one intersection point with $S_{(x,y)}$ in each $\mathfrak{p} \in \Lambda$ version of $\mathbb{R} \times \mathcal{H}_\mathfrak{p}$ and this intersection point lies where $u \in (0, R + \ln\rho)$. As $y$ is increased, the corresponding intersection point will move in a continuous fashion until it leaves $\mathbb{R} \times \mathcal{H}_\mathfrak{p}$. However, this event is precluded by virtue of what is said in the preceding paragraph. This being the case, it follows that each subvariety from the $\mathbb{R} \times (1, 1+\rho^2\varepsilon)$ part of $\mathcal{M}_\Sigma$ intersect C precisely once in each $\mathfrak{p} \in \Lambda$ version of $\mathcal{H}_\mathfrak{p}$, and that the intersection point lies where $u \leq R + \ln\rho_\Diamond$. The analogous argument shows that C intersects each subvariety from the $\mathbb{R} \times (2-\rho^2\varepsilon, 2)$ part of $\mathbb{R} \times \mathcal{H}_\mathfrak{p}$ precisely once in each $\mathfrak{p} \in \Lambda$ version of $\mathcal{H}_\mathfrak{p}$ and that the latter intersection point occurs where $u > -R - \ln\rho$.

Granted the preceding, let $\kappa$ denote the constant that is supplied by Lemma 5.7 for the given values of $\mathcal{I}$ and $\rho$. It follows from what is said in Lemma 5.7 and in the previous two paragraphs that the conclusions of Proposition 5.8 hold with $z_\Diamond < \rho^2\kappa^{-1}$.

## 6. Heegaard Floer curves

This section constitutes a digression to describe in more detail the sorts of pseudoholomorphic subvarieties that arise in Lipshitz' [L] formulation of Heegaard Floer geometry. What is said here summarizes and in some cases elaborates on observations made by Lipshitz. The discussion in this section may be of independent interest to those using the constructions in [L]. In any event, the results here are used in Section 7 and in [KLTIII].

### a) The Heegaard Floer geometry

This first subsection sets up the geometry that is used by Lipshitz to define his version of pseudoholomorphic curves. This subsection has three parts.

*Part 1*: Introduce $M_{(1,2)}$ as short hand to denote the $f^{-1}(1, 2)$ part of M. Define a closed 2-form $w_f$ on $M_{(1,2)}$ as follows: Set, $w_f = w$ on $M_\delta \cap M_{(1,2)}$. Use (1.13) to define $w_f$ on the radius $2\delta$ coordinate ball centered on any given index 1 critical point of $f$. Use the corresponding $(t_-, \theta_-, \varphi_-)$ version to define $w_f$ on the radius $2\delta$ coordinate ball centered on any given index 2 critical point of $f$. The resulting 2-form is annihilated by the Lie derivative of the pseudogradient vector field $\mathfrak{v}$.

Lie transport by $\mathfrak{v}$ identifes $M_{(1,2)}$ with $(1,2) \times \Sigma$ as in Part 5 of Section 1c. By way of review, this identification writes $f$ as the Euclidean coordinate $t$ on the $(1,2)$ factor, and $\mathfrak{v}$ as the dual vector field $\frac{\partial}{\partial t}$. Meanwhile, the 2-form $w_\Sigma$ appears as the area form on $\Sigma$ that is defined in Part 5 of Section 1c. The identification between $M_{(1,2)}$ and $(1,2) \times \Sigma$ is used implicitly in much of what follows.

An almost complex structure on $\mathbb{R} \times M_{(1,2)}$ is defined by the rules listed momentarily. This almost complex structure is denoted in what follows by $J_{HF}$. To set the notation, reintroduce from Part 3 of Section 1c the function $h_+ = 2e^{2t_+}\cos\theta_+ \sin^2\theta_+$ on the radius $8\delta_*$ coordinate ball centered on a given index 1 critical point of $f$. The functions $(f, \varphi_+, h_+)$ define coordinates on the $1 - 3\cos^2\theta_+ > 0$ part of this coordinate ball.

The list refers to the dual vector fields for these coordinates. There is the analogous function $h_- = 2e^{2t_-}\cos\theta_- \sin^2\theta_-$, and coordinates $(f, \varphi_-, h_-)$ and with the corresponding dual vector fields on the $1 - 3\cos^2\theta_- > 0$ part of the radius $8\delta_*$ coordinate ball centered on a given index 2 critical point of $f$. The list also refers to the subsets $T_+$ and $T_-$ of $\Sigma$ that are introduced in Part 4 of Section 1c.

- $J_{HF}\cdot\partial_s = \mathfrak{v}$.
- $J_{HF}$ *is invariant with respect to constant translations along the $\mathbb{R}$ factor of $\mathbb{R}\times M$.*
- $J_{HF}$ *preserves the kernel of* $df$ *and doing so, it define the orientation given by* $w_\Sigma$.
- *Write* $M_{(1,2)}$ *as* $(1,2)\times\Sigma$. *The almost complex structure* $J_{HF}$ *is invariant with respect to translations along the $(1,2)$ factor of $\mathbb{R}\times(1,2)\times(T_+\cap T_-)$.*
- $J_{HF}\frac{\partial}{\partial\varphi_+} = \frac{\partial}{\partial h_+}$ *on $\mathbb{R}$'s product with the radius $8\delta_*$ coordinate ball centered on any index 1 critical point of $f$; and* $J_{HF}\frac{\partial}{\partial\varphi_-} = -\frac{\partial}{\partial h_-}$ *on $\mathbb{R}$'s product with the radius $8\delta_*$ coordinate ball centered on any index 2 critical point of $f$.*

(6.1)

What follows are comments about these rules. With $M_{(1,2)}$ written as $(1,2)\times\Sigma$, the third bullet requires that $J_{HF}$ preserves the $T\Sigma$ factor of $T(\mathbb{R}\times(1,2)\times\Sigma)$ and that it defines the orientation on $T\Sigma$ given by the area form $w_\Sigma$. Compatibility between the rule in the fifth bullet and the rules in the third and fourth bullets is insured by CONSTRAINT 3 in Section 1c. The final comment concerns the 2-form $ds\wedge df + w_\Sigma$ on $\mathbb{R}\times M_{1,2}$. This 2-form is symplectic, and $J_{HF}$ is a compatible, almost complex structure. The identification of $M_{(1,2)}$ with $(1,2)\times\Sigma$ writes this symplectic form as $ds\wedge dt + w_\Sigma$.

*Part 2*: The identification $M_{(1,2)} = (1,2)\times\Sigma$ identifies the geometry just defined with the geometry that Lipshitz introduces in the Section 1 of [L]. Lipshitz introduces in this section of his paper certain sorts of almost complex structure on $\mathbb{R}\times[1,2]\times\Sigma$, and $J_{HF}$ is an example. Meanwhile, The form $ds\wedge df + w_f$ on $\mathbb{R}\times M_{(1,2)}$ when viewed as a 2-form on the manifold $\mathbb{R}\times(1,2)\times\Sigma$ is what Lipshitz calls a *split* symplectic form.

What follows is a converse of sorts to the preceding observations. Let $J_*$ denote an almost complex structure on $\mathbb{R}\times[1,2]\times\Sigma$ of the sort that Lipshitz considers in Section 1 of [L]. This almost complex structure obeys all but possibly the fifth rule in (6.1). The rule given by the fifth bullet in (6.1) can be imposed without compromising either what is said by Lipshitz [L] or by Ozsváth and Szabó in any of their many papers (see [OS2]). The following says more on this point with regards to the Heegaard Floer framework of Ozsváth and Szabó [OS1]. The assumption made in the fifth bullet is equivalent to assuming that the $[0,1]$-parametrized families of almost complex structures that are considered in Section 3.2 of [OS1] are fixed and independent of the $[0,1]$ coordinate near 0 and near 1. The [0,1] dependence of their almost complex structures is used only to insure that certain genericity conditions hold; and their arguments are valid with the constraint just described near {0} and near {1}. In any event, the rule given by the fifth bullet can be weakened considerably without changing the conclusions of the propositions and lemmas that follow. For example, it is most likely sufficient to assume that $J_{HF}$ when viewed as an almost complex structure on $\mathbb{R} \times (1,2) \times \Sigma$ extends smoothly to a neighborhood of $\mathbb{R} \times \{1\} \times C_+$ and to $\mathbb{R} \times \{2\} \times C_-$. This said, bullet five is kept as is.

The rule (**J2**) in Section 1 of [L] requires somewhat more than what is required by that in the fourth bullet of (6.1). However, the extra conditions in Lipshitz' rule (**J2**) are not required for what is done in his paper; nor are they truly necessary for the work of Osváth and Szabó. In any event, the rule (**J2**) from [L] can be imposed here without compromising anything done in what follows or in [KLTIII].

*Part 3*: Let $J_{HF}$ denote an almost complex structure that obeys the rules in (6.1). Lipshitz introduces in Section 1 of [L] certain $J_{HF}$-holomorphic maps into $\mathbb{R} \times [1,2] \times \Sigma$. The domain of such a map is a Riemann surface, S, with punctured boundary. In particular, there is compact, Riemann surface, S̲ with smooth boundary which is such that S is the complement in S̲ of a set of 2G labeled points in the boundary of S̲. Half of these points are labeled as *negative* and half as *positive*. The map, $u$: S → $\mathbb{R} \times [1,2] \times \Sigma$, is contrained so as to satisfy six conditions, these labeled as (**M0**)-(**M6**) in Section 1 of [L]. These conditions are stated in the upcoming (6.2). The statement reintroduces from Part 2 of Section 1c the sets $C_+$ and $C_-$ in $\Sigma$, these the respective union of $\Sigma$'s intersection with the ascending disks from the index 1 critical points of $f$ and the descending disks from the index 2 critical points of $f$. What follows are, in order, Lipshitz's conditions.

- *The map $u$ is smooth on S and it is $J_{HF}$-holomorphic on the interior of S.*
- *The map sends the boundary of S to the union of $\mathbb{R} \times \{2\} \times C_-$ and $\mathbb{R} \times \{1\} \times C_+$.*
- *There are no components of S in the fiber of the projection to $\Sigma$.*
- *Each component of $(\mathbb{R} \times \{2\} \times C_-) \cup (\mathbb{R} \times \{1\} \times C_+)$ contains the image via $u$ of one and only one boundary component of S.*
- *The pull-back via $u$ of the coordinate s limits to $-\infty$ on sequences in S that converge in S̲ to a negative S̲−S point, and s limits to $\infty$ on sequences in S that converge in S̲ to a positive S̲−S point.*
- *The integral of $u^* w_\Sigma$ over S is finite, and there exists $\kappa_u \geq 1$ with the following property: Let $I \subset \mathbb{R}$ denote an interval of length 1. Then the integral of $u^*(ds \wedge dt)$ over $u^{-1}(I \times [1,2] \times \Sigma)$ is bounded by $\kappa_u$.*
- *The map $u$ is an embedding.*

(6.2)

A pair (S, u) as just described is said in [L] and in what follows to be an *admissible* pair. A *singular admissible set* is a set of the form ((S, u), $\vartheta_\Sigma$) where (S, u) is as described above except for the fact that the last bullet in (6.2) need not be obeyed. Even so, u must be almost everywhere 1-1, and its image has but a finite number of singular points with none on the boundary. Meanwhile, $\vartheta_\Sigma$ consists of a finite set of pairs whose first component is a fiber of the projection from $\mathbb{R} \times (1, 2) \times \Sigma$ to $\Sigma$ and whose second component is a positive integer. However, distinct pairs from $\vartheta_\Sigma$ are distinct fibers. The set $\vartheta_\Sigma$ can be empty.

**b) The points in S near its boundary**

This subsection says more about the behavior of an admissible pair (S, u) near the boundary of S. The results below elaborate on various remarks in [L]. Note that what is said here is also valid if u(S) has a finite number of interior singular points. The discussion that follows has three parts.

*Part 1*: Let p denote an index 1 or index 2 critical point of $f$. The lemma that follows describes the function $u^*s$ on the boundary of S. The lemma uses $\partial S$ to denote the boundary of S.

**Lemma 6.1**: *Let $\partial S \subset S$ denote the boundary. The coordinate function s from the $\mathbb{R}$ factor of $\mathbb{R} \times [1, 2] \times \Sigma$ pulls back via u to $\partial S$ as a proper function with no critical points.*

*Proof of Lemma 6.1*: The fact that the $u^*ds$ can not vanish on $T(\partial S)$ owes allegiance to the fact that u embeds $\partial S$ into a level sets of t. The argument that follows explains why this is in the case when the boundary component in question is labeled by an index 1 critical point of $f$. Except for notation, the same argument works for the components with index 2 critical point labels.

To start, let p denote a given index 1 critical point and let $\partial_p S$ denote the corresponding component of the $\partial S$. The map u sends a neighborhood of $\partial_p S$ in S into a neighborhood of $\mathbb{R} \times \{1\} \times C_{p+}$ in $\mathbb{R} \times [1, 2] \times \Sigma$. The latter neighborhood is chosen to have the form $\mathbb{R} \times [1, 1+z_p) \times T_+$ with $z_p > 0$. Given that $z_p$ is small, the functions $(\varphi_+, h_+)$ from Part 3 of Section 1c define coordinates for $T_+$ and so the coordinate $s$ for $\mathbb{R}$, $t$ for $[1, 1+z_p)$ and the pair $(\varphi_+, h_+)$ supply coordinates for this neighborhood. The almost complex structure endomorphism $J_{HF}$ maps $\frac{\partial}{\partial s}$ to $\frac{\partial}{\partial t}$, and by dint of the fifth bullet of (6.1), it maps $\frac{\partial}{\partial \varphi_+}$ to $\frac{\partial}{\partial h_+}$.

Suppose that $u^*ds$ vanishes at some point on $\partial_p S$ so as to derive some nonsense. Given that $u$ is $J_{HF}$-holomorphic, the latter condition requires that $u^*dt$ is also zero at this point. As a consequence, the function $s + it$ pulls back to S with a critical point at this boundary point. The ramifications of this fact are explained momentarily. To set the stage, fix a holomorphic identification between a neighborhood of this boundary point and the half disk $D_+ \subset \mathbb{C}$ whose points have norm less than 1 and nonnegative imaginary part. Use $\zeta$ to denote the complex coordinate on $\mathbb{C}$ and thus on $D_+$.

The pull-back via $u$ of $s + it$ is holomorphic on the interior of S. As explained to the author by Curt McMullen, a version of the Schwarz reflection trick from Theorem 24 in [Ah] can be used to write the pull-back of $s + it$ on $D_+$ as

$$u^*(s + it) = a_0 \zeta^n + \mathcal{O}(|\zeta|^{n+1})$$

(6.3)

with $a_0 \in \mathbb{C}-0$ and $n > 1$. However, (6.3) implies the $t = 1$ level set in S intersects the interior of S, which is assuredly not the case.

*Part 2*: It follows from what is said in Part 1 that the pull-back via $u$ of the coordinates functions $s$ and $t$ on $\mathbb{R} \times [1, 2] \times \Sigma$ give coordinates on a neighborhood of $\partial S$. This part says more about the map $u$ on such a neighborhood.

To this end, let p again denote an index one critical point of $f$. The boundary component $\partial_p S$ can be parametrized by the pull-back of the function $s$. In particular, the image via $u$ of $\partial_p S$ can be parametrized by $\mathbb{R}$ as a graph

$$s \to (s, t = 1, \varphi^S(s))$$

(6.4)

where $\varphi^S: \mathbb{R} \to C_{p+}$ is a smooth map.

**Lemma 6.2**: *The closure of the corresponding boundary component $\partial_p S$ in $\underline{S}$ adds precisely one negative and one positive point of $\underline{S}-S$. Moreover, the map $\varphi^S$ extends continuously to this closure so as to map these points of $\underline{S}-S$ to $C_+ \cap C_-$.*

*Proof of Lemma 6.2*: Let $\zeta$ again denote the standard complex coordinate for $\mathbb{C}$ and $D_+ \subset \mathbb{C}$ denote the half disk where $|\zeta| < 1$ and $\text{im}(\zeta) \geq 0$. Let $q_S \in \underline{S}$ denote a negative boundary point. There is a neighborhood of $q_S$ in $\underline{S}$ whose closure intersects $\underline{S}-S$ only in $q_S$, and is diffeomorphic to the half disk $D_+$ with q corresponding to the origin. This diffeomorphism can be chosen so that the complex structure on $\underline{S}$ near $q_S$ maps to the standard complex structure on $\mathbb{C}$. It proves useful to employ the map $\zeta \to \frac{1}{\pi} \ln(\zeta) + i$ so as to holomorphically identify the $|\zeta| \le e^{-1}$ part of $D_+ - 0$ with the strip $(-\infty, -1] \times [1, 2]$. The boundary component $(-\infty, -1] \times \{1\}$ is mapped by $u$ to a component of $\mathbb{R} \times \{1\} \times C_+$. Let p denote the index 1 critical point of $u$ that labels this component. Meanwhile, the boundary component, $(-\infty, -1] \times \{2\}$, is mapped by $u$ to a component of $\mathbb{R} \times \{1\} \times C_+$ or of $\mathbb{R} \times \{2\} \times C_-$.

Given that the integral of $w_\Sigma$ over S is finite, it follows using Lemma 5.7 that the constant s slices of the images via $u$ of $(-\infty, 0) \times [1, 2]$ converge in $[1, 2] \times \Sigma$ to a curve of the form $[1, 2] \times q$ with q some point in $\Sigma$. The $t = 1$ end of the $u$ image of each such slice must lie in $C_{p+}$ and the $t = 2$ end must therefore lie in $C_-$. Let p´ denote the index 2 critical point of $f$ that labels this component. The $s \to -\infty$ convergence of the image of these slices requires that $q \in C_{p+} \cap C_{p'-}$. This convergence of the slices with what was just said about their limit implies what is asserted for the lemma about the $s \to -\infty$ behavior of $\partial S$ and the map $\varphi^S$. Very much the same argument implies what is said about the $s \to \infty$ behavior of $\partial S$ and $\varphi^S$.

*Part 3*: Let p denote an index 1 or index 2 critical point of $f$. The next lemma says something about the manner in which the restriction of $u$ to the constant $u^*(t)$ slices of S near $\partial_p S$ approximate the map $\varphi^S$. When p is an index 1 critical point, the lemma refers to the annular neighborhood $T_{p+}$ of $C_{p+}$ with its coordinates $(\varphi_+, h_+)$. When p is an index 2 critical point, it refers to $T_{p-}$ and the corresponding coordinates $(\varphi_-, h_-)$. The map $\varphi^S$ has domain $\mathbb{R}$ and so can be viewed as an $\mathbb{R}$ valued map. The lemma views it as such.

**Lemma 6.3**: *Let p denote an index 1 or index 2 critical point of the function $f$. There exists $z_p \in (0, 1)$ and a neighborhood of $\partial_p S$ in S with the following properties: The composition of first $u$ and then the projection from $\mathbb{R} \times [1, 2] \times \Sigma$ to $\Sigma$ maps this neighborhood into $T_{p+}$ or $T_{p-}$ as the case may be. Moreover, the neighborhood has a parametrization by $\mathbb{R} \times [0, z_p)$ whereby $u$ appears as*

- $(s, z) \to (s, t = 1 + z, \varphi_+ = \varphi(s, z), h_+ = \varsigma(s, z))$ *when p has index 1*,
- $(s, z) \to (s, t = 2 - z, \varphi_- = \varphi(s, z), h_- = \varsigma(s, z))$ *when p has index 2*,

*where $\varphi(\cdot)$ and $\varsigma(\cdot)$ are maps from $\mathbb{R} \times [0, z_p)$ to $\mathbb{R}$ that obey*

$$|\varsigma(s, z)| + z^{-1}|\varphi(s, z) - \varphi^S(s)| < \kappa z$$

*with $\kappa$ a positive constant. In addition,*

$$|\partial_s \varsigma(s, z)| + z^{-1}|(\partial_s \varphi)(s, z) - (\partial_s \varphi^S)(s))| < \kappa z .$$

*In general, the absolute values of the derivatives of $\varphi(\cdot)$, $\varsigma(\cdot)$ and $\varphi^S(\cdot)$ to any given specified order enjoy an s and z independent upper bound on $\mathbb{R} \times [0, z_p)$.*

*Proof of Lemma 6.3*: Fix $T \geq 1$ and Lemma 6.2 with fact that $u$ is J-holomorphic implies that the lemma's assertions hold for a neighborhood of the part of $\partial_p S$ where $|u^*s| \leq T$. This understood, it is sufficient to prove that the constants $\kappa$ and $z_p$ can be chosen to be independent of T. The argument that follows takes p to be an index 1 critical point of $f$. But for cosmetics, the same argument works for the index 2 critical points.

To start, let $q^S \in C^p_+ \cap C_-$ denote the $s \to -\infty$ limit on $\partial_p S$. View $q^S$ as a point in $\mathbb{R}/(2\pi\mathbb{Z})$ using the coordinate $\varphi_+$. As noted in the Part 1 above, the constant $s$ slices of S that intersect $u(\partial_p S)$ converge as $s \to \infty$ to the arc $[1, 2] \times q^S$ in $[1, 2] \times \Sigma$. As a consequence, given $\varepsilon > 0$, there exists $T_\varepsilon$ such that the part of $u(S)$ with $s < -T_\varepsilon$ that intersects $u(\partial_p S)$ has distance $\varepsilon$ or less from this arc. In particular, it follows from (6.2) that if $\varepsilon < c_0^{-1}$, then this part of $u(S)$ lies in $[1, 2] \times T_{p+}$; and in particular, it lies in the part of $\mathbb{R} \times [1, 2] \times T_p$ where $|\varphi_+ - q^S| + |\hat{h}_+| < c_0\varepsilon$.

It follows from these last observations and their $s \to \infty$ analogs that a version of $z_p$ can be chosen so as to satisfy the demands of the lemma except possibly for the uniform bounds on $\varphi$ and $\varsigma$. The following says precisely what this means: Given $\varepsilon > 0$, there exists $z_{p,\varepsilon}$ and a neighborhood of $\partial_p S$ with two properties: First, the composition of $u$ and then projection to $\Sigma$ maps this neighborhood into $T_{p+}$. Second, the neighborhood has a parametrization by $\mathbb{R} \times [0, z_p)$ whereby $u$ appears as the indicated graph with $\varphi(\cdot)$ and $\varsigma(\cdot)$ smooth and such that $|\varphi(s, z) - \varphi^S(s)| + |\varsigma(s, z)| < \varepsilon$.

The map $(s, z) \to \zeta(s, z) = \varphi(s, z) + i\varsigma(s, z)$ is holomorphic on $\mathbb{R} \times (0, z_{p,\varepsilon})$ because $u$ is $J_{HF}$-holomorphic. Fix a Euclidean half disk $D_+ \subset \mathbb{R} \times [0, z_{p,\varepsilon})$ of radius $\frac{1}{8} z_{p,\varepsilon}$ with center on any given $(s_0, 0)$ point. The Schwarz reflection argument from Theorem 24 in [Ah] can be repeated here to extend $\zeta$ as a holomorphic function on the disk in $\mathbb{C}$ where the complex coordinate $s + iz$ is such that $|s - s_0|^2 + |z|^2 \leq \frac{1}{64} z_{p,\varepsilon}^2$. The boundary values on the full disk of this extended function $s + iz \to \zeta(s + iz)$ are uniformly bounded. It follows as a consequence of the Cauchy integral formula that its derivatives to any given order are uniformly bounded on the concentric disk with radius $\frac{1}{16} z_{p,\varepsilon}$ with this bound independent of $s_0$. This implies in particular that

$$|\varsigma(s, z)| + z^{-1} |\varphi(s, z) - \varphi^S(s)| \leq c_0 |z|$$

(6.5)

in this smaller disk.  This is the inequality asserted by the lemma for the norms of φ - φ$^S$ and ς.  A similar argument gives the asserted bound for the norms of the respective s-deriviatives; and likewise the asserted bound on the norms of the derivatives to any given order.

**c) Neighborhoods of S̲–S**

This subsection says more about the map $u$ near a point in S̲–S. The discussion that follows talks about a neighborhood of a given negative point. An almost identical story can be told about neighborhoods of the positive points. The discussion here has two parts.

*Part 1*: Let $q_S$ denote a given negative point of S̲–S, and let p and p´ denote the respective index 1 and index 2 critical points of $f$ that label the nearby components of ∂S. Parametrize the S̲–$q_S$ part of a neighborhood of $q_S$ by $(-\infty, 0) \times [1, 2]$ as in the proof of Lemma 6.2. As noted in this same proof, the constant $s$ slices of the $u$ image of domain converge in $[1, 2] \times \Sigma$ to a curve of the form $[1, 2] \times q$ with $q \in C_{p+} \cap C_{p'-}$. This fact has the following consequence:  There exists $s_0 > 1$ such that the that image via $u$ of the complement of some compact set in $(-\infty, -1] \times [1, 2]$ has intersection number 1 with each constant $(s, t)$ slice of $(-\infty, s_0] \times [1,2] \times \Sigma$. What with Lemma 4.7, the preceding observation leads directly to the following assertion: The composition of first the map $u$ and then projection to $(-\infty, 2s_0] \times [1, 2]$ restricts to its inverse image in $(-\infty, 0) \times [1, 2]$ as a diffeomorphism. Moreover, this diffeomorphism is complex analytic by virtue of the fact that $J_{HF}$ maps $\frac{\partial}{\partial s}$ to $\frac{\partial}{\partial t}$. This understood, there exists $s_1 \geq 2s_0$ such that the map $u$ when written using these coordinates can be written as a graph over $(-\infty, s_1] \times [1,2]$ using a map, ψ, from $(-\infty, s_1] \times [1, 2]$ into $T_{p+} \cap T_{p'-}$. This is to say that $q_S$ has a neighborhood in S̲ whose intersection with S can be parametrized by $(-\infty, s_1] \times [1, 2]$ such that the map $u$ appears as the map

$$(s, t) \to (s, t, \psi(s,t)).$$

(6.6)

The map ψ is smooth; it maps $(-\infty, s_1] \times \{1\}$ to $C_{p+}$ and it maps $(-\infty, s_1] \times \{2\}$ to $C_{p-}$; and it is such that $\lim_{s \to \infty} \sup_{t \in [1,2]} \text{dist}(\psi(s,t), q) = 0$.

*Part 2*: It proves useful to holomorphically identify the domain of ψ to the complement of 0 in the $|\zeta| \leq e^{-\frac{\pi}{2} s_1}$ part of the first quadrant of $\mathbb{C}$, the quadrant where both

re($\zeta$) ≥ 0 and im($\zeta$) ≥ 0. This part of $\mathbb{C}$ is denoted by $D_{++}$. This identification of the domain of ψ is done using the map

$$(s, t) \to \zeta = e^{\frac{\pi}{2}(s+i(t-1))}.$$

(6.7)

Use U to denote q's component of $T_{p+} \cap T_{p'-}$. The coordinates ($\varphi_+$, $h_+$) can be used as Euclidean coordinates on U; and no generality lost by assuming that q has coordinate $\varphi_+ = 0$. This parametrization is defined for $|\varphi_+| + |h_+| \le c_0^{-1}$. It follows from what is said in Part 3 of Section 1b that the respective intersections of $C_{p+}$ and $C_{p'-}$ with U appear in terms of the coordinates ($\varphi_+$, $h_+$) as the $h_+ = 0$ axis, and the $\varphi_+ = 0$ axis  Meanwhile, the first and fourth bullets in (6.1) say that $J_{HF}$ on $\mathbb{R} \times [1, 2] \times U$ acts to send $\frac{\partial}{\partial s}$ to $\frac{\partial}{\partial t}$ and $\frac{\partial}{\partial \varphi_+}$ to $\frac{\partial}{\partial h_+}$.

Using this notation, the map ψ takes the form

$$\zeta \to (\varphi_+ = \varphi(\zeta), h_+ = \varsigma(\zeta))$$

(6.8)

where φ and ς are $\mathbb{R}$-valued functions. What was said above about $J_{HF}$ implies that the $\mathbb{C}$-valued function $\lambda = \varphi + i\varsigma$ is holomorphic on $D_{++}$–0. What was said subsequent to (6.6) about ψ says the following about λ: First, $\lim_{\zeta \to 0} \lambda = 0$. Second, the imaginary part of λ is zero on the im($\zeta$) = 0 boundary of $D_{++}$ and the real part of λ is zero on the re($\zeta$) = 0 boundary $D_{++}$.

With the preceding understood, a four-fold application of the Schwarz reflection trick from Theorem 24 in [Ah] extends the domain of the map λ from $D_{++}$ as a holomorphic map whose domain is the complement of the origin in the full $|\zeta| \le e^{-\frac{\pi}{2} s_1}$ disk in $\mathbb{C}$. This extended version of λ also limits to zero as $\zeta$ limits to 0.

It follows from what was just said that λ and thus ψ extend across the origin of $\mathbb{C}$ so as to define a holomorphic map from a disk about the origin in $\mathbb{C}$ to D with the origin going to q.  Such a map has the form

$$\zeta \to w(\zeta) = r\,\zeta^n + \mathcal{O}(\zeta^{n+1})$$

(6.9)

where n ∈ {1, 2, …} and where $r \in \mathbb{C}$–0. This last observation with (6.7) and (6.8) supply the desired picture of *u* near the given point in <u>S</u>–S.

d) **The $\mathbb{R} \times M_{(1,2)}$ image of $u$.**

Of interest in Section 7 and in [KLTIII] is the $u$ image of the interior of S with the latter viewed as a subset of $\mathbb{R} \times M_{(1,2)}$. The $u$ image of the part near $\partial S$ is of particular concern in this subsection, this being the concern of the upcoming Lemma 6.4.

To set notation, introduce $S_0$ to denote $S - \partial S$. The composition of first $u$ and then projection from $\mathbb{R} \times M_{(1,2)}$ to $M_{(1,2)}$ defines a map from $S_0$ to $M_{(1,2)}$ that is denoted by $u_M$. Lemma 6.4 refers to the spherical coordinate $\theta = \theta_+$ or $\theta = \theta_-$ that is introduced in Section 1a for the radius $\delta_*$ coordinate ball centered at a given critical point of $f$ with index respectively 1 or 2.

**Lemma 6.4**: *Let $(S, u)$ denote an admissible pair with it understood that $u(S)$ can have interior singular points. View $u(S_0)$ as a $J_{HF}$-holomorphic subvariety in $\mathbb{R} \times M_{(1,2)}$. There exists $\kappa \geq \delta_*^{-1}$ with the following significance: Fix $z \in (0, \kappa^{-4})$.*
- *The function $1 - 3\cos^2\theta$ is bounded by $\kappa^{-1}$ on the part of $u_M(S_0)$ that lies in the radius $\delta_*$ coordinate ball centered on any given index 1 or index 2 critical point of $f$.*
- *The part of $u_M(S_0)$ where $f < 1 + z$ lies in the union of the radius $\kappa z^{1/2}$ coordinate balls centered on the index 1 critical points of $f$. The part where $f > 2 - z$ lies in the union of the radius $\kappa z^{1/2}$ coordinate balls centered on the index 2 critical points of $f$.*

*Proof of Lemma 6.4*: Fix an index 1 or index 2 critical point of $f$ and introduce the function $h$ as defined in (4.1) for the corresponding version of $\mathcal{H}_p$. Lemma 6.3 bounds $h/f = \frac{\cos\theta \sin^2\theta}{1 - \cos^2\theta}$ on the part of $u_M(S_0)$ that lies in the radius $\delta_*$ coordinate ball centered at p. The existence of such a bound implies what is asserted by the first bullet of the lemma. The bound in the first bullet implies the second bullet.

e) **Linear operators and admissible pairs**

The purpose of this subsection is to summarize some of what is said in [L] about the operator that defines the first order deformations of an admissible pair. The discussion here has six parts. In what follows, $(S, u)$ is a given admissible pair.

*Part 1*: The image via $u$ of $S_0$ is a J-holomorphic submanifold of $\mathbb{R} \times (1, 2) \times \Sigma$ and so its normal bundle inherits from J a complex line bundle structure. Use $N_S \to S_0$ to denote the pull-back of this normal bundle. The map $u$ extends as an embedding of the closed surface S into $\mathbb{R} \times [1, 2] \times \Sigma$. What follows describes the extension of $N_S$ to $\partial S$.

To start, let p denote an index 1 critical point of $f$, and let $\partial_p S$ denote the corresponding boundary component. Lemma 6.3 describes the map $u$ on a certain parametrized neighborhood in S of $\partial_p S$. With respect to this parametrization, the normal bundle appears as the pull-back via $u$ of the restriction to the neighborhood $T_{p+}$ of $T\Sigma$. This bundle has its basis of sections $\{\frac{\partial}{\partial \varphi_+}, \frac{\partial}{\partial h_+}\}$ and the bundle complex structure maps the section $\frac{\partial}{\partial \varphi_+}$ to $\frac{\partial}{\partial h_+}$. There is the analogous basis and description of the complex line bundle structure for $N_S$ on a neighborhood of an index 2 critical point version of $\partial_p S$.

Note that the trivializations given above for respective index 1 and index 2 critical point versions of $\partial_p S$ are compatible on a small neighborhood of a given point in $\underline{S}$−S, this being a consequence of what is said in the fourth bullet of (6.1).

Given that $u$ is an embedding, it maps suitable neighborhoods of distinct negative points in $\underline{S}$−S to distinct points in $T_- \cap T_+$, and it likewise maps neighborhoods of distinct positive points in $\underline{S}$−S to distinct points in $T_- \cap T_+$. These facts about $u$ have an important consequence that is described momentarily. The statement of this consequence refers to a certain hermitian fiber metric on $N_S$. This fiber metric is induced by the Riemannian metric on $\mathbb{R} \times [1, 2] \times \Sigma$ that is defined by the symplectic form $ds \wedge dt + w_\Sigma$ and its compatible almost complex structure $J_{HF}$.

Here is the promised consequence: There exists $\rho_S > 0$ and an *exponential* map $\mathfrak{e}_S$ from $N_S$ to $\mathbb{R} \times [1, 2] \times \Sigma$ that embeds the radius $\rho_S$ disk subbundle as a neighborhood of $u(S)$ that contains the set of points in $\mathbb{R} \times [1, 2] \times \Sigma$ with distance $c_0^{-1}\rho_S$ or less from S. In addition, this number $\rho_S$ and the map $\mathfrak{e}_S$ can be chosen so that the image of each fiber disk of radius $\rho_S$ is a $J_{HF}$-holomorphic disk in $\mathbb{R} \times [1, 2] \times \Sigma$. Such a choice proves useful in subsequent sections. A version of $\mathfrak{e}_S$ with these properties can be constructed by mimicking what is done in Lemma 5.4 of [T3]. Use $N_0 \subset N_S$ to denote this radius $\rho_S$ disk subbundle.

*Part 2*: Suppose that $u'$: $S \to \mathbb{R} \times [1, 2] \times \Sigma$ is a smooth map that obeys Lipshitz' conditions (**M1**)-(**M6**) in [L], which is to say it is described by the bullets in (6.2) except that it need not be $J_{HF}$-holomorphic. If $\sup_{y \in S} |u(y) - u'(y)| < c_0^{-1}\rho_S$, then the image of $u'$ will lie in $\mathfrak{e}_S(N_0)$. If $u'$ also has algebraic and geometric intersection 1 with each fiber of $N_S$, then $u'$ factors as $\mathfrak{e}_S \circ \eta$ where $\eta$ is a smooth section of $N_S$. Assuming that this is the case, the surface $u'(S_0)$ will be $J_{HF}$-holomorphic if and only if $\eta$ obeys an equation that has the schematic form

$$\bar{\partial} \eta + \mathfrak{r}_1(\eta) \cdot \partial \eta + \mathfrak{r}_0(\eta) = 0$$

(6.10)

where the notation is as follows: First, $\bar{\partial}$ signifies the d-bar operator on sections of $N_S$ as defined using the hermitian metric to give the bundle a holomorphic structure. What is written as $\partial$ is the adjoint operator. Meanwhile, the map $\mathfrak{r}_1\colon N_0 \to N_S \otimes \mathrm{Hom}(T^{1,0}S; T^{0,1}S)$ and the map $\mathfrak{r}_0\colon N_0 \to N_S \otimes T^{0,1}S$ are smooth, fiber preserving maps that vanish along the zero section. If $\eta$ obeys (6.10), then the pair $(S, u')$ is admissible if it is understood that S has the complex structure that is induced via $J_{HF}$ by the embedding $u'$.

As explained next, the non-linear differential operator that is depicted on the left hand side of (6.10) supplies an $\mathbb{R}$-linear operator from $C^\infty(S; N_S)$ to $C^\infty(S; N_S \otimes T^{0,1}S)$. Write the left hand side of (6.10) as a map, $\mathcal{F}\colon C^\infty(S; N_0) \to C^\infty(S; N_S \otimes T^{0,1}S)$. Fix $\eta \in C^\infty(S; N_S)$. The operator in question sends $\eta$ to $D_S \eta = \lim_{\varepsilon \to 0} \varepsilon^{-1} \mathcal{F}(\varepsilon \eta)$. This operator has the form

$$D_S \eta = \bar{\partial} \eta + \nu \eta + \mu \bar{\eta}$$

(6.11)

where $\nu$ is a smooth, bounded section of $C^\infty(S; T^{0,1}S)$ and $\mu$ is a smooth, bounded section of $C^\infty(S; N_S^2 \otimes T^{0,1}S)$.

*Part 3*: The operator $D_S$ has a Fredholm version that comes by suitably choosing its domain. This domain is a certain completion of a linear subspace of $C^\infty(S; N_S)$. The latter subspace is denoted by $C^\infty_\partial(S; N_S)$. Sections in $C^\infty_\partial(S; N_S)$ are constrained on $\partial S$ and in their behavior as points in $\underline{S}-S$ are approached. The constraints that specify membership in $C^\infty_\partial(S; N_S)$ are given momentarily. The statement of the constraints refers implicitly to the identification from Part 1 of $N_S$ along an index 1 critical point version of $\partial_p S$ with the pullback via $u$ of $T\Sigma|_{C^p_+}$. Note in particular that $TC^p_+$ sits in this bundle as a real 1-dimensional subbundle. Likewise, $N_S$ along an index 2 critical point version of $\partial_p S$ is identified with the pullback via $u$ of $T\Sigma|_{C^p_-}$; and $TC^p_-$ is likewise viewed as a real line in the latter.

Here are the constraints: A section $\eta \in C^\infty(S; N_S)$ sits in the subspace $C^\infty_\partial(S; N_S)$ if and only if

- $\eta$ *has compact support*.
- $\eta$ *along an index 1 critical point version of $\partial_p S$ sits in $u^* TC^p_+$*.
- $\eta$ *along an index 2 critical point version of $\partial_p S$ sits in $u^* TC^p_-$*.

(6.12)

The completion of $C^\infty_\partial(S; N_S)$ that defines the domain of S uses the Sobolev $L^2_1$ inner product. Likewise, the image Hilbert space uses the $L^2$ inner product on the space of compactly supported sections of $N_S \otimes T^{0,1}S$. To elaborate, use $\nabla$ to denote the covariant derivative on sections of $N_S$, its tensor powers, and tensor products with powers of T*S that is induced by the Riemanian metric on $\mathbb{R} \times [1, 2] \times \Sigma$. The respective norms that define the domain Hilbert space and image Hilbert space for $D_S$ are the square roots of the functions on $C^\infty_\partial(N; S)$ and $C^\infty(N \otimes T^{0,1}S)$ that assign to a given section $\eta$ the value

$$\int_S (|\nabla \eta|^2 + |\eta|^2) \quad and \quad \int_S |\eta|^2 \;.$$

(6.13)

The index of this Fredholm version of $D_S$ is computed by Lipshitz [L]. Note in this regard that Lipshitz uses weighted versions of these norms. The upcoming Part 4 explains why the weights are not necessary.

An admissible pair $(S, u)$ is said in what follows to be *unobstructed* when the corresponding Fredholm version of $D_S$ has trivial cokernel. Lipshitz proves in Section 3 of [L] that there is a $C^\infty$-residual set of allowed complex structures on $\mathbb{R} \times [1, 2] \times \Sigma$ with only unobstructed admissible pairs.

*Part 4*: A bounded linear operator between Banach spaces is Fredholm if it has closed range, finite dimensional kernel and finite dimensional cokernel. Let $\mathbb{H}_S$ and $\mathbb{L}_S$ denote the respective domain and image Hilbert spaces for $D_S$. The fact that $D_S$ has closed range in $\mathbb{L}_S$ and also finite dimensional kernel follows from two facts that are asserted and then proved momentarily. The notation used has $\|\cdot\|$ denoting the $L^2$ norm whose square is depicted in (6.13). What is denoted by $\eta$ refers to an element in $C^\infty_\partial(N; S)$. What is denoted by $c$ is a constant that depends on $S$ and is greater than 1. What follows are the two promised facts.

- $\|D_S \eta\|^2 \geq c^{-1} \|d\eta\|^2 - c \|\eta\|^2$.
- *There exists $s_* > 1$ such that if $\eta$ has support where $|s| > s_*$, then $\|D_S \eta\|^2 \geq c^{-1} \|\eta\|^2$*.

(6.14)

A standard argument using the Rellich lemma uses (6.14) to deduce that $D_S$ has closed range and finite dimensional kernel. Given that the range is closed, the cokernel is isomorphic to the kernel of the adjoint of $D_S$, this a bounded operator that maps $\mathbb{L}_S$ to $\mathbb{H}_S$. Integration by parts and some standard arguments about elliptic regularity identify the kernel of the latter operator with that of a certain bounded, differential operator that also maps $\mathbb{H}_S$ to $\mathbb{L}_S$. The operator in question is denoted by $D_S^\#$; as explained below, what is asserted by (6.14) is true for $D_S^\#$; and thus its kernel is also finite dimensional.

Turn now to (6.14). The first bullet is derived using a standard integration by parts. The conditions from the second and third bullets of (6.12) insure the vanishing of boundary terms. To see about the second bullet in (6.14), parametrize the map $u$ on a neighborhood of a given negative point in $\underline{S}$–$S$ as described in Part 1 of Section 1c. By way of a reminder, the parametrization is by $(-\infty, s_1] \times [1, 2]$ and is such that the map $u$ appears as a graph $(s, t) \to (s, t, \varphi_+ = \varphi(s, t), h_+ = \varsigma(s, t))$ where $(\varphi, \varsigma)$ are $\mathbb{R}$-valued functions that map the domain into the relevant component of $T_+ \cap T_-$.

The coordinate $\varphi_+$ can be assumed to be $\mathbb{R}$ valued on the image of the pair $(\varphi, \varsigma)$, and such that the point $(\varphi_+ = 0, h_+ = 0)$ is the nearby point in $C_+ \cap C_-$. The locus $C_+$ appears as the $h_+ = 0$ axis and $C_-$ as the $\varphi_+ = 0$ axis. This understood, the pair $(\varphi, \varsigma)$ is such that $\varphi + i\varsigma$ is a holomorphic function of $s + it$. In addition, this pair is constrained so that $\lim_{s \to \infty} (|\varphi| + |\varsigma|) = 0$, so that $\varsigma = 0$ where $t = 1$, and so that $\varphi = 0$ where $t = 2$.

A section of the normal bundle $N_S$ over this part of S pulls back via this parametrization as a map from $(-\infty, s_1] \times [1, 2]$ to $\mathbb{C}$ with it understood that the $\mathbb{C}$-valued 1-form $d\varphi_+ + idh_+$ is used to identify the normal bundle with the product complex line bundle. Write a such a map $\eta$ as $\eta_1 + i\eta_2$ with $\eta_1$ and $\eta_2$ both real valued functions. Then the map $\eta$ comes from $C^\infty(S; N)$ if $\eta_2 = 0$ on $(-\infty, s_1] \times \{1\}$ and $\eta_1 = 0$ on $(-\infty, s_1] \times [1, 2]$. Meanwhile, the operator $D_S$ here is the $\frac{1}{2}(\partial_s + i\partial_t)$ and so

$$D_S\eta = \tfrac{1}{2}(\partial_s\eta_1 - \partial_t\eta_2) + \tfrac{i}{2}(\partial_s\eta_2 + \partial_t\eta_1).$$

(6.15)

As a consequence, if $\eta$ has support where $s > s_1$, then

$$\|D_S\eta\|^2 = \tfrac{1}{4}(\|\partial_s\eta_1\|^2 + \|\partial_s\eta_2\|^2 + \|\partial_t\eta_1\|^2 + \|\partial_t\eta_2\|^2).$$

(6.16)

With (6.16) understood, introduce $C^\infty_*$ to denote the subspace of smooth, $\mathbb{R}$-valued functions on $[1, 2]$ that vanish at one of the end points of $[1, 2]$. The function

$$\zeta \to Q_*(\zeta) = \int_1^2 |\partial_t \zeta|^2 \, dt$$

(6.17)

on $C^\infty_*$ is bounded away from zero in the sense that

$$Q_*(\zeta) > \tfrac{\pi^2}{4} \int_1^2 |\zeta|^2 \, dt.$$

(6.18)

This fact with (6.16) and its analog for the positive points in $\underline{S}-S$ imply what is asserted by the second bullet in (6.14).

The operator $D^\#$ is the formal, $L^2$ adjoint of $D_S$; and so it is a first order operator with leading derivative term given by the formal $L^2$ adjoint of $\bar{\partial}$. This the case, an integration by parts can be employed to prove the $D^{\#}$ version of the first bullet in (6.14). The proof of the second bullet is just like that for $D_S$.

*Part 5*: Let $\mathcal{A}_{HF}$ denote the space of admissible pairs (S, $u$) where the topology is defined as follows: The open neighborhoods of a pair (S, $u$) are generated by a basis whose sets are labeled a positive real number, $\varepsilon$. An admissible pair (S´, $u$´) is in the corresponding set when two conditions are met:

- $\sup_{z \in S}(\text{dist}(u(z), u´(S´))) + \sup_{z \in S´} \text{dist}(u(S), u´(z)) < \varepsilon$.
- $|\int_S u^* \mu - \int_{S´} u´^* \mu| < \varepsilon$ *if $\mu$ is any 2-form on* $\mathbb{R} \times [1, 2] \times \Sigma$ *with* $|\mu| \leq 1, |\nabla \mu| \leq \frac{1}{\varepsilon}$ *and compact support on* $|s| < \frac{1}{\varepsilon}$.

(6.19)

Suppose that (S, $u$) is an admissible pair. What follows is a direct consequence of the definition in (6.19). The pair (S, $u$) has a neighborhood in $\mathcal{A}_{HF}$ with the following property: If (S´, $u$´) is in this neighborhood, then $u´(S´)$ is in the radius $\rho_S$ tubular neighborhood of $u(S)$ and it has geometric intersection number 1 with each fiber. In particular, it can be written as the image of S via the map $\mathfrak{e}_S \circ \eta$ where $\eta$ is a section of $N_0$ that obeys (6.10) and the second two bullets of (6.12). Moreover, $|\eta|$ has limit zero as $|u^*s| \to \infty$. The next lemma is a consequence of these last observations.

**Lemma 6.5**: *Let* (S, $u$) *denote an admissible pair. There is a neighborhood of* (S, $u$) *in the space* $\mathcal{A}_{HF}$ *and, for each* $k \in \{0, 1, 2, \ldots\}$, *a constant* $\kappa_k$ *with the following significance: Suppose that* (S´, $u$´) *is in this neighborhood. Then* $u´(S´)$ *can be written as the image of* $\mathfrak{e}_S \circ \eta$ *where* $\eta$ *is a section of* $N_0$ *that obeys (6.10), the second two bullets in (6.12), and is such that* $\lim_{|s| \to \infty} |\eta| = 0$. *In addition,*

$$|\nabla^k \eta| \leq \kappa_k \, e^{-|u^*s|/\kappa_0} \sup_z |\eta| \, .$$

*Proof of Lemma 6.5*: Without the exponential factor, a bound given by the lemma for the derivatives of $\eta$ follows directly using the sorts of elliptic regularity theorems that can be found in Chapter 6 of Morrey's book [Mo]. Note in this regard that it follows from what is said in Section 6b and in Part 2 above that (6.10) near the boundary of S can be written so as to appear as the standard, linear Cauchy-Riemann equations. This being the case, regularity near the boundary can be proved using classical linear techniques.

The bound with the exponential factor can be proved using what is asserted in the second bullet of (6.14). Here is a brief sketch of how this is done: Fix $\varepsilon > 0$ and $T \gg 1$. Let $\chi: \mathbb{R} \to [0, 1]$ denote a smooth function which is equal to 1 where $s < 0$ and equal to zero where $s > 1$. Take the square of the norm of the expression on the left hand side of (6.10), multiply by $e^{-\varepsilon|u*s|}\chi(T - |u*s|)$ and integrate over S. Integrate by parts using the boundary conditions given by the second bullet in (6.12). Use the fact that $|\eta|$ limits to zero as $|u*s| \to \infty$ with the second bullet in (6.12) to see that the $\varepsilon \to 0$ limit of the resulting integral can be taken, and that this limit implies the bound

$$\int_{|u*s|>T} |\eta|^2 \leq c_0 \int_{|u*s|\in[T-1,T]} |\eta|^2 .$$

(6.20)

The sequence $\{T = n\}_{n=1,2,...}$ of such inequalities implies directly that

$$\int_{|u*s|>T} |\eta|^2 \leq c_0 \, e^{-|u*s|/c_0} \sup_S |\eta| .$$

(6.21)

This last bound with the aforementioned elliptic regularity techniques will give the bounds that are asserted by the lemma.

*Part 6*: The next lemma speaks to the relationship between the kernel of the operator D in (6.11) and solutions to (6.10). The lemma refers to the kernel and cokernel of the Fredholm version of $D_S$ that maps $\mathbb{H}_S$ to $\mathbb{L}_S$. The lemma also uses $\mathbb{I}_S$ to denote the tautological inclusion map from kernel($D_S$) into $\mathbb{H}_S$

**Lemma 6.6**: *Let (S, u) denote a given admissible pair. There exists a constant $\kappa \geq 1$, a neighborhood $U \subset$ kernel($D_S$) of the origin, and a smooth maps $\mathfrak{f}: U \to$ cokernel($D_S$) and $\mathfrak{q}: U \to \mathbb{H}_S \cap C^\infty(S; N_0)$ with the following properties:*
- $\mathfrak{f}(0) = 0$, $\mathfrak{q}(0) = 0$ *and* $d\mathfrak{q}|_0 = \mathbb{I}_S$.
- *Let $\mathcal{N}$ denote the space of smooth sections of $N_0$ that obey (6.10), the second and third bullets of (6.12) and have pointwise norm bounded by $\kappa^{-1}\rho_S$. Then $\mathfrak{q}$ maps $\mathfrak{f}^{-1}(0)$ homeomorphically onto $\mathcal{N}$.*
- *The space $\mathcal{N}$ parametrizes a neighborhood of (S, u) in $\mathcal{A}_{HF}$ in the following way: Let $\eta \in \mathcal{N}$. The corresponding admissible pair (S´, u´) is defined by the condition that $u´(S´) = (\mathfrak{e}_S \circ \eta)(S)$. This is to say that S´ = S but with the complex structure that is induced by the $J_{HF}$-holomorphic embedding $\mathfrak{e}_S \circ \eta$.*

*Proof of Lemma 6.6*: The assertions of the first two bullets are straightfoward applications of the inverse function theorem following standard arguments pioneered by Kuranishi. This is done by viewing the expression on the left hand side of (6.10) as defining a map from a ball in a certain Banach space to a second Banach space. The only subtlety is that the spaces $\mathbb{H}_S$ and $\mathbb{L}_S$ can not be used because the $L^2_1$ norm does not control the supremum norm. Even so, this is a standard issue in the theory of J-holomorphic curves and is dealt with by using a domain and range that give slightly stronger control over the derivatives of the section, see for example [MS]. The third bullet follows directly from what is said in Part 5.

The next lemma is a corollary to the preceding lemma. It refers to the $\mathbb{R}$-action on $\mathcal{A}_{HF}$ where the generator $1 \in \mathbb{R}$ acts on any given pair $(S, u)$ to give the pair $(S, u^1)$ where $u_1$ is the composition of first $u$ and then translation by 1 along the $\mathbb{R}$ factor of $\mathbb{R} \times [1, 2] \times \Sigma$.

**Lemma 6.7**: *Let $\mathcal{A}_{smooth}$ denote the subspace of the moduli space of admissible pairs that contains the unobstructed pairs. The subspace $\mathcal{A}_{smooth}$ is an open, $\mathbb{R}$-invariant set with the structure of a smooth manifold. Moreover, if $(S, u) \in \mathcal{A}_{smooth}$, then the corresponding version of Lemma 6.6's map $\mathfrak{f}$ defines a smooth embedding from a small radius ball about the origin in* $\mathrm{kernel}(D_S)$ *onto an open neighborhood of $(S, u)$ in $\mathcal{A}_{smooth}$.*

*Proof of Lemma 6.7*: The assertion that a neighborhood in $\mathcal{A}_{HF}$ of a given pair $(S, u)$ is homeomorphic to a neighborhood of 0 in the kernel of $D_S$ follows directly from Lemma 6.6. The smooth structure on this neighborhood is defined by requiring that this same homeomorphism define a diffeomorphism. The inverse function theorem can be used to verify that the associated transition functions are $C^\infty$.

**7. ech-HF subvarieties and Heegaard Floer curves**
This section explains how the $\mathbb{R} \times M_\delta$ part of any given ech-HF subvariety determines pseudoholomorphic subvarietis of the sort used in Lipshitz's [L] reformulation of Heegaard Floer homology. The observations of this section are summarized by Proposition 7.1 in Section 7a, Proposition 7.2 in Section 7b and Proposition 7.3 in Section 7c.

### a) ech-HF subvarieties and singular admissible sets

Singular admissible sets enter the embedded contact homology story via the upcoming Proposition 7.1. To set things up for this proposition, it is necessary to add one additional constraint to those in Section 3a for the sorts of almost complex structures on $\mathbb{R} \times Y$ that are used to define the notion of an ech-HF subvariety. To this end, fix an almost complex $J_{HF}$ that obeys (6.1). What follows is the extra constraint on Section 3a's almost complex structure J.

$$J = J_{HF} \text{ on the } f^{-1}(1,2) \text{ part of } \mathbb{R} \times M_\delta.$$
(7.1)

Note in this regard that CONSTRAINTS 1 and 2 in Section 3a are equivalent on $\mathbb{R} \times M_\delta$ to what is required by the first and second bullets in (6.1). CONSTRAINT 3 in Section 3a is equivalent to the rule in the third bullet of (6.1). Meanwhile CONSTRAINTS 4 and 5 in Section 3a are compatible with the constraints in the fourth bullet of (6.1).

A given singular admissible set $Z = ((S, u), \vartheta_\Sigma)$ defines an integral, 2-dimensional current on any given compact subset of $\mathbb{R} \times [1,2] \times \Sigma$, this given by integration over the 2-chain $[u(S)] + \sum_{(Z,m)\in\vartheta_\Sigma} m[Z]$. This current is denoted by $[Z]$. The pairing of any given 2-form $\mu$ with $[Z]$ is denoted by $[Z](\mu)$. The support of the current $[Z]$ is the subvariety $u(S) \cup (\cup_{(Z,m)\in\vartheta_\Sigma} Z)$. The latter is denoted by $|Z|$. Note that integration over $|Z|$ does not give $[Z]$ unless all integers from $\vartheta_\Sigma$ are equal to 1.

**Proposition 7.1**: *Fix an almost complex structure $J_{HF}$ that obeys the conditions in (6.1). Given $\mathcal{I} \geq 1$ and $\varepsilon > 0$, there exists $\kappa > \varepsilon^{-1}$ with the following significance: Define the geometry of Y using data $\delta, x_0, R$ with $\delta < \kappa^{-2}$, and fix an almost complex structure, J, on $\mathbb{R} \times Y$ that obeys the constraints in Section 3a plus (7.1). Fix an ech-HF subvariety C with $I_{ech}(C) < \mathcal{I}$ or with $\int_C w < \mathcal{I}$. There is a singular admissible set, $Z = ((S,u), \vartheta_\Sigma)$, with*

- $\sup_{z \in (C \cap ([-1/\varepsilon, 1/\varepsilon] \times M_\varepsilon))} \text{dist}(z, |Z|) + \sup_{z \in (|Z| \cap ([-1/\varepsilon, 1/\varepsilon] \times M_\varepsilon))} \text{dist}(z, C) < \varepsilon$.

- *Let $\mu$ denote a smooth 2-form on $\mathbb{R} \times M_\varepsilon$ with sup-norm bounded by 1, covariant derivative norm bounded by $\varepsilon^{-1}$ and compact support where $|s| < \varepsilon^{-1}$. Then*

$$\left| \int_C \mu - [Z](\mu) \right| \leq \varepsilon.$$

The rest of this subsection contains the

*Proof of Proposition 7.1*

The proof that follows has eight steps.

Step 1: Suppose that no such $\kappa$ exists. There would necessarily exist a sequence $\{(D_n, C_n)\}_{n=1,2,}$ of pairs of the following sort: First, any given $n \in \{1, 2,...\}$ version of $D_n$ is a data set consisting of 4-tuples $(\delta_n, x_{0,n}, R_n, J_n)$ such that $\delta_n < 1/n$; and such that the geometry of Y and $\mathbb{R} \times Y$ can be defined as in Section 1a and so that the results in the previous sections can be invoked using $\delta = \delta_n$, $x_0 = x_{0,n}$, $R = R_n$, and $J = J_n$. Second, for each such integer n, what is denoted by $C_n$ is an ech-HF subvariety as defined by $J_n$. This subvariety is such that $I_{ech}(C_n) \leq \mathcal{I}$ or that $w(C_n) < \mathcal{I}$. However, there is no pair $(u, S)$ that satisfies the $C = C_n$ version of what is required by the two bullets of the lemma.

Step 2: With Proposition 5.1 in mind, use Proposition 5.5 to find a finite set, $\vartheta$, and a subsequence of $\{C_n\}_{n=1,2,...}$, hence renumbered consecutively from 1, with properties that are described in this step and in Steps 3-5.

Consider first the properties of $\{C_n\}_{n=1,2,...}$. The $s \to \infty$ limit of the constant $s$ slices of a given $n \in \{1, 2, ...\}$ version of $C_n$ defines an element in $\mathcal{Z}_{ech,M}$ and thus an HF-cycle of length G. By the same token, the $s \to -\infty$ limit of the constant $s$ slices of $C_n$ also defines an HF-cycle of length G. The sequence $\{C_n\}_{n=1,2,...}$ is chosen so that the same HF-cycle is defined by the $s \to \infty$ limit of the constant $s$ slices of each $n \in \{1,2,...\}$ version of $C_n$, and so that the same HF-cycle is defined by the $s \to -\infty$ limit of the constant $s$ slices of each such $C_n$. These respective cycles are denoted by $\hat{\vartheta}_+$ and $\hat{\vartheta}_-$.

Step 3: Additional properties of $\{C_n\}_{n=1,2,...}$ are stated momentarily. To set the stage, use Proposition 5.8 to find $z_\Diamond < \delta_*^2$ so that the conclusions of Proposition 5.8 hold for each $n \in \{1, 2, ...\}$ version of $C_n$. If n is sufficiently large, then $\delta_n < \frac{1}{100} z_\Diamond^{1/2}$ and so $C_n$ interects each $\mathbb{R} \times \{1+z_\Diamond\}$ member of $\mathcal{M}_\Sigma$ precisely once in the radius $\delta_*$ coordinate ball centered on each index 1 critical point of $f$. The subvariety $C_n$ likewise intersects each $\mathbb{R} \times \{2-z_\Diamond\}$ member of $\mathcal{M}_\Sigma$ precisely once in the radius $\delta_*$ coordinate ball centered on each index 2 critical point of $f$.

What follows is a consequence of the preceding observations: Fix $z \in (0, z_\Diamond]$. Then all sufficiently large n versions of $C_n$ intersect each $\mathbb{R} \times \{1+z\}$ member of $\mathcal{M}_\Sigma$ precisely once in the radius $\delta_*$ coordinate ball centered on each index 1 critical point of $f$. Meanwhile, $C_n$ intersects each $\mathbb{R} \times \{2-z\}$ member of $\mathcal{M}_\Sigma$ precisely once in the radius $\delta_*$ coordinate ball centered on each index 2 critical point of $f$.

To say more about these intersections, fix $z$ as above. The level set in question is an annulus, as can be seen using the coordinates $(t, \varphi_+, \hat{h}_+)$ if p has index 1 or the coordinates $(t, \varphi_-, \hat{h}_-)$ if p has index 2. In the former case, this is the $t = 1+z$ level set, and the annulus is that where $|\hat{h}_+| < \frac{4}{3\sqrt{3}} \delta_*^2$. In the latter case, the level set has $t = 2-z$ and the annulus is that where $|\hat{h}_-| < \frac{4}{3\sqrt{3}} \delta_*^2$. In either case, let $\Sigma^z_p$ denote the part of this level set in the radius $\delta_*$ coordinate ball centered at p. It follows from Lemma 4.7 that the intersection of $C_n$ with $\mathbb{R} \times \Sigma^z_p$ is a smoothly embedded arc that appears as the graph of a map from $\mathbb{R}$ to $\Sigma^z_p$. The $s \to -\infty$ limit of the latter map has distance $c_0 \delta_n$ or less from the intersection between $\Sigma^z_p$ and the relevant integral curve of $\mathfrak{v}$ from $\hat{\mathfrak{v}}_-$; and the $s \to \infty$ limit has distance $c_0 \delta_n$ or less from the intersection between $\Sigma^z_p$ and the relevant integral curve of $\mathfrak{v}$ from $\hat{\mathfrak{v}}_+$.

Proposition 5.8 has the following additional consequence: Fix $\rho \in (0, \delta_*]$ and there exists $z_\rho \in (0, z_0]$ such that all sufficiently large n versions of $C_n$ intersect $\Sigma^z_p$ in the radius $\rho$ coordinate ball centered on p.

<u>Step 4</u>: Given the assumed bound on $\{I_{ech}(C_n)\}_{n=1,2,\ldots}$, the two bullets of Proposition 5.1 can be used to apply Proposition 5.5 with U first $M_{\delta_*}$, then $M_{\delta_*/2}$, and successively for $U = M_{\delta_*/k}$ with $k \in \{1, 2, \ldots\}$. These successive applications of Proposition 5.5 produce the following: First, a finite set $\vartheta$ whose typical element is a pair $(Z, m)$ with $Z \subset \mathbb{R} \times M_{(1,2)}$ an irreducible, $J_{HF}$-holomorphic subvariety. These are such that no two pair from $\vartheta$ share the same subvariety. Moreover,

- $\sum_{(Z,m) \in \vartheta} \int_{Z \cap ([s,s+1] \times M_{(1,2)})} ds \wedge df = G$ *for each* $s \in \mathbb{R}$.
- $\sum_{(Z,m) \in \vartheta} m \int_Z w_\Sigma \leq \mathcal{I}$.

(7.2)

Also produced is a subsequence (hence renumbered consecutively from 1) of $\{C_n\}_{n=1,2,\ldots}$ such that if $\mathcal{K} \subset \mathbb{R} \times M_{(1,2)}$ is any given compact set, then

- $\lim_{n \to \infty} (\sup_{z \in (C_n \cap \mathcal{K})} \text{dist}(z, \cup_{(Z,m) \in \vartheta} Z) + \sup_{z \in \cup_{(Z,n) \in \vartheta} Z \cap \mathcal{K}} \text{dist}(z, C_n)) = 0$.
- *Let $\mu$ denote a smooth 2-form with compact support on $\mathcal{K}$. Then*

$$\lim_{n \to \infty} \left( \int_{C_n} \mu - \sum_{(Z,m) \in \vartheta} m \int_Z \mu \right) = 0.$$

(7.3)

Let $\vartheta_\Sigma \subset \vartheta$ denote the subset of pairs whose subvariety component is mapped to a constant $f$ slice of $M_{(1,2)}$ via the projection from $\mathbb{R} \times M_{(1,2)}$. Use $\vartheta_*$ to denote $\vartheta - \vartheta_\Sigma$.

Step 5: Let p denote an index 1 or index 2 critical point of $f$. Fix $z \in (0, z_{\Diamond}]$. Each sufficiently large n version of $C_n$ has precisely 1 intersection point with $\mathbb{R} \times \Sigma^z_p$. This with (7.3) has two immediate implications: First, there is precisely one pair from $\vartheta_*$ whose image via the projection from $\mathbb{R} \times M_{(1,2)}$ intersects $\Sigma^z_p$. Second, each pair from $\vartheta_*$ has integer component 1. This being the case, view $\vartheta_*$ now as a set of distinct $J_{HF}$ holomorphic subvarieties. What follows is a consequence of Lemma 4.2 and (7.3): View $M_{(1,2)}$ as $(1,2) \times Z$ and $\mathbb{R} \times M_{(1,2)}$ as $\mathbb{R} \times (1,2) \times \Sigma$. The subvariety $Z_0 = \cup_{Z \in \theta_*} Z$ has intersection number G with each constant $\mathbb{R} \times (1, 2)$ slice of $\mathbb{R} \times (1,2) \times \Sigma$.

The preceding observations with Lemma 4.7 imply the following: There exists $\rho_* \in (0, \delta_*]$ such that if p is any given index 1 or index 2 critical point of $f$, then exactly one irreducible component of the subvariety $Z_0$ is mapped via the projection to $M_{(1,2)}$ into the radius $\rho_*$ coordinate ball centered at p, and its intersection with the product of $\mathbb{R}$ with this ball is a submanifold of $\mathbb{R} \times M_{(1,2)}$. Take $\rho_*$ so that no submanifold from $\vartheta_\Sigma$ projects to the union of these radius $\rho_*$ coordinate balls.

The next observation is a consequence both of Lemma 4.7 and what is said at the end of Step 3. Given $\rho \in (0, \rho_*]$, there exists $z_\rho \in (0, z_{\Diamond}]$ such that if $z \in (0, z_\rho]$, then the subvariety from $\vartheta_*$ whose projection to $M_{(1,2)}$ intersects a given version of $\Sigma^z_p$ does so in the radius $\rho$ coordinate ball centered on the critical point p.

View $\mathbb{R} \times M_{(1,2)}$ as $\mathbb{R} \times (1,2) \times \Sigma$ and view $Z_0$ as a subvariety in $\mathbb{R} \times (1,2) \times \Sigma$. Given Lemma 4.7, was said in the preceding two paragraphs has the following reinterpretation: There exists $z_1 \in (0,1)$ such that the $\mathbb{R} \times (1, 1+z_1)$ part of $Z_0$ is a smoothly embedded submanifold of $\mathbb{R} \times (1,2) \times \Sigma$ that consists of G components; and each component is a graph over its image in $\mathbb{R} \times (1,2)$. Moreover, this graph is defined by a map $\lambda \colon \mathbb{R} \times (1, 2) \to \Sigma$ whose image in $T_+$ and is such that

$$\lim_{t \to 1} \sup_{s \in \mathbb{R}} \text{dist}(\lambda(s,t), C_+) = 0 \,.$$

(7.4)

Finally, each circle from $C_+$ contains the limit points of precisely one $\mathbb{R} \times (1, 1+z_1)$ component of $Z_0$.

There is a completely analogous story for the $\mathbb{R} \times (2-z_1, 2)$ part of $Z_0$ with it understood that $T_-$ and $C_-$ appears in lieu of $T_+$ and $C_+$.

Equation (7.2) and what is said above can be used with Lemma 5.6 to deduce the following: There exists $s_1 \geq 1$ such that the $s > s_1$ part of $Z_0$ is a disjoint union of G smooth, open strips. Each such strip is a graph over its image in $(s_1, \infty) \times (1,2)$. Such a graph is defined by a map, $\psi: (s_1, \infty) \times (1,2) \to \Sigma$ with image in $T_+ \cap T_-$. In addition, there exists a point $q \in C_+ \cap C_-$ such that

$$\lim_{s \to \infty} \sup_{t \in (1,2)} \text{dist}(\psi(s, t), q) = 0$$

(7.5)

Finally, the G points so defined in $C_+ \cap C_-$ pair the set of index 1 critical points of $f$ with the set of index 2 critical points of $f$.

The analogous conclusion holds for the $s < -s_1$ portion of $\mathcal{Z}_0$.

<u>Step 6</u>: The union of $\mathcal{Z}_0 \subset \mathbb{R} \times (1,2) \times \Sigma$ with the constant $\mathbb{R} \times (1, 2)$ slices from $\vartheta_\Sigma$ define a $J_{HF}$-holomorphic subvariety in $\mathbb{R} \times (1,2) \times \Sigma$ with at worst a finite set of singular points. The fact that holomorphic structures on punctured surfaces extend over the punctures implies that there is a *model curve* for $\mathcal{Z}_0$. This model curve is a smooth complex curve with a $J_{HF}$-holomorphic, almost every where 1-1 map to $\mathbb{R} \times (1, 2) \times \Sigma$ whose image is $\mathcal{Z}_0$. Let $S_0$ denote this model curve and let $u$ denote the map. The complement of a finite set in $S_0$ is identified by $u$ with the complement of the singular points in $\mathcal{Z}_0$ and this identification defines the complex structure on this part, and hence all of $S_0$. As explained in the next steps, $S_0$ is the interior of a surface with boundary S, to which $u$ extends so that (S, $u$) with $\vartheta_\Sigma$ defines a singular admissible set. Given what is said by (7.3), the existence of (S, $u$) puts the lie to the assumption that the conclusions of Proposition 7.1 do not hold for the sequence $\{C_n\}_{n=1,2,\ldots}$. This contradiction is avoided only if Proposition 7.1 is true.

<u>Step 7</u>: This step describes S and the compact surface <u>S</u>. The complement of $S_0$ in S is the disjoint union of 2G copies of $\mathbb{R}$, each these in 1-1 correspondence with the index 1 critical points of $f$. Let p denote a given index 1 critical point of p. And let $\partial_p S$ denote the corresponding boundary component of S. A neighborhood of $\partial_p S$ in S is described in the next paragraph.

To start the description, let $Z_p$ denote the component of the $t < 1+z_1$ part of $\mathcal{Z}_0$ that lies in $\mathbb{R} \times (1, 1+z_1) \times T_{p+}$ and let $S_{0p} \subset S_0$ denote $u^{-1}(Z_p)$. The map $u$ identifies $S_{0p}$ with $Z_p$. Since $Z_p$ is a graph over $\mathbb{R} \times (1, 1+z_1)$, this identification with $Z_p$ followed by the projection identifies $S_{0p}$ with $\mathbb{R} \times (1, 1+z_1)$. This identification extends to a neighborhood of $S_{0p}$ in S so as to identify the latter with $\mathbb{R} \times [1, 1+z_1)$. This identification writes the component $\partial S_p$ as $\mathbb{R} \times \{1\}$.

Let p now denote an index 2 critical point and let $Z_p$ in this case denote the component of the $t > 2-z_1$ part of $\mathcal{Z}_0$ in $\mathbb{R} \times (2-z_1, 2) \times T_{p-}$. The $u$-inverse image of $Z_p$ is denoted by $S_{0p}$ and it is identified with $\mathbb{R} \times (2-z_1, 2)$ via the composition of $u$ and the projection from $Z_p$ to the $\mathbb{R} \times (2-z_1, 2)$ factor. A neighborhood in S of this part of $S_0$ has the corresponding identification with $\mathbb{R} \times (2-z_1, 2]$ with $\partial_p S$ appearing as $\mathbb{R} \times \{2\}$.

Consider next $\underline{S}$. The complement of S in $\underline{S}$ consists of 2G points, half labeled positive and half labeled negative. What follows describes the positive points. There is a completely analogous description of the negative points which is left to the reader. The set of positive points in $\underline{S}$–S is in 1-1 correspondence with the G points in $C_+ \cap C_-$ that appear in the G versions of (7.5). Let q denote such a point and let $q_S$ denote its partner in $\underline{S}$–S. Let p and p´ denote the respective index 1 and index 2 critical point versions of $C_{p+}$ and $C_{p'-}$ whose intersection contains q. The point $q_S$ is contained in the closure of $\partial_p S$ and $\partial_{p'} S$. To elaborate, let $Z_{0q+}$ denote the component of the $s \geq s_1$ part of $Z_0$ that appears with q in (7.5). Let $S_{0q+} \subset S_0$ denote the $u$-inverse image of $Z_{0q+}$. The map $u$ identifies $S_{0q+}$ with $Z_{0q+}$ and as a consequence, the projection from $Z_{0q+}$ to $(s_1, \infty) \times (1, 2)$ identifies $S_{0q+}$ with $(s_1, \infty) \times (1, 2)$. With this identification understood, it follows that $S_{0q+}$ has a neighborhood in S that is diffeomorphic as a manifold with boundary to $(s_1, \infty) \times [1, 2]$ with the boundary component $(s_1, \infty) \times \{1\}$ in $\partial_p S$ and the component $(s_1, \infty) \times \{2\}$ in $\partial_{p'} S$.

Let $\zeta$ denote the standard complex coordinate on $\mathbb{C}$ and introduce $D_+$ to denote the $|z| < e^{-\pi s_1}$ part the closed upper half plane. Identify the $(s_1, \infty) \times [1, 2]$ with $D_+$–0 via the map

$$(s, t) \to e^{-\pi(s+it)}.$$

(7.6)

This identification extends to a neighborhood of $q_S$ so as to identify the latter with $D_+$.

Step 8: Look at (6.2) to see that the only remaining task is that of extending the map $u$ to the boundary of S. Note in this regard that (7.4) and its $t \to 2$ analog imply the following: Let p denote an index 1 critical point of $f$. If $u$ extends to $\partial_p S$, then it maps the latter to $\mathbb{R} \times C_{p+}$. By the same token, let p now denote an index 2 critical point of $f$. If $u$ extends to $\partial_p S$, then it maps $\partial_p S$ to $\mathbb{R} \times C_{p-}$. In either case, if the extension is smooth, then it is de-facto an embedding of $\partial_p S$ by virtue of what is said in Step 5 about the map $u$ on $S_{0p}$. By way of a reminder, let p denote an index 1 critical point of $f$. Step 5 writes a neighborhood of $\partial S_p$ in S as $\mathbb{R} \times [1, 1+z_1)$ and writes $u$ on the $\mathbb{R} \times (1, 1+z_1)$ part as the graph of the map $\lambda$ from $\mathbb{R} \times (1, 1+z_1)$ to $T_{p+}$. The extension of $u$ to $\mathbb{R} \times \{1\}$ is given by an extension of $\lambda$ which maps $\mathbb{R} \times \{1\}$ to $C_{p+}$.

The assumption in the fifth bullet of (6.1) implies that $J_{HF}$ on $\mathbb{R} \times (1, 1+z_1) \times T_+$ defines an $s$ and $t$ independent complex structure on $T_+$ such that $\varphi_+ + i\hbar_+$ defines a holomorphic parameter on a neighborhood of any given point. Parametrize a neighborhood of $\partial_p S$ in S as $\mathbb{R} \times [1, 1+z_1)$ in the manner described above. Use $s$ as before to denote the Euclidean coordinate on the $\mathbb{R}$ factor of $\mathbb{R} \times [1, 1+z_1)$; but write the Euclidean coordinate $t$ on the $[1, 1+z_1)$ factor as $1+z$ where z is the Euclidean coordinate on $(0, z_1)$. Use the coordinates $(\varphi_+, h_+)$ to write $\lambda$ as the map

$$(x, z) \to (\varphi_+ = \varphi^S(s, z), h_+ = \varsigma^S(s, z))$$
(7.7)

where $\varphi^S$ and $\varsigma^S$ are smooth functions. The fact that $u$ is J-holomorphic implies in particular that $\varphi^S + i\varsigma^S$ must obey the Cauchy-Riemann equations, which is to say that

$$\partial_s \varphi^S - \partial_z \varsigma^S = 0 \quad and \quad \partial_z \varsigma^S + \partial_s \varphi^S = 0.$$
(7.8)

Equation (7.4) asserts that $\lim_{z \to 0} \sup_{s \in \mathbb{R}} |\varsigma^S| = 0$.

The Cauchy-Riemann equations in (7.8) and the fact that $|\varsigma^S| \to 0$ as $z \to 0$ implies via a version of Schwarz reflection trick from Theorem 24 in [Ah] that the pair $(\varphi^S, \varsigma^S)$ extend smoothly to $z = 0$ boundary of the domain.

**b) ech-HF subvarieties and broken, singular admissible sets**

An ech-HF subvariety that obeys the conditions of Proposition 7.1 can define in principle a finite collection of singular admissible sets. The members of such a collection obey certain constraints that are described Part 1 of the subsection. A collection that obeys these constraints is said in what follows to be a *broken, singular admissible set*. The upcoming Proposition 7.2 in Part 2 of the subsection explains how an ech-HF subvariety determines such a collection.

*Part 1*: Let $N \in \{1, 2, ...\}$ and let $\Xi = \{((S_1, u_1), \vartheta_{\Sigma 1}), ..., ((S_N, u_N), \vartheta_{\Sigma N})\}$ denote an ordered collection of N singular admissible sets. This set $\Xi$ is said to be a broken, singular admissible set when the constraints given in the upcoming (7.9) are met. These constraints refer solely to the case when $N > 1$. A broken admissible set with one element is neither more nor less than an admissible set. The notation in (7.9) uses $\hat{\upsilon}_{k+}$ and $\hat{\upsilon}_{k-}$ when $k \in \{1, ..., N\}$ to denote the respective collections of G points in $C_+ \cap C_-$ that are defined by the positive and negative points of $\underline{S}_k - S_k$. Equation (7.9) also refers to two currents in $\mathbb{R} \times [1, 2] \times \Sigma$ that are defined on any given compact set by a singular admissible set. The first of these currents is the current $[Z]$ from Proposition 7.1. The second is given by integration over the $s \leq 0$ part of the 2-chain that defines $[Z]$; it denoted by $[Z_{<0}]$. Note that both $[Z]$ and $[Z_{<0}]$ have finite pairing with the 2-form $w_\Sigma$. These pairings are denoted by $[Z](w)$ and $[Z_{<0}](w)$.

Here are the promised constraints.

- *No element* $((S, u), \vartheta_\Sigma) \in \Xi$ *is such that both* $\vartheta_\Sigma = \emptyset$ *and* $u(S)$ *is $\mathbb{R}$-invariant.*
- *Each* $Z = ((S, u), \vartheta_\Sigma) \in \Xi$ *is such that* $[Z_{<0}](w) = \frac{1}{2}[Z](w)$.
- $\hat{\upsilon}_{(k-1)+} = \hat{\upsilon}_{k-}$ *for each* $k \in \{2, \ldots, N\}$.

(7.9)

To give an idea as to how to view these constraints, fix $k \in \{2, \ldots, N\}$. What is said in (7.9) implies that the constant, $s \gg 1$ slices of the $u_{k-1}$ image of $S_{k-1}$ are very close to the constant, $s \ll -1$ slices of the $u_k$ image of $S_k$. This fact suggests that the union of a suitable set of $k \in \{1, \ldots, N\}$ dependent translates in $\mathbb{R} \times [1, 2] \times \Sigma$ of the subvarieties that are given by $\Xi$ may well be everywhere very close to a subvariety that is defined by a single admissible set of the form $\{(S, u), \vartheta_\Sigma\}$ where $S$ is a obtained from the disjoint union of the surfaces from the set $\{S_k\}_{1 \le k \le N}$ by identifying a suitable neighborhood of the positive points in each $k \in \{2, \ldots, N\}$ version of $\underline{S}_{k-1}$ with a suitable neighorhood of the negative points in $\underline{S}_k$.

*Part 2*: The focus here is the upcoming Proposition 7.2. The statement of this proposition uses the following notation: Let $Z = \{(S, u), \vartheta_\Sigma\}$ denote a singular admissible set and let $s \in \mathbb{R}$ denote a given point. The $J_{HF}$-holomorphic map from $S$ to $\mathbb{R} \times [1,2] \times \Sigma$ that is obtained by composing first $u$ and then translation by $s$ along the $\mathbb{R}$ factor is denoted by $u^s$. The translation of the elements in $\vartheta_\Sigma$ by $s$ along the $\mathbb{R}$ factor defines a new set of constant $(s,t)$ slices of $\mathbb{R} \times [1,2] \times \Sigma$ that is denoted in by $\vartheta_\Sigma^s$. The new singular admissible set $((S, u^s), \vartheta_\Sigma^s)$ is denoted below by $Z^s$.

The proposition also refers to a set of intervals $\{I_1, \ldots, I_N\} \subset \mathbb{R}$ that is defined using a constant $\varepsilon$ and a corresponding set of real numbers $\{s_1, \ldots, s_N\}$ with $s_{k-1} < s_k - \frac{3}{2\varepsilon}$. If $N = 1$, then $I_1 = \mathbb{R}$. If $N$ is bigger than 1, then the interval $I_1$ is $(-\infty, s_2 - \frac{1}{2\varepsilon})$, the interval $I_N$ is $(s_{N-1} + \frac{1}{2\varepsilon}, \infty)$, and each $k \in \{2, \ldots, N-1\}$ version of $I_k$ is $(s_{k-1} + \frac{1}{2\varepsilon}, s_{k+1} - \frac{1}{2\varepsilon})$. Note that the union of these intervals is $\mathbb{R}$.

**Proposition 7.2**: *Fix an almost complex structure $J_{HF}$ that obeys the conditions in (6.1). Given $\mathcal{I} \ge 1$ there exists a constant $\kappa_{\mathcal{I}} \ge 1$, and given also $\varepsilon > 0$, there exists $\kappa > \varepsilon^{-1}$ with the following significance: Define the geometry of $Y$ using data $\delta, x_0, R$ with $\delta < \kappa^{-2}$, and fix an almost complex structure, $J$, on $\mathbb{R} \times Y$ that obeys the constraints in Section 3a plus*

*(7.1)*. *Let C denote an ech-HF subvariety with either* $I_{ech}(C) \leq \mathcal{I}$ *or with* $\int_C w$. *There exists a broken, singular admissible set* $\Xi = \{\mathcal{Z}_1, \ldots, \mathcal{Z}_N\}$ *with at most* $\kappa_{\mathcal{I}}$ *elements and a corresponding set* $\{s_1, s_2, \ldots, s_N\} \in \mathbb{R}$ *with the properties listed below.*

- $s_{k-1} < s_k - \frac{3}{2\varepsilon}$ *for each* $k \in \{2, \ldots, N\}$.
- $\sup_{z \in (C \cap (I_k \times M_\varepsilon))} \mathrm{dist}(z, |\mathcal{Z}_k^{s_k}|) + \sup_{z \in (|\mathcal{Z}_k^{s_k}| \cap (I_k \times M_\varepsilon))} \mathrm{dist}(z, C) < \varepsilon$ *for each* $k \in \{1, \ldots, N\}$.
- *Let* $\mu$ *denote a smooth 2-form on* $\mathbb{R} \times M_\varepsilon$ *with norm bounded by* 1, *covariant derivative norm bounded by* $\varepsilon^{-1}$ *and compact support on an interval of length* $\frac{2}{\varepsilon}$ *in some* $k \in \{1, \ldots, N\}$ *version of* $I_k$. *Then* $|\int_C \mu - [\mathcal{Z}_k^{s_k}](\mu)| \leq \varepsilon$.

Note that this proposition in the case N = 1 says more than what is said by Proposition 7.1 by virtue of the fact that $I_1 = \mathbb{R}$.

***Proof of Proposition 7.2***: The proof is obtained by applying Proposition 7.1 to suitable translates of C along the $\mathbb{R}$ factor of $\mathbb{R} \times M$. The bound on N is obtained using Proposition 5.1 to bound the sum of the various $\mathcal{Z} \in \Xi$ versions of $[\mathcal{Z}](w)$. Note in addition that the set of singular admissible sets with a given bound for this pairing is compact in the sense that any given sequence with this bound converges in the sense given by the bullets in Proposition 7.2 to a broken, singular admissible set. Proposition 5.1 and Lemma 5.7 are used to deal with the non-compact ends of $I_1$ and $I_N$. As there are no hidden subtleties in the arguments, the details are omitted.

**c) Canonical admissible pairs for ech-HF submanifolds**

This section supplies a refinement to the versions of Propositions 7.1 and 7.2 for the case when $\Xi$ has but a single element which is a true admissible pair. This refinement is used in [KLTIII] to choose such a pair in a 'canonical' fashion. The upcoming Proposition 7.3 states the refinement. What follows directly supply some necessary background.

Let $\mathcal{A}_{smooth}$ denote the subspace in the moduli space of unobstructed admissible pairs. Let $S \to \mathcal{A}_{smooth}$ denote the fiber bundle whose fiber over a given element (S, $u$) is the surface S. Sitting over $S$ is the universal normal bundle, this the complex line bundle whose fiber over $S|_{(S,u)}$ is the normal bundle $N \to S$. Let $\mathcal{N} \to S$ denote this bundle.

To continue, suppose that $\mathcal{K} \subset \mathcal{A}_{smooth}$ has compact image in $\mathcal{A}_{smooth}/\mathbb{R}$. The techniques used for Lemma 5.4 of [T3] work for parametrized families of surfaces and in particular, they can be used to obtain the following data: A constant radius disk subbundle $\mathcal{N}_0 \subset \mathcal{N}|_{\mathcal{K}}$, a smooth map, $\mathbb{R}$-equivariant map $\mathfrak{e}^{\mathcal{K}}: \mathcal{N}_0 \to \mathbb{R} \times [1, 2] \times \Sigma$ and a constant $\rho_{\mathcal{K}} > 0$ with the following property: The map $\mathfrak{e}^{\mathcal{K}}$ restricts to a given $(S, u) \in \mathcal{K}$ fiber $S|_{(S,u)}$ as an exponential map of the sort described in Part 1. In particular, the choice $\mathfrak{e}_S = \mathfrak{e}^{\mathcal{K}}|_{(S,u)}$ embeds each fiber as a $J_{HF}$-holomorphic disk such that the corresponding disk bundle $N_0 = \mathcal{N}_0|_{(S,u)}$ is embedded onto a neighborhood of $S$ that contains all points in $\mathbb{R} \times [1, 2] \times \Sigma$ with distance less than $\rho_{\mathcal{K}}$ from $S$. The parts that follow use implicitly $N_0 = \mathcal{N}_0|_{(S,u)}$ and $\mathfrak{e}_S = \mathfrak{e}^{\mathcal{K}}|_{(S,u)}$ when referring to an exponential map for a given pair $(S, u) \in \mathcal{K}$. In particular, this version is used in Proposition 7.3 when defining the $(S, u) \in \mathcal{K}$ version of the operator $D_S$, and it is used to view submanifolds near $u(S)$ as submanifolds of the corresponding version of $N_0$.

**Proposition 7.3**: *Fix an $\mathbb{R}$-invariant compact set $\mathcal{K} \subset \mathcal{A}_{smooth}$, and there exists $\kappa \geq 1$; and given $\mathcal{I} \geq 1$ and $\varepsilon \in (0, \kappa^3)$, there exists $\kappa_* \geq \kappa^4$ with the following significance: Define the geometry of $Y$ and $\mathbb{R} \times Y$ using $\delta < \kappa_*^{-1}$ and an almost complex structure that obeys the constraints in Section 3a and also (7.1). Let $C$ denote an embedded, ech-HF subvariety with either $I_{ech}(C) < \mathcal{I}$ or with $\int_C w < \mathcal{I}$.*

- *Suppose that the conclusions of Proposition 7.1 hold for $C$ using $\Xi = ((S´, u´), \emptyset)$ with $(S´, u´) \in \mathcal{K}$. There exists a unique $(S, u) \in \mathcal{A}_{smooth}$ such that the intersection of $C$ with $[-\frac{1}{\varepsilon}, \frac{1}{\varepsilon}] \times M_\varepsilon$ can be written as $\mathfrak{e}_S \circ \eta$ where $\eta$ is a section of $N_0 \to S$ with the following properties:*
  a) *The conclusions of Proposition 7.1 hold using the given value of $\varepsilon$ and $(S, u)$.*
  b) *The norm of $\eta$ is bounded by $\kappa\varepsilon$.*
  c) *The section $\eta$ obeys (6.17) on the $(u^*s, u^*t) \in [-\frac{1}{\varepsilon}, \frac{1}{\varepsilon}] \times [1+\varepsilon^2, 2-\varepsilon^2]$ part of $S$, and the respective restrictions of $\eta$ and the kernel of $D_S$ are $L^2$-orthogonal on this same part of $S$.*
- *Suppose that the conclusions of Proposition 7.2 hold for $C$ using $\Xi = ((S´, u´), \emptyset)$ with $(S´, u´) \in \mathcal{K}$. There exists a unique $(S, u) \in \mathcal{A}_{smooth}$ such that the intersection of $C$ with $\mathbb{R} \times M_\varepsilon$ can be written as $\mathfrak{e}_S \circ \eta$ where $\eta$ is a section of $N_0 \to S$ with the following properties:*
  a) *The conclusions of Proposition 7.2 hold using the given value of $\varepsilon$ and $(S, u)$.*
  b) *The norm of $\eta$ is bounded by $\kappa\varepsilon$.*
  c) *The section $\eta$ obeys (6.17) on the $u^*t \in [1+\varepsilon^2, 2-\varepsilon^2]$ part of $S$ and the respective restrictions of $\eta$ and the kernel of $D_S$ are $L^2$-orthogonal on this same part of $S$.*

The proof of this proposition is given momentarily. By way of a parenthetical remark, what is said in this proposition is not needed if (S´, u´) is in a 0 or 1 dimensional component of $\mathcal{A}_{smooth}$ because the latter is a single orbit of the $\mathbb{R}$ action. As it turns out, only those cases where the relevant part of $\mathcal{A}_{smooth}$ is a single $\mathbb{R}$ orbit are needed for the constructions in [KLTIII].

The three steps that follow directly provide background for the proof Proposition 7.3. These steps and the subsequent proof use $c$ to denote a constant with value greater than 1 that depends only on $\mathcal{K}$. Its precise value increases between subsequent appearances.

Step 1: Fix an open, $\mathbb{R}$-invariant set in $\mathcal{A}_{smooth}$ with compact closure in $\mathcal{A}_{smooth}/\mathbb{R}$ whose interior contains $\mathcal{K}$. Let $\mathcal{K}_*$ denote the closure of this $\mathbb{R}$-invariant set. Use in what follows the $\mathcal{K}_*$ version of the subbundle $\widetilde{\mathcal{N}_0}$ with the corresponding exponential map $\mathfrak{e}^{\mathcal{K}_*}$ and constant $\rho_{\mathcal{K}_*}$ to define any given $(S, u) \in \mathcal{K}_*$ versions of the bundle $N_0$ and the corresponding exponential map $\mathfrak{e}_S$. Since $\mathcal{K}_*$ is chosen given $\mathcal{K}$, the constant $\rho_{\mathcal{K}_*}$ can be assumed to depend only on $\mathcal{K}$. For the same reason, the derivatives of $\mathfrak{e}_S$ to any given order can be assumed to depend only on $\mathcal{K}$.

The use of a smoothly varying exponential map implies that the family of vector spaces {kernel($D_S$): $(S, u) \in \mathcal{K}_*$} fit together as a smooth vector bundle of $\mathcal{K}_*$.

Step 2: The remarks in this step are for the most part direct consequences of Lemmas 6.6 and 6.7 using the fact that the various $(S, u)$ versions of $N_0$ and $\mathfrak{e}_S$ come from $\widetilde{\mathcal{N}_0}$ and $\mathfrak{e}^{\mathcal{K}_*}$. In particular, the specific remarks that concern the $(S, u) \in \mathcal{K}_*$ versions of $D_S$ and kernel($D_S$) follow from Lemmas 6.6 and 6.7 and what is said in Step 1 using standard techniques for studying small perturbations of Fredholm operators.

The first remark concerns the size of Lemma 6.6's ball B: There exists $\rho_* > c^{-1}$ such that if $(S, u) \in \mathcal{K}$, then $(S, u)$ has a neighborhood in $\mathcal{K}_*$ that is parametrized by a ball B $\subset$ kernel($D_S$) of radius $\rho_*$. Let q denote the map given by Lemma 6.6. Fix $\eta_0 \in$ B and let (S´, u´) denote the admissible pair that corresponds to $\eta_0$. By way of a reminder, the 2-dimensional surface S´ is the same as S, but its complex structure is now induced by $u´$. Meanwhile, $u´ = \mathfrak{e}_S \circ q(\eta_0)$. The surface S´ has its complex normal bundle, N´; and the identification between $u´(S´)$ and the image of S via the map $\mathfrak{e}_S \circ q(\eta_0)$ supplies an $\mathbb{R}$-linear identification of the bundle N´ with the normal bundle N $\to$ S.

The bundle N´ $\to$ S´ has its disk subbundle $N_0´$ and the corresponding exponential map $\mathfrak{e}_{S´}: N_0´ \to \mathbb{R} \times [1, 2] \times \Sigma$. The afore-mentioned identification of N´ with N identifies $N_0´$ with a subbundle in N. The latter is not $N_0$ in general. With $N_0´$ viewed in N, the exponential map $\mathfrak{e}_{S'}$ can be written on a somewhat smaller radius subbundle $N_1' \subset N_0'$ as a composition of first an embedding $\mathfrak{h}_{S'}: N_1' \to N_0$ and then $\mathfrak{e}_S$. Note in this regard that $\mathfrak{h}_{S'}$ is neither linear nor base point preserving in general. In any event, the following can be assumed:

- *Let $\rho$ denote the constant fiber radius of the disk bundle $\mathcal{N}_0$. The subbundle $N_1'$ has radius $\rho - c|\eta_0|$.*
- *The $\mathfrak{h}_{S'}$ image of $N_1'$ contains the subbundle $N_1 \subset N_0$ of radius $\rho - c^2|\eta_0|$.*
- *The map from $B \times N_1$ to $N_0$ that restricts to any given $\eta_0 \in B$ slice of the domain as $(\mathfrak{h}_{S'} \circ \mathfrak{e}_{S'}^{-1}) \circ \mathfrak{e}_S$ is such that*
  a) *Its differential along $B$ at the origin is -1 times the identity.*
  b) *The norms of its derivatives to any given order are bounded by $c$.*

(7.10)

    <u>Step 3</u>: Suppose again that $(S, u) \in \mathcal{K}$ and that $(S', u') \in \mathcal{K}_*$ is parametrized as in Lemmas 6.6 and 6.7 by a given $\eta_0 \in B$. The identification between $N'$ and $N$ identifies the vector space kernel$(D_{S'})$ with a subspace, $H_{S'}$, in $C^\infty(S; N) \cap L^2(S; N)$.

    Let $\Pi': L^2(S'; N') \to$ kernel$(D_{S'})$ denote the $L^2(S'; N')$ orthogonal projection. With $N'$ viewed as $N$ and with kernel$(D_{S'})$ identified with $H_{S'}$, this projection will not in general agree with the $L^2(S; N)$ orthogonal projection. Even so, the following can be said: Let $\Pi: L^2(S; N) \to$ kernel$(D_S)$ denote the $L^2(S; N)$ orthogonal projection. There is a linear isomorphism $H_{S'}:$ kernel$(D_{S'}) \to$ kernel$(D_S)$ such that

$$\|(H_{S'}\Pi' - \Pi)(\cdot)\|_{L^2(S;N)} \le c\|\eta_0\|_{L^2(S;N)} \|\cdot\|_{L^2(S';N')}.$$

(7.11)

Moreover, the linear map $(H_{S'}\Pi' - \Pi)$ varies smoothly as $\eta_0$ varies in $B$; and in particular, its differential and that of its derivatives to any given order have norm bounded by $c$.

***Proof of Proposition 7.3***: The proof of the first bullet is identical but for the introduction of an extra $\mathbb{R}$-dependent cut-off function. This being the case, no more is said here about the first bullet. Note that the roles of $(S', u')$ and $(S, u)$ have been inverted in the proof. This is to avoid a proliferation of 'primed' symbols. This understood, it is assumed in what follows that $(S, u) \in \mathcal{K}$ is such that the assumptions of Proposition 7.2 hold for C with $\Xi = ((S, u), \emptyset)$. Meanwhile, it is assumed that the conclusions of the second bullet of Proposition 7.3 hold for a unique pair $(S', u') \in \mathcal{A}_{smooth}$. The proof has two steps. The first finds the desired $(S', u')$ and the second proves that there is only one such $(S', u')$.

Step 1:  A purely $\mathcal{K}$-dependent upper bound for $\varepsilon$ can be chosen so that C's intersection with the $u*t \in [1+\varepsilon^2, 2-\varepsilon^2]$ part of $\mathbb{R} \times [1,2] \times M_\varepsilon$ is in the image of via $\mathfrak{e}_S$ of the bundle $N_1$.  This being the case, there is a smooth section, $\eta$, of $N_1$ over this same part of S with norm bounded by $c\varepsilon$ and such that this part of C can be written as $\mathfrak{e}_S \circ \eta$. View $\eta$ as a piecewise smooth section of N over the whole of S by declaring it to be zero on the part of S where $u*t$ is not in $[1+\varepsilon^2, 2-\varepsilon^2]$.  This done, use $\eta_C \in \text{kernel}(D_S)$ to denote $\Pi\eta$.

Let $\eta_0 \in B$ and let $(S´, u´)$ denote the corresponding admissible pair.  The subvariety C can also be written $\mathfrak{e}_{S´} \circ \eta´$ where $\eta´$ is a section of $\mathfrak{h}_{S´}^{-1}(N_1) \subset N_1´$ over the $u*t \in [1+\varepsilon^2, 2-\varepsilon^2]$ part of S´.  Extend $\eta´$ over the remainder of S´ as zero.  It follows from (7.10) and (7.11) that

$$H_{S´} \Pi´\eta´ = \eta_C - \eta_0 + \mathfrak{r}(\eta)$$

(7.12)

where $\mathfrak{r}(\eta)$ has $L^2$ norm bounded $c \|\eta_0\|_{L^2(S;N)} \|\eta\|_{L^2(S´;N´)}$.

What with (7.10) and (7.11), the implicit function theorem finds $c \geq 1$ that makes the following true:  If $\varepsilon < c^{-1}$, then there exists a $\eta_0 \in B$ such that the corresponding $\eta´$ obeys $\Pi´\eta´ = 0$.  Moreover, $\eta_0$ differs from $\eta_C$ by at most $c\varepsilon^2$.

Step 2:  Suppose next that $(S_1´, u_1´)$ and $(S_2´, u_2´)$ satisfy the conclusions of the second bullet of the proposition.  Given that $\mathcal{K}_*$ is compact, there exists $c \geq 1$ such that if $\varepsilon < c^{-1}$, then the conclusions of Proposition 7.2 require that $u_1´(S_1´)$ and $u_2´(S_2´)$ be written respectively as $\mathfrak{e}_S(q(\eta_{01}))$ and $\mathfrak{e}_S(q(\eta_{02}))$ where $\eta_{01}$ and $\eta_{02}$ are points in B with norm bounded by $c\varepsilon$.

Let $\eta_1´$ and $\eta_2´$ denote the respective $(S_1´, u_1´)$ and $(S_2´, u_2´)$ versions of $\eta´$.  The left hand sides of the respective $(S_1´, u_1´)$ and $(S_2´, u_2´)$ versions of (7.12) are zero by definition, and as a consequence, subtracting these two versions of (7.12) yields an equation of the form

$$\eta_{02} - \eta_{01} + \mathfrak{r}_{1,2} = 0$$

(7.13)

where $|\mathfrak{r}_{1,2}| \leq c\varepsilon |\eta_{02} - \eta_{01}|$.  An equation of this sort implies that $\eta_{01} = \eta_{02}$ if $\varepsilon \leq c^{-1}$.

**A. Appendix**

The purpose of this appendix is to supply a proof of the assertion that a certain twisted version of embedded contact homology on Y can be defined by the rules laid out by Hutchings in [Hu1]. This assertion is the content of Theorem 2.2 in [KLTI]. The first part of this appendix briefly describes the relevant chain complex, its differential and various other important endomorphisms. The upcoming Theorem A.1 restates the assertion that these endomorphisms are well defined. The remaining parts of this appendix explain why Theorem A.1 is true.

**a) The ech chain complex**

This subsection describes the embedded contact homology chain complex of interest, the differential, and various other relevant endomorphisms. The discussion here has four parts.

*Part 1*: The relevant embedded contact homology chain complex is the free $\mathbb{Z}$ module that is generated by certain principle $\mathbb{Z}$ bundle over Proposition 2.8's set $\mathcal{Z}_{ech,M}$. The definition is defined as follows: Fix a fiducial element $\Theta_0 \in \mathcal{Z}_{ech,M}$. This done, use $\hat{\mathcal{Z}}_{ech,M}$ to denote the set of equivalence classes of pairs of the form $(\Theta, Z)$ where $\Theta \in \mathcal{Z}_{ech,M}$ and where $Z \in H_2(M; [\Theta] - [\Theta_0])$. To state the equivalence relation, introduce by way of notation, $\gamma^{(z_0)}$ to denote the closed integral curve of $v$ that contains the fiducial point $z_0 \in \Sigma$. Pairing with the Poincaré dual of $\gamma^{(z_0)}$ defines a homomorphism from the $\mathbb{Z}$-module of closed 2-cycles to $\mathbb{Z}$. This pairing is denoted by $\langle \gamma^{(z_0)}, \cdot \rangle$. Granted that such is the case, the equivalence relation posits

$$(\Theta, Z) \sim (\Theta', Z') \text{ if and only if both } \Theta = \Theta' \text{ and } \langle \gamma^{(z_0)}, Z - Z' \rangle = 0.$$
(A.1)

The element $1 \in \mathbb{Z}$ acts to send a give element $(\Theta, Z)$ to $(\Theta, Z + [S_0])$ where $[S_0]$ is the cycle given by the $u = 0$ sphere in $\mathcal{H}_0$.

The embedded contact homology chain complex of interest is the free $\mathbb{Z}$-module generated by $\hat{\mathcal{Z}}_{ech,M}$. This module is denoted in what follows by $\mathbb{Z}(\hat{\mathcal{Z}}_{ech,M})$. At the risk of pedantry, any given element in $\mathbb{Z}(\hat{\mathcal{Z}}_{ech,M})$ is a finite, integer weighted linear combination of generators.

The linear functional on $H_2(M; \mathbb{Z})$ given by the pairing with the class $c_{1M}$ maps to a subgroup of $2\mathbb{Z}$. Let $p_M$ denote the divisibility of this subgroup. With (A.1) understood, rules laid out by Hutchings [H2] can be employed in the context at hand to give each generator of $\hat{\mathcal{Z}}_{ech,M}$ a relative $\mathbb{Z}/(p_M\mathbb{Z})$ degree and so give the module $\mathbb{Z}(\hat{\mathcal{Z}}_{ech,M})$ a relative $\mathbb{Z}/(p_M\mathbb{Z})$ grading.

*Part 2*: The differential that defines the homology of this chain complex is given by a certain endomorphism of the module $\mathbb{Z}(\hat{\mathcal{Z}}_{ech,M})$. This endomorphism is defined by its action on each generators. The corresponding homology enjoys an action of the algebra $\mathbb{Z}[\mathbb{U}] \otimes (\wedge^*(H_1(Y; \mathbb{Z})/\text{torsion}))$ whose generators are also defined by endomorphisms of $\mathbb{Z}(\hat{\mathcal{Z}}_{\text{ech},M})$.

In general, a given endomorphism of $\mathbb{Z}(\hat{\mathcal{Z}}_{\text{ech},M})$ is defined by its action on the generating set. This is to say that it is determined by a rule that assigns to any given element $\hat{\Theta} \in \hat{\mathcal{Z}}_{\text{ech},M}$ a formal sum of the form

$$\hat{\Theta} \to \sum_{\hat{\Theta}' \in \hat{\mathcal{Z}}_{\text{ech},M}} N_{\hat{\Theta}',\hat{\Theta}} \, \hat{\Theta}' \, ,$$

(A.2)

where any given coefficient is an integer, and where only finitely many are non-zero. The endomorphism is thus defined by the corresponding set of coefficients. The definition of those that give the differential is described below. The definition of those that give that generate the action of $\mathbb{Z}(\mathbb{U}) \otimes (\wedge^*(H_1(Y; \mathbb{Z})/\text{torsion}))$ are described in the upcoming Part 3.

According to the rules laid out by Hutchings [Hu1], the coefficients that appear in (A.2) for the case of the differential are defined using certain sets of J-holomorphic subvarieties in $\mathbb{R} \times Y$. In particular, each orderd pair $\{\hat{\Theta}', \hat{\Theta}\} \subset \hat{\mathcal{Z}}_{\text{ech},M}$ labels such a set, this denoted in what follows by $\mathcal{M}_1(\hat{\Theta}', \hat{\Theta})$. To give the criteria for membership, write $\hat{\Theta}$ as $(\Theta, Z)$ and $\hat{\Theta}'$ as $(\Theta', Z')$. A given subvariety C is a member of this set when four conditions are satisfied, the first three being

$$\Theta_{C+} = \Theta, \;\; \Theta_{C-} = \Theta' \;\; and \;\; I_{\text{ech}}(C) = 1.$$

(A.3)

To state the fourth, introduce $[C]_Y$ to denote the element in $H_2(Y; [\Theta]-[\Theta'])$ that is defined by the C's image in Y via the projection from $\mathbb{R} \times Y$. The fourth condition asserts that

$$Z = [C]_Y + Z' \, .$$

(A.4)

The group $\mathbb{R}$ acts on $\mathcal{M}_1(\hat{\Theta}', \hat{\Theta})$ via constant translations of its elements along the $\mathbb{R}$ factor of $\mathbb{R} \times Y$. With this in mind, the coefficient $N_{\hat{\Theta}',\hat{\Theta}}$ in the version of (A.2) that gives the differential can be calculated if $\mathcal{M}_1(\hat{\Theta}', \hat{\Theta})/\mathbb{R}$ is a finite set, and if a certain Fredholm operator that is associated to each element in this set has trivial cokernel.

*Part 3*: The other relevant endomorphisms of $\mathbb{Z}(\hat{\mathcal{Z}}_{\text{ech},M})$ generate an action of the algebra $\mathbb{Z}[\mathbb{U}] \otimes (\wedge^*(H_1(Y; \mathbb{Z})/\text{torsion}))$ on the resulting homology. The coefficients that appear in the relevant versions of (A.2) are also computed using certain sets of J-holomorphic subvarieties on $\mathbb{R} \times Y$.

The definition of an endomorphism that defines the action of $\mathbb{U}$ on the homology of requires first the choice of a point in Y that does not lie on any closed integral curve of $v$ from elements in $\mathcal{Z}_{ech,M}$. Choose such a point in $\mathcal{H}_0$ or in $M_\delta$ where $f \notin [1, 2]$. Let y denote this point. According to Hutchings (see Section 11 of [HS]), a given $\hat{\Theta}', \hat{\Theta} \in \hat{\mathcal{Z}}_{ech,M}$ labeled coefficient $N_{\hat{\Theta}',\hat{\Theta}}$ in the corresponding version of (A.2) is computed from data supplied by the set, $\mathcal{M}_{2,y}(\hat{\Theta}',\hat{\Theta})$, of J-holomorphic subvarieties with membership defined as follows: A subvariety C is a member when (A.4) holds and when the $I_{ech}(C) = 2$ version of (A.3) holds. In addition, C must contain the point $(0, y) \in \mathbb{R} \times Y$. Let $\mathcal{M}_0(\hat{\Theta}',\hat{\Theta})$ denote the set of J-holomorphic subvarieties that obey (A.4) and the $I_{ech} = 0$ version of (A.3). It follows from Proposition 4.1 and what is said in Section 3 that the set of subvarieties that comprise $\mathcal{M}_{2,y}(\hat{\Theta}',\hat{\Theta})$ enjoy a 1-1 correspondence with the set that comprises $\mathcal{M}_0(\hat{\Theta}',\hat{\Theta})$ because a subvariety in the former is the union of some subvariety in the latter with the sphere from Proposition 3.1's moduli space that contains the point $(0, y)$. As noted parenthetically after Proposition 3.1, each sphere in the latter's moduli space has $I_{ech} = 2$. This identification of the respective coefficient in question can be calculated if the following conditions hold: The space $\mathcal{M}_0(\hat{\Theta}',\hat{\Theta}) = \emptyset$ when $\hat{\Theta}' \neq \hat{\Theta}$ and $\mathcal{M}_0(\hat{\Theta},\hat{\Theta})$ has but the one element $C = \cup_{\gamma \in \Theta}(\mathbb{R} \times \gamma)$. In addition, the cokernel of a certain Fredholm operator that is associated to each sphere in Proposition 3.1's must have trivial cokernel.

A set of endomorphisms of $\mathbb{Z}(\hat{\mathcal{Z}}_{ech,M})$ that generate the desired action of $\wedge^*(H_1(Y;\mathbb{Z})/\text{torsion})$ on the homology of the chain complex is defined using a chosen, suitably generic set of smooth 1-cycles in Y that represent a basis for $H_1(Y;\mathbb{Z})/\text{torsion}$. Each cycle from this set must be disjoint from the integral curves of $v$ that come from elements in $\mathcal{Z}_{ech,M}$. Let $\gamma$ denote such a cycle. With $\gamma$, chosen, any given $\hat{\Theta}', \hat{\Theta} \in \hat{\mathcal{Z}}_{ech,M}$ version of the coefficient $N_{\hat{\Theta}',\hat{\Theta}}$ for the corresponding endomorphism of $\mathbb{Z}(\hat{\mathcal{Z}}_{ech,M})$ is computed using the subset of subvarieties in $\mathcal{M}_1(\hat{\Theta}',\hat{\Theta})$ that intersect $\{0\} \times \gamma \subset \mathbb{R} \times Y$ and their corresponding intersection points. See also what is said in Section 11 of [HS].

*Part 4*: The theorem that follows asserts that the desired version of embedded contact homology can be defined in the manner outlined above. It also asserts that this version of contact homology admits the desired action of $\mathbb{Z}[\mathbb{U}] \otimes (\wedge^*(H_1(M;\mathbb{Z})/\text{torsion}))$. The theorem uses $\mathcal{J}$ to denote the $C^\infty$-Fréchét space of almost complex structures on $\mathbb{R} \times Y$ with the following two properties: Each $J \in \mathcal{J}$ obeys the constraints given in Section 3a; and each obeys (7.1) for some almost complex structure $J_{HF}$ that obeys (6.1).

**Theorem A.1**: *Fix the data $\delta, x_0, R$ to define the geometry of $Y$. There exists a residual set $\mathcal{J}_{ech} \subset \mathcal{J}$ with the following significance: Fix an almost complex structure $J \in \mathcal{J}_{ech}$ to define the almost complex geometry of $\mathbb{R} \times Y$. Then the endomorphisms that define the desired embedded contact homology differential on $\mathbb{Z}(\hat{\mathcal{Z}}_{ech,M})$ and the endomorphism that define the generators of the desired action of $\mathbb{Z}(\mathbb{U}) \otimes (\wedge^*(H_1(Y;\mathbb{Z})/\text{torsion}))$ on the resulting homology can be defined by Hutchings rules. The differential decreases the relative $\mathbb{Z}/(p_M\mathbb{Z})$ of each generator by 1 and so the resulting homology has a relative $\mathbb{Z}/(p_M\mathbb{Z})$ grading. The action of $\mathbb{Z}(\mathbb{U}) \otimes (\wedge^*(H_1(Y;\mathbb{Z})/\text{torsion}))$ is such that $\mathbb{U}$ decreases the relative grading degree by 2 and the generators of the action of $\wedge^*(H_1(Y;\mathbb{Z})/\text{torsion})$ decrease it by 1. The action of $\mathbb{Z}(\mathbb{U}) \otimes (\wedge^*(H_1(Y;\mathbb{Z})/\text{torsion}))$ on the homology is independent of the chosen point in $\mathcal{H}_0 \cup (M_\delta - f^{-1}([1,2]))$ that is used to define the $\mathbb{U}$-endomorphism, and it is independent of the cycles that are used to define the endomorphisms that give the action of $H_1(Y;\mathbb{Z})/\text{torsion}$.*

The proof of this theorem is given momentarily. By way of a look ahead, the proof has four aspects. The first verifies that (A.5) and (A.6) hold when J comes from a certain residual subset of $\mathcal{J}$. The second verifies that the endomorphism that describes the differential has square zero. The third verifies that the endomorphisms of $\mathbb{Z}(\hat{\mathcal{Z}}_{ech,M})$ that are meant to define the action of $\mathbb{U}$ and of $H_1(Y;\mathbb{Z})$ induce an action of the algebra $\mathbb{Z}[\mathbb{U}] \otimes (\wedge^*(H_1(Y;\mathbb{Z})/\text{torsion}))$. The fourth verifies that this action is independent of the chosen point in $Y$ and the chosen 1-cycles that are used to construct the defining endomorphisms of $\mathbb{Z}(\hat{\mathcal{Z}}_{ech,M})$. The first aspect of the proof constitutes Subappendix Ab, and the remaining aspects constitute Subappendix Ac. What is said about the relative grading follows from the definitions of the endomorphisms.

**b) Proof of Theorem A.1: The endomorphisms**

Fix an almost complex structure $J \in \mathcal{J}$. What is said in Part 2 of Appendix Aa about the definition of the differential on $\mathbb{Z}(\hat{\mathcal{Z}}_{ech,M})$ can be summarized as follows: The endomorphism that defines the differential can computed according to the rules laid out by Hutching when the following conditions are met:

- *Fix $\hat{\Theta} \in \hat{\mathcal{Z}}_{ech,M}$. All but a finite set of $\hat{\Theta}' \in \hat{\mathcal{Z}}_{ech,M}$ versions of $\mathcal{M}_1(\hat{\Theta}', \hat{\Theta})$ are empty.*
- *Each $\hat{\Theta}', \hat{\Theta} \in \hat{\mathcal{Z}}_{ech,M}$ version of $\mathcal{M}_1(\hat{\Theta}', \hat{\Theta})$ is a finite set of $\mathbb{R}$-orbits.*
- *The Fredholm operator associated to each such $\mathbb{R}$-orbit has trivial cokernel.*

(A.5)

What is said in Part 3 of Appendix Aa concerning the $\mathbb{U}$ endomorphism can be summarized as follows: The desired endomorphism can defined using the rules laid out by Hutchings if there exists $y \in Y$ such that the following conditions hold:

- *If $\hat{\Theta}$, $\hat{\Theta}' \in \hat{\mathcal{Z}}_{ech,M}$ are distinct, then $\mathcal{M}_0(\hat{\Theta}', \hat{\Theta}) = \emptyset$.*
- *If $\hat{\Theta} \in \hat{\mathcal{Z}}_{ech,M}$ then $\mathcal{M}_0(\hat{\Theta}, \hat{\Theta})$ has but the one element $C = \cup_{\gamma \in \Theta}(\mathbb{R} \times \gamma)$.*
- *The relevant Fredholm operator for any given sphere from Proposition 3.1's moduli space has trivial cokernel*

(A.6)

What is said in this same part of Appendix Aa also implies the following: The endomorphisms on $\mathbb{Z}(\hat{\mathcal{Z}}_{ech,M})$ that generate the action of $\wedge^*(H_1(Y;\mathbb{Z})/\text{torsion})$ on the homology can be defined by Hutchings' rules from a chosen, suitably generic basis set of 1-cycles when (A.5) holds.

The three parts that follow explain why (A.5) and (A.6) hold for the almost complex structures in a certain residual subset of $\mathcal{J}$.

*Part 1*: This part and Part 2 address the cokernel issue. To this end, let C denote a given J-holomorphic subvariety in $\mathbb{R} \times Y$ with $\Theta_{C+} = \Theta$ and $\Theta_{C-} = \Theta'$. The Fredholm operator in question is the standard deformation operator for J-holomorphic subvarieties. If C is immersed, then it has the same form as the operator that is depicted in (6.11). The operator in this case defines an $\mathbb{R}$-linear, Fredholm operator with domain the $L^2_1$ Hilbert space of sections of C's complex normal bundle and range the $L^2$ Hilbert space of sections of the tensor product of the latter bundle with the complex line bundle $T^{0,1}C$. The operator in the general case is described in Section 4 of [HT].

The set of almost complex structures on $\mathbb{R} \times Y$ is given a $C^\infty$-Fréchét manifold structure by viewing it as a submanifold in $C^\infty(\mathbb{R} \times Y; \text{End}(T(\mathbb{R} \times Y)))$. The set $\mathcal{J}$ likewise inherits the structure of a $C^\infty$-Fréchét manifold. A standard argument (see, e.g. the proof of Theorem 1.8 in [Hu3] or what is said in [HT]) can be used to find a $C^\infty$-residual subset in the space of almost complex structures satisfying (3.4) and the CONSTRAINTS 1 and 2 in Section 3a whose members are characterized as follows: Take any two elements from $\mathcal{Z}_{ech,M}$, and any pseudoholomorphic subvariety whose positive and negative ends are associated to the chosen elements. Then the corresponding Fredholm operator has trivial cokernel.

A close look at the proof (see, eg the appendix in [T2]) gives a slightly weaker statement if the residual set is to sit in $\mathcal{J}$. In particular, arguments that mimic those in Section 3 of [L] can be used to find a $C^\infty$-residual set in $\mathcal{J}$ whose members are characterized as follows:

*Fix an almost complex structure from this set. Take any two elements from $\mathcal{Z}_{ech,M}$, and any pseudoholomorphic subvariety whose positive and negative ends are associated to the chosen elements. Then the corresponding Fredholm operator has trivial cokernel if the subvariety has no irreducible components from Proposition 3.2's space $\mathcal{M}_\Sigma$.*

(A.7)

As explained next, the third bullet in (A.5) and all of (A.6) hold if J is from this residual. Part 2 of this subappendix describes what is involved in the proof that (A.7) characterizes a residual subset of $\mathcal{J}$.

Take $J \in \mathcal{J}$ so that (A.7) holds. Consider first (A.5). It follows from (1.1) in [Hu2] that each subvariety in any given $\hat{\Theta}'$, $\hat{\Theta} \in \hat{\mathcal{Z}}_{ech,M}$ version of $\mathcal{M}_1(\hat{\Theta}', \hat{\Theta})$ is an embedded submanifold. This being the case, it lacks irreducible components from $\mathcal{M}_\Sigma$ and so the conclusions of (A.7) apply. The third bullet in (A.5) restates these conclusions.

To see about (A.6), choose the point y in either $\mathcal{H}_0$ or in $M_\delta$ where $f$ is either between (0, 1) or (2, 3). Hutching's index inequality in (1.1) of [Hu2] gives the first two bullets in (A.6). The assertion given by the third bullet follows using degree arguments along the lines of those used in Section 2.3 of [T5]. These degree arguments do not require (A.7).

*Part 2*: This part of the sub-appendix explains why $\mathcal{J}$ has a $C^\infty$-residual subset whose members are characterized by (A.7). To start, fix $J \in \mathcal{J}$, a pair of elements from $\mathcal{Z}_{ech,M}$ and a J-holomorphic subvariety whose respective postive and negative ends are associated to the chosen elements. Let C denote this subvariety and let $\mathfrak{D}_C$ denote the associated Fredholm operator. The complement in C of a finite set of points will be embedded and so have a complex normal bundle. An endomorphism from $T\mathcal{J}|_J$ on this part of C defines a homorphism from C's anti-canonical bundle to C's complex normal bundle. If the section has support in this part of C, then it defines a section of the domain space of $\mathfrak{D}_C$. With this point in mind, let $\mathcal{T} \subset T\mathcal{J}$ denote a given subbundle. The elements of $\mathcal{T}|_J$ defines in the manner outlined above, a subspace of elements in the range space of $\mathfrak{D}_C$. This subspace is denoted in what follows by $\mathcal{R}_\mathcal{T}|_C$. Let $\Pi_C$ denote the tautological projection from the range Hilbert space of $\mathfrak{D}_C$ to the cokernel of $\mathfrak{D}_C$. A subbundle $\mathcal{T} \subset T\mathcal{J}$ is said to be *large* if the following is true:

*Suppose that $J \in \mathcal{J}$ and that C is an ech-HF subvariety as defined by J. Then $\Pi_C$ is surjective on the subspace $\mathcal{R}_\mathcal{T}|_C$.*

(A.8)

Suppose that $\mathcal{T} \subset T\mathcal{J}$ is large, and also tangent to a foliation of $\mathcal{J}$. The Smale-Sard theorem can be used to prove the following: Fix $J' \in \mathcal{J}$, an integer $k \geq 0$ and $\varepsilon > 0$.

Let $\mathcal{F}_{J'}$ denote the leaf of the foliation through J´. There exist endomorphisms $j_0 \in \mathcal{T}|_{J'}$ and $j_1$ of $T(\mathbb{R} \times Y)$ such that

- $J = J' + j_0 + j_1 \in \mathcal{F}_{J'}$ *is such that (A.7) holds*.
- *The $C^k$-norm of $j_0$ is less than ε and that of $j_1$ is at most ε times the $C^k$-norm of $j_0$*.

(A.9)

Any given J´ ∈ $\mathcal{J}$ can be modified by changing its corresponding $J_{HF}$ subject to (6.1). Modifications of this sort define a foliation of $\mathcal{J}$ whose tangent subbundle in $T\mathcal{J}$ is large. This can be proved using Aronzajn's unique continuation principle [Ar].

*Part 3*: This part of the subappendix explains why the first two bullets of (A.5) hold when J is characterized by (A.7). The first bullet follows directly from (5.5) and Lemma 5.2. To see about the second bullet, note that (A.7) implies that any given $\hat{\Theta}'$, $\hat{\Theta}$ ∈ $\hat{\mathcal{Z}}_{ech,M}$ of $\mathcal{M}_1(\hat{\Theta}', \hat{\Theta})$ consists of a countable set of 1-dimensional manifolds and that each such manifold is an orbit of the action of $\mathbb{R}$. This understood, the second bullet in (A.5) follows if $\mathcal{M}_1(\hat{\Theta}', \hat{\Theta})$ is compact. But for one additional comment, the proof that $\mathcal{M}_1(\hat{\Theta}', \hat{\Theta})$ is compact uses the bounds in Proposition 5.1 as input for arguments that differ only in notation from those that prove Assertion (a)(i) of Theorem 1.8 in [Hu3].

The extra comment is needed to address the fact that the situation here corresponds to the case of Theorem 1.8 in [Hu3] when what Hutchings denotes by d is the same as the genus of Σ. The analog here of Hutching's surface Σ is the Heegaard surface, and the analog here of d is the number of segments that comprise the $M_δ$ part of the union of the integral curves of ν from any given Θ ∈ $\hat{\mathcal{Z}}_{ech,M}$. The versions here of d and genus Σ are equal, but Assertion (a)(i) of Theorem 1.8 of [Hu3] as stated requires d > genus(Σ). Even so, the $I_{ech}$ = 1 and d = genus(Σ) version of Assertion (a)(i) of this Theorem 1.8 is true [Hu4]. In fact, the argument for this case is identical to that given in [Hu3] for the d > genus(Σ) cases. To say slightly more about why this is, note that the proof of Assertion (a)(i) invokes an inequality for the embedded contact homology index; this is the inequality written in Case 2 of the proof of Lemma 9.4 in [Hu3]. The proof that Assertion (a)(i) of Theorem 1.8 in [Hu3] holds requires that this inequality have negative right hand side, and it is this requirement that determines the relative sizes of d and the genus of Σ in the statement of Assertion (a)(i). Meanwhile, the right hand side of this inequality is observedly negative when $I_{ech}$ = 1 and d = genus(Σ).

By way of comparison, there is a Heegaard Floer analog of the $I_{ech}$ = 1 and d = genus(Σ) version of Assertion (i)(a) of Theorem 1.8 in [Hu3], this being Lemma 8.2 in [L]. What is more, the proof of the latter invokes the Heegaard Floer analog of the inequality that appears in Case 2 of the proof of Lemma 9.4 of [Hu3].

**c) The square of the differential and the action of $\mathbb{Z}[\mathbb{U}] \otimes (\wedge^*(H_1(Y;\mathbb{Z})/\text{torsion})$**

Fix an almost complex structure for which (2.7) holds so as to use Hutchings' rules to define the endomorphism of $\mathbb{Z}(\hat{\mathcal{Z}}_{\text{ech},M})$ that is meant to be the desired differential. Part of what follows explains why this endomorphism has square zero. Fix a point, y, in either $\mathcal{H}_0$ or in the part of $M_\delta$ where $f \notin [1,2]$ to define the $\mathbb{U}$ endomorphism. Likewise, fix a basis of smooth cycles for $H_1(Y;\mathbb{Z})/\text{torsion}$ to define the endomorphisms that are meant to generate the $\wedge^*(H_1(Y;\mathbb{Z})/\text{torsion})$ action on the embedded contact homology. What follows also expains why these endomorphisms generate the desired action of $\mathbb{Z}[\mathbb{U}] \otimes (\wedge^*(H_1(Y;\mathbb{Z})/\text{torsion})$, and why this action is independent of chosen point y and of the chosen basis of 1-cycles.

To handle these assigned tasks, introduce the nested exhaustion of $\mathcal{Z}_{\text{ech},M}$ by $\{\mathcal{Z}_{\text{ech},M}^L\}_{L \geq 1}$ of subsets with that are defined here as follows: Define the length of any given $\Theta \in \mathcal{Z}_{\text{ech},M}$ to be the sum of the integrals of $\hat{a}$ over the the closed integral curves of $v$ from $\Theta$. The subset $\mathcal{Z}_{\text{ech},M}^L$ contains the elements with length less than L. Any given $\mathcal{Z}_{\text{ech},M}^L$ is a finite set and the union of all of them is $\mathcal{Z}_{\text{ech},M}$. Moreover, if $\Theta \in \mathcal{Z}_{\text{ech},M}^L$ and if $C \subset \mathbb{R} \times Y$ is a J-holomorphic subvariety with $\Theta_{C+} \in \mathcal{Z}_{\text{ech},M}^L$ and with $\Theta_{C-} \in \mathcal{Z}_{\text{ech},M}$, then $\Theta_{C-}$ is also in $\mathcal{Z}_{\text{ech},M}^L$.

For any given $L \geq 1$, use $\hat{\mathcal{Z}}_{\text{ech},M}^L$ to denote the restriction of $\hat{\mathcal{Z}}_{\text{ech},M}$ to $\mathcal{Z}_{\text{ech},M}$. Each such set generates the corresponding submodule $\mathbb{Z}(\hat{\mathcal{Z}}_{\text{ech},M}^L) \subset \mathbb{Z}(\hat{\mathcal{Z}}_{\text{ech},M})$. It follow from what is at the end of the preceding paragraph that each such submodule is mapped to itself by the differential, the $\mathbb{U}$-endomorphism, and by the endomorphisms that are meant to generate the action of $\wedge^*(H_1(Y;\mathbb{Z})/\text{torsion})$. This being the case, it is sufficient to prove that $\delta^2 = 0$ on each $L \geq 1$ version of $\mathbb{Z}(\hat{\mathcal{Z}}_{\text{ech},M}^L)$. Likewise, it is sufficient to prove the the assertions in Theorem A.1 about the other endomorphisms for their restrictions to each $L \geq 1$ version of $\mathbb{Z}(\hat{\mathcal{Z}}_{\text{ech},M}^L)$.

With the preceding understood, fix $L \geq 1$. What follows explains why the differential has square zero on $\mathbb{Z}(\hat{\mathcal{Z}}_{\text{ech},M}^L)$. The arguments for the assertions that concern the other endomorphisms amount to straightforward variations of that given below and so the latter are omitted.

To start, invoke (A.5) to conclude that the set of $\mathbb{R}$ orbits that comprise $\cup_{\hat{\Theta}',\hat{\Theta} \in \hat{\mathcal{Z}}_{\text{ech},M}^L} \mathcal{M}_1(\hat{\Theta}',\hat{\Theta})$ is finite. What follows is a consequence of this fact. Fix a $C^\infty$ neighborhood of J in the space of almost complex structures that obey (3.4) and also CONSTRAINTS 1 and 2 in Section 3a. If the J´ is in such a neighborhood, and if J - J´ has sufficiently small $C^6$ norm, then there is a J´ analog of any given $\hat{\Theta}', \hat{\Theta} \in \hat{\mathcal{Z}}_{\text{ech},M}^L$ version of $\mathcal{M}_1(\hat{\Theta}',\hat{\Theta})$, and the elements in the latter enjoy a 1-1 correspondence with those in J's version of $\mathcal{M}_1(\hat{\Theta}', \hat{\Theta})$. The correspondence is such that if C is a submanifold from the J version, and C´ is from the J´ version, then C´ sits in a small radius tubular neighborhood of C and is isotopic to C in this neighborhood by an isotopy whose $C^k$ norm for any given k ∈ {1, 2,...} is bounded by a k-dependent multiple of the $C^{k+3}$ norm of J´ - J. This partner C´ can be constructed using perturbation theoretic methods because the cokernel of $\mathfrak{D}_C$ is trivial.

These perturbation methods construct an injective map from the J version of $\mathcal{M}_1(\hat{\Theta}', \hat{\Theta})$ to the J´ version. Meanwhile, a limiting argument that differs only cosmetically from what is said in the proof of Assertion (a)(i) from Theorem 1.8 in [Hu3] proves this injective map is also surjective onto the J´ version of $\mathcal{M}_1(\hat{\Theta}', \hat{\Theta})$.

Let C´ denote a submanifold from the J´ version of $\mathcal{M}_1(\hat{\Theta}', \hat{\Theta})$. The fact that this submanifold is $C^6$ close to its partner in the J version of $\mathcal{M}_1(\hat{\Theta}', \hat{\Theta})$ implies that the operator $\mathfrak{D}_{C'}$ also has trivial cokernel. More to the point, it implies that the coefficient $N_{\hat{\Theta}', \hat{\Theta}}$ that appears in (A.2) in the case of the differential can be computed using the submanifolds from the J´ version of $\mathcal{M}_1(\hat{\Theta}', \hat{\Theta})$. As a consquence, the endomorphism that defines the action of the differential on $\mathbb{Z}(\hat{\mathcal{Z}}_{ech,M}^L)$ can be computed using J´ in lieu or J. This is important by virtue of what is said in the second paragraph of Part 1 of Appendix Ab: The almost complex structure J´ can be chosen so that every J´-holomorphic subvariety whose respective positive and negative ends are associated to elements in $\mathcal{Z}_{ech,M}$ is such that its corresponding Fredholm operator has trivial cokernel. This being the case, the arguments in [HT] can be repeated with only cosmetic changes to see that the J´ version of the differential has square zero on $\mathbb{Z}(\hat{\mathcal{Z}}_{ech,M}^L)$. As noted a moment ago, latter endomorphism is identical to the J version; and so J's differential also has square zero on $\mathbb{Z}(\hat{\mathcal{Z}}_{ech,M}^L)$.

**References**

[Ah]       L. Ahlfors, Complex Analysis; McGraw Hill 1979.
[Ar]       N. Aronszajn, *A unique continuation theorem for elliptic equations or inequalities of the second order*; J. Math. Pures. Appl. **36** (1957), 253–249.
[HWZ]   H. Hofer, K. Wysocki, and E. Zehnder, *Properties of pseudoholomorphic curves in symplectizations. I. Asymptotics*, Ann. Inst. H. Poincare Anal. Non Lineaire **13** (1996), 337-379.
[Hum]   C. Hummel, Gromov's compactness theorem for pseudoholomorphic curves, Birkhauser 1997.

| | |
|---|---|
| [Hu1] | M. Hutchings, *Embedded contact homology and its applications*; Proceedings of the International Congress of Mathematicians, Hyderabad India 2010, and arXiv:1003.3209v1. |
| [Hu2] | M. Hutchings, *The embedded contact homology index revisited*; in <u>New Perspective and Challenges in Symplectic Field Theory</u>, 263-267, CRM Proc. Lecture Notes, **49**, American Mathematical Society 2009. |
| [Hu3] | M. Hutchings, *An index inequality for embedded pseudoholomorphic curves in symplectizations*; J. Eur. Math. Soc. **4** (2002) 313-361. |
| [Hu4] | M. Hutchings, private communications. |
| [HS] | M. Hutchings and M. Sullivan, *Rounding corners of polygons and the embedded contact homology of $T^3$*; Geom. and Topol. **10** (2006) 169-266. |
| [HT] | M. Hutchings and C. H. Taubes, *Gluing pseudoholomorphic cylinders along branched covers II*; J. Symplectic Geom. **7** (2009) 29-133. |
| [KLTI] | C. Kutluhan, Y-J. Lee and C. H. Taubes, *HF=HM I: Heegaard Floer homology and Seiberg-Witten Floer homology;* arXiv . |
| [KLTIII] | C. Kutluhan, Y-J. Lee and C. H. Taubes, *HF=HM III: Holomorphic curves and the differential for the ech and Heegard Floer homology correspondence;* in preparation. |
| [L] | R. Lipshitz, *A cylindrical reformulation of Heegaard Floer homology*; Geom. Topol. **10** (2006) 955-1097. |
| [Mc] | D. McDuff, *The local behavior of pseudoholomorphic curves in almost complex four manifolds*; J. Differential Geom., **34** (1991), 143–164. |
| [MS] | D. McDuff and D. Salamon, <u>J-holomorphic curves and quantum cohomology</u>, University Lecture Series **6**, American Mathematical Society 1994. |
| [Mo] | C. B. Morrey, <u>Multiple Integrals in the Calculus of Variations</u>; Springer-Verlag 1966. |
| [OS1] | P. Ozsváth and Z. Szabó, *Holomorphic disks and topological invariants for closed 3-manifolds*; Annals of Math **159** (2004) 1027-1158. |
| [OS2] | P. Ozsváth and Z. Szabó, *Heegaard diagrams and Floer homology*; pages 1083-1099 <u>Proceedings of the International Congress of Mathematicians, Madrid 2006, Volume 2</u>,American Mathematical Society, 2007. |
| [S] | R. Siefring, *The relative asymptotic behavior of pseudoholomorphic half-cylinders*; Comm. Pure Appl. Math. **61** (2008) 1631-1684. |
| [T1] | C. H. Taubes, *A compendium of pseudoholomorphic beasts in $\mathbb{R}\times(S^1\times S^2)$*; Geom. and Topol, **6** (2002) 657-814. |
| [T2] | C. H. Taubes, *Tamed to compatible: Symplectic forms via moduli space integration*; arXiv:0910.5440. |
132

[T3]     C. H. Taubes, *SW=>Gr: From the Seiberg-Witten equations to pseudoholomorphic curves*; in Seiberg-Witten and Gromov invariants for symplectic 4-manifolds, International Press, Boston MA, USA; 2000.

[T4]     C. H. Taubes, *Embedded contact homology and Seiberg-Witten Floer homology IV*, Geom. and Topol., to appear.

[T5]     C. H. Taubes, *Pseudoholomorphic punctured spheres in $\mathbb{R} \times (S^1 \times S^2)$: Properties and existence*; Geom. and Topol. **10** (2006) 785-928.

[W]      J. Wolfson, *Gromov's compactness of pseudoholomorphic curves and symplectic geometry*; J. Differential Geom. **28** (1988) 383–405.

[Y]      R. Ye, *Gromov's compactness theorem for pseudo holomorphic curves*; Trans. Amer. Math. Soc. 342 (1994), no. 2, 671–694.

Cagatay Kutluhan; Department of Mathematics, Columbia University, New York, NY 10027. *Email address*: kutluhan@umich.edu.

Yi-Jen Lee; Department of Mathematics, Purdue University, West Lafayette, IN 47907 *Email address*: yjlee@math.purdue.edu.

Clifford Henry Taubes; Department of Mathematics, Harvard University, Cambridge, MA 02138. *Email address*: chtaubes@math.harvard.edu.